\documentclass[10pt, a4paper,
oneside, headinclude,footinclude]{scrartcl}
\usepackage[utf8]{inputenc}

\usepackage{silence}
\usepackage[
nochapters, % Turn off chapters since this is an article        
pdfspacing, % Makes use of pdftex’ letter spacing capabilities via the microtype package
dottedtoc % Dotted lines leading to the page numbers in the table of contents
]{classicthesis} % The layout is based on the Classic Thesis style
\usepackage{amsopn}
\usepackage[T1]{fontenc} % Use 8-bit encoding that has 256 glyphs
\usepackage{multirow}
\usepackage[utf8]{inputenc} % Required for including letters with accents
\usepackage{hhline}
\usepackage{graphicx} % Required for including images
\graphicspath{{Figures/}} % Set the default folder for images
\usepackage{indentfirst}
\usepackage{enumitem} % Required for manipulating the whitespace between and within lists
\usepackage{savesym}
\usepackage{mathrsfs}
\usepackage{geometry}
\usepackage{esint}
\usepackage{longtable}

\usepackage[
    backend=biber,
    style=numeric,
    natbib=false,
    url=false, 
    doi=true,
    eprint=true,
    maxnames=50
]{biblatex}

\usepackage{tikz}
\usepackage[automark]{scrlayer-scrpage}

\definecolor{halfgray}{gray}{0.55}

\usepackage{stmaryrd}

\DeclareFontFamily{U}{rsfs}{\skewchar\font127}
\DeclareFontShape{U}{rsfs}{m}{n}{%
  <-6.5>rsfs5 <6.5-7.5>rsfs7 <7.5->rsfs10
}{}
\DeclareFontFamily{U}{stmry}{}
\DeclareFontShape{U}{stmry}{m}{n}{%
  <-5.5>stmary5 <5.5-6.5>stmary6 <6.5-7.5>stmary7
  <7.5-8.5>stmary8 <8.5-9.5>stmary9 <9.5->stmary10
}{}

\usepackage{subfig} % Required for creating figures with multiple parts (subfigures)

\usepackage{amsmath, amssymb,amsthm,amsfonts} % For including math equations, theorems, symbols, etc

\usepackage[english,capitalize]{cleveref}

\usepackage{thmtools}
\usepackage{thm-restate}

\theoremstyle{plain}
\newtheorem{teorema}{Theorem}[section]
\newtheorem{proposizione}[teorema]{Proposition}
\newtheorem{lemma}[teorema]{Lemma}
\newtheorem{corollario}[teorema]{Corollary}
\newtheorem*{theorem*}{Theorem}

\theoremstyle{definition}
\newtheorem{definizione}{Definition}[section]

\theoremstyle{remark}

\makeatletter
\AtBeginDocument{%
  \def\theHtheorem2{\theHsubsection.\arabic{theorem2}}%
  \def\theHproposition2{\theHsubsection.\arabic{proposition2}}%
}
\makeatother

\newcommand{\R}{\mathbb{R}}
\newcommand{\N}{\mathbb{N}}

\newcommand{\Z}{\mathbb{Z}}

\newcommand{\Haus}{\mathscr{H}}
\newcommand{\Leb}{\mathscr{L}}

\newcommand{\Lip}{\mathrm{Lip}}

\newcommand{\dist}{\mathrm{dist}}
\newcommand{\supp}{\mathrm{supp}}

\newcommand{\dV}{d_V\kern-1pt}
\newcommand{\dW}{d_W\kern-1pt}

\def\dist{\mathop\mathrm{dist}} 
\def\supp{\mathop\mathrm{supp}}

\DeclareFontFamily{U}{rcjhbltx}{}
\DeclareFontShape{U}{rcjhbltx}{m}{n}{<->rcjhbltx}{}
\DeclareSymbolFont{hebrewletters}{U}{rcjhbltx}{m}{n}

\DeclareMathSymbol{\lamed}{\mathord}{hebrewletters}{108}
\DeclareMathSymbol{\mem}{\mathord}{hebrewletters}{109}
\DeclareMathSymbol{\ayin}{\mathord}{hebrewletters}{96}
\DeclareMathSymbol{\tsadi}{\mathord}{hebrewletters}{118}
\DeclareMathSymbol{\qof}{\mathord}{hebrewletters}{113}
\DeclareMathSymbol{\shin}{\mathord}{hebrewletters}{152}
\DeclareMathSymbol{\tav}{\mathord}{hebrewletters}{90}

\def\XXint#1#2#3{{\setbox0=\hbox{$#1{#2#3}{\int}$ }
		\vcenter{\hbox{$#2#3$ }}\kern-.6\wd0}}

\def\R{\mathbb{R}}

\def\Lip{\mathop\mathrm{Lip}} 						%lipschitz constant \ck
\def\dist{\mathop\mathrm{dist}} 						%distance
\def\supp{\mathop\mathrm{supp}}					%support
\newcounter{const}
\makeatletter
\AtBeginDocument{\let\constant@label\label}
\newcommand{\newC}{%
  \begingroup
  \refstepcounter{const}\ensuremath{C_{\theconst}}%
  \@ifnextchar\label{\newC@label}{\endgroup}%
}
\def\newC@label\label#1{%
  \constant@label{#1}%
  \endgroup
}
\makeatother
\newcommand{\oldC}[1]{\ensuremath{C_{\ref{#1}}}}

\newcounter{eps}

\newcommand{\vertiii}[1]{{\left\vert\kern-0.25ex\left\vert\kern-0.25ex\left\vert #1 
    \right\vert\kern-0.25ex\right\vert\kern-0.25ex\right\vert}}

\hypersetup{
colorlinks=true, breaklinks=true,bookmarksnumbered,
urlcolor=webbrown, linkcolor=RoyalBlue, citecolor=webgreen, % Link colors
pdftitle={}, % PDF title
pdfauthor={\textcopyright}, % PDF Author
pdfsubject={}, % PDF Subject
pdfkeywords={}, % PDF Keywords
pdfcreator={pdfLaTeX}, % PDF Creator
pdfproducer={LaTeX with hyperref and ClassicThesis} % PDF producer
}

\usepackage{tikz}
\ifcsname teorema\endcsname\else\theoremstyle{plain}\ifcsname c@proposizione\endcsname\newtheorem{teorema}[proposizione]{Theorem}\else\newtheorem{teorema}{Theorem}[section]\fi\fi

\title{\normalfont\spacedallcaps{Sharp decay of the spherical averages of the Fourier transform of open sets}} 
\author{%
\spacedlowsmallcaps{Alex Iosevich}\thanks{Alex Iosevich was partially supported by NSF grant DMS-2506858 and a 2026 Simons Fellowship in Mathematics.},
\spacedlowsmallcaps{Andrea Merlo}\thanks{Andrea Merlo was supported by IKERBASQUE; by project PID2024-157724NB-I00, funded by MICIU/AEI/10.13039/501100011033/FEDER, EU; by the Europa Excelencia grant EUR2025-165042, funded by the Ministerio de Ciencia, Innovaci\'on y Universidades and Agencia Estatal de Investigaci\'on; and by grant IT1875-26 of the Basque Government.}\hspace{0.3em}
\spacedlowsmallcaps{and Mihalis Mourgoglou}\thanks{Mihalis Mourgoglou was supported by IKERBASQUE and partially supported by project PID2024-157724NB-I00, funded by MICIU/AEI/10.13039/501100011033/FEDER, EU, and by grants IT-1615-22 and IT1875-26 of the Basque Government.}.
}

\date{} 

\begin{document}

\renewcommand{\sectionmark}[1]{\markright{\spacedlowsmallcaps{#1}}} 
\lehead{\mbox{\llap{\small\thepage\kern1em\color{halfgray} \vline}\color{halfgray}\hspace{0.5em}\rightmark\hfil}} 
\pagestyle{scrheadings}
\maketitle 
\setcounter{tocdepth}{2}

\begin{abstract}
For a bounded domain $\Omega\subset\mathbb R^d$, the spherical average $E_\Omega(R)=\int_{\mathbb S^{d-1}}|\widehat{\mathbb{1}_\Omega}(R\omega)|^2d\omega$ is never $O(R^{-a})$ with $a>d+1$, and for smooth boundaries it is $O(R^{-d-1})$, as it is for the ball, whether or not the boundary is curved. In 2003, Brandolini, Hofmann and Iosevich extended this to every convex domain and to every domain with $C^{1,\frac{1}{2}}$ boundary, and this regularity condition has not been improved since. How much regularity does the decay $R^{-d-1}$ require? We show that among H\"older classes the answer is $C^{1,\frac{1}{4}}$. Every bounded domain with $C^{1,\frac{1}{4}}$ boundary satisfies $E_\Omega(R)\leq C_\Omega R^{-d-1}$, while for every $0<s<\frac{1}{4}$ any smooth domain can be altered near a single boundary point to produce a $C^{1,s}$ domain for which this bound fails along a sequence of radii. More generally, when the boundary has fractional Sobolev regularity $W^{1+s,p}$, we determine the optimal decay exponent under a Dini condition on the normal, and in a large range of $s$ and $p$ without it. The upper bounds follow from $L^2$ estimates for nondegenerate oscillatory integral operators whose phases are distances between points of two graphs with Dini or fractional Sobolev gradients, and hence barely better than $C^1$. These are far too rough for the integrations by parts techniques from the smooth theory. We replace derivatives by finite differences adapted to the phase, a device that we expect to be useful in a variety of related problems. 
\end{abstract}

\tableofcontents

\section{Introduction}

One of the most time-honored lines of investigation in harmonic analysis is the study of the decay of the Fourier transform of measures. Let $\mu$ be a compactly supported finite Borel measure, and let
$$ \widehat{\mu}(\xi)=\int e^{-i x \cdot \xi} d\mu(x)$$ denote its Fourier transform. In general, $|\widehat{\mu}(\xi)|$ need not decay as $|\xi|$ goes to infinity. For example, if $\mu$ is the natural measure on the triadic Cantor set, then $\widehat{\mu}(\xi)$ does not tend to $0$ as $|\xi| \to \infty$ due to the self-similar structure of the triadic Cantor set. See \cite[Chapter~8]{Mattila2015Fourier}. Similarly, if $\mu$ is a nonzero measure supported on a compact piece of a $(d-1)$-dimensional hyperplane in ${\mathbb R}^d$, the Fourier transform does not decay in the direction orthogonal to the hyperplane.

On the other hand, if $\mu$ is the surface measure on a smooth compact convex hypersurface in ${\mathbb R}^d$ with everywhere nonvanishing Gaussian curvature, then the classical method of stationary phase gives
$$ |\widehat{\mu}(\xi)| \leq C(1+|\xi|)^{-\frac{d-1}{2}}. $$
See \cite[Chapter~VIII, Section~3]{Stein1993Harmonic}.

The converse is also true. Iosevich and Lu \cite{IosevichLu2000} show that this uniform decay estimate for the surface measure on a smooth compact hypersurface forces the Gaussian curvature to be nonzero everywhere. Their argument passes through the sharp $L^2$ restriction theorem. For smooth compact convex hypersurfaces of finite type, Iosevich \cite{Iosevich1999Scaling} proves that $L^2$ restriction estimates imply the corresponding isotropic Fourier decay.

If the Gaussian curvature is allowed to vanish, the optimal decay above need not hold. For example, if the surface is given by the equation $x_1^{2m}+x_2^{2m}+\dots+x_d^{2m}=1$, where $m\geq2$ is an integer, and $\mu$ is the surface measure, then
$$ |\widehat{\mu}(\xi)| \leq C(1+|\xi|)^{-\frac{d-1}{2m}}.$$
This follows from the finite type estimates of Bruna, Nagel and Wainger \cite{BrunaNagelWainger1988}. The exponent cannot be improved, as can be seen along the coordinate axes.

While the optimal pointwise decay need not hold if the Gaussian curvature is allowed to vanish, the average spherical decay is sufficient for a variety of applications. More precisely, we may consider
$$ {\left( \int_{S^{d-1}} {|\widehat{\mu}(R \omega)|}^2 d\omega \right)}^{\frac{1}{2}},$$ where $\mu$ is, once again, the surface measure. For general finite measures, such spherical averages play a central role in distance set problems. See \cite{Mattila1987Spherical}. In lattice point counting, it is convenient to consider the Fourier transform of the indicator of the set bounded by the surface. Spherical average estimates for indicators of planar convex bodies give bounds for lattice point discrepancy averaged over rotations and translations. See \cite{BrandoliniRigoliTravaglini1998}. To this end, let $\Omega\subset\mathbb R^d$ be a bounded domain, with $d\geq2$. If $\partial\Omega$ is smooth, its spherical Fourier average satisfies
$$
E_\Omega(R):=\int_{\mathbb S^{d-1}}|\widehat{\mathbb{1}_\Omega}(R\omega)|^2\,d\omega\leq C_\Omega R^{-d-1},\qquad R\geq2;
$$
see \cite{BCT2001}. This is the same power as for the ball, although no curvature assumption is imposed on $\partial\Omega$. No bounded domain of finite perimeter does better. By a result of Beretti and Gennaioli, see \cref{perimeter-frequency-average}, the limit superior of $R^{d+1}E_\Omega(R)$ is at least $2^d\pi^{d-1}$ times the perimeter of $\Omega$, so a bound $E_\Omega(R)\leq CR^{-a}$ with $a>d+1$ is impossible. This concerns upper bounds only. The average may vanish at individual radii, and for the unit ball it vanishes whenever $R$ is a zero of the Bessel function $J_{\frac{d}{2}}$. Brandolini, Hofmann and Iosevich proved in \cite{bhi} that $E_\Omega(R) \leq C_{\Omega}R^{-d-1}$ for every bounded convex domain, without any regularity assumption on its boundary, and for every bounded domain with $C^{1,\frac{1}{2}}$ boundary. More precisely, their second result assumes that $|\langle x-y,\nu_\Omega(y)\rangle|\leq C|x-y|^{\frac{3}{2}}$ locally on $\partial\Omega$, where $\nu_\Omega$ is the outward unit normal.

The question addressed in this paper is how much regularity of the boundary is needed when convexity is not assumed. We prove that, among H\"older classes, the answer is $C^{1,\frac{1}{4}}$.

\begin{teorema}\label{thm:holder}
Let $d\geq2$ and $0<s<1$. Every bounded domain $\Omega\subset\mathbb R^d$ with $C^{1,s}$ boundary satisfies
$$E_\Omega(R)\leq C_\Omega(R^{-d-1}+R^{-d-\frac{1}{2}-2s})$$
for every $R\geq2$. If $0<s<\frac{1}{4}$, then every bounded $C^\infty$ domain can be changed inside an arbitrarily small neighborhood of any prescribed point of its boundary into a bounded domain $\Omega_*$ with $C^{1,s}$ boundary such that
$$\limsup_{R\to\infty}\frac{\log E_{\Omega_*}(R)}{\log R}\geq-d-\frac{1}{2}-2s.$$
In particular, the bound $E_\Omega(R)\leq C_\Omega R^{-d-1}$ holds for every bounded domain with $C^{1,s}$ boundary if and only if $s\geq\frac{1}{4}$.
\end{teorema}

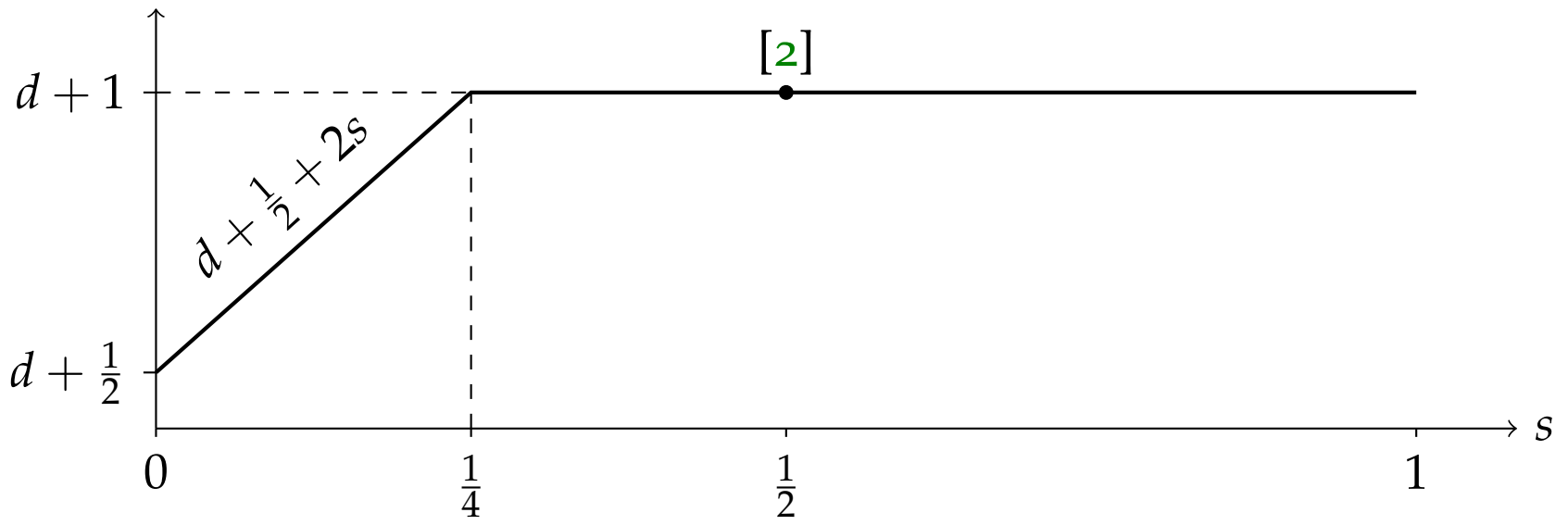
\begin{figure}[ht]
\centering
\begin{tikzpicture}[x=9cm,y=4cm]
\draw[->] (0,0.4) -- (1.08,0.4) node[right] {$s$};
\draw[->] (0,0.4) -- (0,1.15);
\draw (0,0.5) -- (-0.01,0.5) node[left] {$d+\frac{1}{2}$};
\draw (0,1) -- (-0.01,1) node[left] {$d+1$};
\draw (0,0.4) -- (0,0.385) node[below] {$0$};
\draw (0.25,0.4) -- (0.25,0.385) node[below] {$\frac{1}{4}$};
\draw (0.5,0.4) -- (0.5,0.385) node[below] {$\frac{1}{2}$};
\draw (1,0.4) -- (1,0.385) node[below] {$1$};
\draw[dashed] (0.25,0.4) -- (0.25,1);
\draw[dashed] (0,1) -- (0.25,1);
\draw[thick] (0,0.5) -- node[sloped,above] {$d+\frac{1}{2}+2s$} (0.25,1) -- (1,1);
\fill (0.5,1) circle (1.5pt) node[above] {\cite{bhi}};
\end{tikzpicture}
\caption{The supremum of the numbers $a$ such that every bounded domain with $C^{1,s}$ boundary satisfies $E_\Omega(R)=O(R^{-a})$, as a function of $s$, $0<s<1$. It equals $d+\frac{1}{2}+2s$ for $s<\frac{1}{4}$ and $d+1$ for $s\geq\frac{1}{4}$. The dot marks the condition $C^{1,\frac{1}{2}}$ of \cite{bhi}.}
\label{fig:threshold}
\end{figure}

For $s<\frac{1}{4}$ the two halves of the theorem match, so the power $R^{-d-\frac{1}{2}-2s}$ cannot be improved in the class $C^{1,s}$, and the best exponent depends on $s$ as shown in Figure~\ref{fig:threshold}. The theorem is a statement about the class $C^{1,s}$ as a whole. It does not say that a domain with $E_\Omega(R)\leq C_\Omega R^{-d-1}$ must have a $C^{1,\frac{1}{4}}$ boundary, and indeed every convex domain has this decay by \cite{bhi}, however rough its boundary may be. Theorem~\ref{thm:holder} is the case $p=\infty$ of \cref{prop:upper} and \cref{sharp-counterexamples} below, because the normal of a $C^{1,s}$ boundary is H\"older continuous, and hence Dini continuous.

More generally, we determine the optimal decay in classes where the graph gradients have fractional Sobolev regularity, under a Dini assumption on the normal. The counterexamples satisfy this additional assumption as well, and the sharpness is obtained on one fixed domain along a sequence of radii. Without the Dini assumption we obtain the same optimal powers in the ranges specified below.

Fix $0<s<1$ and $1\leq p\leq\infty$. We write $\partial\Omega\in C^1\cap W^{1+s,p}$ if $\partial\Omega$ is a compact embedded $C^1$ hypersurface and the gradients of its graph representations belong to $W^{s,p}$. For $p=\infty$ this is the class $C^{1,s}$. In order to fix notations, we let
$$
\gamma=\begin{cases}
2s,&p=\infty,\\
s\min\{\frac{3p}{2},2\},&d=2,\ p<\infty,\\
s\min\{p,2\},&d\geq3,\ p<\infty.
\end{cases}
$$

In this generality, the first result shows that, when $\gamma<\frac{1}{2}$, the decay $R^{-d-1}$ can fail after perturbing any smooth bounded domain in an arbitrarily small neighborhood of one boundary point. The modified boundary has the prescribed $C^1\cap W^{1+s,p}$ regularity and its normal is Dini continuous. 

\begin{proposizione}
Let $d\geq2$, $0<s<1$ and $1\leq p\leq\infty$. There exists a bounded domain $\Omega_*\subset\mathbb R^d$ with $\partial\Omega_*\in C^1\cap W^{1+s,p}$ and Dini continuous outward unit normal such that
$$
\limsup_{R\to\infty}\frac{\log E_{\Omega_*}(R)}{\log R}\geq-d-\min\left\{1,\frac{1}{2}+\gamma\right\}.
$$
If $\gamma<\frac{1}{2}$, the domain $\Omega_*$ can be obtained by changing any prescribed bounded $C^\infty$ domain inside an arbitrarily small neighborhood of any prescribed point of its boundary.
\end{proposizione}

Since the lower bound holds on one fixed domain, it rules out every decay exponent larger than $d+\min\{1,\frac{1}{2}+\gamma\}$.
The upper bounds proved in \cref{prop:upper-intro} give the converse conclusion under the Dini assumption and in the further ranges stated below. The constants $C_\Omega$ are independent of $R$ and may depend on the graph charts, their regularity bounds and, when assumed, the Dini modulus of the normal, as well as on $d,s,p$.

\begin{proposizione}\label{prop:upper-intro}
Let $\Omega$ be a bounded domain with $\partial\Omega\in C^1\cap W^{1+s,p}$ and let $R\geq2$. If the normal has a Dini modulus, then
$$
E_\Omega(R)\leq C_\Omega (R^{-d-1}+R^{-d-\frac{1}{2}-\gamma}).
$$
Without this additional assumption, if $p<\infty$, then
$$
E_\Omega(R)\leq C_\Omega R^{-d-1}+C_\Omega R^{-\frac{d+3}{2}-s\min\{\frac{3p}{2},2\}}.
$$
If $d\geq3$ and $p\in[1,\infty]$, one also has
$$
E_\Omega(R)\leq C_\Omega R^{-d-1}+C_\Omega
\begin{cases}
R^{-d-\frac{1}{2}-s\min\{p,2\}},&sp>d-2,\\
(\log(2+R))^{\frac{p}{2}}R^{-d-\frac{1}{2}-s\min\{p,2\}},&sp=d-2,\\
(\log(2+R))^{\frac{1}{2}}R^{-\frac{d+3+sp}{2}-s\min\{p,2\}},&sp<d-2.
\end{cases}
$$
\end{proposizione}

\begin{figure}[ht]
\centering
\begin{tikzpicture}[scale=0.92]
\filldraw[fill=gray!12] plot[smooth cycle,tension=0.8] coordinates {(0,0) (1.2,-0.4) (2.4,0.1) (2.6,1.2) (1.6,2) (0.3,1.7) (-0.4,0.8)};
\node at (1.1,0.8) {$\Omega$};
\fill (1.6,2) circle (1.2pt);
\draw[dashed] (1.6,2) circle (0.35);
\node[above] at (1.6,2.35) {$x_0$};
\draw[->] (2.05,2.15) to[bend left=20] (3.7,1.65);
\begin{scope}[shift={(4,1)}]
\fill[gray!12] (-0.2,-0.9) -- (-0.2,0) -- (0.4,0)
 -- plot[domain=0.4:0.65,samples=120] (\x,{0.09*(1-((\x-0.525)/0.125)^2)^2*cos(deg(150.8*(\x-0.525)))})
 -- (0.8,0) -- plot[domain=0.8:1.3,samples=160] (\x,{0.16*(1-((\x-1.05)/0.25)^2)^2*cos(deg(75.4*(\x-1.05)))})
 -- (1.6,0) -- plot[domain=1.6:2.6,samples=160] (\x,{0.27*(1-((\x-2.1)/0.5)^2)^2*cos(deg(31.42*(\x-2.1)))})
 -- (3.2,0) -- plot[domain=3.2:5.2,samples=160] (\x,{0.45*(1-((\x-4.2)/1)^2)^2*cos(deg(12.57*(\x-4.2)))})
 -- (6.4,0) -- (6.4,-0.9) -- cycle;
\draw (-0.2,0) -- (0.4,0)
 -- plot[domain=0.4:0.65,samples=120] (\x,{0.09*(1-((\x-0.525)/0.125)^2)^2*cos(deg(150.8*(\x-0.525)))})
 -- (0.8,0) -- plot[domain=0.8:1.3,samples=160] (\x,{0.16*(1-((\x-1.05)/0.25)^2)^2*cos(deg(75.4*(\x-1.05)))})
 -- (1.6,0) -- plot[domain=1.6:2.6,samples=160] (\x,{0.27*(1-((\x-2.1)/0.5)^2)^2*cos(deg(31.42*(\x-2.1)))})
 -- (3.2,0) -- plot[domain=3.2:5.2,samples=160] (\x,{0.45*(1-((\x-4.2)/1)^2)^2*cos(deg(12.57*(\x-4.2)))})
 -- (6.4,0) node[right] {$\partial\Omega_*$};
\fill (0,0) circle (1.2pt);
\node[below] at (0,-0.05) {$x_0$};
\node[scale=0.6] at (0.2,0.1) {$\cdots$};
\node[below] at (1.05,-0.2) {$x_3$};
\node[below] at (2.1,-0.3) {$x_2$};
\node[below] at (4.2,-0.45) {$x_1$};
\node at (5.85,-0.55) {$\Omega_*$};
\end{scope}
\end{tikzpicture}
\caption{The construction behind the lower bounds, in the plane and not to scale. Near a point $x_0$ of the boundary of a smooth domain $\Omega$, the boundary is first made flat, and the flat piece is then replaced by the graph of a sum of oscillating bumps with disjoint supports accumulating at $x_0$. The bump centered at $x_k$ oscillates with frequency $R_k$ and makes $E_{\Omega_*}(R_k)$ large, while its amplitude and width are small enough to keep $\partial\Omega_*$ in the prescribed regularity class.}
\label{fig:construction}
\end{figure}
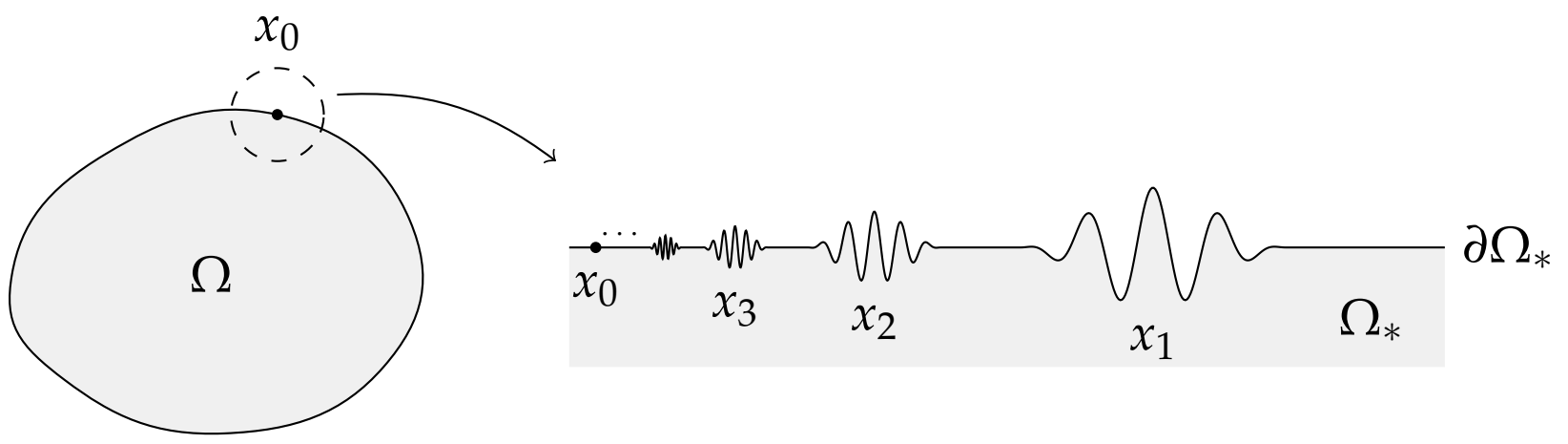

Our starting point was the construction of the counterexamples. When $\gamma<\frac{1}{2}$, we first make a flat patch on a smooth boundary and then perturb it by an oscillating graph, choosing the amplitude and the lengths of its support according to the regularity that we want to preserve. The perturbation produces a large Fourier transform on a small set of directions close to the tangent plane. Placing such perturbations on disjoint supports and separating their frequencies, we obtain one fixed domain with the required limiting exponent, see \cref{one-perturbation}, \cref{sequence-perturbations} and Figure~\ref{fig:construction}. The possibility of obtaining slow decay while leaving the domain unchanged outside an arbitrarily small neighborhood was the first indication that the upper bounds should be sought on individual graph patches.

This also suggests where the difficulty of the upper bounds lies. For a fixed $\omega\in\mathbb S^{d-1}$ the phase $x\mapsto\langle\omega,x\rangle$ is constant on hyperplanes orthogonal to $\omega$, and therefore gives no cancellation in those directions. After the divergence theorem and a partition of unity, the boundary is reduced to graphs $X(u)=(u,A(u))$ with small gradient, and we have to estimate integrals of the form
$$
J_R(\omega)=\int_Ue^{-iR\langle\omega,X(u)\rangle}\psi(u)\langle\omega,(-\nabla A(u),1)\rangle\,du,
$$
where $\psi$ is a $C^1$ cutoff. The angular caps around the normal to the reference plane already give the contribution $R^{-d-1}$ to $E_\Omega$. Indeed, if $u,v$ belong to the graph patch, the chord $X(u)-X(v)$ is almost parallel to the reference plane and cannot be parallel to any $\omega$ in those caps. Hence the angular phase in the expansion of $|J_R|^2$ has no critical points there, and integration by parts on the sphere gives the required bound. We are thus left with the directions close to the reference plane.

Even for those directions there is one graph variable in which the oscillation can be controlled. For fixed $\omega$, this is the direction of its projection onto the reference plane, and on a sufficiently small angular patch we may choose a fixed first coordinate so that $\omega_1+\omega_d\partial_1A$ is bounded away from zero. Mollifying at scale $R^{-1}$ and integrating by parts in that coordinate replaces the rough coefficient of $J_R$ by first differences of $\nabla A$ at distances at most $R^{-1}$, see \eqref{regularized-graph-integral}. This use of first differences already appears in \cite{bhi}. The remaining problem concerns the $d-2$ graph variables perpendicular to the chosen first coordinate, and this explains why a further argument is needed only when $d\geq3$. More precisely, the estimate to be recovered is the power $R^{-d+\frac{3}{2}}$ for the squared $L^2$ norm of the graph integral with an arbitrary $L^2$ coefficient, as in \cref{spherical-graph-estimate}. For smooth graphs this is the bound associated with a one-sided fold, see \cite[Theorem 2.1]{GS}, where the fold condition is supplied by the curvature of the sphere and does not require any curvature of the graph.

To see how the missing decay can be recovered, we expand the squared norm and integrate first in $\omega$. Before restricting to the directions close to the reference plane, the kernel is
$$
\int_{\mathbb S^{d-1}}e^{-iR\langle\omega,X(u)-X(v)\rangle}\,d\omega,
$$
which is a Bessel kernel depending only on $|X(u)-X(v)|$. Estimating this kernel in absolute value gives only the power $R^{-\frac{d-1}{2}}$, while the graph estimate requires $R^{-d+\frac{3}{2}}$. Thus taking absolute values loses a factor $R^{-\frac{d-2}{2}}$ which must come from the oscillation of the kernel itself. We retain its two distance phases $e^{iR|X(u)-X(v)|}$ and $e^{-iR|X(u)-X(v)|}$, and the same phases occur after inserting the angular cutoffs, see \cref{spherical-integral-expansion}. Fixing the first graph coordinate at the two points and rescaling the remaining variables gives an oscillatory operator in $d-2$ variables whose phase is the distance between the two rescaled graph points,
$$\Psi(x,y)=\sqrt{1+|v+x-y|^2+(P(x)-Q(y))^2},$$
where $v\in\mathbb R^{d-2}$ is fixed and $P,Q$ are rescaled slices of the graph function, see \cref{estimate-in-d-minus-two-variables}. The mixed derivative of $\Psi$ is invertible, and the $L^2$ norm of this operator supplies the missing factor. The issue is therefore to prove the usual nondegenerate operator estimate for this particular phase at a regularity where the repeated integrations by parts of the smooth theory cannot be performed.

In the classical arguments, integration by parts is applied to the kernel of the composition of the operator with its adjoint and differentiates coefficients built from the phase, see \cite[Lemma 2.3 and Proposition 2.5]{GS} and \cite[Section 2]{Comech}. Here these coefficients depend on first derivatives of the graph functions, and their higher derivatives need not exist. We replace differentiation by translation of the integration variable, choosing the translation so that the exponential of the phase increment stays a fixed distance from one. We divide the cutoff by this exponential minus one, and subtracting the translated integrals puts a finite difference on the whole quotient, see \eqref{finite-difference-identity}. Iterating the identity leaves high differences of the smooth cutoff and remainders controlled by first differences of the graph gradients. Differences already replace derivatives in \cite{bhi}, but there they are not iterated inside an oscillatory operator. Here they are, and the iteration works because of the special form of $\Psi$. After localization to small cubes and a linear change of variables, the mixed part of $\Psi$ becomes $\langle x,w\rangle$ plus a sum of products $f_\mu(x)g_\mu(w)$ of Lipschitz functions with $\sum_\mu\operatorname{Lip}(f_\mu)\operatorname{Lip}(g_\mu)$ small, see \eqref{target-local-phase}. The phase increment, the reciprocal of the exponential minus one and the annular cutoffs therefore expand into absolutely summable series in which the variables are separated. Apart from a piece on which the differences of the graph gradients are large, every remainder is again an operator of the kind being estimated, multiplied by a small coefficient, and its norm is absorbed into the left-hand side of \eqref{operator-norm-before-absorption}. The Dini assumption makes those remainders summable, whereas the fractional Sobolev assumption gives the bounds that depend on $sp$ and on the number $d-2$ of remaining variables. This proves the operator estimate in \cref{estimate-in-d-minus-two-variables} and, in turn, the graph estimate needed above. Returning to the coefficient of $J_R$, the first differences produced by mollification contribute a factor $R^{-2s}$ to the squared norm for $C^{1,s}$ boundaries. Together with $R^{-d+\frac{3}{2}}$ and the factor $R^{-2}$ from the divergence theorem, this gives $R^{-d-\frac{1}{2}-2s}$ and accounts for the threshold $s=\frac{1}{4}$.

\section{Preliminaries and main results}

Let $d\geq2$, $0<s<1$, and $1\leq p\leq\infty$. In every Euclidean dimension $m\geq1$ we use the Fourier-transform convention
$$
\mathscr Ff(\xi)=\widehat f(\xi)=\int_{\mathbb R^m}e^{-i\xi\cdot x}f(x) dx.
$$
For a bounded domain $\Omega\subset\mathbb R^d$ define
$$
E_\Omega(R)=\int_{\mathbb S^{d-1}}|\widehat{\mathbb{1}_\Omega}(R\omega)|^2 d\omega.
$$
We write $\partial\Omega\in C^1\cap W^{1+s,p}$ if $\partial\Omega$ is a compact embedded $C^1$ hypersurface and the gradients of its graph representations belong to $W^{s,p}$. When $p=\infty$ this means $C^{1,s}$. For the outward unit normal $\nu_\Omega$, let
$$
\omega_\nu(t)=\sup\{|\nu_\Omega(x)-\nu_\Omega(y)|:x,y\in\partial\Omega,\ |x-y|\leq t\},
$$
and put
$$
\gamma=
\begin{cases}
2s,&p=\infty,\\
s\min\{\frac{3p}{2},2\},&d=2,\ p<\infty,\\
s\min\{p,2\},&d\geq3,\ p<\infty.
\end{cases}
$$
The Dini condition used below is
$$
\int_0^1\omega_\nu(t) \frac{dt}{t}<\infty.
$$

\subsection{Results}

\begin{corollario}\label{sharp}
Under the Dini assumption, the supremum of the numbers $a$ such that every domain in the stated class satisfies $E_\Omega(R)=O_\Omega(R^{-a})$, with the implicit constant allowed to depend on $\Omega$, is
$$
d+\min\left\{1,\frac{1}{2}+\gamma\right\}.
$$
The same assertion holds without the Dini assumption if $p=\infty$, or if $p<\infty$ and $sp\geq d-2$. If $p=\infty$, or if $p<\infty$ and $sp>d-2$, every such domain satisfies
$$
E_\Omega(R)\leq C_\Omega (R^{-d-1}+R^{-d-\frac{1}{2}-\gamma}).
$$
If $p<\infty$ and $sp=d-2$, then
$$
E_\Omega(R)\leq C_\Omega R^{-d-1}+C_\Omega(\log(2+R))^{\frac{p}{2}}R^{-d-\frac{1}{2}-\gamma}.
$$
If $sp=d-2$ and $p<4(d-2)$, the full second term, including its logarithmic factor, is bounded by $C_\Omega R^{-d-1}$, and may therefore be included in the first term. At $sp=d-2$ and $p\geq4(d-2)$, the critical power is sharp in this supremum sense, although we do not prove that the logarithm can be removed at the critical power. In particular, every bounded $C^{1,s}$ domain satisfies $E_\Omega(R)\leq C_\Omega R^{-d-1}$ if and only if $s\geq\frac{1}{4}$.
\end{corollario}

Brandolini, Hofmann, and Iosevich proved the power $R^{-d-1}$ for every bounded convex domain and, without convexity, under $|(x-y)\cdot\nu_\Omega(y)|\leq C|x-y|^{\frac{3}{2}}$ for all $x,y\in\partial\Omega$ \cite{bhi}. Lebedev studied the different question of global $L^q(\mathbb R^d)$ integrability of $\widehat{\mathbb{1}_\Omega}$ for domains with $C^1$ boundary \cite{lebedev}. The assertions above concern the spherical $L^2$ mean at each radius, include nonconvex domains, and prove lower bounds in the same boundary classes.

\section{Elementary facts used below}

We use Fourier inversion, Plancherel's identity, the Hardy-Littlewood maximal theorem, and the corresponding coefficient identity for Fourier series on a torus in the forms stated in \cite[Chapters 2 and 3]{grafakos}. The normalization is displayed whenever it matters. The remaining estimates specific to the argument are proved below.

We first record the translation estimate with all constants needed in the scaling arguments. If $m\geq1$, $0<s<1$, and $1\leq p<\infty$, we use
\begin{equation}
    [g]_{W^{s,p}(\mathbb R^m)}^p=\int_{\mathbb R^m}\int_{\mathbb R^m}\frac{|g(u)-g(v)|^p}{|u-v|^{m+sp}} du dv
=\int_{\mathbb R^m}\frac{\|g( \cdot+h)-g\|_{L^p(\mathbb R^m)}^p}{|h|^{m+sp}} dh.
\label{normaWsp}
\end{equation}

The second identity follows from the change of variables $h=v-u$ and Tonelli's theorem.

\begin{lemma}\label{translation}
Let $g\in W^{1,p}(\mathbb R^m)$ have compact support and let $0<s<1$. Then
$$
\|g( \cdot+h)-g\|_p\leq\min\{2\|g\|_p,|h|\|\nabla g\|_p\}
$$
and
$$
[g]_{W^{s,p}(\mathbb R^m)}^p\leq C_{m,s,p}\|g\|_p^{p(1-s)}\|\nabla g\|_p^{ps}.
$$
If instead $g$ is bounded and Lipschitz, then
$$
[g]_{C^{0,s}(\mathbb R^m)}\leq 2^{1-s}\|g\|_\infty^{1-s}\operatorname{Lip}(g)^s.
$$
\end{lemma}

\begin{proof}
The bound by $2\|g\|_p$ is the triangle inequality. If $g\in C_c^\infty(\mathbb R^m)$, then $g(u+h)-g(u)=\int_0^1\langle\nabla g(u+th),h\rangle dt$, and Minkowski's inequality bounds the $L^p$ norm of this difference by $|h|\|\nabla g\|_p$. Smooth compactly supported functions are dense in $W^{1,p}(\mathbb R^m)$, see \cite[Chapter 4]{evansgariepy}, and both sides of this inequality are continuous in the $W^{1,p}$ norm, so it holds for every $g$ in the lemma. If $\|g\|_p$ or $\|\nabla g\|_p$ vanishes, then $g=0$, because a function in $W^{1,p}(\mathbb R^m)$ with vanishing gradient is constant and $g$ has compact support. Otherwise we insert the minimum into \eqref{normaWsp}, pass to polar coordinates in $h$ and split the radial integral at $r_0=\frac{\|g\|_p}{\|\nabla g\|_p}$, which yields
$$
\begin{aligned}
\relax[g]_{W^{s,p}(\mathbb R^m)}^p&\leq|\mathbb S^{m-1}|\|\nabla g\|_p^p\int_0^{r_0}r^{p-sp-1}dr+2^p|\mathbb S^{m-1}|\|g\|_p^p\int_{r_0}^\infty r^{-sp-1}dr\\
&=|\mathbb S^{m-1}|\left(\frac{1}{p(1-s)}+\frac{2^p}{sp}\right)\|g\|_p^{p(1-s)}\|\nabla g\|_p^{ps}.
\end{aligned}
$$
The H\"older estimate is proved in the same way, starting from $|g(u+h)-g(u)|\leq\min\{2\|g\|_\infty,\operatorname{Lip}(g)|h|\}$ and comparing $|h|$ with $\frac{2\|g\|_\infty}{\operatorname{Lip}(g)}$, the cases $\|g\|_\infty=0$ and $\operatorname{Lip}(g)=0$ being trivial.
\end{proof}

The next lemma collects the cutoff, periodization, and Fourier-series argument used several times below.

\begin{lemma}\label{fourier-expansions}
Let $\mathcal R_0$ and $\mathcal R_1$ be bounded rectangles in $\mathbb R^q$, with the closure of $\mathcal R_0$ contained in the interior of $\mathcal R_1$, and let $N$ be an integer such that $N>\frac{q}{2}$. If $F\in C^{2N}(\mathcal R_1)$, then on $\mathcal R_0$ one has an absolutely convergent expansion
$$
F(z)=\sum_{k\in\mathbb Z^q}c_kE_k(z).
$$
The functions $E_k$ may be chosen in the form $E_k(z):=e^{i(\langle\xi_k,z\rangle+\theta_k)}$, where $\xi_k\in\mathbb R^q$, $\theta_k\in\mathbb R$, and
$|\xi_k|\leq C_{\mathcal R_0,\mathcal R_1,q}(1+|k|)$. 
More generally, let $q=q_1+\cdots+q_m$ and write $z=(z_1,\ldots,z_m)$ with $z_j\in\mathbb R^{q_j}$. Then
$$
E_k(z)=E_{k,1}(z_1)\cdots E_{k,m}(z_m),
$$
where $|E_{k,j}|=1$ and
$\Lip(E_{k,j})\leq C_{\mathcal R_0,\mathcal R_1,q}(1+|k|)$. Consequently, if $z_j=z_j(x_j)$ is Lipschitz, then
$\Lip(E_{k,j}\circ z_j)\leq C_{\mathcal R_0,\mathcal R_1,q}(1+|k|)\Lip(z_j)$.
For every $r\geq0$ such that $2N>q+r$, the coefficients satisfy
$$
\sum_{k\in\mathbb Z^q}(1+|k|)^r|c_k|\leq C_{\mathcal R_0,\mathcal R_1,q,N,r}\sum_{|\alpha|\leq2N}\|\partial^\alpha F\|_{L^\infty(\mathcal R_1)}.
$$
\end{lemma}

\begin{proof}
Choose a smooth function $\chi$ supported in $\mathcal R_1$ and equal to one on $\mathcal R_0$, and fix an invertible diagonal matrix $L$ and a vector $a$, depending only on $\mathcal R_0$ and $\mathcal R_1$, such that $z\mapsto Lz+a$ maps the support of $\chi$ into the open cube $(0,2\pi)^q$. The $2\pi$-periodic extension of $\zeta\mapsto(\chi F)(L^{-1}(\zeta-a))$ is of class $C^{2N}$, and integrating by parts with the identity $(1-\Delta)^Ne^{-i\langle k,\zeta\rangle}=(1+|k|^2)^Ne^{-i\langle k,\zeta\rangle}$ shows that its Fourier coefficients satisfy
$$
|c_k|\leq C_{\mathcal R_0,\mathcal R_1,q,N}(1+|k|^2)^{-N}\sum_{|\alpha|\leq2N}\|\partial^\alpha F\|_{L^\infty(\mathcal R_1)}.
$$
There are at most $C_q(1+j)^{q-1}$ points $k\in\mathbb Z^q$ with $j\leq|k|<j+1$, so the sum $\sum_k(1+|k|)^r(1+|k|^2)^{-N}$ is finite whenever $2N>q+r$. This proves the weighted bound for the coefficients and, with $r=0$, the absolute and uniform convergence of the Fourier series, whose sum is the periodic extension itself by the uniqueness of Fourier coefficients. Returning to the variable $z$, we obtain the expansion on $\mathcal R_0$ with $E_k(z)=e^{i\langle k,Lz+a\rangle}$, that is, with $\xi_k=Lk$ and $\theta_k=\langle k,a\rangle$. Since $L$ is diagonal, $E_k$ is the constant $e^{i\langle k,a\rangle}$ times the product of the exponentials $e^{i\langle k_j,L_jz_j\rangle}$, where $k_j$ and $L_j$ are the blocks of $k$ and $L$ corresponding to $z_j$, and each factor has modulus one and Lipschitz constant at most $\|L\||k|$. The bound for $\Lip(E_{k,j}\circ z_j)$ then follows from the chain rule.
\end{proof}

\begin{lemma}\label{L2-fourier-expansion}
Let $n,q\geq1$. Let $\mathcal R_0$ and $\mathcal R_1$ be bounded open rectangles in $\mathbb R^q$ such that $\mathcal R_0$ is compactly contained in $\mathcal R_1$. Let $N$ be an integer such that $N>\frac{q}{2}$ and suppose $\mathcal U\subset\mathbb R^n$ is measurable, and let
$B:\mathcal U\times\mathcal R_1\to\mathbb C$
is a measurable function. Assume that $B(u,z)=0$ for almost every $(u,z)\in\mathcal U\times(\mathcal R_1\setminus\mathcal R_0)$ and that, for every multi-index $\alpha$ with $|\alpha|\leq N$, the weak derivative of $B$ with respect to $z$ satisfies
$$
\partial_z^\alpha B\in L^2(\mathcal U\times\mathcal R_1).
$$
Then, there are functions $B_k\in L^2(\mathcal U)$, an invertible matrix $A\in\mathbb R^{q\times q}$ and a vector $a\in\mathbb R^q$, where $A$ and $a$ depend only on $\mathcal R_0$ and $\mathcal R_1$, such that, upon setting
$$
E_k(z):=e^{i\langle k,Az+a\rangle},
\qquad k\in\mathbb Z^q,\qquad\text{one has}\qquad
B(u,z)=\sum_{k\in\mathbb Z^q}B_k(u)E_k(z)\qquad\text{in $L^2(\mathcal U\times\mathcal R_0)$,}
$$
and
$$
\sum_{k\in\mathbb Z^q}\|B_k\|_{L^2(\mathcal U)}
\leq C_{\mathcal R_0,\mathcal R_1,q,N}
\Big(
\sum_{|\alpha|\leq N}
\|\partial_z^\alpha B\|_{L^2(\mathcal U\times\mathcal R_1)}^2
\Big)^{\frac{1}{2}}.
$$
Let us further assume 
$K:\mathcal R_0\times\mathcal U\to\mathbb C$ is a  measurable function such that
$$
\operatorname*{ess\,sup}_{z\in\mathcal R_0}
\|K(z, \cdot )\|_{L^2(\mathcal U)}<\infty
\qquad\text{and that there is $L>0$ such that}\qquad
\Big\|
\int_{\mathcal U}K(z,u)b(u) du
\Big\|_{L^2(\mathcal R_0)}
\leq L\|b\|_{L^2(\mathcal U)},
$$
for every $b\in L^2(\mathcal U)$. Then
$$
\Big\|
\int_{\mathcal U}K(z,u)B(u,z) du
\Big\|_{L^2(\mathcal R_0)}
\leq C_{\mathcal R_0,\mathcal R_1,q,N}L
\Big(
\sum_{|\alpha|\leq N}
\|\partial_z^\alpha B\|_{L^2(\mathcal U\times\mathcal R_1)}^2
\Big)^{\frac{1}{2}}.
$$
\end{lemma}

\begin{proof}
Fix an invertible matrix $A$ and a vector $a$, depending only on $\mathcal R_0$ and $\mathcal R_1$, such that $z\mapsto Az+a$ maps $\mathcal R_1$ into the open cube $(0,2\pi)^q$. Since $B$ vanishes on $\mathcal U\times(\mathcal R_1\setminus\mathcal R_0)$ and $\mathcal R_0$ is compactly contained in $\mathcal R_1$, for almost every $u$ the function $\zeta\mapsto B(u,A^{-1}(\zeta-a))$, extended by zero to $[0,2\pi]^q$ and then periodically, keeps its weak derivatives of order at most $N$. We let $B_k(u)$ be its $k$-th Fourier coefficient. Parseval's identity, integrated in $u$, shows that $\sum_kB_kE_k$ converges to $B$ in $L^2(\mathcal U\times\mathcal R_0)$. Applied to the derivatives in $\zeta$, whose Fourier coefficients are $(ik)^\alpha B_k(u)$, it also shows that
$$
\sum_{k\in\mathbb Z^q}(1+|k|^2)^N\|B_k\|_{L^2(\mathcal U)}^2\leq C_{\mathcal R_0,\mathcal R_1,q,N}\sum_{|\alpha|\leq N}\|\partial_z^\alpha B\|_{L^2(\mathcal U\times\mathcal R_1)}^2,
$$
because $(1+|k|^2)^N\leq C_{q,N}\sum_{|\alpha|\leq N}|k^\alpha|^2$ and each derivative in $\zeta$ is a combination of derivatives in $z$ with coefficients depending only on $A$. Since $N>\frac{q}{2}$, the series $\sum_k(1+|k|^2)^{-N}$ converges, and the Cauchy-Schwarz inequality
$$
\sum_{k\in\mathbb Z^q}\|B_k\|_{L^2(\mathcal U)}\leq\left(\sum_{k\in\mathbb Z^q}(1+|k|^2)^{-N}\right)^{\frac{1}{2}}\left(\sum_{k\in\mathbb Z^q}(1+|k|^2)^N\|B_k\|_{L^2(\mathcal U)}^2\right)^{\frac{1}{2}}
$$
proves the estimate for the coefficients.

For the last assertion, let $J_1\subset J_2\subset\cdots$ be finite sets whose union is $\mathbb Z^q$. Since $|E_k|=1$, the assumed operator bound and the triangle inequality show that
$$
\left\|\int_{\mathcal U}K(z,u)\sum_{k\in J_m}B_k(u)E_k(z)du\right\|_{L^2(\mathcal R_0)}\leq L\sum_{k\in J_m}\|B_k\|_{L^2(\mathcal U)}\leq C_{\mathcal R_0,\mathcal R_1,q,N}L\left(\sum_{|\alpha|\leq N}\|\partial_z^\alpha B\|_{L^2(\mathcal U\times\mathcal R_1)}^2\right)^{\frac{1}{2}}.
$$
For almost every $z\in\mathcal R_0$, the Cauchy-Schwarz inequality in $u$ bounds the difference between $\int_{\mathcal U}K(z,u)B(u,z)du$ and the integral on the left by $\|K(z,\cdot)\|_{L^2(\mathcal U)}$ times the $L^2(\mathcal U)$ norm of $B(\cdot,z)-\sum_{k\in J_m}B_k(\cdot)E_k(z)$. Since $\|K(z,\cdot)\|_{L^2(\mathcal U)}$ is essentially bounded in $z$ and the expansion converges in $L^2(\mathcal U\times\mathcal R_0)$, the $L^2(\mathcal R_0)$ norm of this difference tends to zero as $m\to\infty$, and the estimate passes to the limit.
\end{proof}

The finite difference argument in the proof of \cref{estimate-in-d-minus-two-variables} manipulates series of products of three functions, each depending on one variable. The next lemma records once and for all the rules that make these manipulations legitimate.

\begin{lemma}\label{separated-expansions}
Let $m\geq1$ and let $E,E'\subset\mathbb R^m$ be measurable sets of positive measure. We say that a function $c$ on $E\times E'\times\mathbb R^m$ has a separated expansion of size at most $S$ if there are bounded measurable functions $u_\nu$ on $E$, $v_\nu$ on $E'$ and $b_\nu$ on $\mathbb R^m$, indexed by a countable set, such that $\sum_\nu\|u_\nu\|_{L^\infty}\|v_\nu\|_{L^\infty}\|b_\nu\|_{L^\infty}\leq S$ and
$$c(x,x',w)=\sum_\nu u_\nu(x)v_\nu(x')b_\nu(w)$$
for almost every $(x,x',w)$. Such a series converges absolutely almost everywhere, and $|c|\leq S$ almost everywhere. If $c_1$ and $c_2$ have separated expansions of sizes at most $S_1$ and $S_2$, then multiplying the two series term by term produces a separated expansion of $c_1c_2$, indexed by pairs of indices, of size at most $S_1S_2$. If $c$ has a separated expansion of size at most $S$ and $\sum_{k\geq0}|a_k|S^k<\infty$, then multiplying out the powers produces a separated expansion of $\sum_{k\geq0}a_kc^k$, indexed by the tuples $(k,\nu_1,\ldots,\nu_k)$, of size at most $\sum_{k\geq0}|a_k|S^k$. Finally, for every $h\in\mathbb R^m$, replacing each $b_\nu$ by $b_\nu(\cdot-h)$ produces a separated expansion of $c(x,x',w-h)$ of the same size.
\end{lemma}

\begin{proof}
For almost every $(x,x',w)$, the inequalities $|u_\nu(x)|\leq\|u_\nu\|_{L^\infty}$, $|v_\nu(x')|\leq\|v_\nu\|_{L^\infty}$ and $|b_\nu(w)|\leq\|b_\nu\|_{L^\infty}$ hold for all the countably many indices at once, and at such points the series is dominated by $S$. Now let $c_1$ and $c_2$ be as in the statement. The norm of a product of two bounded functions is at most the product of their norms, and Tonelli's theorem for double series of nonnegative terms shows that the sum, over all pairs of indices, of the products of the six norms is the product of the two sizes. The multiplied series therefore has size at most $S_1S_2$ and converges absolutely almost everywhere. Its sum is $c_1c_2$, because an absolutely convergent double series may be summed in any order, in particular as the product of its two factors. Applying this $k-1$ times expands $c^k$ into a series indexed by $(\nu_1,\ldots,\nu_k)$ of size at most $S^k$, and summing these expansions with the weights $a_k$ proves the assertion about power series. The last assertion holds because translations do not change $L^\infty$ norms.
\end{proof}

\begin{lemma}\label{normal}
Let $U\subset\mathbb R^{d-1}$ be bounded and open, and let
$f\in C^1(U)$ be $1$-Lipschitz. The upward unit normal to the
graph parametrized by $X(u)=(u,f(u))$ is
$$
\nu(u)=\frac{(-\nabla f(u),1)}{\sqrt{1+|\nabla f(u)|^2}}.
$$
For every $p\in[1,\infty]$, we have $[\nu]_{W^{s,p}(U)}\leq2[\nabla f]_{W^{s,p}(U)}$ and
$$
|u-v|\leq|X(u)-X(v)|\leq\sqrt2\,|u-v|.
$$
The seminorm bound also holds on the graph, using Euclidean
distance and surface measure, with a constant depending only
on $d,s,p$.
\end{lemma}

\begin{proof}
The derivative of the map $v\mapsto(1+|v|^2)^{-\frac{1}{2}}(-v,1)$ has norm at most $(1+|v|^2)^{-\frac{1}{2}}\leq1$, so that $|\nu(u)-\nu(v)|\leq|\nabla f(u)-\nabla f(v)|$, and the bound for $[\nu]_{W^{s,p}(U)}$ follows. Since $f$ is $1$-Lipschitz, we have $|u-v|^2\leq|X(u)-X(v)|^2=|u-v|^2+|f(u)-f(v)|^2\leq2|u-v|^2$. Finally, the Jacobian $\sqrt{1+|\nabla f(u)|^2}$ in the area formula lies between $1$ and $\sqrt2$, so integrals over the graph with respect to $\mathcal H^{d-1}$ are comparable to the corresponding integrals over $U$, and together with the lower bound $|X(u)-X(v)|\geq|u-v|$ this transfers the seminorm bound to the graph.
\end{proof}

For a scalar function $h$ on $\mathbb S^{d-1}$, extend it to $\mathbb R^d\setminus\{0\}$ by $h(\frac{x}{|x|})$ and let $\nabla_{\mathbb S}^jh(\omega)$ mean the $j$-th Euclidean derivative of this radial extension at $\omega$, with $\nabla_{\mathbb S}^0h=h$.

\begin{lemma}\label{spherical-kernel}
There is a constant $C_d$ such that, for every $c\in C^{2d+5}(\mathbb S^{d-1})$, every $z\in\mathbb R^d$, and every $R\geq1$,
$$
\Big|\int_{\mathbb S^{d-1}}e^{-iR\omega\cdot z}c(\omega) d\omega\Big|
\leq C_d(1+R|z|)^{-\frac{d-1}{2}}\sum_{j=0}^{2d+5}\|\nabla_{\mathbb S}^jc\|_\infty.
$$
If $c$ is supported in a fixed spherical coordinate patch on which
$$
|z-(z\cdot\omega)\omega|\geq c_0|z|
$$
then
$$
\Big|\int_{\mathbb S^{d-1}}e^{-iR\omega\cdot z}c(\omega) d\omega\Big|
\leq C_{d,c_0}(1+R|z|)^{-d}\sum_{j=0}^d\|\nabla_{\mathbb S}^jc\|_\infty.
$$
\end{lemma}

\begin{proof}
If $R|z|\leq1$, both integrals are bounded by $|\mathbb S^{d-1}|\|c\|_\infty$, and both estimates hold because $(1+R|z|)^{-d}\geq2^{-d}$. We may therefore assume that $R|z|>1$. The two estimates are then the classical bounds for oscillatory integrals with the smooth phase $\omega\mapsto-\langle\omega,\frac{z}{|z|}\rangle$ and the large parameter $\mu=R|z|$, and we only have to keep track of the number of derivatives of $c$ that enter. In a fixed coordinate patch $\omega=\varphi(y)$ the integral becomes $\int e^{i\mu f(y)}u(y)dy$, where $f(y)=-\langle\varphi(y),\frac{z}{|z|}\rangle$ stays in a bounded subset of every $C^k$ space as $z$ varies, and $u=(c\circ\varphi)J$, with $J$ the density of the surface measure in these coordinates. By the chain rule, the derivatives of $u$ of order at most $k$ are bounded by a constant multiple of $\sum_{j\leq k}\|\nabla_{\mathbb S}^jc\|_\infty$.

For the second assertion, recall that $z-\langle z,\omega\rangle\omega$ is the tangential gradient at $\omega$ of the function $\omega\mapsto\langle\omega,z\rangle$. The hypothesis therefore bounds $|\nabla f|$ from below on the support of $u$ by a constant multiple of $c_0$, the constant depending only on the patch, and \cite[Theorem 7.7.1]{HormanderI}, applied with $k=d$, shows that
$$
\left|\int e^{i\mu f(y)}u(y)dy\right|\leq C_{d,c_0}\mu^{-d}\sum_{|\alpha|\leq d}\sup|\partial^\alpha u|,
$$
which is the second estimate.

For the first assertion we rotate $\frac{z}{|z|}$ to $e_d$, which changes neither the integral nor the norms of the derivatives of $c$, so that the phase becomes $-\mu\omega_d$. A fixed smooth partition of unity splits $c$ into a piece supported in $\{\omega_d>\frac{1}{2}\}$, a piece supported in $\{\omega_d<-\frac{1}{2}\}$ and a piece supported in $\{|\omega_d|<\frac{3}{4}\}$. On the last set the tangential gradient of $\omega\mapsto\omega_d$ has length $(1-\omega_d^2)^{\frac{1}{2}}>\frac{\sqrt7}{4}$, so after splitting that piece among finitely many fixed coordinate patches, the second assertion bounds its contribution by $C_d\mu^{-d}\sum_{j\leq d}\|\nabla_{\mathbb S}^jc\|_\infty$. In the cap $\{\omega_d>\frac{1}{2}\}$ we use the coordinates $\omega=(y,\sqrt{1-|y|^2})$, in which the phase is $\mu f(y)$ with $f(y)=-\sqrt{1-|y|^2}$. The only critical point of $f$ is $y=0$, the Hessian of $f$ there is the identity, and $\frac{|y|}{|\nabla f(y)|}=\sqrt{1-|y|^2}\leq1$. We now apply \cite[Theorem 7.7.5]{HormanderI} with $k=\lceil\frac{d-1}{2}\rceil$. Each term of the resulting expansion is $\mu^{-\frac{d-1}{2}-j}$ times a combination of derivatives of $u$ of order $2j<2k$ at the origin, and the remainder is bounded by a constant multiple of $\mu^{-k}\sum_{|\alpha|\leq2k}\sup|\partial^\alpha u|$. Since $k\geq\frac{d-1}{2}$ and $2k\leq d$, the contribution of the cap is at most $C_d\mu^{-\frac{d-1}{2}}\sum_{j\leq d}\|\nabla_{\mathbb S}^jc\|_\infty$. The cap $\{\omega_d<-\frac{1}{2}\}$ is treated in the same way, and the first assertion follows, since $d\leq2d+5$.
\end{proof}

\begin{lemma}\label{flat-patch}
Let $\Omega\subset\mathbb R^d$ be a bounded $C^\infty$ domain. There are numbers $r_\Omega>0$ and $C_\Omega\geq1$ such that, for every $x_0\in\partial\Omega$ and every $0<\rho\leq r_\Omega$, there is a bounded $C^\infty$ domain $\Omega_{x_0,\rho}$ satisfying
$$\Omega_{x_0,\rho}\cap B(x_0,3\rho)^c=\Omega\cap B(x_0,3\rho)^c,\qquad \partial\Omega_{x_0,\rho}\cap B(x_0,3\rho)^c=\partial\Omega\cap B(x_0,3\rho)^c.$$
Moreover, $d_\Haus(\partial\Omega,\partial\Omega_{x_0,\rho})\leq C_\Omega\rho^2$, and, if $\nu_\Omega(x_0)$ is the outward unit normal to $\partial\Omega$ at $x_0$, then
$$\partial\Omega_{x_0,\rho}\cap B(x_0,2\rho)=x_0+\nu_\Omega(x_0)^\perp\cap B(x_0,2\rho),\quad \Omega_{x_0,\rho}\cap B(x_0,2\rho)=\{x\in B(x_0,2\rho):(x-x_0)\cdot\nu_\Omega(x_0)<0\}.$$
\end{lemma}

\begin{proof}
The compactness of $\partial\Omega$ gives numbers $r_\Omega>0$ and $M_\Omega\geq1$ such that $25M_\Omega r_\Omega \leq1$ and 
with the following property. For every $x_0\in\partial\Omega$, up to isometries, there is a smooth function $h$ defined on $U_{x_0}:=\{v\in\mathbb R^{d-1}:|v|<3r_\Omega\}$ such that
$$e_d=\nu_\Omega(x_0),\qquad h(0)=0,\qquad \nabla h(0)=0,\qquad \|\nabla^2h\|_\infty\leq M_\Omega,$$
$$\Omega\cap B(x_0,3r_\Omega)=\{x_0+v+te_d:v\in U_{x_0},\ t<h(v)\}\cap B(x_0,3r_\Omega).$$
Fix $0<\rho\leq r_\Omega$ and choose a positive smooth function $\varphi$ such that $0\leq\varphi\leq1$, $\varphi(v)=1$ if $|v|\leq2\rho$, and $\varphi(v)=0$ if $|v|\geq\frac{5\rho}{2}$ and define $h'(v)=h(v)(1-\varphi(v))$.
Inside $B(x_0,3\rho)$ set
$$\Omega_{x_0,\rho}\cap B(x_0,3\rho)=\{x_0+v+te_d:|v|<3\rho,\ t<h'(v)\}\cap B(x_0,3\rho),$$
and let $\Omega_{x_0,\rho}$ agree with $\Omega$ outside this ball.

Taylor's formula implies that  $|h(v)|\leq \frac{M_\Omega|v|^2}{2}$, and the graphs of $h$ and $h'$ can differ only when $|v|<\frac{5\rho}{2}$. In this region, since $|h'|\leq|h|$ and $25M_\Omega r_\Omega\leq1$, we have $\max\{|h(v)|,|h'(v)|\}\leq\frac{25M_\Omega\rho^2}{8}\leq\frac{\rho}{2}$, so that the modified graph stays inside $B(x_0,3\rho)$, and hence 
$$d_\Haus(\partial \Omega,\partial \Omega_{x_0,\rho})\leq\lVert h-h'\rVert_\infty\leq 9M_\Omega\rho^2,$$
where the first inequality is obtained pairing the points $x_0+v+h(v)e_d$ and $x_0+v+h'(v)e_d$.
The conclusion follows with $C_\Omega=\max\{1,9M_\Omega\}$.
\end{proof}

\begin{lemma}\label{perturb-flat-patch}
Let $\Omega\subset\mathbb R^d$ be a bounded $C^\infty$ domain, and let $r_\Omega$ be the number given by \cref{flat-patch}. For every $x_0\in\partial\Omega$, every $0<\rho<r_\Omega$, and every $f\in C_c^1(B^{d-1}(0,\rho))$ satisfying $\|f\|_\infty+\|\nabla f\|_\infty\leq\rho$, there is a bounded $C^1$ domain $\Omega_{x_0,\rho,f}$ which agrees with $\Omega_{x_0,\rho}$ outside $B(x_0,2\rho)$ and consequently agrees with $\Omega$ outside $B(x_0,3\rho)$. Up to an isometry, one has
$$\partial\Omega_{x_0,\rho,f}\cap B(0,2\rho)=\{x_0+u+f(u)e_d:|u|<2\rho\}.$$
Moreover, for every integrable function $G$ supported in $B(0,2\rho)$, one has
$$\int_{\mathbb R^d}(\mathbb1_{\Omega_{x_0,\rho,f}}-\mathbb1_{\Omega_{x_0,\rho}})(x)G(x) dx=\int_{B^{d-1}(0,\rho)}\int_0^{f(u)}G(u+ze_d)dzdu.$$
\end{lemma}

\begin{proof}
Without loss of generality we may assume that $x_0=0$ and $\nu_\Omega(x_0)=e_d$. By \cref{flat-patch}, we have 
$$\partial\Omega_{x_0,\rho}\cap B(0,2\rho)=e_d^\perp\cap B(0,2\rho)\qquad\text{and}\qquad\Omega_{x_0,\rho}\cap B(0,2\rho)=\{u+ze_d\in B(0,2\rho):z<0\}.$$
Inside $B(0,2\rho)$ define
$$\Omega_{x_0,\rho,f}\cap B(0,2\rho):=\{u+ze_d\in B(0,2\rho):z<f(u)\},$$
and let $\Omega_{x_0,\rho,f}$ agree with $\Omega_{x_0,\rho}$ outside this ball.

The two domains can differ only at points $u+ze_d$ for which $|u|<\rho$ and $|z|\leq|f(u)|\leq\rho$. Every such point belongs to $B(0,\sqrt2\rho)$, which is compactly contained in $B(0,2\rho)$. Therefore the two definitions agree in a neighbourhood of $\partial B(0,2\rho)$ and determine a bounded $C^1$ domain. Its boundary in $B(0,2\rho)$ is the graph of the extension of $f$, which proves the first identity.

For every fixed $u\in B^{d-1}(0,\rho)$, the vertical section of $\mathbb1_{\Omega_{x_0,\rho,f}}-\mathbb1_{\Omega_{x_0,\rho}}$ is equal to $1$ on $(0,f(u))$ if $f(u)>0$, is equal to $-1$ on $(f(u),0)$ if $f(u)<0$, and vanishes otherwise. Fubini's theorem gives the asserted identity, where the integral from $0$ to $f(u)$ is understood as an oriented integral. 
\end{proof}

\section{Lower bounds}

In this section we prove that every domain with smooth boundary can be perturbed at any point at a sufficiently small scale in such a way that the rate of decay of the average Fourier transform is dictated by the roughness of the surface of the domain. 

\subsection{The decay at bad scales}

In order to do so, we fix a bounded $C^\infty$ domain $\Omega$, a point $x_0\in\partial\Omega$, and $0<\rho<\min\{r_\Omega,1\}$, where $r_\Omega$ is the number given by \cref{flat-patch} and let $\Omega_0=\Omega_{x_0,\rho}$ be the domain given by \cref{flat-patch}. Since $E_\Omega$ is invariant under isometries, we may assume without loss of generality that $x_0=0$ and $\nu_\Omega(x_0)=e_d$. 

We identify $e_d^\perp$ with $\R^{d-1}$ and we write $u=(x,y)\in\mathbb R\times\mathbb R^{d-2}$, using only $x$ when $d=2$. Finally, in what follows, we fix a nonnegative function $\psi\in C_c^\infty((-\frac{1}{16},\frac{1}{16}))$ with positive integral and when $d\geq3$, we also choose a nonnegative function $\eta\in C_c^\infty(B^{d-2}(0,\frac{1}{16}))$ with positive integral and $\|\eta\|_\infty=1$. After multiplying $\psi$ by a positive constant, we assume that $\|\psi'\|_\infty+\frac{3\|\psi\|_\infty}{2}+\|\psi\|_\infty\|\nabla\eta\|_\infty\leq1$. When $d=2$ there is no variable $y$. In that case the cutoff $\eta$, the parameters $D$ and $D_k$ below and all the conditions involving them are absent, the factors $D^{d-2}$ and $D_k^{d-2}$ in the estimates are equal to one, and the normalization above reads $\|\psi'\|_\infty+\frac{3\|\psi\|_\infty}{2}\leq1$.

\begin{proposizione}\label{one-perturbation}
Let us fix $R\geq2$, $a>0$, $0<L\leq\rho$ and $0<D\leq\rho$
in such a way that 
$$RL\geq1,\qquad aR\leq\rho,\qquad RD\geq1,\qquad\text{and}\qquad\frac{L}{RD^2}\leq1.$$
When $d=2$, the parameter $D$ and the conditions involving it are omitted. For every $\phi\in[0,2\pi]$, we let $\Omega_{R,\phi}:=\Omega_{x_0,\rho,f_{R,\phi}}$ be the smooth domain given by \cref{perturb-flat-patch}, where
$$f_{R,\phi}(x,y)=\begin{cases}a\psi(\frac{x}{L})\eta(\frac{y}{D})\cos(Rx+\phi),&d\geq3,\\
a\psi(\frac{x}{L})\cos(Rx+\phi),&d=2.\end{cases}$$
For every $0<s<1$ and $1\leq p<\infty$, its unit normal satisfies
$$[\nu_{\Omega_{R,\phi}}]_{W^{s,p}(\partial\Omega_{R,\phi})}^p\leq C+Ca^pR^{p(1+s)}LD^{d-2},$$
while for $p=\infty$ one has $[\nu_{\Omega_{R,\phi}}]_{C^{0,s}(\partial\Omega_{R,\phi})}\leq C+CaR^{1+s}$. In addition, for at least one value of $\phi$, there holds
$$E_{\Omega_{R,\phi}}(R)\geq ca^2R^{-d+\frac{3}{2}}L^{\frac{3}{2}}D^{d-2}.$$
\end{proposizione}

\begin{proof}
Thanks to the choice of our parameters, the support of $\psi(\frac{x}{L})$ lies in $(-\frac{\rho}{16},\frac{\rho}{16})$ and for $d\geq3$, the support of $\eta(\frac{y}{D})$ lies in $B^{d-2}(0,\frac{\rho}{16})$. Hence $f_{R,\phi}$ is supported in $B^{d-1}(0,\rho)$. Differentiating and using $L^{-1}\leq R$ and $D^{-1}\leq R$, we obtain
$$\|f_{R,\phi}\|_\infty+\|\nabla f_{R,\phi}\|_\infty\leq aR\leq\rho,\qquad \|\nabla^2f_{R,\phi}\|_\infty\leq CaR^2.$$
Thus \cref{perturb-flat-patch} defines the domain $\Omega_{R,\phi}$.
The support of $f_{R,\phi}$ has $(d-1)$-dimensional measure at most $LD^{d-2}$ when $d\geq3$ and at most $L$ when $d=2$. It follows that
$$\|\nabla f_{R,\phi}\|_p+R^{-1}\|\nabla^2f_{R,\phi}\|_p\leq CaR(LD^{d-2})^{\frac{1}{p}}.$$
Let us observe that if $p=\infty$, we have 
$$[\nabla f_{R,\phi}]_{C^{0,s}}\leq \sup_{\substack{x\neq y\\|x-y|\leq R^{-1}}}\frac{|\nabla f_{R,\phi}(x)-\nabla f_{R,\phi}(y)|}{|x-y|^s}+\sup_{\substack{|x-y|\geq R^{-1}}}\frac{|\nabla f_{R,\phi}(x)-\nabla f_{R,\phi}(y)|}{|x-y|^s}\leq CaR^{1+s}.$$
Thus, we can assume that $1\leq p<\infty$. 
The first estimate in \cref{translation}, applied to $\nabla f_{R,\phi}$, therefore gives
$$\|\nabla f_{R,\phi}( \cdot+h)-\nabla f_{R,\phi}\|_p\leq CaR(LD^{d-2})^{\frac{1}{p}}\min\{1,R|h|\}.$$
Thanks to \eqref{normaWsp}, we see that 
\begin{equation}
    \begin{split}
        [\nabla f_{R,\phi}]_{W^{s,p}}^p\leq&\int_{|h|\geq R^{-1}}\frac{\|\nabla f_{R,\phi}(\cdot+h)-\nabla f_{R,\phi}\|_{L^p(\mathbb R^{d-1})}^p}{|h|^{d-1+sp}}dh+\int_{|h|\leq R^{-1}}\frac{\|\nabla f_{R,\phi}(\cdot+h)-\nabla f_{R,\phi}\|_{L^p(\mathbb R^{d-1})}^p}{|h|^{d-1+sp}}dh.
    \end{split}
\end{equation}
Passing to polar coordinates the integral with $|h|\leq R^{-1}$ contains $R^p\int_0^{R^{-1}}r^{p-sp-1} dr=\frac{R^{sp}}{p(1-s)}$, while the integral with $|h|\geq R^{-1}$ contains $\int_{R^{-1}}^\infty r^{-sp-1} dr=\frac{R^{sp}}{sp}$. Consequently
$$[\nabla f_{R,\phi}]_{W^{s,p}}^p\leq Ca^pR^{p(1+s)}LD^{d-2},$$
where $C$ depends only on $d,p$ and $s$, and on the fixed functions $\psi$ and $\eta$. By \cref{normal}, the preceding estimates control the normal on the perturbed graph. The remaining part of the boundary is smooth and independent of $R,a,L,D$ and $\phi$, while the support of the perturbation has positive distance from the nonflat part of $\partial\Omega_0$. This proves the stated estimates for $\nu_{\Omega_{R,\phi}}$.

We next prove the lower bound for the decay of the average Fourier transform. Write $\xi=(\xi_1,\xi',\xi_d)$, with all terms containing $\xi'$ or $y$ omitted when $d=2$. Choose a smooth cutoff supported in $B(0,2\rho)$ and equal to one on the support of $\mathbb1_{\Omega_{R,\phi}}-\mathbb1_{\Omega_0}$. Applying the integral identity in \cref{perturb-flat-patch} to the product of this cutoff and $e^{-i\xi\cdot(u,z)}$ gives
$$\widehat{\mathbb1_{\Omega_{R,\phi}}}(\xi)=\widehat{\mathbb1_{\Omega_0}}(\xi)+T_{f_{R,\phi}}(\xi),\qquad\text{where}\qquad T_f(\xi):=\int e^{-i(\xi_1x+\langle \xi', y\rangle)}\int_0^{f(x,y)}e^{-i\xi_dz} dz dx dy.$$
For $\omega\in\mathbb S^{d-1}$, let us define
$$q_R(\omega):=\frac{1}{2\pi}\int_0^{2\pi}\widehat{\mathbb1_{\Omega_{R,\phi}}}(R\omega)e^{-i\phi} d\phi=\frac{1}{2\pi}\int_0^{2\pi}T_{f_{R,\phi}}(R\omega)e^{-i\phi} d\phi,$$
where the identity comes from the fact that $\widehat{\mathbb1_{\Omega_0}}(\xi)$ is independent of $\phi$. For real $t,v$, one has
$$\Big|\int_0^ve^{-itz} dz-v\Big|=\Big|\int_0^v(e^{-itz}-1) dz\Big|\leq\frac{|t|v^2}{2}.$$
Using this inequality with $t=R\omega_d$ and $v=f_{R,\phi}(x,y)$ in the definition of $T_{f_{R,\phi}}$, we obtain
$$\Big|q_R(\omega)-\int e^{-iR(\omega_1x+\omega'\cdot y)}\frac1{2\pi}\int_0^{2\pi}f_{R,\phi}(x,y)e^{-i\phi} d\phi dx dy\Big|\leq\frac{R|\omega_d|}{4\pi}\int_0^{2\pi}\int|f_{R,\phi}(x,y)|^2 dx dy d\phi.$$
The integral with respect to $\phi$ can be computed explicitly. Indeed,
$$\frac1{2\pi}\int_0^{2\pi}\cos(Rx+\phi)e^{-i\phi} d\phi=\frac12e^{iRx},$$
and therefore
$$\frac1{2\pi}\int_0^{2\pi}f_{R,\phi}(x,y)e^{-i\phi} d\phi=\begin{cases}\frac a2\psi(\frac{x}{L})\eta(\frac{y}{D})e^{iRx},&d\geq3,\\[2mm]
\frac a2\psi(\frac{x}{L})e^{iRx},&d=2.\end{cases}$$
A change of variables now gives
$$\int e^{-iR(\omega_1x+\omega'\cdot y)}\frac1{2\pi}\int_0^{2\pi}f_{R,\phi}(x,y)e^{-i\phi} d\phi dx dy=\begin{cases}\frac{aLD^{d-2}}2\widehat\psi (RL(\omega_1-1))\widehat\eta(RD\omega'),&d\geq3,\\[2mm]
\frac{aL}2\widehat\psi (RL(\omega_1-1)),&d=2.\end{cases}$$
Finally,
$$\int_0^{2\pi}\int|f_{R,\phi}(x,y)|^2 dx dy d\phi=\begin{cases}\pi a^2LD^{d-2}\|\psi\|_2^2\|\eta\|_2^2,&d\geq3,\\
\pi a^2L\|\psi\|_2^2,&d=2.\end{cases}$$
Combining the preceding identities yields
\begin{equation}
\begin{cases} |q_R(\omega)-\frac{aLD^{d-2}}2\widehat\psi (RL(\omega_1-1))\widehat\eta(RD\omega')|\leq Ca^2R|\omega_d|LD^{d-2},&d\geq3,\\[2mm]
 |q_R(\omega)-\frac{aL}2\widehat\psi (RL(\omega_1-1))|\leq Ca^2R|\omega_d|L,&d=2,\end{cases}
\label{stimeqR}
\end{equation}
where the constant $C$ above depends on $\psi$ and $\eta$. 
The number $\widehat\psi(0)$ is positive and, when $d\geq3$, so is $\widehat\eta(0)$ since $\psi$ and $\eta$ are nonnegative and nonzero. Fix once and for all $0<\varepsilon\leq\frac{1}{4}$ such that $|\widehat\psi(t)|\geq\frac{\widehat\psi(0)}{2}$ for $|t|\leq10\varepsilon$ and, when $d\geq3$, $|\widehat\eta(v)|\geq\frac{\widehat\eta(0)}{2}$ for $|v|\leq10\varepsilon$. We also require $C\varepsilon\leq\frac{\widehat\psi(0)\widehat\eta(0)}{16}$ when $d\geq3$ and $C\varepsilon\leq\frac{\widehat\psi(0)}{8}$ when $d=2$, where $C$ is the constant in \eqref{stimeqR}. Such a choice is possible by the continuity of the Fourier transforms at the origin.

Parametrize the sphere near $(1,0,\ldots,0)$ by setting $\omega_1=(1-|\omega'|^2-\omega_d^2)^{\frac{1}{2}}$ and let us define the set $\Gamma_R$ as the family of those $(\omega_1,\omega',\omega_d)\in \mathbb S^{d-1}$ such that $|\omega'|\leq\frac{\varepsilon}{RD}$ and $|\omega_d|\leq\frac{\varepsilon}{(RL)^{\frac{1}{2}}}$. Observe that if $d=2$, only the second inequality is meaningful. Since $RL\geq1$ and $RD\geq1$ the area formula implies
$$\Haus^{d-1}(\Gamma_R)\geq c(RD)^{-d+2}(RL)^{-\frac{1}{2}},$$
where the constant $c$ depends only on $d$ and $\varepsilon$ and hence only on $d$, $\psi$ and $\eta$ for $d\geq 3$. 
The inequality $1-(1-t)^{\frac{1}{2}}\leq t$ for $0\leq t\leq1$ implies the following estimates
$$RL|1-\omega_1|\leq\begin{cases}\varepsilon^2(\frac{L}{RD^2}+1),&d\geq3,\\
\varepsilon^2,&d=2.\end{cases}$$

On $\Gamma_R$, one also has $|RD\omega'|\leq\varepsilon$ for $d\geq3$. Moreover, the above estimate and the assumption $\frac{L}{RD^2}\leq1$ implies that 
$$RL|1-\omega_1|\leq2\varepsilon^2\qquad(d\geq3).$$
It follows that
$$|RL(\omega_1-1)|\leq2\varepsilon^2\leq10\varepsilon,\qquad |RD\omega'|\leq\varepsilon\leq10\varepsilon,$$
where the second inequality is omitted when $d=2$, and thanks to the choice of $\varepsilon$, we infer that 
$$ |\widehat\psi (RL(\omega_1-1))|\geq\frac{\widehat\psi(0)}2\qquad\text{and, when }d\geq3,\qquad|\widehat\eta(RD\omega')|\geq\frac{\widehat\eta(0)}2,\qquad\text{for every }\omega\in \Gamma_R.$$
This implies, thanks to \eqref{stimeqR}, that for every $\omega\in \Gamma_R$ the following estimate for $q_R(\omega)$ holds 
\begin{equation}
    \begin{split}
        |q_R(\omega)|\geq &\begin{cases}\frac{aLD^{d-2}}{8}\widehat\psi (0)\widehat\eta(0)- Ca^2R|\omega_d|LD^{d-2},&d\geq3,\\
        \frac{aL}{4}\widehat\psi (0)-Ca^2R|\omega_d|L,&d=2,\end{cases}\\
        \geq &\begin{cases}\frac{aLD^{d-2}}{8}\widehat\psi (0)\widehat\eta(0)- C\varepsilon a^2(RL)^{\frac{1}{2}}D^{d-2},&d\geq3,\\
        \frac{aL}{4}\widehat\psi (0)-C\varepsilon a^2(LR)^{\frac{1}{2}},&d=2,\end{cases}
    \end{split}
\end{equation}
Since $aR\leq\rho\leq1$ and $RL\geq1$, we have 
$$a^2(RL)^{\frac{1}{2}}=aL\frac{aR}{(RL)^{\frac{1}{2}}}\leq aL.$$
The choice of $\varepsilon$ therefore allows us to absorb the error terms in the preceding estimate, and we conclude that
$$|q_R(\omega)|\geq \begin{cases}\frac{aLD^{d-2}}{16}\widehat\psi (0)\widehat\eta(0),&d\geq3,\\
\frac{aL}{8}\widehat\psi (0),&d=2,\end{cases}.$$
In particular absorbing $\widehat\psi(0)$ and $\widehat\eta(0)$ into the constant $c$, we have
$$|q_R(\omega)|\geq caLD^{d-2}\qquad \text{for }\omega\in\Gamma_R.$$
Applying Cauchy-Schwarz in $\phi$ gives $(2\pi)^{-1}\int_0^{2\pi}|\widehat{\mathbb1_{\Omega_{R,\phi}}}(R\omega)|^2 d\phi\geq|q_R(\omega)|^2$ and integrating first in $\omega$ and then in $\phi$ shows that some value of $\phi$ satisfies
$$E_{\Omega_{R,\phi}}(R)\geq ca^2R^{-d+\frac{3}{2}}L^{\frac{3}{2}}D^{d-2}.$$
This concludes the proof. 
\end{proof}

\begin{proposizione}\label{sequence-perturbations}
Fix $0<s<1$ and $1\leq p\leq\infty$ and let $R_k\geq2$ tend to infinity. Let $b_k>0$, $0<L_k\leq\rho$, $0<D_k\leq\rho$ and $x_k\in(-\frac{\rho}{4},\frac{\rho}{4})$ be such that the supports in the $x$ variable of the functions below are pairwise disjoint and their union has compact closure in $(-\frac{\rho}{2},\frac{\rho}{2})$. For every $k\geq1$ and every $\phi\in[0,2\pi]$, define
$$f_{k,\phi}(x,y)=\begin{cases}b_k\psi(\frac{x-x_k}{L_k})\eta(\frac{y}{D_k})\cos(R_k(x-x_k)+\phi),&d\geq3,\\
b_k\psi(\frac{x-x_k}{L_k})\cos(R_k(x-x_k)+\phi),&d=2.\end{cases}$$
Assume further that the following algebraic constraints are satisfied, with the conditions involving $D_k$ omitted when $d=2$,
\begin{equation}
    R_kL_k\geq1,\qquad b_kR_k\leq\rho, \qquad R_kD_k\geq1\qquad\text{and}\qquad\frac{L_k}{R_kD_k^2}\leq1,
    \label{equazionestimealgebriche}
\end{equation}
Suppose also that
$\sum_k b_kR_k(1+\log R_k)<\infty$ and that
$$\sum_k b_k^pR_k^{p(1+s)}L_kD_k^{d-2}<\infty \qquad\text{when $1\leq p<\infty$, and}\qquad \sup_k  b_kR_k^{1+s}<\infty\qquad\text{when $p=\infty$}.$$
Finally, we suppose that, for every $k\geq1$, there holds
\begin{equation}
    \begin{split}
        \sum_{j>k}b_jL_jD_j^{d-2}\leq\frac1{32}b_kL_kD_k^{d-2}.
        \label{tailestimate}
    \end{split}
\end{equation}
Under these assumptions there exist phases $\phi_k\in[0,2\pi]$ such that $F=\sum_kf_{k,\phi_k}$ belongs to $C_c^1(B^{d-1}(0,\rho))$ and the domain $\Omega_F:=\Omega_{x_0,\rho,F}$ has $\partial\Omega_F\in C^1\cap W^{1+s,p}$, with the usual interpretation $C^{1,s}$ when $p=\infty$, and further if $\omega_F$ is the modulus of continuity of the exterior normal of $\Omega_F$ we have that 
$$\int_0^1\omega_F(r)\frac{dr}{r}<\infty.$$
Moreover,
$$E_{\Omega_F}(R_k)\geq cb_k^2R_k^{-d+\frac{3}{2}}L_k^{\frac{3}{2}}D_k^{d-2},$$
and the domain $\Omega_F$ agrees with the prescribed domain $\Omega$ outside $B(x_0,3\rho)$.
\end{proposizione}

\begin{proof}
All the estimates in the first part of the proof are uniform with respect to the phases. We may therefore write temporarily $f_k=f_{k,\phi_k}$ and $F=\sum_kf_k$ for an arbitrary choice of the phases, which will be selected below.
Arguing as in \cref{one-perturbation}, with $a,R,L,D$ replaced by $b_k,R_k,L_k,D_k$, it is immediate to see that
$$\|f_k\|_\infty+\|\nabla f_k\|_\infty\leq Cb_kR_k,\qquad \|\nabla^2f_k\|_\infty\leq Cb_kR_k^2.$$ 
The first summability assumption gives uniform convergence of $F=\sum_kf_k$ and $\sum_k\nabla f_k$ and this implies that $F$ is of class $C^1$ and $\nabla F=\sum_k\nabla f_k$. The supports of $f_k$ are disjoint and $b_kR_k\leq\rho$, so the fixed normalization of $\psi$ and $\eta$ give $\|F\|_\infty+\|\nabla F\|_\infty\leq\rho$. Hence \cref{perturb-flat-patch} defines the domain $\Omega_F$.

Suppose first that $p<\infty$. The regularity estimate in \cref{one-perturbation} gives
$$[\nabla f_k]_{W^{s,p}}^p\leq Cb_k^pR_k^{p(1+s)}L_kD_k^{d-2}.$$
For every $u,v\in\mathbb R^{d-1}$, at most two of the vectors $\nabla f_k(u)-\nabla f_k(v)$ are nonzero. Therefore
$$|\nabla F(u)-\nabla F(v)|^p\leq2^{p-1}\sum_k|\nabla f_k(u)-\nabla f_k(v)|^p.$$
Integrating this inequality and using Tonelli's theorem, we obtain
\begin{equation*}
    \begin{split}
        [\nabla F]_{W^{s,p}(\mathbb R^{d-1})}^p=&\int_{\mathbb R^{d-1}}\int_{\mathbb R^{d-1}}\frac{|\nabla F(u)-\nabla F(v)|^p}{|u-v|^{d-1+sp}} du dv\leq2^{p-1}\sum_k\int_{\mathbb R^{d-1}}\int_{\mathbb R^{d-1}}\frac{|\nabla f_k(u)-\nabla f_k(v)|^p}{|u-v|^{d-1+sp}}dudv\\
        \leq &2^{p-1}\sum_k[\nabla f_k]_{W^{s,p}(\mathbb R^{d-1})}^p\leq C\sum_{k\in \N}b_k^pR_k^{p(1+s)}L_kD_k^{d-2}<\infty.
        \nonumber
    \end{split}
\end{equation*}
When $p=\infty$, the same observation gives
$$[\nabla F]_{C^{0,s}(\mathbb R^{d-1})}=\sup_{u\neq v}\frac{|\nabla F(u)-\nabla F(v)|}{|u-v|^s}\leq2\sup_{k\in \N}[\nabla f_k]_{C^{0,s}(\mathbb R^{d-1})}\leq C\sup_{k\in \N}b_kR_k^{1+s}.$$
The pointwise derivative estimates also give $\sup_{|u-v|\leq t}|\nabla F(u)-\nabla F(v)|\leq C\sum_kb_kR_k\min\{1,R_kt\}$ and since all terms are nonnegative, Tonelli's theorem yields
$$\int_0^1\sup_{|u-v|\leq t}|\nabla F(u)-\nabla F(v)| \frac{dt}{t}\leq C\sum_{k\in \N}b_kR_k\int_0^1\min\{1,R_kt\}\frac{dt}{t}\leq C\sum_{k\in \N}b_kR_k(1+\log R_k)<\infty,$$
which shows that the modulus of continuity of $\nabla F$ satisfies a Dini condition.

When $p<\infty$, since $\nabla F\in W^{s,p}$, \cref{normal} controls the $W^{s,p}$ seminorm of $\nu_{\Omega_F}$ on the perturbed graph. When $p=\infty$, the same formula controls its $C^{0,s}$ seminorm. The remaining part of the boundary is smooth. Choose finitely many boundary charts and $t_0>0$ such that any two boundary points at distance less than $t_0$ belong to a common chart. The preceding estimates control all such pairs, while the contribution of pairs at distance at least $t_0$ is finite because the normal is bounded and the boundary is compact. Hence the normal has the asserted fractional regularity on $\partial\Omega_F$. The same argument and \cref{normal} give
$$\omega_\nu(t)\leq C\Big(t+\sup_{|u-v|\leq Ct}|\nabla F(u)-\nabla F(v)|\Big),\qquad 0<t<t_0,$$
and the preceding estimate proves the Dini condition for $\nu_{\Omega_F}$.
It remains to choose the phases. For every $k$, define
$$q_k(\omega)=\frac1{2\pi}\int_0^{2\pi}T_{f_{k,\phi}}(R_k\omega)e^{-i\phi} d\phi\qquad\text{where}\qquad T_f(\xi):=\int e^{-i(\xi_1x+\langle \xi', y\rangle)}\int_0^{f(x,y)}e^{-i\xi_dz}dzdxdy.$$
The computation in the proof of \cref{one-perturbation}, with $a,R,L,D$ replaced by $b_k,R_k,L_k,D_k$ and with the translation $x\mapsto x-x_k$ yields a compact set $\Gamma_k$ in the sphere $\mathbb S^{d-1}$ such that 
$$|q_k(\omega)|\geq\begin{cases}\frac{\widehat\psi(0)\widehat\eta(0)}{16}b_kL_kD_k^{d-2},&d\geq3,\\[2mm]
\frac{\widehat\psi(0)}{8}b_kL_k,&d=2,\end{cases}\qquad \text{for every }\omega\in\Gamma_k,$$
and
$$\int_{\Gamma_k}|q_k(\omega)|^2 d\omega\geq cb_k^2R_k^{-d+\frac{3}{2}}L_k^{\frac{3}{2}}D_k^{d-2},$$
where $c$ depends only on $\psi$ and $\eta$. We now choose the phases inductively. Suppose that $\phi_1,\ldots,\phi_{k-1}$ have already been fixed, with no condition when $k=1$. For $\phi\in[0,2\pi]$, set
$$F_{k,\phi}=\sum_{j<k}f_{j,\phi_j}+f_{k,\phi},\qquad \text{and}\qquad\Omega_{k,\phi}=\Omega_{x_0,\rho,F_{k,\phi}}.$$
Since the supports of the functions are disjoint, the integral identity in \cref{perturb-flat-patch} gives
$$\widehat{\mathbb1_{\Omega_{k,\phi}}}(\xi)=\widehat{\mathbb1_{\Omega_0}}(\xi)+\sum_{j<k}T_{f_{j,\phi_j}}(\xi)+T_{f_{k,\phi}}(\xi).$$
The first two terms on the right-hand side are independent of $\phi$. Their integrals against $e^{-i\phi}$ vanish, and therefore
$$\frac1{2\pi}\int_0^{2\pi}\widehat{\mathbb1_{\Omega_{k,\phi}}}(R_k\omega)e^{-i\phi} d\phi=q_k(\omega).$$
The Cauchy-Schwarz inequality and integration on $\Gamma_k$ give
$$\|q_k\|_{L^2(\Gamma_k)}^2\leq\frac1{2\pi}\int_0^{2\pi}\|\widehat{\mathbb1_{\Omega_{k,\phi}}}(R_k \cdot)\|_{L^2(\Gamma_k)}^2 d\phi.$$
Hence there exists $\phi_k\in[0,2\pi]$ such that
$\|\widehat{\mathbb1_{\Omega_{k,\phi_k}}}(R_k \cdot)\|_{L^2(\Gamma_k)}\geq\|q_k\|_{L^2(\Gamma_k)}$.
We fix this value of $\phi_k$ and continue the induction. After all the phases have been chosen, set
$$f_k=f_{k,\phi_k},\qquad F=\sum_kf_k,\qquad \Omega^{(k)}=\Omega_{x_0,\rho,\sum_{j\leq k}f_j}.$$

For every $\xi\in\mathbb R^d$, the definition of $T_{f_j}$ implies 
$$|T_{f_j}(\xi)|\leq\|f_j\|_1\leq\begin{cases}\widehat\psi(0)\widehat\eta(0)b_jL_jD_j^{d-2},&d\geq3,\\
\widehat\psi(0)b_jL_j,&d=2.\end{cases}$$
The tail assumption with $k=1$ therefore implies that $\sum_jT_{f_j}(\xi)$ converges absolutely and uniformly in $\xi$. The disjointness of the supports and the integral identity in \cref{perturb-flat-patch} therefore give
$$\widehat{\mathbb1_{\Omega_F}}(\xi)=\widehat{\mathbb1_{\Omega^{(k)}}}(\xi)+\sum_{j>k}T_{f_j}(\xi).$$
Moreover, for every $\omega\in\Gamma_k$, we have
$$\Big|\sum_{j>k}T_{f_j}(R_k\omega)\Big|\leq\begin{cases}\frac{\widehat\psi(0)\widehat\eta(0)}{32}b_kL_kD_k^{d-2},&d\geq3,\\[2mm]
\frac{\widehat\psi(0)}{32}b_kL_k,&d=2,\end{cases}\leq\frac12|q_k(\omega)|.$$
Consequently, for every $k\in \N$ we have 
$$\|\widehat{\mathbb1_{\Omega_F}}(R_k \cdot)\|_{L^2(\Gamma_k)}\geq\|\widehat{\mathbb1_{\Omega^{(k)}}}(R_k \cdot)\|_{L^2(\Gamma_k)}-\Big\|\sum_{j>k}T_{f_j}(R_k \cdot)\Big\|_{L^2(\Gamma_k)}\geq\frac12\|q_k\|_{L^2(\Gamma_k)}.$$
It follows that
$$E_{\Omega_F}(R_k)\geq\int_{\Gamma_k}|\widehat{\mathbb1_{\Omega_F}}(R_k\omega)|^2 d\omega\geq\frac14\int_{\Gamma_k}|q_k(\omega)|^2 d\omega\geq cb_k^2R_k^{-d+\frac{3}{2}}L_k^{\frac{3}{2}}D_k^{d-2},$$
concluding the proof of the proposition. 
\end{proof}

\begin{definizione}
    Let us fix $d\geq 2$, $0<s<1$ and $1\leq p\leq \infty$. We let 
$$\gamma(d,s,p):=\begin{cases}
s\min\{\frac{3p}{2},2\},&d=2,\\
s\min\{p,2\},&d\geq3.\end{cases}$$
\end{definizione}

\begin{proposizione}\label{sharp-counterexamples}
Let $d\geq2$, $0<s<1$ and $1\leq p\leq\infty$. There exists a bounded domain $\Omega_*$ such that $\partial\Omega_*\in C^1\cap W^{1+s,p}$, whose exterior unit normal satisfies
$$\int_0^1\omega_{\nu_{\Omega_*}}(t)\frac{dt}{t}<\infty,\qquad\text{where}\qquad \omega_{\nu_{\Omega_*}}(t):=\sup_{|x-y|\leq t}|\nu(x)-\nu(y)|.$$
and such that
$$\limsup_{R\to\infty}\frac{\log E_{\Omega_*}(R)}{\log R}\geq -d-\min\left\{1,\frac{1}{2}+\gamma(d,s,p)\right\}.$$
If $\gamma(d,s,p)<\frac{1}{2}$, the domain $\Omega_*$ can be obtained by changing any prescribed bounded $C^\infty$ domain inside an arbitrarily small neighbourhood of any prescribed point of its boundary.
\end{proposizione}

\begin{proof}
Suppose first that $\gamma\geq\frac{1}{2}$. This is the least interesting case since the decay of the Fourier transform of the ball has such decay rate. Let $\Omega_*=B(0,1)$. This domain is smooth, and its exterior normal is Lipschitz, hence it has all the required fractional regularity and satisfies the Dini condition. The asymptotic formula for the Fourier transform of the ball gives
$$\widehat{\mathbb1_{B(0,1)}}(R\omega)=c_dR^{-\frac{d+1}{2}}\cos\Big(R-\frac{(d+1)\pi}{4}\Big)+O_d(R^{-\frac{d+3}{2}}).$$
Taking $R_k=2\pi k+\frac{(d+1)\pi}{4}$, we easily obtain thanks to radial symmetry of the ball that 
$$E_{B(0,1)}(R_k)\geq cR_k^{-d-1}.$$
The upper bound for domains with Dini normal and the inequality $\gamma\geq\frac{1}{2}$ gives
$$E_{B(0,1)}(R)\leq C\bigl(R^{-d-1}+R^{-d-\frac{1}{2}-\gamma}\bigr)\leq2CR^{-d-1}.$$
It follows that
$$\limsup_{R\to\infty}\frac{\log E_{B(0,1)}(R)}{\log R}=-d-1.$$

Suppose now that $\gamma<\frac{1}{2}$. Let $\Omega$ be a prescribed bounded $C^\infty$ domain, let $x_0\in\partial\Omega$, and fix $\rho<\min\{\frac{r_\Omega}{10},1\}$. We apply \cref{flat-patch} at $x_0$ and then use the notation introduced before \cref{one-perturbation}. For every $k\geq1$, set
$$x_k=\rho2^{-k-2},\qquad c_k=\rho2^{-k-4},\qquad \ell_k=\rho2^{-k-6},\qquad R_k=e^{Mk^2},$$
where $M>0$ will be fixed below. When $d\geq3$, define
$$
(b_k,L_k,D_k)=\begin{cases}
(c_kR_k^{-1},\ell_k,\rho R_k^{-\frac{sp}{d-2}}),&1\leq p<2,\\
(c_kR_k^{-1-s},\ell_k,\rho),&2\leq p\leq\infty.
\end{cases}
$$
When $d=2$, define
$$
(b_k,L_k)=\begin{cases}
(c_kR_k^{-1},\ell_kR_k^{-sp}),&1\leq p<\frac{4}{3},\\
(c_kR_k^{-1-s},\ell_k),&\frac{4}{3}\leq p\leq\infty.
\end{cases}
$$

We now verify separately that the four choices of the parameters satisfy all the assumptions of \cref{sequence-perturbations}. We fix $M$ once and for all so that
$$M\geq\frac{\log(\frac{33}{4})}{1+s}+\log(128\rho^{-1})\begin{cases}
\frac{1}{1-\frac{2sp}{d-2}},&d\geq3,\ 1\leq p<2,\\
1,&d\geq3,\ 2\leq p\leq\infty,\\
\frac{1}{1-sp},&d=2,\ 1\leq p<\frac{4}{3},\\
1,&d=2,\ \frac{4}{3}\leq p\leq\infty.
\end{cases}$$

First, let us observe that the functions $f_k$ have disjoint supports. Since $L_k\leq\ell_k$ and $\psi$ is supported in $(-\frac{1}{16},\frac{1}{16})$, the support in the $x$ variable of $f_{k,\phi}$ is contained in the interval $(x_k-\frac{\ell_k}{16},x_k+\frac{\ell_k}{16})$.
For every $k\geq1$, our definitions imply that
$$x_k-x_{k+1}=\rho2^{-k-3},\qquad\frac{\ell_k+\ell_{k+1}}{16}=3\rho2^{-k-11}<\rho2^{-k-3}.$$
Thus these intervals are pairwise disjoint. Since $x_k\to0$ and $x_1+\frac{\ell_1}{16}<\frac{\rho}{2}$, their union has compact closure in $(-\frac{\rho}{2},\frac{\rho}{2})$. We now check \eqref{equazionestimealgebriche} in the four cases. 

\smallskip

Suppose first that $d\geq3$ and $1\leq p<2$. In this case $0<L_k\leq\rho$, $0<D_k\leq\rho$, and
$b_kR_k=c_k=\rho2^{-k-4}\leq\rho$ and thanks to the definition of $\gamma$, we have that $sp<\frac{1}{2}$.  This implies in particular that
$$M\geq\frac{\log(\frac{128}{\rho})}{1-\frac{2sp}{d-2}}\geq\log\left(\frac{128}{\rho}\right).$$
Therefore, the choice of our parameters imply 
$$R_kL_k=\rho2^{-k-6}e^{Mk^2}\geq\rho2^{-k-6}\Big(\frac{128}{\rho}\Big)^k=2^{6k-6}\rho^{1-k}\geq1.$$
$$R_kD_k=\rho e^{M(1-\frac{sp}{d-2})k^2}\geq\rho e^{M(1-\frac{2sp}{d-2})}\geq128.$$
$$\frac{L_k}{R_kD_k^2}=\rho^{-1}2^{-k-6}e^{-M(1-\frac{2sp}{d-2})k^2}\leq\rho^{-1}2^{-k-6}\frac{\rho}{128}=2^{-k-13}\leq1.$$
Thus all the algebraic assumptions are satisfied. For the fractional regularity, direct substitution gives
$$b_k^pR_k^{p(1+s)}L_kD_k^{d-2}=\rho^{d-2}c_k^p\ell_k.$$

\smallskip

Suppose next that $d\geq3$ and $2\leq p\leq\infty$. As in the previous case,
$$R_kL_k=\rho2^{-k-6}e^{Mk^2}\geq1,\qquad R_kD_k=\rho e^{Mk^2}\geq128\geq1,$$
$$\frac{L_k}{R_kD_k^2}=\rho^{-1}2^{-k-6}e^{-Mk^2}\leq2^{-k-13}\leq1,\qquad\text{and}\qquad b_kR_k=c_kR_k^{-s}\leq c_k\leq\rho.$$
We further have 
$$b_k^pR_k^{p(1+s)}L_kD_k^{d-2}=\rho^{d-2}c_k^p\ell_k,\qquad\text{when $1\leq p<\infty$, while }\qquad b_kR_k^{1+s}=c_k\qquad \text{for }p=\infty.$$

\smallskip

Now let $d=2$ and $1\leq p<\frac{4}{3}$. In this case, by construction, $0<L_k\leq \rho$,
$b_kR_k=c_k\leq\rho$, and
$$R_kL_k=\rho2^{-k-6}e^{M(1-sp)k^2}\geq\rho2^{-k-6}\Big(\frac{128}{\rho}\Big)^k=2^{6k-6}\rho^{1-k}\geq1\quad\text{and}\quad b_k^pR_k^{p(1+s)}L_k=c_k^pR_k^{sp}\ell_kR_k^{-sp}=c_k^p\ell_k.$$

Finally, suppose that $d=2$ and $\frac{4}{3}\leq p\leq\infty$. In this case we have
$$R_kL_k=\rho2^{-k-6}e^{Mk^2}\geq1,\qquad b_kR_k=c_kR_k^{-s}\leq c_k\leq\rho.$$
$$b_k^pR_k^{p(1+s)}L_k=c_k^p\ell_k\qquad\text{for }1\leq p<\infty\text{ and}\qquad b_kR_k^{1+s}=c_k\qquad\text{for }p=\infty.$$

We can now verify the summability assumptions simultaneously. When $p<\infty$, the preceding identities give
$$\sum_kb_k^pR_k^{p(1+s)}L_kD_k^{d-2}= \rho^{d-2}\sum_kc_k^p\ell_k<\infty,$$
while for $p=\infty$ we have
$$\sup_kb_kR_k^{1+s}=\sup_kc_k<\infty.$$
Moreover, in each of the four cases $b_kR_k\leq c_k$, and hence
$$\sum_kb_kR_k(1+\log R_k)\leq\sum_kc_k(1+Mk^2)<\infty.$$

It remains to check the tail assumption \eqref{tailestimate}. In the four cases we have
$$
b_kL_kD_k^{d-2}=\begin{cases}
\rho^{d-2}c_k\ell_kR_k^{-1-sp},&d\geq3,\ 1\leq p<2,\\
\rho^{d-2}c_k\ell_kR_k^{-1-s},&d\geq3,\ 2\leq p\leq\infty,\\
c_k\ell_kR_k^{-1-sp},&d=2,\ 1\leq p<\frac{4}{3},\\
c_k\ell_kR_k^{-1-s},&d=2,\ \frac{4}{3}\leq p\leq\infty.
\end{cases}
$$
Let us fix $j<k$ and observe that 
$$\frac{b_kL_kD_k^{d-2}}{b_jL_jD_j^{d-2}}\leq \frac{c_k\ell_k}{c_j\ell_j}\Big(\frac{R_j}{R_k}\Big)^{1+sq(p,s,d)},\qquad\text{where }q(p,s,d):=\begin{cases}
p,&d\geq3,\ 1\leq p<2,\\
1,&d\geq3,\ 2\leq p\leq\infty,\\
p,&d=2,\ 1\leq p<\frac{4}{3},\\
1,&d=2,\ \frac{4}{3}\leq p\leq\infty.\end{cases}$$
In these notations we have 
$$\frac{b_kL_kD_k^{d-2}}{b_jL_jD_j^{d-2}}=2^{-2(k-j)}e^{-M(1+sq(p,s,d))(k^2-j^2)}\leq \Big(\frac{e^{-M(1+s)}}{4}\Big)^{k-j}.$$
The choice of $M$ gives $\frac{e^{-M(1+s)}}{4}\leq\frac{1}{33}$ and as a consequence
$$\frac{\sum_{j>k}b_jL_jD_j^{d-2}}{b_kL_kD_k^{d-2}}\leq\sum_{r=1}^\infty\Big(\frac1{33}\Big)^r=\frac1{32}.$$
We have now verified every assumption of \cref{sequence-perturbations}. Thus, 
\cref{sequence-perturbations} gives phases $\phi_k$ and a single function $F=\sum_kf_{k,\phi_k}$. The domain
$$\Omega_*=\Omega_{x_0,\rho,F}$$
has the required fractional regularity and agrees with the prescribed domain $\Omega$ outside $B(x_0,3\rho)$ and the modulus of continuity of the normal satisfies the required Dini integrability condition.

We are left to prove the lower bounds for $E_{\Omega_*}(R_k)$. Thanks to \cref{sequence-perturbations}, we have 
\begin{equation*}
\begin{split}
E_{\Omega_F}(R_k)&\geq cb_k^2R_k^{-d+\frac{3}{2}}L_k^{\frac{3}{2}}D_k^{d-2}=\begin{cases}
c\rho^{d-2}c_k^2\ell_k^{\frac{3}{2}}R_k^{-d-\frac{1}{2}-sp},&d\geq3,\ 1\leq p<2,\\
c\rho^{d-2}c_k^2\ell_k^{\frac{3}{2}}R_k^{-d-\frac{1}{2}-2s},&d\geq3,\ 2\leq p\leq\infty,\\
cc_k^2\ell_k^{\frac{3}{2}}R_k^{-\frac{5}{2}-\frac{3sp}{2}},&d=2,\ 1\leq p<\frac{4}{3},\\
cc_k^2\ell_k^{\frac{3}{2}}R_k^{-\frac{5}{2}-2s},&d=2,\ \frac{4}{3}\leq p\leq\infty
\end{cases}  \geq cc_k^2\ell_k^{\frac{3}{2}}R_k^{-d-\frac{1}{2}-\gamma}.
    \nonumber
\end{split}
\end{equation*}
Since
$c_k^2\ell_k^{\frac{3}{2}}=\rho^{\frac{7}{2}}2^{-17}2^{-\frac{7k}{2}}$ and $\log R_k=Mk^2$ we obtain
$$\limsup_{R\to\infty}\frac{\log E_{\Omega_*}(R)}{\log R}\geq-d-\frac12-\gamma,$$
and this concludes the proof. 
\end{proof}

\subsection{Scales of optimal decay}

The following consequence of a result of Beretti and Gennaioli gives arbitrarily large scales at which every set of finite perimeter has decay of order $R^{-d-1}$.

\begin{corollario}\label{perimeter-frequency-average}
Let $B\subset\mathbb R^d$ be measurable, with finite measure and finite perimeter $\operatorname{Per}(B)=\mathcal H^{d-1}(\partial^*B)$, where $\partial^*B$ denotes the reduced boundary of $B$. Then
$$\lim_{R\to\infty}\frac1R\int_0^Rr^{d+1}E_B(r) dr=2^d\pi^{d-1}\operatorname{Per}(B).$$
Consequently,
$$\liminf_{R\to\infty}R^{d+1}E_B(R)\leq2^d\pi^{d-1}\operatorname{Per}(B)\leq\limsup_{R\to\infty}R^{d+1}E_B(R).$$
\end{corollario}

\begin{proof}
Beretti and Gennaioli prove in \cite[Corollary 1.4]{berettigennaioli} that
$$
\lim_{T\to\infty}\frac{2\pi^2}{T}\int_{|\eta|\leq T}|\eta|^2\Big|\int_Be^{-2\pi i\langle \eta, x\rangle } dx\Big|^2 d\eta=\operatorname{Per}(B).
$$
Their Fourier transform differs from ours by a factor $2\pi$ in the frequency variable, since $\int_Be^{-2\pi i\langle\eta,x\rangle}dx=\widehat{\mathbb1_B}(2\pi\eta)$. Integrating in polar coordinates and substituting $\xi=2\pi\eta$ and $R=2\pi T$, we find
$$
\frac{1}{R}\int_0^Rr^{d+1}E_B(r)dr=\frac{1}{R}\int_{|\xi|\leq R}|\xi|^2|\widehat{\mathbb1_B}(\xi)|^2d\xi=\frac{(2\pi)^{d+1}}{T}\int_{|\eta|\leq T}|\eta|^2|\widehat{\mathbb1_B}(2\pi\eta)|^2d\eta,
$$
and the right-hand side tends to $\frac{(2\pi)^{d+1}}{2\pi^2}\operatorname{Per}(B)=2^d\pi^{d-1}\operatorname{Per}(B)$. This is the first assertion. The second follows at once, because the function $r\mapsto r^{d+1}E_B(r)$ cannot stay above $2^d\pi^{d-1}\operatorname{Per}(B)+\varepsilon$, or below $2^d\pi^{d-1}\operatorname{Per}(B)-\varepsilon$, for all large $r$ while its averages converge to $2^d\pi^{d-1}\operatorname{Per}(B)$.
\end{proof}

Thus every set of finite perimeter has arbitrarily large scales at which its spherical Fourier average is bounded from above by a fixed multiple of $R^{-d-1}$. No regularity of its boundary is needed.

\section{Upper bounds}

The main result of this section consists of the following upper bounds for the spherical average of the Fourier transform.

\begin{proposizione}\label{prop:upper}
Let $d\geq2$, $0<s<1$ and $1\leq p\leq\infty$, and let $\Omega\subset\mathbb R^d$ be a bounded domain with $\partial\Omega\in C^1\cap W^{1+s,p}$. If $1\leq p<\infty$, then there exists a constant such that for every $R\geq2$, we have
\begin{equation}
    E_\Omega(R)\leq C_\Omega R^{-d-1}+C_\Omega R^{-\frac{d+3}{2}-s\min\{\frac{3p}{2},2\}}.
    \label{stimafinaleI}
\end{equation}
In addition, if $d\geq3$ and $1\leq p<\infty$, one also has
$$E_\Omega(R)\leq C_\Omega R^{-d-1}+C_\Omega\begin{cases}R^{-d-\frac{1}{2}-s\min\{p,2\}},&sp>d-2,\\
(\log(2+R))^{\frac{p}{2}}R^{-d-\frac{1}{2}-s\min\{p,2\}},&sp=d-2,\\
(\log(2+R))^{\frac{1}{2}}R^{-\frac{d+3+sp}{2}-s\min\{p,2\}},&sp<d-2.\end{cases}$$
In addition, if the modulus of continuity of the outward normal satisfies $\int_0^1\omega_{\nu_\Omega}(t) \frac{dt}{t}<\infty$, then
$$E_\Omega(R)\leq C_\Omega R^{-d-1}+C_\Omega R^{-d-\frac{1}{2}-\gamma}.$$
\end{proposizione}

\subsection{The reduction to the graph case}

In this subsection we prove that \cref{prop:upper} can be reduced to the proof of the estimates for a graph. 
This was actually the starting point of our argument. We realized that proving
the bounds for domains essentially reduces to proving suitable
estimates for graphs. The local perturbations producing the lower
bounds first suggested this reduction. The same viewpoint is
suggested by the oscillatory integrals themselves: for each fixed
$\omega$, the phase $x\mapsto\langle\omega,x\rangle$ is constant
on hyperplanes orthogonal to $\omega$, and the boundary contributions
that may cancel can be grouped along these hyperplanes. This suggests
that the relevant pieces can, at least heuristically, be assembled
into a single graph. We therefore first establish the upper bounds
assuming \cref{spherical-graph-estimate}, whose proof is given
at the end of this section.

\begin{proof}[Proof of \cref{prop:upper}, assuming \cref{spherical-graph-estimate}] We divide the proof into several parts.

\smallskip

\spacedlowsmallcaps{Reduction to graph integrals.} We first cover the boundary by graphs and choose a partition of unity. We use fixed cubes and impose the gradient bound required later by \cref{spherical-graph-estimate}. Set
$$
U_0:=(-1,1)^{d-1},\qquad U_1:=(-2,2)^{d-1},\qquad U:=(-3,3)^{d-1},\qquad c_0:=\frac{1}{8\sqrt{d-1}}.
$$
For $d\geq3$, let $\varepsilon_0$ be given by \cref{spherical-graph-estimate} for these cubes and $c_0$. By compactness, we can choose points $x_1,\ldots,x_m\in\partial\Omega$, numbers $r_j>0$, orthogonal transformations $Q_j$ with $Q_je_d=\nu_\Omega(x_j)$, and functions $A_j\in C^1(U)$ such that
$$
X_j(u):=x_j+r_jQ_j(u,A_j(u)),\qquad u\in U,
$$
parametrize portions of $\partial\Omega$, the sets $X_j(U_0)$ cover $\partial\Omega$, and
$$
A_j(0)=0,\qquad \nabla A_j(0)=0,\qquad \|\nabla A_j\|_{L^\infty(U)}\leq\min\left\{\frac{1}{10},c_0\right\}.
$$
When $d\geq3$, we also require $\|\nabla A_j\|_{L^\infty(U)}\leq\varepsilon_0$. 
Choose a nonnegative $C^1$ partition of unity $\{\chi_j\}_{j=1}^m$ subordinate to the cover $\{X_j(U_0)\}_{j=1}^m$, and define $\psi_j:=\chi_j\circ X_j\in C_c^1(U_0)$.  Define
$$
J_{R,j}(\omega):=\int_{\partial^*\Omega}\chi_j(x)e^{-iR\langle\omega,x\rangle}\langle\omega,\nu_\Omega(x)\rangle\,d\mathcal H^{d-1}(x),
$$
and notice that the divergence theorem and the Cauchy-Schwarz inequality imply
\begin{equation}
E_\Omega(R)=R^{-2}\int_{\mathbb S^{d-1}}\Big|\sum_{j=1}^mJ_{R,j}(\omega)\Big|^2\,d\omega\leq mR^{-2}\sum_{j=1}^m\int_{\mathbb S^{d-1}}|J_{R,j}(\omega)|^2\,d\omega.
    \label{stimolatu}
\end{equation}
It remains to estimate one integral over a fixed graph. Let us fix $j$, drop the dependence of all the objects above on $j$, let $X(u):=(u,A(u))$ and let us introduce the operator
\begin{equation}
J_R(\omega):=\int_Ue^{-iR\langle\omega,X(u)\rangle}\psi(u)\langle\omega,(-\nabla A(u),1)\rangle\,du.
\label{graph-integral}
\end{equation}
The area formula gives the identity
$$J_{R,j}(\omega)=r_j^{d-1}e^{-iR\langle\omega,x_j\rangle}J_{r_jR}(Q_j^T\omega).$$
Since orthogonal maps preserve spherical measure, bounds for \eqref{graph-integral} give bounds for $J_{R,j}$ after replacing the radius by $r_jR$ and multiplying by $r_j^{2d-2}$. There are only finitely many fixed values of $r_j$, so these factors enter only the final constant. By the assumption on $\partial\Omega$, each $\nabla A_j$ belongs to $W^{s,p}(U)$ when $p<\infty$ and is $s$-H\"older continuous when $p=\infty$.

\smallskip

\spacedlowsmallcaps{Decay along the good directions.} We first estimate the contribution of $\{|\pi_H\omega|\leq\frac{1}{2}\}$, where $\pi_H$ is the projection onto the first $d-1$ coordinates. Choose $\eta\in C^\infty(\mathbb S^{d-1})$, with $0\leq\eta\leq1$, equal to one on this set and supported in $\{|\pi_H\omega|<\frac{3}{4}\}$. For fixed $u,v\in U$, we apply the second estimate in \cref{spherical-kernel} with geometric parameter $\frac{2}{5}$ and the choices
$$
c(\omega):=\eta(\omega)\langle\omega,(-\nabla A(u),1)\rangle\langle\omega,(-\nabla A(v),1)\rangle,\qquad z:=X(u)-X(v).
$$
We verify the hypotheses separately on the upper and lower hemispheres. Since $|\omega_d|>\frac{\sqrt7}{4}$ on $\supp\eta$, the corresponding restrictions of $c$ are smooth and supported in the respective coordinate charts. Here $u,v$ are fixed, so angular differentiation acts only on $\eta$ and on the two linear functions of $\omega$. Since $|(-\nabla A(u),1)|\leq2$, the angular derivatives of $c$ through order $d$ are bounded uniformly in $u,v$ by a constant depending only on $d$ and $\eta$. To check the geometric hypothesis, convexity of $U$ and the slope bound give
$$
|\langle z,\omega\rangle|\leq|\pi_H\omega||u-v|+|\omega_d||A(u)-A(v)|\leq\Big(\frac34+\frac1{10}\Big)|u-v|\leq\frac9{10}|z|.
$$
Consequently, an elementary computation shows  that $|z-\langle z,\omega\rangle\omega|\geq \frac{2}{5}|z|$. Thus, \cref{spherical-kernel} applies and we deduce that 
$$\Big|\int_{\mathbb S^{d-1}}e^{-iR\omega\cdot z}c(\omega) d\omega\Big|
\leq 4C_{d,\frac{2}{5}}(1+R|z|)^{-d}\sum_{j=0}^d\|\nabla_{\mathbb S}^j\eta\|_\infty\leq 4C_{d,\frac{2}{5}}(1+R|u-v|)^{-d}\sum_{j=0}^d\|\nabla_{\mathbb S}^j\eta\|_\infty,$$
where in the last inequality we used the fact that $\|\nabla A\|\leq \frac{1}{10}$. This implies in particular
\begin{equation}
\begin{split}
&\int_{\mathbb S^{d-1}}\eta(\omega)|J_R(\omega)|^2\,d\omega
=\int_U\int_U\psi(u)\psi(v)\int_{\mathbb S^{d-1}}e^{-iR\langle\omega,X(u)-X(v)\rangle}\eta(\omega)\langle\omega,(-\nabla A(u),1)\rangle\langle\omega,(-\nabla A(v),1)\rangle\,d\omega\,du\,dv\\
&\leq  4d^2C_{d,\frac{2}{5}}\sum_{j=0}^d\|\nabla_{\mathbb S}^j\eta\|_\infty \int_U\int_U\frac{\psi(u)\psi(v)}{(1+R|u-v|)^{d}}dudv \leq4d^2C_{d,\frac{2}{5}}\sum_{j=0}^d\|\nabla_{\mathbb S}^j\eta\|_\infty\int_U\int_U\frac{|\psi(u)|^2+|\psi(v)|^2}{(1+R|u-v|)^d}\,du\,dv\\
&\leq4d^2C_{d,\frac{2}{5}}\sum_{j=0}^d\|\nabla_{\mathbb S}^j\eta\|_\infty\|\psi\|_{L^2(U)}^2\int_{\mathbb R^{d-1}}(1+R|z|)^{-d}\,dz=4d^2\Leb^{d-1}(B(0,1))C_{d,\frac{2}{5}}R^{-(d-1)}\sum_{j=0}^d\|\nabla_{\mathbb S}^j\eta\|_\infty \|\psi\|_{L^2(U)}^2.
\label{stimaconeta1}
\end{split}
\end{equation}

\smallskip

\spacedlowsmallcaps{The contribution of the bad directions.} We shall integrate by parts to replace the coefficient in \eqref{graph-integral} by differences of $\psi$ and $\nabla A$ over distances at most $R^{-1}$. The bad directions are those directions that are more or less aligned with the tangent plane on which we are constructing the graph. We will split the strip of the sphere close to the tangent and we will divide the strip into patches where all the elements in $\mathbb{S}^{d-1}$ on which we are integrating have approximately the same projection $\tilde \omega$ on the tangent. This is where \cref{spherical-graph-estimate} will come into play. 

Write $1-\eta$ as a finite sum of nonnegative smooth functions $\{\vartheta\in \Theta\}$, each supported where one of the first $d-1$ coordinates has absolute value greater than $2c_0$. We will prove the estimate for a function $\vartheta$ supported where $\omega_1>2c_0$ since the estimates we will obtain will be essentially invariant by rotations of the plane $\omega_d=0$. Since we are assuming $\omega_1>2c_0$, then 
$$
q(u,\omega):=\partial_1\langle\omega,X(u)\rangle=\omega_1+\omega_d\partial_1A(u)\geq2c_0-\|\nabla A\|_\infty\geq c_0,\qquad\text{for every } u\in U_1.
$$
Let us notice that the quotient $\frac{\psi(u)\langle\omega,(-\nabla A(u),1)\rangle}{q(u,\omega)}$ is continuous. Choose a nonnegative $\rho\in C_c^\infty(B^{d-1}(0,1))$ with integral one. Take $R$ large enough that $\supp\psi+\mathrm{cl}(B(0,R^{-1}))\Subset U_0$ and $U_0+B(0,R^{-1})\subset U_1$. Observe that 
$$
q(u,\omega)e^{-iR\langle\omega,X(u)\rangle}=-\frac1{iR}\partial_{u_1}e^{-iR\langle\omega,X(u)\rangle},
$$
therefore integrating by parts
\begin{equation*}
\begin{split}
&\int_Ue^{-iR\langle\omega,X(u)\rangle}q(u,\omega)\int_{U_1}R^{d-1}\rho(R(u-v))\frac{\psi(v)\langle\omega,(-\nabla A(v),1)\rangle}{q(v,\omega)}\,dv\,du\\
&\qquad=-\frac1{iR}\int_U\partial_{u_1}(e^{-iR\langle\omega,X(u)\rangle})\int_{U_1}R^{d-1}\rho(R(u-v))\frac{\psi(v)\langle\omega,(-\nabla A(v),1)\rangle}{q(v,\omega)}\,dv\,du\\
&\qquad=\frac1{iR}\int_Ue^{-iR\langle\omega,X(u)\rangle}\int_{U_1}R^{d-1}\partial_{u_1}(\rho(R(u-v)))\frac{\psi(v)\langle\omega,(-\nabla A(v),1)\rangle}{q(v,\omega)}\,dv\,du\\
&\qquad=\frac1i\int_Ue^{-iR\langle\omega,X(u)\rangle}\int_{U_1}R^{d-1}(\partial_1\rho)(R(u-v))\frac{\psi(v)\langle\omega,(-\nabla A(v),1)\rangle}{q(v,\omega)}\,dv\,du.
\nonumber
\end{split}
\end{equation*}
For fixed $u\in U_0$, the change of variables $y=R(u-v)$ implies that
\begin{equation}
\int_{U_1}R^{d-1}\rho(R(u-v))\,dv=1,\qquad \int_{U_1}R^{d-1}(\partial_1\rho)(R(u-v))\,dv=0.
    \label{propovviemapotenti}
\end{equation}
This shows that the following identity holds true
\begin{equation}
\begin{split}
&\int_{U_0}e^{-iR\langle\omega,X(u)\rangle}\frac{\psi(u)\langle\omega,(-\nabla A(u),1)\rangle}{q(u,\omega)}\int_{U_1}R^{d-1}(q(u,\omega)\rho(R(u-v))+i(\partial_1\rho)(R(u-v)))dvdu\\
&\qquad\qquad\qquad\qquad\qquad=\int_{U_0}e^{-iR\langle\omega,X(u)\rangle}\psi(u)\langle\omega,(-\nabla A(u),1)\rangle\,du=J_R(\omega).
\label{espressioneperJR}
\end{split}
\end{equation}
On the other hand, the integration by parts just proved gives
\begin{equation}
\begin{split}
&\int_{U_0}e^{-iR\langle\omega,X(u)\rangle}\int_{U_1}R^{d-1}\Big(q(u,\omega)\rho(R(u-v))+i(\partial_1\rho)(R(u-v))\Big)\frac{\psi(v)\langle\omega,(-\nabla A(v),1)\rangle}{q(v,\omega)}\,dv\,du\\
&\qquad=\Big(\frac1i+i\Big)\int_{U_0}e^{-iR\langle\omega,X(u)\rangle}\int_{U_1}R^{d-1}(\partial_1\rho)(R(u-v))\frac{\psi(v)\langle\omega,(-\nabla A(v),1)\rangle}{q(v,\omega)}\,dv\,du=0.
\label{riscritturadi0masimileaJR}
\end{split}
\end{equation}
Subtracting \eqref{espressioneperJR} and \eqref{riscritturadi0masimileaJR}
\begin{equation}
\begin{split}
\frac{J_R(\omega)}{R^{d-1}}&=\int_{U_0}e^{-iR\langle\omega,X(u)\rangle}a_R(u,\omega)du,\\
a_R(u,\omega)&:=\int_{U_1}(q(u,\omega)\rho(R(u-v))+i(\partial_1\rho)(R(u-v)))\\
&\times\left(\frac{\psi(u)\langle\omega,(-\nabla A(u),1)\rangle}{q(u,\omega)}-\frac{\psi(v)\langle\omega,(-\nabla A(v),1)\rangle}{q(v,\omega)}\right)dv.
\end{split}
\label{regularized-graph-integral}
\end{equation}
For $u,v\in U_1$ and since we are supposing $\omega_1>2c_0$, we have
\begin{equation*}
\begin{split}
&\Big|\frac{\psi(u)\langle\omega,(-\nabla A(u),1)\rangle}{q(u,\omega)}-\frac{\psi(v)\langle\omega,(-\nabla A(v),1)\rangle}{q(v,\omega)}\Big|\\
&\qquad\qquad\leq\frac{2}{c_0}|\psi(u)-\psi(v)|+\frac{\|\psi\|_{L^\infty(U)}}{c_0}|\nabla A(u)-\nabla A(v)|+\frac{2\|\psi\|_{L^\infty(U)}}{c_0^2}|q(u,\omega)-q(v,\omega)|\\
&\qquad\qquad\leq\frac{2}{c_0}\|\nabla\psi\|_{L^\infty(U)}|u-v|+\Big(\frac1{c_0}+\frac{2}{c_0^2}\Big)\|\psi\|_{L^\infty(U)}|\nabla A(u)-\nabla A(v)|.
\nonumber
\end{split}
\end{equation*}

We now insert the preceding estimate into the integral in $v$. Since $|q(u,\omega)|\leq2$ and both $\rho(R(u-v))$ and $(\partial_1\rho)(R(u-v))$ vanish unless $|u-v|\leq R^{-1}$, we obtain
\begin{equation}
\begin{split}
&R^{d-1}\left|\int_{U_1}(q(u,\omega)\rho(R(u-v))+i(\partial_1\rho)(R(u-v)))\left(\frac{\psi(u)\langle\omega,(-\nabla A(u),1)\rangle}{q(u,\omega)}-\frac{\psi(v)\langle\omega,(-\nabla A(v),1)\rangle}{q(v,\omega)}\right)\,dv\right|\\
&\quad\leq(2\|\rho\|_\infty+\|\partial_1\rho\|_\infty)R^{d-1}\int_{B(u,R^{-1})}\Big(\frac{2}{c_0}\|\nabla\psi\|_\infty|u-v|+\Big(\frac1{c_0}+\frac2{c_0^2}\Big)\|\psi\|_\infty|\nabla A(u)-\nabla A(v)|\Big)\,dv\\
&\quad\leq\Leb^{d-1}(B^{d-1}(0,1))(2\|\rho\|_\infty+\|\partial_1\rho\|_\infty)\Big(\frac{2\|\nabla\psi\|_\infty}{c_0R}+\Big(\frac1{c_0}+\frac2{c_0^2}\Big)\|\psi\|_\infty\fint_{B(u,R^{-1})}|\nabla A(u)-\nabla A(v)|\,dv\Big).
\label{stimaaRpd}
\end{split}
\end{equation}
To write these estimates, denote by $a_R(u,\omega)$ the integral in $v$ in \eqref{regularized-graph-integral}. Thus there holds
$$
J_R(\omega)=R^{d-1}\int_{U_0}e^{-iR\langle\omega,X(u)\rangle}a_R(u,\omega)du
$$
and \eqref{stimaaRpd} implies that
\begin{equation*}
\begin{split}
|J_R(\omega)|&\leq\alpha_{d-1}(2\|\rho\|_\infty+\|\partial_1\rho\|_\infty)\Big(\frac{2|U_0|\|\nabla\psi\|_\infty}{c_0R}+\Big(\frac1{c_0}+\frac2{c_0^2}\Big)\|\psi\|_\infty\int_{U_0}\fint_{B(u,R^{-1})}|\nabla A(u)-\nabla A(v)|dvdu\Big),
\nonumber
\end{split}
\end{equation*}
where $\alpha_{d-1}:=\Leb^{d-1}(B^{d-1}(0,1))$. Let us observe that 
\begin{equation}
\begin{split}
\int_{\mathbb S^{d-1}}\vartheta(\omega)|J_R(\omega)|^2\,d\omega&=R^{2d-2}\int_{U_0}\int_{U_0}\int_{\mathbb S^{d-1}}e^{-iR\langle\omega,X(u)-X(v)\rangle}\vartheta(\omega)a_R(u,\omega)\overline{a_R(v,\omega)}\,d\omega\,dv\,du,
\end{split}
\label{upper-square}
\end{equation}
where the function $a_{R}(u,\omega)$ was introduced in \eqref{regularized-graph-integral}. 
\smallskip

\spacedlowsmallcaps{The reduction for the bad directions.} We shall prove that there is a constant $\newC\label{C:upper-first}$, depending only on $d,c_0,\rho,\vartheta,\|\psi\|_{C^1(U)}$, such that
\begin{equation}
\int_{\mathbb S^{d-1}}\vartheta(\omega)|J_R(\omega)|^2\,d\omega\leq\oldC{C:upper-first}R^{-\frac{d-1}{2}}\big\|R^{-1}+\fint_{B(\cdot,R^{-1})}|\nabla A(\cdot)-\nabla A(v)|dv\big\|_{L^{\frac{4}{3}}(U_0)}^2.
\label{upper-first-reduction}
\end{equation}
Let us first show how this estimate gives the first bound in the proposition. Assuming \eqref{upper-first-reduction}, it remains to estimate the norm on its right using $\nabla A\in W^{s,p}(U)$.
Observe that, for $p<\infty$, Jensen's inequality gives
\begin{equation}
\begin{split}
&\frac{R^{d-1}}{\alpha_{d-1}}\int_{U_0}\int_{B(u,R^{-1})}|\nabla A(u)-\nabla A(v)|^p\,dv\,du
=\frac{R^{d-1}}{\alpha_{d-1}}\int_{U_0}\int_{B(u,R^{-1})}\frac{|\nabla A(u)-\nabla A(v)|^p}{|u-v|^{d-1+sp}}|u-v|^{d-1+sp}\,dv\,du\\
\leq&\frac{R^{-sp}}{\alpha_{d-1}}\int_{U_0}\int_{B(u,R^{-1})}\frac{|\nabla A(u)-\nabla A(v)|^p}{|u-v|^{d-1+sp}}\,dv\,du\leq\frac{R^{-sp}}{\alpha_{d-1}}\int_U\int_U\frac{|\nabla A(u)-\nabla A(v)|^p}{|u-v|^{d-1+sp}}\,dv\,du=\frac{R^{-sp}}{\alpha_{d-1}}[\nabla A]_{W^{s,p}(U)}^p.
\label{stimaWspnablaA}
\end{split}
\end{equation}
The second inequality uses the fact that $B(u,R^{-1})\subset U_1$. Consequently, by Young's inequality
$$
\int_{U_0}\Big(R^{-1}+\fint_{B(u,R^{-1})}|\nabla A(u)-\nabla A(v)|\,dv\Big)^p\,du\leq2^{p-1}\Big(|U_0|R^{-p}+\alpha_{d-1}^{-1}R^{-sp}[\nabla A]_{W^{s,p}(U)}^p\Big).
$$
Further notice that since $R\geq2$ and $\|\nabla A\|_\infty\leq\frac{1}{10}$, for every $u\in U_0$ we also have the bound
\begin{equation}
    R^{-1}+\fint_{B(u,R^{-1})}|\nabla A(u)-\nabla A(v)|\,dv\leq R^{-1}+\fint_{B(u,R^{-1})}(|\nabla A(u)|+|\nabla A(v)|)\,dv\leq\frac12+2\|\nabla A\|_\infty\leq1.
    \label{upperboundcon1}
\end{equation}
Let $1\leq p<\infty$ and $1\leq r<\infty$. We distinguish the two cases. 
If $p\leq r$, the pointwise bound by one gives
\begin{equation*}
\begin{split}
&\Big\|R^{-1}+\fint_{B(\cdot,R^{-1})}|\nabla A(\cdot)-\nabla A(v)|\,dv\Big\|_{L^r(U_0)}^2\leq\Big(\int_{U_0}\Big(R^{-1}+\fint_{B(u,R^{-1})}|\nabla A(u)-\nabla A(v)|\,dv\Big)^p\,du\Big)^{\frac{2}{r}}\\
&\qquad\qquad\qquad\leq 2^{\frac{2(p-1)}{r}}\Big(\Leb^{d-1}(U_0)R^{-p}+\int_{U_0}\Big(\fint_{B(u,R^{-1})}|\nabla A(u)-\nabla A(v)|\,dv\Big)^p\,du\Big)^{\frac{2}{r}}.
\nonumber
\end{split}
\end{equation*}
Thanks to Jensen's inequality and \eqref{upperboundcon1}, we infer that 
\begin{equation}
    \begin{split}
        \Big\|R^{-1}+\fint_{B(\cdot,R^{-1})}|\nabla A(\cdot)-\nabla A(v)|\,dv\Big\|_{L^r(U_0)}^2\leq& 2^{\frac{2(p-1)}{r}}(\Leb^{d-1}(U_0)R^{-p}+\alpha_{d-1}^{-1}R^{-sp}[\nabla A]_{W^{s,p}(U)}^p)^{\frac{2}{r}}\\
        =&2^{\frac{2(p-1)}{r}}(\Leb^{d-1}(U_0)R^{-(1-s)p}+\alpha_{d-1}^{-1}[\nabla A]_{W^{s,p}(U)}^p)^{\frac{2}{r}}R^{-\frac{2sp}{r}}.
    \end{split}
\end{equation}
If on the other hand $p>r$, then H\"older's and Young's inequalities gives
\begin{equation*}
\begin{split}
&\Big\|R^{-1}+\fint_{B(\cdot,R^{-1})}|\nabla A(\cdot)-\nabla A(v)|\,dv\Big\|_{L^r(U_0)}^2\leq\Leb^{d-1}(U_0)^{\frac{2}{r}-\frac{2}{p}}\Big(\int_{U_0}\Big(R^{-1}+\fint_{B(u,R^{-1})}|\nabla A(u)-\nabla A(v)|\,dv\Big)^p\,du\Big)^{\frac{2}{p}}\\
&\qquad\qquad\qquad\leq2^{2-\frac{2}{p}}\Leb^{d-1}(U_0)^{\frac{2}{r}-\frac{2}{p}}\Big(\Leb^{d-1}(U_0)R^{-p}+\int_{U_0}\fint_{B(u,R^{-1})}|\nabla A(u)-\nabla A(v)|^pdvdu\Big)^{\frac{2}{p}}\\
&\qquad\qquad\qquad\qquad\qquad\leq2^{2-\frac{2}{p}}\Leb^{d-1}(U_0)^{\frac{2}{r}-\frac{2}{p}}\Big(\Leb^{d-1}(U_0)R^{-(1-s)p}+\alpha_{d-1}^{-1}[\nabla A]_{W^{s,p}(U)}^p\Big)^{\frac{2}{p}}R^{-2s}.
\nonumber
\end{split}
\end{equation*}
By construction, we know that $\Leb^{d-1}(U_0)\geq1$ and that $r\geq1$. Therefore, we have
\begin{equation}
\begin{split}
&\Big\|R^{-1}+\fint_{B(\cdot,R^{-1})}|\nabla A(\cdot)-\nabla A(v)|\,dv\Big\|_{L^r(U_0)}^2\leq(2\Leb^{d-1}(U_0))^{2-\frac{2}{p}}\Big(\Leb^{d-1}(U_0)+\frac{[\nabla A]_{W^{s,p}(U)}^p}{\alpha_{d-1}}\Big)^{\frac{2}{p}}\begin{cases}
R^{-\frac{2sp}{r}},&p\leq r,\\
R^{-2s},&p>r.
\end{cases}
\end{split}
\label{upper-increment-norms}
\end{equation}
Taking $r=\frac{4}{3}$ in \eqref{upper-increment-norms} and substituting into \eqref{upper-first-reduction}, for $p<\infty$ we obtain
\begin{equation}
\int_{\mathbb S^{d-1}}\vartheta(\omega)|J_R(\omega)|^2\,d\omega\leq C\Big(\Leb^{d-1}(U_0)+\frac{[\nabla A]_{W^{s,p}(U)}^p}{\alpha_{d-1}}\Big)^{\frac{2}{p}}R^{-\frac{d-1}{2}-s\min\{\frac{3p}{2},2\}},
\label{upper-first-graph}
\end{equation}
where $C$ is a constant depending only on $U_0$ and $p$. 
If on the other hand $p=\infty$, an elementary computation with H\"older continuity implies
$$
\Big\|R^{-1}+\fint_{B(\cdot,R^{-1})}|\nabla A(\cdot)-\nabla A(v)|\,dv\Big\|_{L^r(U_0)}^2\leq C\Leb^{d-1}(U_0)^{\frac{2}{r}}(1+[\nabla A]_{C^{0,s}(U)})^2R^{-2s},\qquad \text{for every }1\leq r<\infty,
$$
where $C$ is a constant depending only on $d$. In this case, taking $r=\frac{4}{3}$ and substituting into \eqref{upper-first-reduction} gives
\begin{equation}
    \int_{\mathbb S^{d-1}}\vartheta(\omega)|J_R(\omega)|^2\,d\omega\leq C \Leb^{d-1}(U_0)^{\frac{3}{2}}(1+[\nabla A]_{C^{0,s}(U)})^2R^{-\frac{d-1}{2}-2s}.
    \label{stimaconp=infty}
\end{equation}

From these formulas, we now prove that the conclusion of the proposition follows. 
Since the functions chosen above satisfy $\eta+\sum_{\vartheta\in \Theta}\vartheta=1$, we have
$$
\int_{\mathbb S^{d-1}}|J_R(\omega)|^2\,d\omega=\int_{\mathbb S^{d-1}}\eta(\omega)|J_R(\omega)|^2\,d\omega+\sum_{\vartheta\in \Theta}\int_{\mathbb S^{d-1}}\vartheta(\omega)|J_R(\omega)|^2\,d\omega.
$$
Thus, we infer that for every $1\leq p\leq \infty$ and every $s\in (0,1)$ and every $d\geq 2$ we have 
\begin{equation}
    \int_{\mathbb S^{d-1}}|J_R(\omega)|^2\,d\omega\overset{\eqref{stimaconeta1},\eqref{upper-first-graph},\eqref{stimaconp=infty}}{\leq} CR^{-d+1}+CR^{-\frac{d-1}{2}-s\min\{\frac{3p}{2},2\}},
\end{equation}
where $C$ depends only on the graph patch. Inserting this bound into \eqref{stimolatu}, which carries the factor $R^{-2}$ coming from the divergence theorem, and summing over the finitely many charts after the change of radius described above, we obtain
$$
E_\Omega(R)\leq C_\Omega R^{-d-1}+C_\Omega R^{-\frac{d+3}{2}-s\min\{\frac{3p}{2},2\}}.
$$
Here $C_\Omega$ depends only on $\Omega$ and its geometry. This proves \eqref{stimafinaleI}. 

\smallskip

\spacedlowsmallcaps{The further bounds for $d\geq3$.} For $p<\infty$, we shall prove that there is a constant $\newC\label{C:upper-graph}$, independent of $R$, such that
\begin{equation}
\begin{split}
&\int_{\mathbb S^{d-1}}\vartheta(\omega)|J_R(\omega)|^2\,d\omega\\
&\qquad\qquad\leq\oldC{C:upper-graph}\Big\|R^{-1}+\fint_{B(\cdot,R^{-1})}|\nabla A(\cdot)-\nabla A(v)|\,dv\Big\|_{L^2(U_0)}^2\begin{cases}
R^{-d+\frac{3}{2}},&sp>d-2,\\
R^{-d+\frac{3}{2}}(\log(2+R))^{\frac{p}{2}},&sp=d-2,\\
R^{-\frac{d-1+sp}{2}}(\log(2+R))^{\frac{1}{2}},&sp<d-2.
\end{cases}
\end{split}
\label{upper-refined-graph}
\end{equation}
The constant depends only on $d,s,p$, the cubes $U_0,U_1,U$, the number $c_0$, the functions $\rho,\vartheta$, and the norms $\|\psi\|_{C^1(U)}$ and $[\nabla A]_{W^{s,p}(U)}$. Let us first show why this estimate gives the further bounds in the proposition. Taking $r=2$ in \eqref{upper-increment-norms} gives
$$
\Big\|R^{-1}+\fint_{B(\cdot,R^{-1})}|\nabla A(\cdot)-\nabla A(v)|\,dv\Big\|_{L^2(U_0)}^2\leq(2\Leb^{d-1}(U_0))^{2-\frac{2}{p}}\Big(\Leb^{d-1}(U_0)+\frac{[\nabla A]_{W^{s,p}(U)}^p}{\alpha_{d-1}}\Big)^{\frac{2}{p}}R^{-s\min\{p,2\}}.
$$
Assuming \eqref{upper-refined-graph}, the above computation \eqref{stimaconeta1} and \eqref{stimolatu}, finally imply
\begin{equation}
E_\Omega(R)\leq C_\Omega R^{-d-1}+C_\Omega\begin{cases}
R^{-d-\frac{1}{2}-s\min\{p,2\}},&sp>d-2,\\
(\log(2+R))^{\frac{p}{2}}R^{-d-\frac{1}{2}-s\min\{p,2\}},&sp=d-2,\\
(\log(2+R))^{\frac{1}{2}}R^{-\frac{d+3+sp}{2}-s\min\{p,2\}},&sp<d-2.
\end{cases}
\label{upper-boundary-sum}
\end{equation}
The fixed changes of radius and the finitely many graph estimates are included in $C_\Omega$, which is independent of $R$. This concludes the proof of the proposition. We are therefore left to prove \eqref{upper-first-reduction} and \eqref{upper-refined-graph}, and to establish the Dini bound.

\smallskip

\spacedlowsmallcaps{Proof of }\eqref{upper-first-reduction}. We start from \eqref{upper-square} and apply the first estimate in \cref{spherical-kernel} with
$$
z:=X(u)-X(v),\qquad c(\omega):=R^{2d-2}\vartheta(\omega)a_R(u,\omega)\overline{a_R(v,\omega)}.
$$
Here $u,v$ are fixed. The quotients defining $a_R$ are smooth in $\omega$ because $q\geq c_0$, and $\supp\vartheta\Subset\{\omega_1>2c_0\}$. Thus $c$ extends smoothly by zero to the sphere. The conclusion is
\begin{equation}
    \Big|\int_{\mathbb S^{d-1}}e^{-iR\langle\omega,X(u)-X(v)\rangle}c(\omega)\,d\omega\Big|\leq C_d(1+R|X(u)-X(v)|)^{-\frac{d-1}{2}}\sum_{k=0}^{2d+5}\|\nabla_{\mathbb S}^kc\|_\infty.
    \label{stimaconcccc}
\end{equation}
The proof of \eqref{upper-first-reduction} reduces to showing that there is a constant $\newC\label{C:upper-kernel}$, depending only on $d,c_0,\rho,\vartheta,\|\psi\|_{C^1(U)}$, such that
\begin{equation}
\begin{split}
\sum_{k=0}^{2d+5}\|\nabla_{\mathbb S}^kc\|_\infty&\leq\oldC{C:upper-kernel}\Big(R^{-1}+\fint_{B(u,R^{-1})}|\nabla A(u)-\nabla A(y)|\,dy\Big)\Big(R^{-1}+\fint_{B(v,R^{-1})}|\nabla A(v)-\nabla A(y)|\,dy\Big).
\end{split}
\label{upper-kernel-derivatives}
\end{equation}
Indeed, substituting this bound \eqref{stimaconcccc} and using $|X(u)-X(v)|\geq|u-v|$, \eqref{upper-square} gives
\begin{equation*}
\begin{split}
\int_{\mathbb S^{d-1}}\vartheta(\omega)|J_R(\omega)|^2\,d\omega=\int_{U_0}\int_{U_0}\int_{\mathbb S^{d-1}}&e^{-iR\langle\omega,X(u)-X(v)\rangle}c(\omega)\,d\omega\,dv\,du\\
\leq\int_{U_0}\int_{U_0}\Big|\int_{\mathbb S^{d-1}}e^{-iR\langle\omega,X(u)-X(v)\rangle}c(\omega)\,d\omega\Big|\,dv\,du&\leq C_d\int_{U_0}\int_{U_0}\frac{\sum_{k=0}^{2d+5}\|\nabla_{\mathbb S}^kc\|_\infty}{(1+R|X(u)-X(v)|)^{\frac{d-1}{2}}}\,dv\,du.
\nonumber
\end{split}
\end{equation*}
Using \eqref{upper-kernel-derivatives} to estimate the derivatives of $c$ on the sphere gives
\begin{equation*}
\begin{split}
&\int_{\mathbb S^{d-1}}\vartheta(\omega)|J_R(\omega)|^2\,d\omega\\
&\qquad\qquad\leq C_d\oldC{C:upper-kernel}\int_{U_0}\int_{U_0}\frac{(R^{-1}+\fint_{B(u,R^{-1})}|\nabla A(u)-\nabla A(y)|\,dy)(R^{-1}+\fint_{B(v,R^{-1})}|\nabla A(v)-\nabla A(y)|\,dy)}{(1+R|X(u)-X(v)|)^{\frac{d-1}{2}}}\,dv\,du.
\nonumber
\end{split}
\end{equation*}
We now replace the distance between the graph points by $|u-v|$. Since $X(u)=(u,A(u))$, for every $u\neq v$, we have
$$
\frac{1}{(1+R|X(u)-X(v)|)^{\frac{d-1}{2}}}\leq\frac{1}{(R|X(u)-X(v)|)^{\frac{d-1}{2}}}\leq\frac{R^{-\frac{d-1}{2}}}{|u-v|^{\frac{d-1}{2}}}.
$$
This yields in particular
\begin{equation}
\begin{split}
&\int_{\mathbb S^{d-1}}\vartheta(\omega)|J_R(\omega)|^2\,d\omega\\
&\qquad\leq \frac{C_d\oldC{C:upper-kernel}}{R^{\frac{d-1}{2}}}\int_{U_0}\int_{U_0}\frac{(R^{-1}+\fint_{B(u,R^{-1})}|\nabla A(u)-\nabla A(y)|dy)(R^{-1}+\fint_{B(v,R^{-1})}|\nabla A(v)-\nabla A(y)|dy)}{|u-v|^{\frac{d-1}{2}}}
dvdu.
\label{stima1}
\end{split}
\end{equation}
Now we show that there is a constant $\newC\label{C:upper-potential}$, depending only on $d$, such that, for every nonnegative $f\in L^{\frac{4}{3}}(\mathbb R^{d-1})$, we have
\begin{equation}
\int_{\mathbb R^{d-1}}\int_{\mathbb R^{d-1}}\frac{f(u)f(v)}{|u-v|^{\frac{d-1}{2}}}\,dv\,du\leq\oldC{C:upper-potential}\|f\|_{L^{\frac{4}{3}}(\mathbb R^{d-1})}^2.
\label{upper-potential}
\end{equation}
Fixing $t>0$ and splitting $B(u,t)$ into annuli and using the Hardy-Littlewood maximal function $\mathcal Mf(u):=\sup_{r>0}\frac{1}{|B(u,r)|}\int_{B(u,r)}|f(y)|dy$, we have 
\begin{equation*}
\begin{split}
\int_{|u-v|\leq t}\frac{f(v)}{|u-v|^{\frac{d-1}{2}}}\,dv&\leq\sum_{j=0}^\infty(2^{-j-1}t)^{-\frac{d-1}{2}}\int_{B(u,2^{-j}t)}f(v)\,dv\leq\frac{2^{\frac{d-1}{2}}\alpha_{d-1}}{1-2^{-\frac{d-1}{2}}}t^{\frac{d-1}{2}}\mathcal Mf(u).
\nonumber
\end{split}
\end{equation*}
On the other hand, H\"older's inequality implies that
$$
\int_{|u-v|>t}\frac{f(v)}{|u-v|^{\frac{d-1}{2}}}\,dv\leq\|f\|_{L^{\frac{4}{3}}}\Big(\int_{|z|>t}|z|^{-2(d-1)}\,dz\Big)^{\frac{1}{4}}=\alpha_{d-1}^{\frac{1}{4}}t^{-\frac{d-1}{4}}\|f\|_{L^{\frac{4}{3}}}.
$$
For $f\neq0$, take $\tilde t=(\frac{\|f\|_{L^{\frac{4}{3}}}}{\mathcal Mf(u)})^{\frac{4}{3(d-1)}}$. Thus, we have
\begin{equation}
    \begin{split}
        &\int_{\mathbb R^{d-1}}\int_{\mathbb R^{d-1}}\frac{f(u)f(v)}{|u-v|^{\frac{d-1}{2}}}\,dv\,du=\int_{\mathbb R^{d-1}}f(u)\Big(\int_{|u-v|\leq \tilde t}\frac{f(v)}{|u-v|^{\frac{d-1}{2}}}\,dv+\int_{|u-v|> \tilde t}\frac{f(v)}{|u-v|^{\frac{d-1}{2}}}\,dv\Big)du\\
        &\qquad\qquad\qquad\leq\int_{\mathbb R^{d-1}}f(u)\Big(\frac{2^{\frac{d-1}{2}}\alpha_{d-1}}{1-2^{-\frac{d-1}{2}}}\tilde t^{\frac{d-1}{2}}\mathcal Mf(u)+\alpha_{d-1}^{\frac{1}{4}}\tilde t^{-\frac{d-1}{4}}\|f\|_{L^{\frac{4}{3}}}\Big) du\\
&=\Big(\frac{2^{\frac{d-1}{2}}\alpha_{d-1}}{1-2^{-\frac{d-1}{2}}}+\alpha_{d-1}^{\frac{1}{4}}\Big)\|f\|_{L^{\frac{4}{3}}}^{\frac{2}{3}}\int_{\mathbb R^{d-1}}f(u)(\mathcal Mf(u))^{\frac{1}{3}}\,du\leq\Big(\frac{2^{\frac{d-1}{2}}\alpha_{d-1}}{1-2^{-\frac{d-1}{2}}}+\alpha_{d-1}^{\frac{1}{4}}\Big)\|f\|_{L^{\frac{4}{3}}}^{\frac{5}{3}}\|\mathcal Mf\|_{L^{\frac{4}{3}}}^{\frac{1}{3}}.
\label{stigvjgvjhgvjghv}
    \end{split}
\end{equation}
Thus, by \eqref{stima1}, \eqref{stigvjgvjhgvjghv} and the Hardy-Littlewood maximal theorem \cite[Section 2.1]{grafakos}, we have
$$
\int_{\mathbb S^{d-1}}\vartheta(\omega)|J_R(\omega)|^2\,d\omega\leq CR^{-\frac{d-1}{2}}\Big\|R^{-1}+\fint_{B(\cdot,R^{-1})}|\nabla A(\cdot)-\nabla A(y)|\,dy\Big\|_{L^{\frac{4}{3}}(U_0)}^2,
$$
Thus \eqref{upper-kernel-derivatives} implies \eqref{upper-first-reduction} where $C$ is a constant depending only on $d,c_0,\rho,\vartheta,\|\psi\|_{C^1(U)}$.

We now prove \eqref{upper-kernel-derivatives}. Since $c(\omega)=R^{2d-2}\vartheta(\omega)a_R(u,\omega)\overline{a_R(v,\omega)}$, we have 
\begin{equation*}
\begin{split}
\|\nabla_{\mathbb S}^kc\|_\infty&\leq R^{2d-2}\sum_{\alpha+\beta+\gamma=k}\frac{k!}{\alpha!\beta!\gamma!}\|\nabla_{\mathbb S}^\alpha\vartheta\|_\infty\sup_{\omega_1>2c_0}|\nabla_{\mathbb S}^\beta a_R(u,\omega)|\sup_{\omega_1>2c_0}|\nabla_{\mathbb S}^\gamma a_R(v,\omega)|,\qquad \text{for every }0\leq k\leq2d+5.
\nonumber
\end{split}
\end{equation*}
Thus, to prove \eqref{upper-kernel-derivatives}, it suffices to show that there is a constant $\newC\label{C:upper-amplitude}$, depending only on $d,c_0,\rho,\|\psi\|_{C^1(U)}$, such that
\begin{equation}
\sup_{\omega_1>2c_0}|\nabla_{\mathbb S}^ja_R(u,\omega)|\leq\oldC{C:upper-amplitude}R^{1-d}\Big(R^{-1}+\fint_{B(u,R^{-1})}|\nabla A(u)-\nabla A(y)|\,dy\Big),\qquad \text{for every }0\leq j\leq2d+5.
\label{graph-amplitude-bound}
\end{equation}
Indeed, applying this estimate at $u$ and $v$ in the preceding product contributes $R^{1-d}$ from each factor, so that $R^{2d-2}R^{1-d}R^{1-d}=1$. Summing over $k$ then gives \eqref{upper-kernel-derivatives}, with
$$
\oldC{C:upper-kernel}:=\oldC{C:upper-amplitude}^2\sum_{k=0}^{2d+5}3^k\sum_{a=0}^k\|\nabla_{\mathbb S}^a\vartheta\|_\infty.
$$
We recall that
\begin{equation*}
\begin{split}
a_R(u,\omega)&=\int_{U_1}\Big(q(u,\omega)\rho(R(u-v))+i(\partial_1\rho)(R(u-v))\Big)\Big(\frac{\psi(u)\langle\omega,(-\nabla A(u),1)\rangle}{q(u,\omega)}-\frac{\psi(v)\langle\omega,(-\nabla A(v),1)\rangle}{q(v,\omega)}\Big)\,dv.
\nonumber
\end{split}
\end{equation*}
Fix $u\in U_0$. The integrand is smooth with respect to $\omega$ on $\{\omega_1>2c_0\}$, because its denominators are at least $c_0$. For each integer $j\geq0$, its angular derivatives of order $j$ are bounded by a constant depending only on $d,c_0,j,\rho,\|\psi\|_\infty$ times $\mathbb{1}_{B(u,R^{-1})}(v)$. Differentiation under the integral is therefore justified at every order, and $a_R(u,\cdot)$ is smooth on $\{\omega_1>2c_0\}$.

We claim that, for every $0\leq k\leq2d+5$, differentiating $k$ times with respect to $\omega$ under the integral yields
$\nabla_{\mathbb S}^ka_R(u,\omega)=\mathrm{I}+\mathrm{II}+\mathrm{III}+\mathrm{IV}$, where
\begin{equation*}
\begin{split}
\mathrm{I}&:=\int_{U_1}(\psi(u)-\psi(v))\Big(\rho(R(u-v))\nabla_{\mathbb S}^k\langle\omega,(-\nabla A(u),1)\rangle+i(\partial_1\rho)(R(u-v))\nabla_{\mathbb S}^k\frac{\langle\omega,(-\nabla A(u),1)\rangle}{q(u,\omega)}\Big)\,dv,
\nonumber
\end{split}
\end{equation*}
\begin{equation*}
\begin{split}
\mathrm{II}&:=-\int_{U_1}\psi(v)\Big(\rho(R(u-v))\nabla_{\mathbb S}^k\langle\pi_H\omega,\nabla A(u)-\nabla A(v)\rangle+i(\partial_1\rho)(R(u-v))\nabla_{\mathbb S}^k\frac{\langle\pi_H\omega,\nabla A(u)-\nabla A(v)\rangle}{q(u,\omega)}\Big)\,dv,
\nonumber
\end{split}
\end{equation*}
\begin{equation*}
\begin{split}
\mathrm{III}&:=-\int_{U_1}\psi(v)(\partial_1A(u)-\partial_1A(v))\Big(\rho(R(u-v))\nabla_{\mathbb S}^k\frac{\omega_d\langle\omega,(-\nabla A(v),1)\rangle}{q(v,\omega)}\Big)\,dv.
\nonumber
\end{split}
\end{equation*}
\begin{equation*}
\begin{split}
\mathrm{IV}&:=-\int_{U_1}\psi(v)(\partial_1A(u)-\partial_1A(v))\Big(i(\partial_1\rho)(R(u-v))\nabla_{\mathbb S}^k\frac{\omega_d\langle\omega,(-\nabla A(v),1)\rangle}{q(u,\omega)q(v,\omega)}\Big)\,dv.
\nonumber
\end{split}
\end{equation*}

Let us show that this is indeed the case. We compute the derivatives of $a_R$ with respect to $\omega$ and then collect the terms containing the same differences of $\psi$ and $\nabla A$. Since $\rho(R(u-v))$ and $(\partial_1\rho)(R(u-v))$ are independent of $\omega$, cancelling $q(u,\omega)$ in the first term gives
\begin{equation*}
\begin{split}
\nabla_{\mathbb S}^ka_R(u,\omega)&=\int_{U_1}\rho(R(u-v))\nabla_{\mathbb S}^k\Big(\psi(u)\langle\omega,(-\nabla A(u),1)\rangle-\psi(v)\frac{q(u,\omega)\langle\omega,(-\nabla A(v),1)\rangle}{q(v,\omega)}\Big)\,dv\\
&\qquad\qquad+i\int_{U_1}(\partial_1\rho)(R(u-v))\nabla_{\mathbb S}^k\Big(\frac{\psi(u)\langle\omega,(-\nabla A(u),1)\rangle}{q(u,\omega)}-\frac{\psi(v)\langle\omega,(-\nabla A(v),1)\rangle}{q(v,\omega)}\Big)\,dv.
\nonumber
\end{split}
\end{equation*}
We now rewrite the two expressions inside the derivatives. In the first one, writing $q(u,\omega)=q(v,\omega)+(q(u,\omega)-q(v,\omega))$ and using the definition of $q$, we have 
\begin{equation*}
\begin{split}
&\quad\psi(u)\langle\omega,(-\nabla A(u),1)\rangle-\psi(v)\frac{q(u,\omega)\langle\omega,(-\nabla A(v),1)\rangle}{q(v,\omega)}\\
&=(\psi(u)-\psi(v))\langle\omega,(-\nabla A(u),1)\rangle-\psi(v)\langle\pi_H\omega,\nabla A(u)-\nabla A(v)\rangle\\&\qquad-\psi(v)(\partial_1A(u)-\partial_1A(v))\frac{\omega_d\langle\omega,(-\nabla A(v),1)\rangle}{q(v,\omega)}.
\nonumber
\end{split}
\end{equation*}
The second integrand can be on the other hand rewritten as follows
\begin{equation*}
\begin{split}
&\quad\frac{\psi(u)\langle\omega,(-\nabla A(u),1)\rangle}{q(u,\omega)}-\frac{\psi(v)\langle\omega,(-\nabla A(v),1)\rangle}{q(v,\omega)}\\
&=\frac{\psi(u)\langle\omega,(-\nabla A(u),1)\rangle-\psi(v)\langle\omega,(-\nabla A(v),1)\rangle}{q(u,\omega)}+\psi(v)\langle\omega,(-\nabla A(v),1)\rangle\Big(\frac1{q(u,\omega)}-\frac1{q(v,\omega)}\Big)\\
&=(\psi(u)-\psi(v))\frac{\langle\omega,(-\nabla A(u),1)\rangle}{q(u,\omega)}-\psi(v)\frac{\langle\pi_H\omega,\nabla A(u)-\nabla A(v)\rangle}{q(u,\omega)}\\&\qquad-\psi(v)(\partial_1A(u)-\partial_1A(v))\frac{\omega_d\langle\omega,(-\nabla A(v),1)\rangle}{q(u,\omega)q(v,\omega)}.
\nonumber
\end{split}
\end{equation*}
Differentiating these two identities $k$ times with respect to $\omega$ and collecting the terms, we obtain the decomposition of $\nabla_{\mathbb S}^ka_R$ claimed above.

We first estimate $\mathrm{I}$. Its integrand vanishes unless $|u-v|<R^{-1}$, so $|\psi(u)-\psi(v)|\leq R^{-1}\|\nabla\psi\|_\infty$. Since $B(u,R^{-1})$ has volume $\alpha_{d-1}R^{1-d}$, these two factors give $R^{-d}$ once we estimate the derivatives in $\omega$ independently on $u$. 
For fixed $u$, the function $\omega\mapsto\langle\omega,(-\nabla A(u),1)\rangle$ and its first derivative have norm at most $2$, and all its higher derivatives vanish. Moreover, since $q(u,\omega)=\langle e_1+\partial_1A(u)e_d,\omega\rangle$, we have
$$
\Big|D_\omega^j\frac1{q(u,\omega)}\Big|=j!\frac{|e_1+\partial_1A(u)e_d|^j}{q(u,\omega)^{j+1}}\leq2^jj!c_0^{-j-1}.
$$
For $j\geq1$, the product rule therefore gives
\begin{equation*}
\begin{split}
\Big|D_\omega^j\frac{\langle\omega,(-\nabla A(u),1)\rangle}{q(u,\omega)}\Big|&\leq\sum_{\ell=0}^j\binom{j}{\ell}\Big|D_\omega^\ell\langle\omega,(-\nabla A(u),1)\rangle\Big|\Big|D_\omega^{j-\ell}\frac1{q(u,\omega)}\Big|\\
&=|\langle\omega,(-\nabla A(u),1)\rangle|\Big|D_\omega^j\frac1{q(u,\omega)}\Big|+j\Big|D_\omega\langle\omega,(-\nabla A(u),1)\rangle\Big|\Big|D_\omega^{j-1}\frac1{q(u,\omega)}\Big|\\
&\leq2\cdot2^jj!c_0^{-j-1}+2j\cdot2^{j-1}(j-1)!c_0^{-j}=2^jj!c_0^{-j-1}(2+c_0)\leq2^{j+2}j!c_0^{-j-1}.
\nonumber
\end{split}
\end{equation*}
Notice that since the Euclidean derivatives of order up to $k$ control the tangential derivatives of order $k$, there must be a constant $\newC\label{C:upper-angular}$, depending only on $d,c_0$, such that
\begin{equation*}
\begin{split}
|\mathrm{I}|&\leq\oldC{C:upper-angular}\int_{B(u,R^{-1})}|\psi(u)-\psi(v)|\Big(|\rho(R(u-v))|+|(\partial_1\rho)(R(u-v))|\Big)\,dv\\
&\leq\oldC{C:upper-angular}(\|\rho\|_\infty+\|\partial_1\rho\|_\infty)R^{-1}\|\nabla\psi\|_\infty\Leb^{d-1}(B(u,R^{-1}))=\alpha_{d-1}\oldC{C:upper-angular}(\|\rho\|_\infty+\|\partial_1\rho\|_\infty)\|\nabla\psi\|_\infty R^{-d}.
\nonumber
\end{split}
\end{equation*}

We next estimate $\mathrm{II}$ by bounding the derivatives with respect to $\omega$ in its integrand by a constant times $|\nabla A(u)-\nabla A(v)|$. For fixed $u,v$, the function $\omega\mapsto\langle\pi_H\omega,\nabla A(u)-\nabla A(v)\rangle$ is linear. Its value and first derivative have norm at most $|\nabla A(u)-\nabla A(v)|$, and all its higher derivatives vanish. Thus, for $j\geq1$, the product rule gives
\begin{equation*}
\begin{split}
&\Big|D_\omega^j\frac{\langle\pi_H\omega,\nabla A(u)-\nabla A(v)\rangle}{q(u,\omega)}\Big|\leq\sum_{\ell=0}^j\binom{j}{\ell}\Big|D_\omega^\ell\langle\pi_H\omega,\nabla A(u)-\nabla A(v)\rangle\Big|\Big|D_\omega^{j-\ell}\frac1{q(u,\omega)}\Big|\\
&=|\langle\pi_H\omega,\nabla A(u)-\nabla A(v)\rangle|\Big|D_\omega^j\frac1{q(u,\omega)}\Big|+j\Big|D_\omega\langle\pi_H\omega,\nabla A(u)-\nabla A(v)\rangle\Big|\Big|D_\omega^{j-1}\frac1{q(u,\omega)}\Big|\\
&\leq|\nabla A(u)-\nabla A(v)|\Big(2^jj!c_0^{-j-1}+j\,2^{j-1}(j-1)!c_0^{-j}\Big)\leq2^{j+1}j!c_0^{-j-1}|\nabla A(u)-\nabla A(v)|.
\nonumber
\end{split}
\end{equation*}
As above, we observe that these estimates control the tangential derivatives through order $2d+5$. Thus, we obtain 
$$
\begin{aligned}
|\mathrm{II}|&\leq\newC\label{C:upper-angularII}\|\psi\|_\infty\int_{B(u,R^{-1})}\Big(|\rho(R(u-v))|+|(\partial_1\rho)(R(u-v))|\Big)|\nabla A(u)-\nabla A(v)|\,dv\\
&\leq\oldC{C:upper-angularII}(\|\rho\|_\infty+\|\partial_1\rho\|_\infty)\|\psi\|_\infty\int_{B(u,R^{-1})}|\nabla A(u)-\nabla A(v)|\,dv.
\end{aligned}
$$

For $\mathrm{III}$, we use the estimate proved in $\mathrm{I}$ for $\frac{\langle\omega,(-\nabla A(v),1)\rangle}{q(v,\omega)}$ and differentiate its product with $\omega_d$. Since $D_\omega\omega_d=e_d$ and $D_\omega^\ell\omega_d=0$ for $\ell\geq2$, for $j\geq1$ we have
\begin{equation*}
\begin{split}
&\Big|D_\omega^j\frac{\omega_d\langle\omega,(-\nabla A(v),1)\rangle}{q(v,\omega)}\Big|\leq\sum_{\ell=0}^j\binom{j}{\ell}|D_\omega^\ell\omega_d|\Big|D_\omega^{j-\ell}\frac{\langle\omega,(-\nabla A(v),1)\rangle}{q(v,\omega)}\Big|\\
&\quad=|\omega_d|\Big|D_\omega^j\frac{\langle\omega,(-\nabla A(v),1)\rangle}{q(v,\omega)}\Big|+j|e_d|\Big|D_\omega^{j-1}\frac{\langle\omega,(-\nabla A(v),1)\rangle}{q(v,\omega)}\Big|\\
&\leq2^{j+2}j!c_0^{-j-1}+j\,2^{j+1}(j-1)!c_0^{-j}=2^{j+1}j!c_0^{-j-1}(2+c_0)\leq2^{j+3}j!c_0^{-j-1}.
\nonumber
\end{split}
\end{equation*}
Since $|\partial_1A(u)-\partial_1A(v)|\leq|\nabla A(u)-\nabla A(v)|$, we obtain
$$
\begin{aligned}
|\mathrm{III}|&\leq\newC\label{C:upper-angularIII}\|\rho\|_\infty\|\psi\|_\infty\int_{B(u,R^{-1})}|\nabla A(u)-\nabla A(v)|\,dv.
\end{aligned}
$$
Finally, in order to estimate $\mathrm{IV}$, we observe that 
$$
\frac{\omega_d\langle\omega,(-\nabla A(v),1)\rangle}{q(u,\omega)q(v,\omega)}=\frac{\omega_d\langle\omega,(-\nabla A(v),1)\rangle}{q(v,\omega)}\frac1{q(u,\omega)}.
$$
The derivatives of the first factor have just been estimated in $\mathrm{III}$, and those of the second were estimated in $\mathrm{I}$. Applying the product rule and substituting these two bounds gives, for every $j\geq0$,
\begin{equation*}
\begin{split}
&\Big|D_\omega^j\frac{\omega_d\langle\omega,(-\nabla A(v),1)\rangle}{q(u,\omega)q(v,\omega)}\Big|\leq\sum_{\ell=0}^j\binom{j}{\ell}\Big|D_\omega^\ell\frac{\omega_d\langle\omega,(-\nabla A(v),1)\rangle}{q(v,\omega)}\Big|\Big|D_\omega^{j-\ell}\frac1{q(u,\omega)}\Big|\\
&\quad\leq\sum_{\ell=0}^j\binom{j}{\ell}2^{\ell+3}\ell!c_0^{-\ell-1}2^{j-\ell}(j-\ell)!c_0^{-(j-\ell+1)}=\sum_{\ell=0}^j2^{j+3}j!c_0^{-j-2}=2^{j+3}(j+1)!c_0^{-j-2}.
\nonumber
\end{split}
\end{equation*}
As above, we therefore get 
$$
\begin{aligned}
|\mathrm{IV}|&\leq\newC\label{C:upper-angularIV}\|\partial_1\rho\|_\infty\|\psi\|_\infty\int_{B(u,R^{-1})}|\nabla A(u)-\nabla A(v)|\,dv.
\end{aligned}
$$
We now insert the four estimates into $\nabla_{\mathbb S}^ka_R=\mathrm{I}+\mathrm{II}+\mathrm{III}+\mathrm{IV}$. The estimates for $\mathrm{III}$ and $\mathrm{IV}$ together have the same coefficient $\|\rho\|_\infty+\|\partial_1\rho\|_\infty$ as the estimate for $\mathrm{II}$. Hence
$$
\begin{aligned}
\sup_{\omega_1>2c_0}|\nabla_{\mathbb S}^ka_R(u,\omega)|&\leq\newC\label{noncelafacciopiu}(\|\rho\|_\infty+\|\partial_1\rho\|_\infty)\Big(\alpha_{d-1}\|\nabla\psi\|_\infty R^{-d}+2\|\psi\|_\infty\int_{B(u,R^{-1})}|\nabla A(u)-\nabla A(v)|\,dv\Big),
\end{aligned}
$$
with $\oldC{noncelafacciopiu}:=(\oldC{C:upper-angular}+\oldC{C:upper-angularII}+\oldC{C:upper-angularIII}+\oldC{C:upper-angularIV})$
and this proves \eqref{graph-amplitude-bound}.

\smallskip

\spacedlowsmallcaps{Proof of }\eqref{upper-refined-graph}. Let $d\geq3$ and $p<\infty$. We first apply \cref{spherical-graph-estimate} with $\vartheta_0=\vartheta$ and an arbitrary $b\in L^2(U)$ supported in $U_0$. The cubes and $\varepsilon_0$ were fixed at the beginning. We have $\|\nabla A\|_\infty\leq\varepsilon_0$, $\nabla A\in W^{s,p}(U)$ and $q(u,\omega)\geq c_0$ on $U_1\times\supp\vartheta$. Thus, all its hypotheses hold. This implies that there is a constant $\newC\label{C:upper-scalar}$, depending only on $d,s,p,U_0,U_1,U,c_0,\vartheta,[\nabla A]_{W^{s,p}(U)}$ and independent of $R,b$, for which
\begin{equation}
\int_{\mathbb S^{d-1}}\vartheta(\omega)\Big|\int_Ue^{-iR\langle\omega,X(u)\rangle}b(u)\,du\Big|^2\,d\omega\leq L^2\|b\|_{L^2(U)}^2,\qquad L^2:=\oldC{C:upper-scalar}\begin{cases}
R^{-d+\frac{3}{2}},&sp>d-2,\\
R^{-d+\frac{3}{2}}(\log(2+R))^{\frac{p}{2}},&sp=d-2,\\
R^{-\frac{d-1+sp}{2}}(\log(2+R))^{\frac{1}{2}},&sp<d-2.
\end{cases}
\label{primadifourier}
\end{equation}

We now write the integral involving $J_R$ in the form required by \cref{L2-fourier-expansion}. Since $\omega_1>2c_0$ on $\supp\vartheta$, we can use $w=(\omega_2,\ldots,\omega_d)$ as the variable of integration. Then $\omega(w)=(\sqrt{1-|w|^2},w)$ and
\begin{equation*}
\begin{split}
\int_{\mathbb S^{d-1}}\vartheta(\omega)|J_R(\omega)|^2\,d\omega
&\overset{\eqref{regularized-graph-integral}}{=}\int_{|w|<\sqrt{1-4c_0^2}}\frac{\vartheta(\omega(w))}{\omega_1(w)}\Big|\int_{U_0}e^{-iR\langle\omega(w),X(u)\rangle}R^{d-1}a_R(u,\omega(w))\,du\Big|^2\,dw.
\nonumber
\end{split}
\end{equation*}
To apply \cref{L2-fourier-expansion}, we need a coefficient defined on $U_1$, vanishing outside $U_0$, whose derivatives with respect to $w$ belong to $L^2$. We therefore choose $\zeta\in C^\infty(\mathbb S^{d-1})$, equal to one on $\supp\vartheta$ and compactly supported in $\{\omega_1>2c_0\}$, and replace $a_R(u,\omega(w))$ by $\zeta(\omega(w))a_R(u,\omega(w))$. This leaves the displayed integral unchanged, since its weight vanishes wherever $\vartheta(\omega(w))=0$.

We can now identify the two functions to which the lemma will be applied. For $|w|<\sqrt{1-4c_0^2}$, set
\begin{equation*}
\begin{split}
K(w,u)&:=\Big(\frac{\vartheta(\omega(w))}{\omega_1(w)}\Big)^{\frac{1}{2}}e^{-iR\langle\omega(w),X(u)\rangle},\qquad\text{and}\qquad B(u,w):=R^{d-1}\zeta(\omega(w))a_R(u,\omega(w)),
\nonumber
\end{split}
\end{equation*}
and extend both functions by zero elsewhere. The support of $\zeta\circ\omega$ is compactly contained in this ball, so the extension of $B(u,\cdot)$ is smooth. With these notations we can thus rewrite
$$
\int_{\mathbb S^{d-1}}\vartheta(\omega)|J_R(\omega)|^2\,d\omega=\Big\|\int_{U_0}K(w,u)B(u,w)\,du\Big\|_{L^2(U_0,dw)}^2.
$$
Thus, in \cref{L2-fourier-expansion}, we take $\mathcal U=U_0$, $\mathcal R_0=U_0$, $\mathcal R_1=U_1$, $n=q=d-1$ and $N=d$. 
We verify the hypotheses concerning $K$. Since $\omega_1(w)>2c_0$ wherever $K(w,u)\neq0$, we have
$$
\mathop{\mathrm{ess\,sup}}_{w\in U_0}\|K(w,\cdot)\|_{L^2(U_0)}^2\leq\frac{\Leb^{d-1}(U_0)\|\vartheta\|_\infty}{2c_0}.
$$
For every $b\in L^2(U_0)$, extended by zero to $U$, the same change of variables gives
\begin{equation*}
\begin{split}
\Big\|\int_{U_0}K(\cdot,u)b(u)\,du\Big\|_{L^2(U_0)}^2&=\int_{\mathbb S^{d-1}}\vartheta(\omega)\Big|\int_Ue^{-iR\langle\omega,X(u)\rangle}b(u)\,du\Big|^2\,d\omega\overset{\eqref{primadifourier}}{\leq}L^2\|b\|_{L^2(U_0)}^2.
\nonumber
\end{split}
\end{equation*}
Thus, once we verify that the derivatives of $B$ with respect to $w$ through order $d$ belong to $L^2(U_0\times U_1)$, \cref{L2-fourier-expansion} gives
$$
\int_{\mathbb S^{d-1}}\vartheta(\omega)|J_R(\omega)|^2\,d\omega\leq C_{U_0,U_1,d-1,d}^2L^2\sum_{|\alpha|\leq d}\|\partial_w^\alpha B\|_{L^2(U_0\times U_1)}^2.
$$
We now estimate the sum on the right using \eqref{graph-amplitude-bound}.  Where $\omega_1(w)>2c_0$, the product rule and the chain rule give
\begin{equation*}
\begin{split}
\partial_{w_\ell}B(u,w)&=R^{d-1}\partial_{w_\ell}(\zeta\circ\omega)(w)a_R(u,\omega(w))+R^{d-1}\zeta(\omega(w))\nabla_{\mathbb S}a_R(u,\omega(w))[\partial_{w_\ell}\omega(w)],\\
&\qquad\qquad\qquad\partial_{w_\ell}\omega(w)=e_{\ell+1}-\frac{w_\ell}{\omega_1(w)}e_1.
\nonumber
\end{split}
\end{equation*}
Inductively, each term contains an angular derivative of $a_R$ of order at most $d$, multiplied by derivatives of $\zeta\circ\omega$ and $\omega$. Since $\omega_1(w)>2c_0$ on $\supp(\zeta\circ\omega)$, the derivatives of $w\mapsto(\sqrt{1-|w|^2},w)$ through order $d$ are bounded there in terms of $d,c_0$. Hence there is a constant $\newC\label{C:upper-parameters}$, depending only on $d,c_0,\zeta$, such that
\begin{equation*}
\begin{split}
\sum_{|\alpha|\leq d}|\partial_w^\alpha B(u,w)|^2&\leq\oldC{C:upper-parameters}R^{2d-2}\sum_{k=0}^d\sup_{\omega_1>2c_0}|\nabla_{\mathbb S}^ka_R(u,\omega)|^2\\
&\overset{\eqref{graph-amplitude-bound}}{\leq}\oldC{C:upper-parameters}R^{2d-2}\sum_{k=0}^d\oldC{C:upper-amplitude}^2R^{2-2d}\Big(R^{-1}+\fint_{B(u,R^{-1})}|\nabla A(u)-\nabla A(v)|\,dv\Big)^2\\
&=(d+1)\oldC{C:upper-parameters}\oldC{C:upper-amplitude}^2\Big(R^{-1}+\fint_{B(u,R^{-1})}|\nabla A(u)-\nabla A(v)|\,dv\Big)^2.
\nonumber
\end{split}
\end{equation*}
Outside the support of $\zeta\circ\omega$, all derivatives of $B(u,\cdot)$ vanish, so this estimate holds on $U_1$. Integrating in $u\in U_0$ and $w\in U_1$, we obtain
\begin{equation}
\sum_{|\alpha|\leq d}\|\partial_w^\alpha B\|_{L^2(U_0\times U_1)}^2\leq(d+1)\Leb^{d-1}(U_1)\oldC{C:upper-parameters}\oldC{C:upper-amplitude}^2\Big\|R^{-1}+\fint_{B(\cdot,R^{-1})}|\nabla A(\cdot)-\nabla A(v)|\,dv\Big\|_{L^2(U_0)}^2.
\label{upper-parameter-norms}
\end{equation}
The right-hand side is finite because $\nabla A$ is bounded and $U_0$ has finite measure. The function $B(u,\cdot)$ is smooth, and \eqref{upper-parameter-norms} shows that its derivatives with respect to $w$ through order $d$ belong to $L^2(U_0\times U_1)$. All the hypotheses of \cref{L2-fourier-expansion} are therefore satisfied.

We can therefore apply the lemma. Substituting \eqref{upper-parameter-norms} into its conclusion gives
\begin{equation*}
\begin{split}
\int_{\mathbb S^{d-1}}\vartheta(\omega)|J_R(\omega)|^2\,d\omega&\leq C_{U_0,U_1,d-1,d}^2(d+1)\Leb^{d-1}(U_1)\oldC{C:upper-parameters}\oldC{C:upper-amplitude}^2L^2\Big\|R^{-1}+\fint_{B(\cdot,R^{-1})}|\nabla A(\cdot)-\nabla A(v)|\,dv\Big\|_{L^2(U_0)}^2.
\nonumber
\end{split}
\end{equation*}
This concludes the proof of \eqref{upper-refined-graph}.

\smallskip

\spacedlowsmallcaps{The Dini bound.} Suppose first that $d\geq3$. We apply the Dini estimate in \cref{spherical-graph-estimate} to obtain the factor $R^{-d+\frac{3}{2}}$ in place of the three factors in \eqref{primadifourier}. We first verify its additional hypothesis, namely the Dini continuity of the graph gradients.
For each graph chosen at the beginning, the formula for the normal gives
$$
\nabla A_j(u)=-\frac{\pi_HQ_j^T\nu_\Omega(X_j(u))}{\langle Q_j^T\nu_\Omega(X_j(u)),e_d\rangle},\qquad \langle Q_j^T\nu_\Omega(X_j(u)),e_d\rangle=\frac1{\sqrt{1+|\nabla A_j(u)|^2}}\geq\frac12.
$$
Subtracting the quotients at $u,v$ and using $|\nabla A_j(v)|\leq\frac{1}{10}$ and a few algebraic computations, we obtain
\begin{equation*}
\begin{split}
|\nabla A_j(u)-\nabla A_j(v)|\leq 2|\nu_\Omega(X_j(u))-\nu_\Omega(X_j(v))|.
\nonumber
\end{split}
\end{equation*}
Since $|X_j(u)-X_j(v)|\leq2r_j|u-v|$, the preceding estimate gives
$\omega_{\nabla A_j}(t)\leq2\omega_{\nu_\Omega}(2r_jt)$.
Consequently, using $\omega_{\nu_\Omega}\leq2$ to bound the integral over $\tau>1$, we have
\begin{equation*}
\begin{split}
\int_0^1\omega_{\nabla A_j}(t)\,\frac{dt}{t}\leq2\int_0^{2r_j}\omega_{\nu_\Omega}(\tau)\,\frac{d\tau}{\tau}\leq2\int_0^1\omega_{\nu_\Omega}(\tau)\,\frac{d\tau}{\tau}+4\log(\max\{1,2r_j\})<\infty.
\nonumber
\end{split}
\end{equation*}

 We apply the Dini estimate in \cref{spherical-graph-estimate} with $\vartheta_0=\vartheta$ and $b\in L^2(U)$ supported in $U_0$. Besides the Dini condition just verified, its hypotheses are the bounds $\|\nabla A\|_\infty\leq\varepsilon_0$ and $q(u,\omega)\geq c_0$ on $U_1\times\supp\vartheta$, already checked above. There is therefore a constant $\newC\label{C:upper-dini-scalar}$, depending only on $d,U_0,U_1,U,c_0,\vartheta$ and the Dini modulus of $\nabla A$, such that
$$
\int_{\mathbb S^{d-1}}\vartheta(\omega)\Big|\int_Ue^{-iR\langle\omega,X(u)\rangle}b(u)\,du\Big|^2\,d\omega\leq\oldC{C:upper-dini-scalar}R^{-d+\frac{3}{2}}\|b\|_{L^2(U)}^2.
$$
We now use the Dini estimate in place of \eqref{primadifourier} in the application of \cref{L2-fourier-expansion}. With the same functions $K,B$ and the other hypotheses already verified above, the lemma gives
\begin{equation*}
\begin{split}
\int_{\mathbb S^{d-1}}\vartheta(\omega)|J_R(\omega)|^2\,d\omega&=\Big\|\int_{U_0}K(w,u)B(u,w)\,du\Big\|_{L^2(U_0,dw)}^2\leq C_{U_0,U_1,d-1,d}^2\oldC{C:upper-dini-scalar}R^{-d+\frac{3}{2}}\sum_{|\alpha|\leq d}\|\partial_w^\alpha B\|_{L^2(U_0\times U_1)}^2.
\nonumber
\end{split}
\end{equation*}
The right-hand side can be estimated by \eqref{upper-parameter-norms}, which gives
\begin{equation*}
\begin{split}
\int_{\mathbb S^{d-1}}\vartheta(\omega)|J_R(\omega)|^2\,d\omega&\leq C_{U_0,U_1,d-1,d}^2(d+1)\Leb^{d-1}(U_1)\oldC{C:upper-parameters}\oldC{C:upper-amplitude}^2\oldC{C:upper-dini-scalar}R^{-d+\frac{3}{2}}\Big\|R^{-1}+\fint_{B(\cdot,R^{-1})}|\nabla A(\cdot)-\nabla A(v)|\,dv\Big\|_{L^2(U_0)}^2.
\nonumber
\end{split}
\end{equation*}
This implies in particular that 
\begin{equation}
\begin{split}
\int_{\mathbb S^{d-1}}\vartheta(\omega)|J_R(\omega)|^2\,d\omega&\leq C_{U_0,U_1,d-1,d}^2\oldC{C:upper-dini-scalar}R^{-d+\frac{3}{2}}\sum_{|\alpha|\leq d}\|\partial_w^\alpha B\|_{L^2(U_0\times U_1)}^2\\
&\overset{\eqref{upper-parameter-norms}}{\leq}\oldC{C:upper-dini-graphsoeing}R^{-d+\frac{3}{2}}\Big\|R^{-1}+\fint_{B(\cdot,R^{-1})}|\nabla A(\cdot)-\nabla A(v)|\,dv\Big\|_{L^2(U_0)}^2.
\end{split}
\label{upper-dini-graph}
\end{equation}
The constant $\newC\label{C:upper-dini-graphsoeing}$ depends only on $d,U_0,U_1,U,c_0,\rho,\vartheta,\|\psi\|_{C^1(U)}$ and the Dini modulus of $\nabla A$, and is independent of $R$.
Hence, for $p<\infty$, we substitute \eqref{upper-increment-norms} with $r=2$ into \eqref{upper-dini-graph}, to get
\begin{equation*}
\begin{split}
\int_{\mathbb S^{d-1}}\vartheta(\omega)|J_R(\omega)|^2\,d\omega&\leq\oldC{C:upper-dini-graphsoeing}(2\Leb^{d-1}(U_0))^{2-\frac{2}{p}}\Big(\Leb^{d-1}(U_0)+\frac{[\nabla A]_{W^{s,p}(U)}^p}{\alpha_{d-1}}\Big)^{\frac{2}{p}} R^{-d+\frac{3}{2}-s\min\{p,2\}}.
\nonumber
\end{split}
\end{equation*}
For $p=\infty$, H\"older continuity instead gives
$$R^{-1}+\fint_{B(u,R^{-1})}|\nabla A(u)-\nabla A(v)|\,dv\leq R^{-1}+[\nabla A]_{C^{0,s}(U)}R^{-s}\leq(1+[\nabla A]_{C^{0,s}(U)})R^{-s}.$$
In particular, $\omega_{\nabla A}(t)\leq[\nabla A]_{C^{0,s}(U)}t^s$, so the Dini condition needed above holds automatically. Substitution into \eqref{upper-dini-graph} yields
$$\int_{\mathbb S^{d-1}}\vartheta(\omega)|J_R(\omega)|^2\,d\omega\leq\oldC{C:upper-dini-graphsoeing}\Leb^{d-1}(U_0)(1+[\nabla A]_{C^{0,s}(U)})^2R^{-d+\frac{3}{2}-2s}.$$
These are the desired graph bounds with exponent $-d+\frac{3}{2}-\gamma$, since $\gamma=s\min\{p,2\}$ for finite $p$ and $\gamma=2s$ for $p=\infty$.

For $d=2$, \eqref{upper-first-graph} and \eqref{stimaconp=infty} already give the graph bound with exponent $-\frac{1}{2}-\gamma=-d+\frac{3}{2}-\gamma$. Indeed, in this dimension $\gamma=s\min\{\frac{3p}{2},2\}$ for finite $p$ and $\gamma=2s$ for $p=\infty$.

This concludes the proof of the proposition.
\end{proof}

\subsection{An estimate in \texorpdfstring{$d-2$}{d-2} variables}

In this subsection we prove the technical core of this paper, an $L^2$ estimate for an oscillatory integral operator in $d-2$ variables whose phase is the distance between points of two graphs. The proof uses the special form of this phase in an essential way, so the result should not be read as an estimate for arbitrary nondegenerate phases of the same regularity.

\begin{proposizione}\label{estimate-in-d-minus-two-variables}
Let $d\geq3$, set $I=(-\frac{1}{2},\frac{3}{2})^{d-2}$, and fix $\chi_0,\chi_1\in C_c^\infty(I)$.
For every $C_0>0$ there is $\varepsilon_*=\varepsilon_*(d,C_0)>0$ with the following property. Let $v\in\mathbb R^{d-2}$ and $P,Q\in C^1(\mathbb R^{d-2})$ satisfy
\begin{equation}
    \|\nabla P\|_{L^\infty(\mathbb R^{d-2})}+\|\nabla Q\|_{L^\infty(\mathbb R^{d-2})}\leq\varepsilon_*,\qquad \text{and}\qquad|v|+\sup_{x,y\in I}|P(x)-Q(y)|\leq C_0.
    \label{stimenecessarioes}
\end{equation}
For every $\lambda\geq1$, we let
$$V_\lambda f(x):=\chi_0(x)\int_Ie^{i\lambda\Psi(x,y)}\chi_1(y)f(y) dy,$$
where
$$\Psi(x,y):=\sqrt{1+|v+x-y|^2+(P(x)-Q(y))^2},$$
and we let $\omega_P(t):=\sup_{|x-y|\leq t}|\nabla P(x)-\nabla P(y)|$ and $\omega_Q(t):=\sup_{|x-y|\leq t}|\nabla Q(x)-\nabla Q(y)|$ be the moduli of continuity of $\nabla P$ and $\nabla Q$. Then the following assertions hold.
\begin{itemize}
\item[(i)] If $\int_0^1\omega_P(t)\frac{dt}{t}<\infty$ or $\int_0^1\omega_Q(t)\frac{dt}{t}<\infty$, then, for every $\lambda\geq1$, we have
$$\|V_\lambda\|_{L^2(I)\to L^2(I)}^2\leq C\lambda^{-(d-2)}.$$
\item[(ii)] Let $0<s<1$ and $1\leq p<\infty$. If $\nabla P,\nabla Q\in W^{s,p}(\mathbb R^{d-2})$, then, for every $\lambda\geq1$, the norm $\|V_\lambda\|_{L^2(I)\to L^2(I)}^2$ is at most
$$\begin{cases}C\lambda^{-(d-2)}(1+[\nabla P]_{W^{s,p}}^p+[\nabla Q]_{W^{s,p}}^p),&sp>d-2,\\
C\lambda^{-(d-2)}(1+(\log(2+\lambda))^{p-1}([\nabla P]_{W^{s,p}}^p+[\nabla Q]_{W^{s,p}}^p)),&sp=d-2,\\
C\lambda^{-(d-2)}+C\lambda^{-sp}([\nabla P]_{W^{s,p}}^p+[\nabla Q]_{W^{s,p}}^p),&sp<d-2.\end{cases}$$
\end{itemize}
In both assertions, $C$ is independent of $\lambda$. In (i), it depends only on $d,C_0,\chi_0,\chi_1,\omega_P$ and $\omega_Q$. In (ii), it depends only on $d,C_0,\chi_0,\chi_1,s$ and $p$, and the gradients of $P$ and $Q$ enter the bound only through the seminorms displayed above.
\end{proposizione}

The proof of \cref{estimate-in-d-minus-two-variables} occupies the rest of this subsection. We first show that the mixed derivative of $\Psi$ is uniformly invertible, and we use this to bring the phase into a normal form on small cubes. We then estimate a model operator whose phase has this form, and this is where the finite differences enter. Finally, we verify that the phase has the required form and choose the parameters according to the regularity of $P$ and $Q$.

\subsubsection{Uniform invertibility of the mixed derivative}\label{sss:mixed}

Differentiating $\Psi$ gives 
\begin{equation}
\begin{split}
    \partial_x\partial_y\Psi&=-\frac{\mathrm{Id}_{\mathbb R^{d-2}}+\nabla P(x)\otimes\nabla Q(y)}{\Psi}\\
&\qquad\qquad+\frac{1}{\Psi^3}\Big((v+x-y+(P(x)-Q(y))\nabla P(x))\otimes
(v+x-y+(P(x)-Q(y))\nabla Q(y))\Big).
\end{split}
\label{mixed-derivative}
\end{equation}
Since $x,y\in I$, one has $|x-y|\leq2\sqrt{d-2}$. The assumptions on $v$, $P$ and $Q$ therefore give
$$1\leq\Psi(x,y)\leq A:=\sqrt{1+(C_0+2\sqrt{d-2})^2+C_0^2}.$$
Let
$$M_0(x,y):=-\frac{\mathrm{Id}_{\mathbb R^{d-2}}}{\Psi}+\frac{(v+x-y)\otimes(v+x-y)}{\Psi^3}$$
be the matrix obtained from the right-hand side of \eqref{mixed-derivative} by removing the terms containing $\nabla P(x)$ and $\nabla Q(y)$. The eigenvalues of $M_0(x,y)$ have absolute values $\Psi^{-1}$ and $\frac{1+(P(x)-Q(y))^2}{\Psi^3}$,
and both are at least $A^{-3}$. Consequently,
$$|M_0(x,y)\xi|\geq A^{-3}|\xi|\qquad\text{for every }\xi\in\mathbb R^{d-2}.$$

We now compare $\partial_x\partial_y\Psi$ with $M_0(x,y)$. Expanding the tensor product in \eqref{mixed-derivative} gives
\begin{equation}
    \begin{split}
        \partial_x\partial_y\Psi-M_0=-&\frac{\nabla P(x)\otimes\nabla Q(y)}{\Psi}+\frac{(P(x)-Q(y))^2}{\Psi^3}\nabla P(x)\otimes\nabla Q(y)\\
        +&\frac{P(x)-Q(y)}{\Psi^3}\Big(\nabla P(x)\otimes(v+x-y)+(v+x-y)\otimes\nabla Q(y)\Big).
    \end{split}
\end{equation}
Since the operator norm of $a\otimes b$ is $|a||b|$, if $\|\nabla P\|_\infty+\|\nabla Q\|_\infty\leq1$, then
$$\|\partial_x\partial_y\Psi-M_0\|\leq C_1(|\nabla P(x)|+|\nabla Q(y)|),$$
where $C_1=1+C_0(C_0+2\sqrt{d-2})+C_0^2$.
Choose
$0<\varepsilon_*\leq\min\{1,(2A^3C_1)^{-1}\}$ and observe that if $\|\nabla P\|_\infty+\|\nabla Q\|_\infty\leq\varepsilon_*$, then
$$|\partial_x\partial_y\Psi \xi|\geq|M_0(x,y)\xi|-\|\partial_x\partial_y\Psi-M_0\||\xi|\geq\frac{A^{-3}}2|\xi|$$
for every $\xi\in\mathbb R^{d-2}$. Thus $\partial_x\partial_y\Psi$ is invertible and
$$\|(\partial_x\partial_y\Psi)^{-1}\|\leq2A^3.$$
The formula in \eqref{mixed-derivative} also gives $\|D_xD_y\Psi\|\leq C(d,C_0)$. Fix $C_{\mathrm{lin}}\geq1$, depending only on $d$ and $C_0$, such that
$$
\|D_xD_y\Psi(x,y)\|+\|(D_xD_y\Psi(x,y))^{-1}\|\leq C_{\mathrm{lin}}
$$
for all $x,y\in I$ and all $v,P,Q$ satisfying the assumptions of the proposition.

\subsubsection{A local normal form for the phase}\label{sss:normal}

The finite-difference argument below needs an exact representation of the mixed part of the phase as a sum of products of a function of $x$ and a function of $y$.

Let $L>1$ and $0<\delta\leq1$ satisfy
$$L>C_{\mathrm{lin}}\sqrt{d-2}+1\qquad\text{and}\qquad2\delta\sqrt{d-2}(1+C_{\mathrm{lin}}L)<\min(\dist(\supp\chi_0,\partial I),\dist(\supp\chi_1,\partial I)).$$
Choose finite families of open cubes $\{C_{0,a}\}_a$ and $\{C_{1,b}\}_b$ of side $\delta$, with sides parallel to the coordinate axes, compactly contained in $I$, which cover $\supp\chi_0$ and $\supp\chi_1$, respectively, and such that every cube meets the corresponding support. Choose smooth decompositions subordinate to these families,
$$\chi_0=\sum_a\chi_{0,a},\qquad\chi_1=\sum_b\chi_{1,b},$$
so that $\supp\chi_{0,a}\subset C_{0,a}$ and $\supp\chi_{1,b}\subset C_{1,b}$. Thus $V_\lambda$ is the finite sum, over $a$ and $b$, of the operators
$$V_\lambda^{a,b}f(x):=\chi_{0,a}(x)\int_Ie^{i\lambda\Psi(x,y)}\chi_{1,b}(y)f(y)dy.$$
Fix $a$ and $b$, and denote the centres of $C_{0,a}$ and $C_{1,b}$ by $x_0$ and $y_0$. Define
$$P_\delta(x):=\frac{P(x_0+\delta x)-P(x_0)}{\delta},\qquad Q_\delta(y):=\frac{Q(y_0+\delta y)-Q(y_0)}{\delta},$$
and let
$$z_0:=(v+x_0-y_0,P(x_0)-Q(y_0)),\quad U(x):=(x,P_\delta(x)),\quad W(y):=(y,Q_\delta(y)),\quad H(z):=(1+|z|^2)^{\frac{1}{2}}.$$
In these coordinates, the phase becomes
$$\Psi_\delta(x,y):=\Psi(x_0+\delta x,y_0+\delta y)=H(z_0+\delta(U(x)-W(y))).$$
Notice that for every $\xi,\eta\in\mathbb R^{d-2}$, the chain rule gives
\begin{equation}
D_xD_y\Psi(x_0,y_0)[\xi,\eta]=-D^2H(z_0)[(\xi,\langle\nabla P(x_0),\xi\rangle),(\eta,\langle\nabla Q(y_0),\eta\rangle)].
\label{id:doppiaderivata}\end{equation}
Indeed, at $(x_0,y_0)$, the derivative of $(v+x-y,P(x)-Q(y))$ with respect to $x$ in the direction $\xi$ is $(\xi,\langle\nabla P(x_0),\xi\rangle)$, while its derivative with respect to $y$ in the direction $\eta$ is $-(\eta,\langle\nabla Q(y_0),\eta\rangle)$.

Now, let us define the normalized mixed part of the rescaled phase by
$$\Phi_\delta(x,y):=\delta^{-2}(\Psi_\delta(x,y)-\Psi_\delta(x,0)-\Psi_\delta(0,y)+\Psi_\delta(0,0)).$$
By definition,
\begin{equation}
\Psi_\delta(x,y)=\Psi_\delta(x,0)+\Psi_\delta(0,y)-\Psi_\delta(0,0)+\delta^2\Phi_\delta(x,y).
\label{splittingfigo}\end{equation}
We now compare $V_\lambda^{a,b}$ with the operator having phase $\Phi_\delta$ and the rescaled cutoffs:
$$\widetilde V_\lambda^{a,b}f(x):=\chi_{0,a}(x_0+\delta x)\int_Ie^{i\lambda\delta^2\Phi_\delta(x,y)}\chi_{1,b}(y_0+\delta y)f(y)dy.$$
Both rescaled cutoffs are supported in $(-\frac{1}{2},\frac{1}{2})^{d-2}\subset I$. Evaluating the original operator at $X=x_0+\delta x$ and $Y=y_0+\delta y$ and changing variables in its integral, we obtain
\begin{equation}\begin{split}
(V_\lambda^{a,b}f)(x_0+\delta x)&=\chi_{0,a}(x_0+\delta x)\int_Ie^{i\lambda\Psi(x_0+\delta x,Y)}\chi_{1,b}(Y)f(Y)dY\\
&=\delta^{d-2}\chi_{0,a}(x_0+\delta x)\int_{\frac{I-y_0}{\delta}}e^{i\lambda\Psi_\delta(x,y)}\chi_{1,b}(y_0+\delta y)f(y_0+\delta y)dy\\
&=\delta^{d-2}e^{i\lambda\Psi_\delta(x,0)}\widetilde V_\lambda^{a,b}(e^{i\lambda(\Psi_\delta(0,\cdot)-\Psi_\delta(0,0))}f(y_0+\delta\,\cdot))(x).
\end{split}\end{equation}
where the last equality follows from \eqref{splittingfigo}. The changes of variables in the input and output norms give
\begin{equation}\begin{split}
\|f(y_0+\delta\,\cdot)\|_2^2=\delta^{-(d-2)}\|f\|_2^2,\qquad\text{and}\qquad
\|(V_\lambda^{a,b}f)(x_0+\delta\,\cdot)\|_2^2=\delta^{-(d-2)}\|V_\lambda^{a,b}f\|_2^2.
\end{split}\end{equation}
Taking the supremum over the $f\in L^2$ with $\|f\|_2\leq 1$ we have
\begin{equation}
\|V_\lambda^{a,b}\|_{2\to2}=\delta^{\frac{d-2}{2}}\delta^{d-2}\delta^{-\frac{d-2}{2}}\|\widetilde V_\lambda^{a,b}\|_{2\to2}=\delta^{d-2}\|\widetilde V_\lambda^{a,b}\|_{2\to2}.
\label{identitadiV}
\end{equation}
In Section~\ref{sss:representation} we will put $\Phi_\delta$ in the form required in Section~\ref{sss:model}. For this purpose, since $U(0)=W(0)=0$, we observe that the following identity holds
\begin{equation}
\Phi_\delta(x,y)=\delta^{-2}\Big(H(z_0+\delta(U(x)-W(y)))-H(z_0+\delta U(x))-H(z_0-\delta W(y))+H(z_0)\Big).
\label{strutturaphidelta}\end{equation}

\subsubsection{The norm of a model operator}\label{sss:model}

Let us introduce some real-valued function $\phi$ that has the following properties
\begin{equation}
\phi(x,w)=\langle x,w\rangle+\sum_\mu f_\mu(x)g_\mu(w),\qquad \text{with}\qquad \sum_\mu\Lip(f_\mu)\Lip(g_\mu)\leq\varepsilon,
\label{target-local-phase}
\end{equation}
where $\varepsilon>0$ is a sufficiently small fixed number. 
We use fixed cutoffs $\eta_0,\eta_1\in C_c^\infty(\mathbb R^{d-2})$.
 We will prove in Section~\ref{sss:representation} that the function $\Phi_\delta$ can be written in this form after a linear change of the integration variable and thus all the estimates that we obtain here will be used to estimate the operator norm of $\widetilde V_\lambda^{a,b}$. Let us define the operator 
$$T_\lambda b(x):=\eta_0(x)\int e^{i\lambda\phi(x,w)}\eta_1(w)b(w)dw.$$
Finite differences of $|\eta_1|^2$ are sums of translates of the same function. We therefore define
\begin{equation}
T_{\lambda,z}b(x):=\eta_0(x)\int e^{i\lambda\phi(x,w)}\eta_1(w-z)b(w)dw,\qquad N_\lambda:=\sup_{z\in\mathbb R^{d-2}}\|T_{\lambda,z}\|_{2\to2}^2.
\label{translated-operator-family}
\end{equation}
The original operator is $T_{\lambda,0}=T_\lambda$.
Cauchy-Schwarz gives
$N_\lambda\leq\|\eta_0\|_{L^2}^2\|\eta_1\|_{L^2}^2<\infty$.
In the following, we fix $z$. Every estimate will be uniform in $z$, and the norms of the cutoff and of its derivatives are unaffected by this translation.
If we let  $b,u\in C_c^\infty(\mathbb R^{d-2})$, Fubini's theorem implies that 
$$\int T_{\lambda,z} b(x)\overline{u(x)}\,dx=\int b(w)\overline{\overline{\eta_1(w-z)}\int e^{-i\lambda\phi(x',w)}\overline{\eta_0(x')}u(x')\,dx'}\,dw.$$
Therefore, the adjoint operator gives
$$T_{\lambda,z}^*u(w)=\overline{\eta_1(w-z)}\int e^{-i\lambda\phi(x',w)}\overline{\eta_0(x')}u(x')\,dx'.$$
Substituting this formula into $T_{\lambda,z} T_{\lambda,z}^*u$ and applying Fubini once more, we obtain
\begin{equation}
T_{\lambda,z} T_{\lambda,z}^*u(x)=\int K_{\lambda,z}(x,x')u(x')dx' \quad\text{where}\quad K_{\lambda,z}(x,x'):=\eta_0(x)\overline{\eta_0(x')}\int e^{i\lambda(\phi(x,w)-\phi(x',w))}|\eta_1(w-z)|^2\,dw.
\label{kernel-of-local-square}
\end{equation}
Since $\|T_{\lambda,z} T_{\lambda,z}^*\|_{2\to2}\leq N_\lambda$, and $z$ is arbitrary,
it is enough to estimate the integral operator with kernel $K_{\lambda,z}$ 
uniformly in $z$.
This computation will take place dividing, as usual, the scales. This means that we will cut the kernel $K_{\lambda,z}$ into two parts. One will be the kernel close to the origin and the other will be the kernel at infinity.

Fix $M\geq1$ and choose, depending only on $d$, a radial, nonincreasing function $\rho\in C_c^\infty(\mathbb R^{d-2})$, equal to one on $B(0,1)$ and zero outside $B(0,2)$. Set $r_j:=\frac{2^jM}{\lambda}$ for $j\geq0$. We first decompose $T_{\lambda,z} T_{\lambda,z}^*$ into the operators
$$\mathcal K_{\lambda,z}^{\mathrm{near}}F(x):=\int\rho\Big(\frac{x-x'}{r_0}\Big)K_{\lambda,z}(x,x')F(x')\,dx',$$
$$\mathcal K_{j,z}F(x):=\int\Big(\rho\Big(\frac{x-x'}{2r_j}\Big)-\rho\Big(\frac{x-x'}{r_j}\Big)\Big)K_{\lambda,z}(x,x')F(x')\,dx'.$$
Since $r_{j+1}=2r_j$, telescoping gives
$$\rho\Big(\frac{x-x'}{r_0}\Big)+\sum_{j=0}^{J}\Big(\rho\Big(\frac{x-x'}{2r_j}\Big)-\rho\Big(\frac{x-x'}{r_j}\Big)\Big)=\rho\Big(\frac{x-x'}{r_{J+1}}\Big).$$
For sufficiently large $J$, the right-hand side equals one whenever $x,x'\in\supp\eta_0$. Therefore
\begin{equation}
T_{\lambda,z} T_{\lambda,z}^*=\mathcal K_{\lambda,z}^{\mathrm{near}}+\sum_{j\geq0}\mathcal K_{j,z}.
\label{initial-operator-decomposition}
\end{equation}
For fixed $\lambda$, only finitely many $\mathcal K_{j,z}$ are nonzero, since their kernels vanish unless
$r_j\leq|x-x'|\leq4r_j$ and $ x,x'\in\supp\eta_0$.
In particular,
\begin{equation}
\|T_{\lambda,z}\|_{2\to2}^2\leq\|\mathcal K_{\lambda,z}^{\mathrm{near}}\|_{2\to2}+\sum_{j\geq0}\|\mathcal K_{j,z}\|_{2\to2}.
\label{initial-operator-norm-decomposition}
\end{equation}
We estimate the first term directly. For each remaining term, we first apply finite differences to its kernel and identify the resulting operators before estimating their norms.

We first bound the norm of a variant of $T_{\lambda,z}T_{\lambda,z}^*$ whose kernel carries an additional coefficient with a separated expansion, in the sense of \cref{separated-expansions}. This variant will be used to estimate the norms of the operators $\mathcal K_{j,z}$.

Let $u_\nu,v_\nu,b_\nu\in L^\infty(\mathbb R^{d-2})$, for $\nu\geq1$, assume that
$$\sum_{\nu=1}^{\infty}\|u_\nu\|_{L^\infty}\|v_\nu\|_{L^\infty}\|b_\nu\|_{L^\infty}<\infty\qquad\text{set}\qquad c(x,x',w):=\sum_{\nu=1}^{\infty}u_\nu(x)v_\nu(x')b_\nu(w)\in L^\infty(\R^{3(d-2)}).$$
and define
\begin{equation}
\mathcal R_{z}F(x):=\eta_0(x)\int\overline{\eta_0(x')}F(x')\int e^{i\lambda(\phi(x,w)-\phi(x',w))}|\eta_1(w-z)|^2c(x,x',w)\,dw\,dx'.
    \label{operatorprototype}
\end{equation}

For $m\geq1$, let $\mathcal R_{m,z}$ be the operator obtained by replacing $c$ with its first $m$ summands. Linearity of the integrals and the formula for $T_{\lambda,z}^*$ give
\begin{equation}
    \mathcal R_{m,z}F=\sum_{\nu=1}^{m}u_\nu T_{\lambda,z}\big(b_\nu T_{\lambda,z}^*(v_\nu F)\big).
    \label{idnetitaRcoefficentec}
\end{equation}
Therefore, for every $F\in L^2(\mathbb R^{d-2})$, the triangle inequality gives
\begin{equation}
\begin{split}
\|\mathcal R_{m,z}F\|_{L^2}&\leq N_\lambda\Big(\sum_{\nu=1}^{m}\|u_\nu\|_{L^\infty}\|v_\nu\|_{L^\infty}\|b_\nu\|_{L^\infty}\Big)\|F\|_{L^2},
\end{split}
\label{stimaoperatoreperora}
\end{equation}
where we used $\|T_{\lambda,z}^*\|_{2\to2}=\|T_{\lambda,z}\|_{2\to2}$ and $\|T_{\lambda,z}\|_{2\to2}^2\leq N_\lambda$.
To pass to the infinite sum, observe that the tail of the coefficient is bounded in absolute value by
$$\sum_{\nu>m}\|u_\nu\|_{L^\infty}\|v_\nu\|_{L^\infty}\|b_\nu\|_{L^\infty}.$$
Taking absolute values in the kernel of $\mathcal R_{z}-\mathcal R_{m,z}$ and applying Cauchy-Schwarz consequently gives
$$\lim_{m\to \infty}\|\mathcal R_{z}-\mathcal R_{m,z}\|_{2\to2}\leq \lim_{m\to \infty}\|\eta_0\|_{L^2}^2\|\eta_1\|_{L^2}^2\sum_{\nu>m}\|u_\nu\|_{L^\infty}\|v_\nu\|_{L^\infty}\|b_\nu\|_{L^\infty}=0.$$
This implies in particular that 
\begin{equation}
\|\mathcal R_{z}\|_{2\to2}\leq N_\lambda\sum_{\nu=1}^{\infty}\|u_\nu\|_{L^\infty}\|v_\nu\|_{L^\infty}\|b_\nu\|_{L^\infty}.
\label{separated-coefficient-operator-bound}
\end{equation}

We proceed with the computation of the operator norms of $\mathcal K_{\lambda,z}^{near}$ and of $\mathcal K_{j,z}$. 

\smallskip

\spacedlowsmallcaps{The contribution of $K_{\lambda,z}$ close to the origin.}
Let us estimate the operator norm of $\mathcal K_{\lambda,z}^{near}$. 
Since the exponential has absolute value one, \eqref{kernel-of-local-square} gives
$$|K_{\lambda,z}^{\mathrm{near}}(x,x')|\leq\mathbb{1}_{\{|x-x'|\leq \frac{2M}{\lambda}\}}|\eta_0(x)||\eta_0(x')|\int|\eta_1(w-z)|^2\,dw\leq\mathbb{1}_{\{|x-x'|\leq \frac{2M}{\lambda}\}}\|\eta_0\|_{L^\infty}^2\|\eta_1\|_{L^2}^2.$$
Taking the supremum over $x$ gives
\begin{equation}
\sup_x\int|K_{\lambda,z}^{\mathrm{near}}(x,x')|dx'\leq\|\eta_0\|_{L^\infty}^2\|\eta_1\|_{L^2}^2\Leb^{d-2}(B(0,1))(2M)^{d-2}\lambda^{-(d-2)}.
\label{near-row-bound}
\end{equation}
The same calculation, with $x$ and $x'$ interchanged, gives
\begin{equation}
\sup_{x'}\int|K_{\lambda,z}^{\mathrm{near}}(x,x')|dx\leq\|\eta_0\|_{L^\infty}^2\|\eta_1\|_{L^2}^2\Leb^{d-2}(B(0,1))(2M)^{d-2}\lambda^{-(d-2)}.
\label{near-column-bound}
\end{equation}
For every fixed $x$, we have
$$\Big|\int K_{\lambda,z}^{\mathrm{near}}(x,x')u(x')dx'\Big|^2\leq\int|K_{\lambda,z}^{\mathrm{near}}(x,x')|dx'\int|K_{\lambda,z}^{\mathrm{near}}(x,x')||u(x')|^2dx'.$$
Integrating this inequality with respect to $x$ and using \eqref{near-row-bound}, we obtain
\begin{equation*}
\begin{split}
    \int  \Big|\int K_{\lambda,z}^{\mathrm{near}}(x,x')u(x')dx'\Big|^2 dx\leq \frac{\|\eta_0\|_{L^\infty}^2\|\eta_1\|_{L^2}^2\Leb^{d-2}(B(0,1))(2M)^{d-2}}{\lambda^{d-2}}\int \int|K_{\lambda,z}^{\mathrm{near}}(x,x')||u(x')|^2dx'dx. 
    \nonumber
\end{split}
\end{equation*}
Applying Fubini's theorem, we infer that 
$$\|\mathcal K_{\lambda,z}^{near}u\|_{L^2}^2=\int  \Big|\int K_{\lambda,z}^{\mathrm{near}}(x,x')u(x')dx'\Big|^2 dx\leq \Big(\frac{\|\eta_0\|_{L^\infty}^2\|\eta_1\|_{L^2}^2\Leb^{d-2}(B(0,1))(2M)^{d-2}}{\lambda^{d-2}}\Big)^2\lVert u\rVert_{L^2(\R^{d-2})}^2.$$

\smallskip

\spacedlowsmallcaps{The contribution of $K_{\lambda,z}$ at infinity.}
It remains to estimate the smooth annular contributions $\mathcal K_{j,z}$, whose kernels vanish unless $r_j\leq|x-x'|\leq4r_j$. Recall that $r_j=\frac{2^jM}{\lambda}$, and set $\kappa:=\frac{\pi}{128}$ and $h_j:=\frac{\kappa}{\lambda r_j}$ for every $j\geq0$.

\spacedlowsmallcaps{Splitting the operator in annuli.} 
For each $j$, partition $\mathbb R^{d-2}$ into half-open cubes of side $\frac{r_j}{64\sqrt{d-2}}$, and keep only the cubes that meet $\supp\eta_0$. For each such cube $Q$, consider all such cubes $Q'$ for which there exist $x\in Q$ and $x'\in Q'$ satisfying $r_j\leq|x-x'|\leq4r_j$. Observe that the number of such cubes is bounded by some constant $C_d$ depending only on $d$. 
For each ordered pair of such cubes, with centers $c,c'$, choose a signed coordinate vector $e\in \mathscr E^\pm:=\mathscr E\cup(-\mathscr E)$, where $\mathscr E=\{e_1,\ldots,e_{d-2}\}$ is the standard basis of $\mathbb R^{d-2}$, such that there holds
\begin{equation}
   \langle c-c',e\rangle\geq\frac{|c-c'|}{\sqrt{d-2}}.
    \label{stimacentric}
\end{equation}

For each $e\in\mathscr E^\pm$, let $\mathcal P_{j,e}\label{familyPje}$ be the set of ordered pairs $(Q,Q')$ considered above for which the chosen vector is $e$. Each pair belongs to exactly one of these sets. Define
$$\rho_{j,e}(x,x'):=\Big(\rho\Big(\frac{x-x'}{2r_j}\Big)-\rho\Big(\frac{x-x'}{r_j}\Big)\Big)\sum_{(Q,Q')\in\mathcal P_{j,e}}\mathbb 1_Q(x)\mathbb 1_{Q'}(x').$$
Since the cubes are disjoint, the sum is either zero or one, and it is one exactly when $x$ and $x'$ lie in two of the cubes that form a pair in $\mathcal P_{j,e}$. Moreover, if $x,x'\in\supp\eta_0$ and the difference of the two cutoffs is nonzero, the cubes containing $x$ and $x'$ form a pair which belongs to exactly one of the sets $\mathcal P_{j,e}$. Therefore, for all $x,x'\in\supp\eta_0$,
\begin{equation}
\sum_{e\in\mathscr E^\pm}\rho_{j,e}(x,x')=\rho\Big(\frac{x-x'}{2r_j}\Big)-\rho\Big(\frac{x-x'}{r_j}\Big).
\label{id:split:cubi:culo}
\end{equation}
For each $e\in\mathscr E^\pm$, define
$$\mathcal K_{e,j,z}F(x):=\int\rho_{j,e}(x,x')K_{\lambda,z}(x,x')F(x')\,dx',$$
and notice that \eqref{id:split:cubi:culo}, together with the fact that $K_{\lambda,z}(x,x')$ vanishes unless $x,x'\in\supp\eta_0$, shows that
$\mathcal K_{j,z}=\sum_{e\in\mathscr E^\pm}\mathcal K_{e,j,z}$ and so 
$$\|T_{\lambda,z}\|_{2\to2}^2\leq \|\mathcal K_{\lambda,z}^{near}\|_{2\to 2}+\sum_{j\geq0}\sum_{e \in \mathscr E^\pm}\lVert \mathcal K_{e,j,z}\rVert_{2\to 2}.$$

\spacedlowsmallcaps{The finite difference trick.} In this section we perform the main technical step of the proof. This is how we substitute derivatives in the classical smooth fold estimates with differences, in this rough context. 

Set $a_z(w):=|\eta_1(w-z)|^2$. In the following we will decompose the operator $\mathcal K_{e,j,z}$ as follows
\begin{equation}
\begin{split}
&\mathcal K_{e,j,z}F(x)\\
&=\sum_{(Q,Q')\in\mathcal P_{j,e}}\underbrace{\begin{aligned}[t]&\mathbb1_Q(x)\eta_0(x)\int_{Q'}\overline{\eta_0(x')}\rho_{j,e}(x,x')F(x')\Big[\\&\quad\int e^{i\lambda(\phi(x,w)-\phi(x',w))}\delta_{h_j,e}^{d-1}a_{z}(w)\prod_{k=1}^{d-1}q_j(x,x',w-kh_je)\,dw\Big]dx'\end{aligned}}_{\mathcal K_{e,j,z}^{Q,Q',d-1}F(x)}\\
&+\sum_{(Q,Q')\in\mathcal P_{j,e}}\sum_{\ell=0}^{d-2}\underbrace{\begin{aligned}[t]&\mathbb1_Q(x)\eta_0(x)\int_{Q'}\overline{\eta_0(x')}\rho_{j,e}(x,x')F(x')\Big[\\&\quad\int e^{i\lambda(\phi(x,w)-\phi(x',w))}\delta_{h_j,e}^{\ell}a_{z}(w)\delta_{h_j,e}\Big(\prod_{k=0}^{\ell}q_j(x,x',w-kh_je)\Big)\,dw\Big]dx'.\end{aligned}}_{{\mathcal K_{e,j,z}^{Q,Q',\ell}F(x)}}
\label{splitKe,j}
\end{split}
\end{equation}
where
$$q_j(w):=\frac{1}{e^{i\lambda\Delta_{h_j,e}\phi(x,x',w)}-1},\qquad\delta_{h_j,e}b(w):=b(w-h_je)-b(w),$$
and where we use the convention $\delta_{h_j,e}^ma_{z}=\delta_{h_j,e}(\delta_{h_j,e}^{m-1}a_{z})$, where we set the convention $\delta_{h_j,e}^0a_{z}=a_{z}$. For bookkeeping's sake, we recall that if we are able to prove the above decomposition, then the operator $T_{\lambda,z}T^*_{\lambda,z}$ can be rewritten as
$$T_{\lambda,z}T_{\lambda,z}^* =\mathcal K_{\lambda,z}^{near}+\sum_{j\geq 0}\sum_{e\in \mathscr E^\pm} \sum_{\ell=0}^{d-1}\sum_{(Q,Q')\in \mathcal P_{j,e}}\mathcal K_{e,j,z}^{Q,Q',\ell}.$$

We now explain what estimates we need for the operators in \eqref{splitKe,j}, and how these estimates will control their sum. The representations and coefficient bounds stated below will be established in the following paragraphs. 
Recall that, after \eqref{translated-operator-family}, we fixed $z$. The constants and the bounds below must be independent of the parameter $z$, so that taking the supremum over $z$ gives the required estimate for $N_\lambda$.

We first sharpen the summation over the cube pairs. Each cube occurs in at most $C_d$ pairs in $\mathcal P_{j,e}$ in either position. Moreover, $\mathcal K_{e,j,z}^{Q,Q',\ell}F$ is supported in $Q$ and depends only on $\mathbb1_{Q'}F$. Consequently,
\begin{equation}
\begin{split}
\Big\|\sum_{(Q,Q')\in\mathcal P_{j,e}}\mathcal K_{e,j,z}^{Q,Q',\ell}F\Big\|_{L^2}^2
&=\sum_Q\Big\|\sum_{Q':(Q,Q')\in\mathcal P_{j,e}}\mathcal K_{e,j,z}^{Q,Q',\ell}F\Big\|_{L^2}^2\\
\leq C_d\sum_{(Q,Q')\in\mathcal P_{j,e}}\|\mathcal K_{e,j,z}^{Q,Q',\ell}\|_{2\to2}^2\|\mathbb1_{Q'}F\|_{L^2}^2&\leq C_d^2\sup_{(Q,Q')\in\mathcal P_{j,e}}\|\mathcal K_{e,j,z}^{Q,Q',\ell}\|_{2\to2}^2\|F\|_{L^2}^2,
\end{split}
\label{cube-pair-summation}
\end{equation}
where in the first identity we  use the disjointness of the cubes $Q$, the first inequality is Cauchy-Schwarz, and the last inequality uses the bound on the number of pairs $(Q,Q')\in \mathcal P_{j,e}$ where $Q'$ is fixed. Thus
\begin{equation*}
    \begin{split}
        \|T_{\lambda,z}\|_{2\to2}^2\leq &\|\mathcal K_{\lambda,z}^{near}\|_{2\to 2}+C_d^2\sum_{j\geq 0}\sum_{e\in \mathscr E^\pm} \sum_{\ell=0}^{d-1}\sup_{(Q,Q')\in\mathcal P_{j,e}}\|\mathcal K_{e,j,z}^{Q,Q',\ell}\|_{2\to2},
        \nonumber
    \end{split}
\end{equation*}
When $\ell=d-1$, we shall prove the direct estimate
\begin{equation}
\|\mathcal K_{e,j,z}^{Q,Q',d-1}\|_{2\to2}\leq \frac{C_d}{2^jM}\|\eta_0\|_{L^\infty}^2\|\partial_e^{d-1}(|\eta_1|^2)\|_{L^1}\lambda^{-(d-2)}.
\label{required-principal-estimate}
\end{equation}
Its sum over $j$ is therefore bounded by $\frac{C}{M}\|\eta_0\|_{L^\infty}^2\|\partial_e^{d-1}(|\eta_1|^2)\|_{L^1}\lambda^{-(d-2)}$. In particular 
$$\|T_{\lambda,z}\|_{2\to2}^2
\leq C\|\eta_0\|_{L^\infty}^2\Big(\|\eta_1\|_{L^2}^2M^{d-2}+\frac{\sum_{e\in\mathscr E^\pm}\|\partial_e^{d-1}(|\eta_1|^2)\|_{L^1}}{M}\Big)\lambda^{-(d-2)}+C\sum_{j\geq 0}\sum_{e\in \mathscr E^\pm} \sum_{\ell=0}^{d-2}\sup_{(Q,Q')\in\mathcal P_{j,e}}\|\mathcal K_{e,j,z}^{Q,Q',\ell}\|_{2\to2}.$$
For each remaining term with $0\leq\ell\leq d-2$, which are the terms that dictate, at the end of the day, the rate of average decay of the Fourier transform at infinity, we shall construct bounded functions $u_\nu,v_\nu,b_\nu$
and, for each term, a vector $z_\nu\in\{z+rh_je:r=0,\ldots,\ell\}$. These vectors arise from the differences of the amplitude:
$$\delta_{h_j,e}^{\ell}a_z(w)=\sum_{r=0}^{\ell}(-1)^{\ell-r}\binom{\ell}{r}|\eta_1(w-(z+rh_je))|^2.$$
The term indexed by $r$ contains the squared cutoff appearing in the kernel of $T_{\lambda,z+rh_je}T_{\lambda,z+rh_je}^*$. This will be carefully explained below. When we expand its remaining coefficient into products of functions of $x$, $x'$ and $w$, every resulting operator term therefore uses $T_{\lambda,z+rh_je}$ and its adjoint. In the series below, $z_\nu$ denotes this vector for the $\nu$th term, and the binomial coefficients are included in $b_\nu$. We will prove that there are sequences $\{A_j\}_{j\in\N}$ and $\{B_j\}_{j\in\N}$ such that 
\begin{equation}
\sum_{\nu\in \Xi}\|u_\nu\|_{L^\infty}\|v_\nu\|_{L^\infty}\|b_\nu\|_{L^\infty}\leq A_j
\label{required-coefficient-estimate}
\end{equation}
and, for every $F\in L^2$, there holds
\begin{equation}
\|\mathcal K_{e,j,z}^{Q,Q',\ell}F-\sum_{\nu\in \Xi} u_\nu T_{\lambda,z_\nu}(b_\nu T_{\lambda,z_\nu}^*(v_\nu F))\|_{L^2}\leq B_j\|F\|_{L^2}.
\label{required-remainder-representation}
\end{equation}
Here $u_\nu$ and $v_\nu$ are supported in $Q$ and $Q'$, respectively. The functions and translations may depend on $j,e,\ell,Q,Q'$ and on the original cutoff translation, but these dependences are omitted from the notation. The positive numbers $A_j$ and $B_j$ may depend on $\lambda,M$ and the functions in the phase, but must be uniform in the remaining indices and in the original cutoff translation.
By \eqref{separated-coefficient-operator-bound}, applied with the corresponding cutoff translation, we infer that
\begin{equation}
    \Big\|u_\nu T_{\lambda,z_\nu}\Big(b_\nu T_{\lambda,z_\nu}^*(v_\nu F)\Big)\Big\|_{L^2}\leq N_\lambda\|u_\nu\|_{L^\infty}\|v_\nu\|_{L^\infty}\|b_\nu\|_{L^\infty}\|F\|_{L^2}.
    \label{rappresentazionedeiKQQ'etc}
\end{equation}
Thus, the summability in \eqref{required-coefficient-estimate} therefore gives convergence of the series in operator norm. Summing this estimate and using \eqref{required-remainder-representation}, we obtain
\begin{equation}
\|\mathcal K_{e,j,z}^{Q,Q',\ell}\|_{2\to2}\leq A_jN_\lambda+B_j,\qquad 0\leq\ell\leq d-2.
\label{consequence-remainder-representation}
\end{equation}
Combining \eqref{cube-pair-summation}, \eqref{required-principal-estimate}, \eqref{consequence-remainder-representation}, and then taking the supremum over the original cutoff translation, we have
\begin{equation}
N_\lambda\leq C\Big(M^{d-2}+\frac1M\Big)\lambda^{-(d-2)}+CN_\lambda\sum_{j\geq0}A_j+C\sum_{j\geq0}B_j,
\label{operator-norm-before-absorption}
\end{equation}
where $C$ depends only on $d$ and the fixed cutoffs. If the sum of the $A_j$s is sufficiently small, say 
\begin{equation}
    C\sum_{j\geq0}A_j\leq\frac12,
    \label{sommaA_j}
\end{equation}
we reabsorb the $A_j$ sum on the left obtaining
\begin{equation}
N_\lambda\leq2C\Big(M^{d-2}+\frac1M\Big)\lambda^{-(d-2)}+2C\sum_{j\geq0}B_j.
\label{operator-norm-after-absorption}
\end{equation}
The remaining tasks are therefore to establish the decomposition \eqref{splitKe,j}, prove \eqref{required-principal-estimate}, and construct the representations \eqref{required-remainder-representation} with a sufficiently small sum of the $A_j$ and a bound for $\sum_jB_j$ supplied by the stated regularity assumptions.

\smallskip

\spacedlowsmallcaps{Proof of the decomposition in }\eqref{splitKe,j}.
Now, fix $e\in\mathscr E^\pm$. We next rewrite the integral in $w$ in the kernel of $\mathcal K_{e,j,z}$ by finite differences. Pick two cubes $Q$ and $Q'$ with $(Q,Q')\in\mathcal P_{j,e}$, so that \eqref{stimacentric} holds for this signed direction. 
The next identities are applied on this pair, with this direction fixed.
For every $h>0$, define
$$\Delta_{h,e}\phi(x,x',w):=\phi(x,w+he)-\phi(x',w+he)-\phi(x,w)+\phi(x',w).$$
Let us define $\theta_0:=\kappa r_j^{-1}\langle c-c',e\rangle$.
The side-length of the cubes and their distance imply that 
\begin{equation}
    \frac{3\kappa}{4\sqrt{d-2}}\leq\theta_0\leq5\kappa=\frac{5\pi}{128}<\frac{\pi}{2}.
    \label{stimatheta0}
\end{equation}
Using the properties of $\phi$ stated in \eqref{target-local-phase}, we obtain
\begin{equation}\begin{split}
&\qquad\qquad\qquad\qquad\qquad\lambda\Delta_{h_j,e}\phi(x,x',w)-\theta_0=\frac{\kappa}{r_j}(\langle x-c,e\rangle-\langle x'-c',e\rangle)\\
&\quad+\kappa\sum_\mu\Big(\frac{f_\mu(x)-f_\mu(c)}{r_j}+\frac{f_\mu(c)-f_\mu(c')}{r_j}-\frac{f_\mu(x')-f_\mu(c')}{r_j}\Big)\frac{g_\mu(w+h_je)-g_\mu(w)}{h_j}.
\end{split}\label{separated-translated-phase-difference}\end{equation}
We now put $\lambda\Delta_{h_j,e}\phi(x,x',w)-\theta_0$ in the form \eqref{target-local-phase}. Let us define $\Xi:=\{1,2\}\sqcup\{(\mu,a):a=1,2,3\}$, where $\mu$ ranges over the indices in \eqref{target-local-phase}, and define
$$u_1(x):=\frac{\kappa}{r_j}\mathbb1_Q(x)\langle x-c,e\rangle,\qquad v_1(x'):=\mathbb1_{Q'}(x'),\qquad b_1(w):=1,$$
$$u_2(x):=\mathbb1_Q(x),\qquad v_2(x'):=-\frac{\kappa}{r_j}\mathbb1_{Q'}(x')\langle x'-c',e\rangle,\qquad b_2(w):=1.$$
For every $\mu$, let
$$u_{\mu,1}(x):=\mathbb1_Q(x)\frac{f_\mu(x)-f_\mu(c)}{r_j},\qquad v_{\mu,1}(x'):=\mathbb1_{Q'}(x'),$$
$$u_{\mu,2}(x):=\mathbb1_Q(x)\frac{f_\mu(c)-f_\mu(c')}{r_j},\qquad v_{\mu,2}(x'):=\mathbb1_{Q'}(x'),$$
$$u_{\mu,3}(x):=\mathbb1_Q(x),\qquad v_{\mu,3}(x'):=-\mathbb1_{Q'}(x')\frac{f_\mu(x')-f_\mu(c')}{r_j},$$
and, for $a=1,2,3$, set
$$b_{\mu,a}(w):=\frac{\kappa}{h_j}\big(g_\mu(w+h_je)-g_\mu(w)\big).$$

Observe that in \eqref{separated-translated-phase-difference} we were able to reduce $\Delta_{h_j,e}\phi$ to the following form on 
$Q\times Q'\times\mathbb R^{d-2}$:
\begin{equation}
\lambda\Delta_{h_j,e}\phi(x,x',w)-\theta_0=\sum_{\nu\in \Xi} u_\nu(x)v_\nu(x')b_\nu(w).
\label{separated-phase-series}
\end{equation}

Now that we have expressed $\Delta_{h_j,e}\phi$ in such a way, we estimate $\sum_{\nu\in \Xi}\|u_\nu\|_{L^\infty}\|v_\nu\|_{L^\infty}\|b_\nu\|_{L^\infty}$. Recall that thanks to \eqref{separated-translated-phase-difference} and to the definition of $u_\nu$, $v_\nu$ and $b_\nu$ we have
\begin{equation*}
\begin{split}
&\qquad\qquad\qquad\qquad\sum_{\nu\in \Xi}\|u_\nu\|_{L^\infty}\|v_\nu\|_{L^\infty}\|b_\nu\|_{L^\infty}
=\frac{\kappa}{r_j}\Big(\sup_{x\in Q}|\langle x-c,e\rangle|+\sup_{x'\in Q'}|\langle x'-c',e\rangle|\Big)\\
&+\frac{\kappa}{r_jh_j}\sum_\mu\Big(\|f_\mu-f_\mu(c)\|_{L^\infty(Q)}+|f_\mu(c)-f_\mu(c')|+\|f_\mu-f_\mu(c')\|_{L^\infty(Q')}\Big)\|g_\mu(\cdot+h_je)-g_\mu\|_{L^\infty(\mathbb R^{d-2})}.
\end{split}
\nonumber
\end{equation*}
By construction of the cubes and the choice of $x$ and $x'$ we have 
$$\frac{\kappa}{r_j}\Big(\sup_{x\in Q}|\langle x-c,e\rangle|+\sup_{x'\in Q'}|\langle x'-c',e\rangle|\Big)\leq \frac{\kappa}{32\sqrt{d-2}}.$$

Also, the cube geometry and the Lipschitz bounds give
$$\|f_\mu-f_\mu(c)\|_{L^\infty(Q)}\leq\frac{r_j}{128}\Lip(f_\mu),\qquad \|f_\mu-f_\mu(c')\|_{L^\infty(Q')}\leq\frac{r_j}{128}\Lip(f_\mu),$$
$$|f_\mu(c)-f_\mu(c')|\leq5r_j\Lip(f_\mu),\qquad \|g_\mu(\cdot+h_je)-g_\mu\|_{L^\infty(\mathbb R^{d-2})}\leq h_j\Lip(g_\mu).$$
Substitution therefore yields
\begin{equation}
\begin{split}
\sum_{\nu\in \Xi}\|u_\nu\|_{L^\infty}\|v_\nu\|_{L^\infty}\|b_\nu\|_{L^\infty}
&\leq\frac{\kappa}{32\sqrt{d-2}}+\kappa\Big(\frac1{128}+5+\frac1{128}\Big)\sum_\mu\Lip(f_\mu)\Lip(g_\mu)\\
&\leq\frac{\kappa}{32\sqrt{d-2}}+6\kappa\varepsilon\leq \frac{\kappa}{16\sqrt{d-2}},
\end{split}
\label{residual-phase-coefficient-bound}
\end{equation}
where the last inequality uses $\varepsilon\leq\frac{1}{384\sqrt{d-2}}$.
This estimate implies in particular that 
$$|\lambda\Delta_{h_j,e}\phi(x,x',w)-\theta_0|\leq\sum_{\nu\in \Xi}\|u_\nu\|_{L^\infty}\|v_\nu\|_{L^\infty}\|b_\nu\|_{L^\infty}\leq\frac{\kappa}{16\sqrt{d-2}}.$$
Hence, thanks to \eqref{stimatheta0}, we have that 
\begin{equation}
    \frac{\kappa}{2\sqrt{d-2}}\leq\lambda\Delta_{h_j,e}\phi(x,x',w)\leq6\kappa<\frac{\pi}{2}\quad\text{and thus}\quad|e^{i\lambda\Delta_{h_j,e}\phi(x,x',w)}-1|\geq\frac{2}{\pi}\lambda\Delta_{h_j,e}\phi(x,x',w)\geq\frac{\kappa}{\pi\sqrt{d-2}}.\label{stimaqjpippo}
\end{equation}

Before moving on, we observe that the series in \eqref{separated-phase-series} converges absolutely and uniformly, since
$$\sum_{\nu\in\Xi}|u_\nu(x)v_\nu(x')b_\nu(w)|\leq\sum_{\nu\in\Xi}\|u_\nu\|_{L^\infty}\|v_\nu\|_{L^\infty}\|b_\nu\|_{L^\infty}\leq\frac{\kappa}{16\sqrt{d-2}}.$$
In the language of \cref{separated-expansions}, applied with $E=Q$ and $E'=Q'$, it is a separated expansion of size at most $\frac{\kappa}{16\sqrt{d-2}}$. By that lemma, its powers may be multiplied out term by term, and so may the exponential series used below. The terms of the latter are indexed by the tuples $\gamma=(k,\nu_1,\ldots,\nu_k)$.

Thanks to the computation in \eqref{stimaqjpippo}, the quantity $e^{i\lambda\Delta_{h_j,e}\phi(x,x',w)}-1$ does not vanish. Therefore for every
integrable function
 $\varphi$, we have
$$\int e^{i\lambda(\phi(x,w)-\phi(x',w))}\varphi(w)\,dw=\int e^{i\lambda(\phi(x,w)-\phi(x',w))}\frac{e^{i\lambda\Delta_{h_j,e}\phi(x,x',w)}-1}{e^{i\lambda\Delta_{h_j,e}\phi(x,x',w)}-1}\varphi(w)\,dw.$$
By the definition of $\Delta_{h_j,e}\phi$, there holds
$$e^{i\lambda(\phi(x,w)-\phi(x',w))}e^{i\lambda\Delta_{h_j,e}\phi(x,x',w)}=e^{i\lambda(\phi(x,w+h_je)-\phi(x',w+h_je))}$$
and consequently, we have
\begin{equation}
\begin{split}
    \int e^{i\lambda(\phi(x,w)-\phi(x',w))}\varphi(w)dw&=\int e^{i\lambda(\phi(x,w+h_je)-\phi(x',w+h_je))}\frac{\varphi(w)}{e^{i\lambda\Delta_{h_j,e}\phi(x,x',w)}-1}dw\\
&-\int e^{i\lambda(\phi(x,w)-\phi(x',w))}\frac{\varphi(w)}{e^{i\lambda\Delta_{h_j,e}\phi(x,x',w)}-1}dw.
\end{split}
\end{equation}
Changing coordinates in the first integral, we have
$$\int e^{i\lambda(\phi(x,u)-\phi(x',u))}\frac{\varphi(u-h_je)}{e^{i\lambda\Delta_{h_j,e}\phi(x,x',u-h_je)}-1}\,du.$$
Renaming $u$ as $w$ and combining the two integrals, we obtain
\begin{equation}
\int e^{i\lambda(\phi(x,w)-\phi(x',w))}\varphi(w)\,dw=\int e^{i\lambda(\phi(x,w)-\phi(x',w))}\Big(\underbrace{\frac{\varphi(w-h_je)}{e^{i\lambda\Delta_{h_j,e}\phi(x,x',w-h_je)}-1}-\frac{\varphi(w)}{e^{i\lambda\Delta_{h_j,e}\phi(x,x',w)}-1}}_{(*)}\Big)dw.
\label{finite-difference-identity}
\end{equation}
Now, we will inductively apply \eqref{finite-difference-identity} using the notation already introduced. In the following finite-difference calculation, $x,x'$ and $e$ remain fixed:
$$q_j(w):=\frac{1}{e^{i\lambda\Delta_{h_j,e}\phi(x,x',w)}-1},\qquad\delta_{h_j,e}b(w):=b(w-h_je)-b(w).$$
Notice that
\begin{equation}
\begin{split}
     (*)=&\frac{\varphi(w-h_je)-\varphi(w)}{e^{i\lambda\Delta_{h_j,e}\phi(x,x',w-h_je)}-1}+\varphi(w)\Big(\frac{1}{e^{i\lambda\Delta_{h_j,e}\phi(x,x',w-h_je)}-1}-\frac{1}{e^{i\lambda\Delta_{h_j,e}\phi(x,x',w)}-1}\Big)\\
    &\qquad\qquad\qquad=q_j(w-h_je)\delta_{h_j,e}\varphi(w)+\varphi(w)\delta_{h_j,e}q_j(w).
\end{split}
    \label{identitastepbase}
\end{equation}

Let us prove by induction that, for every integer
$m\geq1$ and every integrable function $a$, there holds
\begin{equation}
\begin{split}
\int e^{i\lambda(\phi(x,w)-\phi(x',w))}a(w)dw&=\int e^{i\lambda(\phi(x,w)-\phi(x',w))}\delta_{h_j,e}^{m}a(w)\prod_{k=1}^{m}q_j(w-kh_je)dw\\
&+\sum_{\ell=0}^{m-1}\int e^{i\lambda(\phi(x,w)-\phi(x',w))}\delta_{h_j,e}^{\ell}a(w)\delta_{h_j,e}\Big(\prod_{k=0}^{\ell}q_j(w-kh_je)\Big)dw.
\end{split}
\label{finite-difference-induction}
\end{equation}

The base step is \eqref{identitastepbase} with $\varphi=a$. So, let us prove the formula for $m+1$ assuming we have obtained the formula for $m$, that is \eqref{finite-difference-induction}. We apply \eqref{finite-difference-identity} and \eqref{identitastepbase} to the first term in the right-hand side of \eqref{finite-difference-induction}, with
$$\varphi(w):=\delta_{h_j,e}^{m}a(w)\prod_{k=1}^{m}q_j(w-kh_je),$$
obtaining
\begin{equation}
\begin{split}
\int e^{i\lambda(\phi(x,w)-\phi(x',w))}a(w)\,dw&=\int e^{i\lambda(\phi(x,w)-\phi(x',w))}\delta_{h_j,e}^{m}a(w)\prod_{k=1}^{m}q_j(w-kh_je)\,dw\\
&\qquad+\sum_{\ell=0}^{m-1}\int e^{i\lambda(\phi(x,w)-\phi(x',w))}\delta_{h_j,e}^{\ell}a(w)\delta_{h_j,e}\Big(\prod_{k=0}^{\ell}q_j(w-kh_je)\Big)\,dw\\
&=\int e^{i\lambda(\phi(x,w)-\phi(x',w))}\delta_{h_j,e}\Big(\delta_{h_j,e}^{m}a(w)\prod_{k=1}^{m}q_j(w-kh_je)\Big)q_j(w-h_je)\,dw\\
&\qquad+\int e^{i\lambda(\phi(x,w)-\phi(x',w))}\delta_{h_j,e}^{m}a(w)\prod_{k=1}^{m}q_j(w-kh_je)\delta_{h_j,e}q_j(w)\,dw\\
&\qquad+\sum_{\ell=0}^{m-1}\int e^{i\lambda(\phi(x,w)-\phi(x',w))}\delta_{h_j,e}^{\ell}a(w)\delta_{h_j,e}\Big(\prod_{k=0}^{\ell}q_j(w-kh_je)\Big)\,dw.
\end{split}
\end{equation}
Expanding the difference and adding and subtracting $\delta_{h_j,e}^{m}a(w)\prod_{k=2}^{m+1}q_j(w-kh_je)$, we obtain
\begin{equation}
\begin{split}
&\delta_{h_j,e}\Big(\delta_{h_j,e}^{m}a(w)\prod_{k=1}^{m}q_j(w-kh_je)\Big)=\delta_{h_j,e}^{m}a(w-h_je)\prod_{k=2}^{m+1}q_j(w-kh_je)-\delta_{h_j,e}^{m}a(w)\prod_{k=1}^{m}q_j(w-kh_je)\\
=&\Big(\delta_{h_j,e}^{m}a(w-h_je)-\delta_{h_j,e}^{m}a(w)\Big)\prod_{k=2}^{m+1}q_j(w-kh_je)+\delta_{h_j,e}^{m}a(w)\Big(\prod_{k=2}^{m+1}q_j(w-kh_je)-\prod_{k=1}^{m}q_j(w-kh_je)\Big)\\
&\qquad\qquad\qquad=\delta_{h_j,e}^{m+1}a(w)\prod_{k=2}^{m+1}q_j(w-kh_je)+\delta_{h_j,e}^{m}a(w)\delta_{h_j,e}\Big(\prod_{k=1}^{m}q_j(w-kh_je)\Big).
\end{split}
\end{equation}
This implies in particular that
\begin{equation}
    \begin{split}
        \int e^{i\lambda(\phi(x,w)-\phi(x',w))}a(w)dw&=\int e^{i\lambda(\phi(x,w)-\phi(x',w))}\delta_{h_j,e}^{m+1}a(w)\prod_{k=1}^{m+1}q_j(w-kh_je)dw\\
        &\qquad+\int e^{i\lambda(\phi(x,w)-\phi(x',w))}\delta_{h_j,e}^{m}a(w)\delta_{h_j,e}\Big(\prod_{k=1}^{m}q_j(w-kh_je)\Big)q_j(w-h_je)dw\\
        &\qquad+\int e^{i\lambda(\phi(x,w)-\phi(x',w))}\delta_{h_j,e}^{m}a(w)\prod_{k=1}^{m}q_j(w-kh_je)\delta_{h_j,e}q_j(w)\,dw\\
&\qquad+\sum_{\ell=0}^{m-1}\int e^{i\lambda(\phi(x,w)-\phi(x',w))}\delta_{h_j,e}^{\ell}a(w)\delta_{h_j,e}\Big(\prod_{k=0}^{\ell}q_j(w-kh_je)\Big)\,dw.
    \end{split}
\end{equation}
Now notice that the second and third line can be rewritten as follows
\begin{equation}
\begin{split}
q_j(w-h_je)\delta_{h_j,e}\Big(\prod_{k=1}^{m}q_j(w-kh_je)\Big)+\prod_{k=1}^{m}q_j(w-kh_je)\delta_{h_j,e}q_j(w)&=\prod_{k=1}^{m+1}q_j(w-kh_je)-\prod_{k=0}^{m}q_j(w-kh_je)\\
&=\delta_{h_j,e}\Big(\prod_{k=0}^{m}q_j(w-kh_je)\Big).
\end{split}
\end{equation}
This implies in particular that
\begin{equation*}
    \begin{split}
         \int e^{i\lambda(\phi(x,w)-\phi(x',w))}a(w)dw&=\int e^{i\lambda(\phi(x,w)-\phi(x',w))}\delta_{h_j,e}^{m+1}a(w)\prod_{k=1}^{m+1}q_j(w-kh_je)dw\\
         &\qquad+\sum_{\ell=0}^{m}\int e^{i\lambda(\phi(x,w)-\phi(x',w))}\delta_{h_j,e}^{\ell}a(w)\delta_{h_j,e}\Big(\prod_{k=0}^{\ell}q_j(w-kh_je)\Big)dw.
         \nonumber
    \end{split}
\end{equation*}
Taking $m=d-1$ and $a=a_z$ in \eqref{finite-difference-induction}, substituting into the kernel of $\mathcal K_{e,j,z}$, and summing over $(Q,Q')\in\mathcal P_{j,e}$ proves \eqref{splitKe,j}.

\medskip

\spacedlowsmallcaps{Proof of }\eqref{required-principal-estimate}.
We first estimate the integral containing $\delta_{h_j,e}^{d-1}a_{z}$. The fundamental theorem of calculus gives
$$\delta_{h_j,e}a_{z}(w)=-h_j\int_0^1\partial_ea_{z}(w-th_je)dt.$$
Applying this formula repeatedly yields, for every $1\leq m\leq d-1$,
\begin{equation}
\delta_{h_j,e}^{m}a_{z}(w)=(-h_j)^m\int_{[0,1]^m}\partial_e^ma_{z}\big(w-h_j(t_1+\cdots+t_m)e\big)dt_1\cdots dt_m.
\label{finite-difference-of-cutoff}
\end{equation}
This identity implies in particular that
\begin{equation*}
\begin{split}
\|\delta_{h_j,e}^{m}a_z\|_{L^1}
&\leq h_j^m\int_{[0,1]^m}\int|\partial_e^ma_z(w-h_j(t_1+\cdots+t_m)e)|\,dw\,dt_1\cdots dt_m\\
&=h_j^m\|\partial_e^ma_z\|_{L^1}=h_j^m\|\partial_e^m(|\eta_1|^2)\|_{L^1}.
\end{split}
\nonumber
\end{equation*}
The last equality follows from $a_z(w)=|\eta_1(w-z)|^2$, so the bound is independent of $z$.
By \eqref{stimaqjpippo}, for every $(x,x',w)\in Q\times Q'\times\mathbb R^{d-2}$,
$$\prod_{k=1}^{d-1}|q_j(x,x',w-kh_je)|\leq\Big(\frac{\pi\sqrt{d-2}}{\kappa}\Big)^{d-1}.$$
Consequently, we have
\begin{equation}
\begin{split}
&\Big|\int e^{i\lambda(\phi(x,w)-\phi(x',w))}\delta_{h_j,e}^{d-1}a_z(w)\prod_{k=1}^{d-1}q_j(x,x',w-kh_je)\,dw\Big|\leq\Big(\frac{\pi\sqrt{d-2}}{\kappa}\Big)^{d-1}h_j^{d-1}\|\partial_e^{d-1}(|\eta_1|^2)\|_{L^1}.
\end{split}
\label{first-integral-estimate}
\end{equation}
Since $0\leq\rho_{j,e}\leq1$, the definition of $\mathcal K_{e,j,z}^{Q,Q',d-1}$ in \eqref{splitKe,j} gives
\begin{equation}
    \begin{split}
      |\mathcal K_{e,j,z}^{Q,Q',d-1}F(x)|\leq&\mathbb1_Q(x)\|\eta_0\|_{L^\infty}^2\Big(\frac{\pi\sqrt{d-2}}{\kappa}\Big)^{d-1}h_j^{d-1}\|\partial_e^{d-1}(|\eta_1|^2)\|_{L^1}\int_{Q'}|F(x')|\,dx'\\
      \leq &\Leb^{d-2}(Q')^{\frac{1}{2}}\mathbb1_Q(x)\|\eta_0\|_{L^\infty}^2\Big(\frac{\pi\sqrt{d-2}}{\kappa}\Big)^{d-1}h_j^{d-1}\|\partial_e^{d-1}(|\eta_1|^2)\|_{L^1} \|\mathbb1_{Q'}F\|_{L^2},
    \end{split}
\end{equation}
for every $F\in L^2$ and where the last inequality comes from Cauchy-Schwarz. This implies in particular  that 
$$\|\mathcal K_{e,j,z}^{Q,Q',d-1}F\|_{L^2}\leq \Leb^{d-2}(Q)^{\frac{1}{2}}\Leb^{d-2}(Q')^{\frac{1}{2}}\|\eta_0\|_{L^\infty}^2\Big(\frac{\pi\sqrt{d-2}}{\kappa}\Big)^{d-1}h_j^{d-1}\|\partial_e^{d-1}(|\eta_1|^2)\|_{L^1} \|\mathbb1_{Q'}F\|_{L^2}.$$
Since by construction we have that  $h_j=\frac{\kappa}{\lambda r_j}$, that $r_j=\frac{2^jM}{\lambda}$, that $Q\subseteq B(c,\frac{r_j}{64})$ and $Q'\subseteq B(c',\frac{r_j}{64})$  we have 
\begin{equation*}
\begin{split}
\|\mathcal K_{e,j,z}^{Q,Q',d-1}\|_{2\to2}
&\leq C_d\|\eta_0\|_{L^\infty}^2\Big(\frac{\pi\sqrt{d-2}}{\kappa}\Big)^{d-1}r_j^{d-2}h_j^{d-1}\|\partial_e^{d-1}(|\eta_1|^2)\|_{L^1}=\frac{C_d}{2^jM}\|\eta_0\|_{L^\infty}^2\|\partial_e^{d-1}(|\eta_1|^2)\|_{L^1}\lambda^{-(d-2)},
\end{split}
\nonumber
\end{equation*}
This proves \eqref{required-principal-estimate}, uniformly in $z,j,e,Q,Q'$.

\medskip
\spacedlowsmallcaps{Structure of the operators $\mathcal K_{e,j,z}^{Q,Q',\ell}$ for $0\leq \ell\leq d-2$.}
Fix $j\geq0$, $e\in\mathscr E^\pm$, $(Q,Q')\in\mathcal P_{j,e}$ and $0\leq\ell\leq d-2$. First of all let us observe that, as anticipated in the previous paragraph, we have
$$\delta_{h_j,e}^{\ell}a_z(w)=\sum_{r=0}^{\ell}(-1)^{\ell-r}\binom{\ell}{r}|\eta_1(w-(z+rh_je))|^2.$$
Substituting this identity into \eqref{splitKe,j} gives
$$\mathcal K_{e,j,z}^{Q,Q',\ell}F(x)=\sum_{r=0}^{\ell}(-1)^{\ell-r}\binom{\ell}{r}\mathcal K_{e,j,z}^{Q,Q',\ell,r}F(x),$$
where 
\begin{equation*}
\begin{split}
\mathcal K_{e,j,z}^{Q,Q',\ell,r}F(x):=\eta_0(x)\int\overline{\eta_0(x')}F(x')\int e^{i\lambda(\phi(x,w)-\phi(x',w))}|\eta_1(w-(z+rh_je))|^2c_{Q,Q',j,e,\ell}(x,x',w)dwdx',
\end{split}
\nonumber
\end{equation*}
and 
\begin{equation}
    c_{Q,Q',j,e,\ell}(x,x',w):=\mathbb1_Q(x)\mathbb1_{Q'}(x')\rho_{j,e}(x,x')\delta_{h_j,e}\Big(\prod_{k=0}^{\ell}q_j(x,x',w-kh_je)\Big).
    \label{defcdime}
\end{equation}
In the following we shall construct a decomposition
\begin{equation}
c_{Q,Q',j,e,\ell}(x,x',w)=c^{\mathrm g}_{Q,Q',j,e,\ell}(x,x',w)+c^{\mathrm t}_{Q,Q',j,e,\ell}(x,x',w)=\sum_{\nu\in \Xi} u_\nu(x)v_\nu(x')b_\nu(w)+c^{\mathrm t}_{Q,Q',j,e,\ell}(x,x',w),
\label{remainder-coefficient-splitting}
\end{equation}
where $u_\nu$ and $v_\nu$ are supported in $Q$ and $Q'$, respectively. We will find sequences $A_j$ and $B_j$, that will be independent of $e,\ell,Q,Q'$ and $z$, with the properties claimed in \eqref{required-remainder-representation} and \eqref{sommaA_j} such that the coefficients $A_j$ will therefore satisfy the following bound
\begin{equation}
\sum_{\nu\in \Xi}\|u_\nu\|_{L^\infty}\|v_\nu\|_{L^\infty}\|b_\nu\|_{L^\infty}\leq 2^{-\ell}A_j,
\label{bounded-coefficient-condition}
\end{equation}
while the $B_j$s bounds the following quantity
\begin{equation}
2^\ell\|\eta_0\|_{L^\infty}^2\|\eta_1\|_{L^\infty}^2\Leb^{d-2}(
Q)^{\frac{1}{2}}\Leb^{d-2}(Q')^{\frac{1}{2}}\operatorname*{ess\,sup}_{(x,x')\in Q\times Q'}\int|c^{\mathrm t}_{Q,Q',j,e,\ell}(x,x',w)|dw\leq B_j.
\label{integrable-coefficient-condition}
\end{equation}
Let us show what we can obtain with assuming that \eqref{remainder-coefficient-splitting}-\eqref{integrable-coefficient-condition} hold. Since the operators 
$$G_{e,j,z}^{Q,Q',\ell,r}F(x):=\eta_0(x)\int\overline{\eta_0(x')}F(x')\int e^{i\lambda(\phi(x,w)-\phi(x',w))}|\eta_1(w-(z+rh_je))|^2c_{Q,Q',j,e,\ell}^\mathrm{g}(x,x',w)dwdx'$$
are in the form of \eqref{operatorprototype}, we know by \eqref{idnetitaRcoefficentec} that under the assumption \eqref{remainder-coefficient-splitting}, there holds the identity
$$G_{e,j,z}^{Q,Q',\ell,r}F=\sum_{\nu\in \Xi} u_\nu T_{\lambda,z+rh_je}\Big(b_\nu T_{\lambda,z+rh_je}^*(v_\nu F)\Big).$$
By \eqref{separated-coefficient-operator-bound} and \eqref{bounded-coefficient-condition}, this series converges in operator norm and satisfies
\begin{equation*}
\Big\|\sum_{\nu\in \Xi} u_\nu T_{\lambda,z+rh_je}\Big(b_\nu T_{\lambda,z+rh_je}^*(v_\nu F)\Big)\Big\|_{L^2}\leq2^{-\ell}A_jN_\lambda\|F\|_{L^2}.
\nonumber
\end{equation*}
On the other hand, we also know that 
\begin{equation}
    \begin{split}
      &  \|\mathcal K_{e,j,z}^{Q,Q',\ell,r}F-G_{e,j,z}^{Q,Q',\ell,r}F\|_{L^2}\\
&\qquad\qquad\leq \|\eta_0(\cdot)\int\overline{\eta_0(x')}F(x')\int e^{i\lambda(\phi(\cdot,w)-\phi(x',w))}|\eta_1(w-(z+rh_je))|^2c_{Q,Q',j,e,\ell}^{\mathrm t}(x,x',w)dwdx'\|_{L^2}\\
&\qquad\qquad\leq \|\eta_1\|^2_{L^\infty}\|\eta_0\|^2_{L^\infty}\|\mathbb 1_{Q}\|_{L^2}\int_{Q'}|F(x')|dx'\operatorname*{ess\,sup}_{(x,x')\in Q\times Q'}\int |c_{Q,Q',j,e,\ell}^{\mathrm t}|(x,x',w)dw\\
&\qquad\qquad\leq \|\eta_0\|_{L^\infty}^2\|\eta_1\|_{L^\infty}^2\Leb^{d-2}(
Q)^{\frac{1}{2}}\Leb^{d-2}(Q')^{\frac{1}{2}}\operatorname*{ess\,sup}_{(x,x')\in Q\times Q'}\int|c^{\mathrm t}_{Q,Q',j,e,\ell}(x,x',w)|dw\|F\|_{L^2}\leq 2^{-\ell}B_j\|F\|_{L^2}.
    \end{split}
\end{equation}
Consequently, $\|\mathcal K_{e,j,z}^{Q,Q',\ell,r}\|_{2\to2}\leq2^{-\ell}(A_jN_\lambda+B_j)$ and hence summing over $r$, we infer
\begin{equation*}
\|\mathcal K_{e,j,z}^{Q,Q',\ell}\|_{2\to2}\leq\sum_{r=0}^{\ell}\binom{\ell}{r}2^{-\ell}(A_jN_\lambda+B_j)=A_jN_\lambda+B_j.
\nonumber
\end{equation*}
Notice further that the series of operators that approximates $\mathcal K_{e,j,z}^{Q,Q',\ell}$ can be written explicitly as
$$\sum_{r=0}^{\ell}\sum_{\nu\in \Xi}(-1)^{\ell-r}\binom{\ell}{r}u_\nu T_{\lambda,z+rh_je}\Big(b_\nu T_{\lambda,z+rh_je}^*(v_\nu F)\Big).$$
Enumerating the pairs $(r,\nu)$ therefore gives precisely \eqref{required-coefficient-estimate} and \eqref{required-remainder-representation}.

We now construct the splitting of $c_{Q,Q',j,e,\ell}(x,x',w)$ claimed in \eqref{remainder-coefficient-splitting}. Let us go back to the definition of $c_{Q,Q',j,e,\ell}(x,x',w)$ in \eqref{defcdime}.

\spacedlowsmallcaps{The structure of the product of $q_j$s}
In the following, in order to simplify the notation, when the points $x$ and $x'$ are fixed we will drop them and write $q_j(w)$ for $q_j(x,x',w)$ and $\Delta_{h_j,e}\phi(w+\tau e)$ for $\Delta_{h_j,e}\phi(x,x',w+\tau e)$. 
We have to rewrite the variation of the product of the $q_j$s and the annular cutoff $\rho$ as series of products, up to an error satisfying the hypothesis \eqref{remainder-coefficient-splitting}-\eqref{integrable-coefficient-condition}.  
First of all, by the definition of $\delta_{h_j,e}$, we have that
\begin{equation}
\begin{split}
\delta_{h_j,e}\Big(\prod_{k=0}^{\ell}q_j(w-kh_je)\Big)&=\prod_{k=1}^{\ell+1}q_j(w-kh_je)-\prod_{k=0}^{\ell}q_j(w-kh_je)\\
&=\Big(\prod_{k=1}^{\ell}q_j(w-kh_je)\Big)\Big(q_j(w-(\ell+1)h_je)-q_j(w)\Big).
\end{split}
\label{difference-of-translated-reciprocal-products}
\end{equation}
In addition, for every $\tau,\tau'\in\mathbb R$, by definition of $q_j$ we have
\begin{equation}
\begin{split}
&\qquad q_j(w+\tau e)-q_j(w+\tau'e)=q_j(w+\tau e)q_j(w+\tau'e)(e^{i\lambda\Delta_{h_j,e}\phi(w+\tau'e)}-e^{i\lambda\Delta_{h_j,e}\phi(w+\tau e)})\\
=&i\lambda(\Delta_{h_j,e}\phi(w+\tau'e)-\Delta_{h_j,e}\phi(w+\tau e))q_j(w+\tau e)q_j(w+\tau'e)
\int_0^1e^{i\lambda((1-s)\Delta_{h_j,e}\phi(w+\tau e)+s\Delta_{h_j,e}\phi(w+\tau'e))}ds,
\end{split}
\label{reciprocal-difference-representation}
\end{equation}
where in the last identity we used the fundamental theorem of calculus. Notice that the denominators that are contained in $q_j$ do not vanish by \eqref{stimaqjpippo}. Combining \eqref{difference-of-translated-reciprocal-products} and \eqref{reciprocal-difference-representation}, with $\tau=-(\ell+1)h_j$ and $\tau'=0$, gives
\begin{equation}
\begin{split}
&\delta_{h_j,e}\Big(\prod_{k=0}^{\ell}q_j(w-kh_je)\Big)\\
&\qquad=i\lambda(\Delta_{h_j,e}\phi(w)-\Delta_{h_j,e}\phi(w-(\ell+1)h_je))\Big(\prod_{k=0}^{\ell+1}q_j(w-kh_je)\Big)\int_0^1e^{i\lambda((1-s)\Delta_{h_j,e}\phi(w-(\ell+1)h_je)+s\Delta_{h_j,e}\phi(w))}ds.
\label{rappresentazione-con-cui-lavoriamo}
\end{split}
\end{equation}

Let us observe that by \eqref{stimatheta0}, we have 
\begin{equation}
|e^{i\theta_0}-1|=2\sin\left(\frac{\theta_0}{2}\right)\geq\frac{2\theta_0}{\pi}\geq\frac{\kappa}{\pi\sqrt{d-2}}.
   \label{stimadenominatoretheta0} 
\end{equation}
For each $\gamma=(k,\nu_1,\ldots,\nu_k)$, where $k\geq1$ and the $\nu_a$s ranging in $\Xi$, define
$$U_\gamma(x):=\prod_{a=1}^k u_{\nu_a}(x),\qquad V_\gamma(x'):=\prod_{a=1}^k v_{\nu_a}(x'),\qquad B_\gamma(w):=-\frac{e^{i\theta_0}}{e^{i\theta_0}-1}\frac{i^k}{k!}\prod_{a=1}^k b_{\nu_a}(w).$$
We denote the set of these tuples $\gamma$ by $\Xi'$.
The estimate \eqref{residual-phase-coefficient-bound} implies
\begin{equation}
\begin{split}
&\sum_{\gamma\in \Xi'}\|U_\gamma\|_{L^\infty}\|V_\gamma\|_{L^\infty}\|B_\gamma\|_{L^\infty}\leq\frac{1}{|e^{i\theta_0}-1|}\sum_{k=1}^{\infty}\frac{1}{k!}\Big(\sum_{\nu\in \Xi}\|u_\nu\|_{L^\infty}\|v_\nu\|_{L^\infty}\|b_\nu\|_{L^\infty}\Big)^k\leq\frac{\pi\sqrt{d-2}}{\kappa}(e^{\frac{\kappa}{16\sqrt{d-2}}}-1)<\frac12,
\label{stimazionecool}
\end{split}
\end{equation}
where the last inequality uses $e^t-1\leq2t$ for $0\leq t\leq\frac{1}{2}$. Hence, by \eqref{separated-phase-series} and \cref{separated-expansions}, we deduce that 
\begin{equation}
    \begin{split}
        -\frac{e^{i\theta_0}}{e^{i\theta_0}-1}(e^{i(\lambda \Delta_{h_j,e}\phi(x,x',w)-\theta_0)}-1)=-\frac{e^{i\theta_0}}{e^{i\theta_0}-1}\sum_{k=1}^\infty\frac{i^k}{k!}\Big(\lambda \Delta_{h_j,e}\phi(x,x',w)-\theta_0\Big)^k\\
        \overset{\eqref{separated-phase-series}}{=}-\frac{e^{i\theta_0}}{e^{i\theta_0}-1}\sum_{k=1}^\infty\frac{i^k}{k!}\Big(\sum_{\nu\in \Xi} u_\nu(x)v_\nu(x')b_\nu(w)\Big)^k=\sum_{\gamma\in \Xi'}U_\gamma(x)V_\gamma(x')B_\gamma(w),
    \end{split}
\end{equation}
Rearranging terms, we see that 
\begin{equation}
e^{i\lambda\Delta_{h_j,e}\phi(x,x',w)}-1=(e^{i\theta_0}-1)\Big(1-\sum_{\gamma\in \Xi'} U_\gamma(x)V_\gamma(x')B_\gamma(w)\Big).
\label{stimaespansioneesponenx}
\end{equation}
Observe that thanks to \eqref{stimazionecool}, we know that 
$\|\sum_{\gamma\in \Xi'} U_\gamma(x)V_\gamma(x')B_\gamma(w)\|\leq \frac{1}{2}$. Recalling \eqref{stimadenominatoretheta0}  holds, we have that the following chain of identities is well defined
\begin{equation}
    \begin{split}
        q_j(x,x',w)=\frac{1}{e^{i\theta_0}-1}\frac{1}{1-\sum_\gamma U_\gamma(x)V_\gamma(x')B_\gamma(w)}=\frac{1}{e^{i\theta_0}-1}\sum_{m=0}^\infty\Big(\sum_{\gamma\in \Xi'} U_\gamma(x)V_\gamma(x')B_\gamma(w)\Big)^m
    \end{split}
\end{equation}
By \eqref{stimazionecool}, the series $\sum_{\gamma\in\Xi'}U_\gamma V_\gamma B_\gamma$ is a separated expansion of size at most $\frac{1}{2}$, so \cref{separated-expansions}, applied to the geometric series, multiplies out the powers and turns the right-hand side into a single absolutely convergent series of separated products.

So far we have been able to express the function $q_j$ as a series of triple products of functions that depend separately on $x,x',w$. For each tuple $\alpha=(m,\gamma_1,\ldots,\gamma_m)$, define its three factors by
$$\widetilde u_\alpha(x):=\prod_{a=1}^mU_{\gamma_a}(x),\qquad \widetilde v_\alpha(x'):=\prod_{a=1}^mV_{\gamma_a}(x'),\qquad \widetilde b_\alpha(w):=\frac{1}{e^{i\theta_0}-1}\prod_{a=1}^mB_{\gamma_a}(w).$$
These formulas make sense when $m\geq 1$. For $m=0$ we let $\widetilde u_0:=\mathbb 1_Q$, $\widetilde v_0(x'):=\mathbb 1_{Q'}$ and $\widetilde b_0:=(e^{i\theta_0}-1)^{-1}$.
We denote by $\Xi''$ the set of the tuples $\alpha$ with $m\geq1$, together with the index $0$.
The preceding computation proves
\begin{equation}
q_j(x,x',w)=\sum_{\alpha\in \Xi''}\widetilde u_\alpha(x)\widetilde v_\alpha(x')\widetilde b_\alpha(w),\qquad (x,x',w)\in Q\times Q'\times\mathbb R^{d-2},
\label{separated-reciprocal-expansion}
\end{equation}
with
\begin{equation}
\sum_{\alpha\in \Xi''}\|\widetilde u_\alpha\|_{L^\infty}\|\widetilde v_\alpha\|_{L^\infty}\|\widetilde b_\alpha\|_{L^\infty}\leq\frac{1}{|e^{i\theta_0}-1|}\sum_{m=0}^{\infty}2^{-m}\leq\frac{2\pi\sqrt{d-2}}{\kappa}.
\label{reciprocal-coefficient-bound}
\end{equation}
By the last assertion of \cref{separated-expansions}, replacing $w$ by $w-kh_je$ in \eqref{separated-reciprocal-expansion} changes only the factors $\widetilde b_\alpha$ and preserves the size. Applying \eqref{separated-reciprocal-expansion} at $w-kh_je$, for $k=0,\ldots,\ell+1$, and multiplying the resulting series by the rule for products in \cref{separated-expansions}, we obtain
\begin{equation}
\begin{split}
&\prod_{k=0}^{\ell+1}q_j(x,x',w-kh_je)=\sum_{\alpha_0,\ldots,\alpha_{\ell+1}\in \Xi''}\Big(\prod_{k=0}^{\ell+1}\widetilde u_{\alpha_k}(x)\Big)\Big(\prod_{k=0}^{\ell+1}\widetilde v_{\alpha_k}(x')\Big)\Big(\prod_{k=0}^{\ell+1}\widetilde b_{\alpha_k}(w-kh_je)\Big).
\end{split}
\label{translated-reciprocal-products}
\end{equation}
By \eqref{reciprocal-coefficient-bound}, this is a separated expansion whose size satisfies
\begin{equation}
\sum_{\alpha_0,\ldots,\alpha_{\ell+1}\in\Xi''}\prod_{k=0}^{\ell+1}\|\widetilde u_{\alpha_k}\|_{L^\infty}\|\widetilde v_{\alpha_k}\|_{L^\infty}\|\widetilde b_{\alpha_k}\|_{L^\infty}\leq\left(\frac{2\pi\sqrt{d-2}}{\kappa}\right)^{\ell+2}.
\label{translated-reciprocal-size}
\end{equation}

\spacedlowsmallcaps{The expansion of the integral.} Now, we move on with the study of the coefficient $c_{Q,Q',j,e,\ell}$. The term we will study in this paragraph is the integral in the expansion for $c_{Q,Q',j,e,\ell}$ given in \eqref{rappresentazione-con-cui-lavoriamo}.   
Thanks to \eqref{separated-phase-series}, we infer that 
\begin{equation*}
\begin{split}
&\lambda((1-s)\Delta_{h_j,e}\phi(x,x',w-(\ell+1)h_je)+s\Delta_{h_j,e}\phi(x,x',w))-\theta_0=\sum_\nu u_\nu(x)v_\nu(x')((1-s)b_\nu(w-(\ell+1)h_je)+s b_\nu(w)).
\nonumber
\end{split}
\end{equation*}
For every $\gamma=(k,\nu_1,\ldots, \nu_k)\in {\Xi'}$, let us introduce the function
$$C_\gamma(w):=\frac{e^{i\theta_0}i^k}{k!}\int_0^1\prod_{a=1}^k\Big((1-s)b_{\nu_a}(w-(\ell+1)h_je)+s b_{\nu_a}(w)\Big)ds,$$
and notice that since $|(1-s)b_\nu(w-(\ell+1)h_je)+s b_\nu(w)|\leq\|b_\nu\|_{L^\infty}$ for every $0\leq s\leq1$, and hence
$$\|C_\gamma\|_{L^\infty}\leq \frac{1}{k!}\prod_{a=1}^k\|b_{\nu_a}\|_{L^\infty}.$$
For the constant term set $U_0(x):=\mathbb1_Q(x)$, $V_0(x'):=\mathbb1_{Q'}(x')$ and $C_0(w):=e^{i\theta_0}$.
Using the definitions of $U_\gamma,V_\gamma$ and the preceding estimate for $C_\gamma$, we obtain
\begin{equation}
\begin{split}
&\sum_{\gamma\in\{0\}\cup\Xi'}\|U_\gamma\|_{L^\infty}\|V_\gamma\|_{L^\infty}\|C_\gamma\|_{L^\infty}\leq\sum_{\gamma\in\{0\}\cup\Xi'}\Big(\prod_{a=1}^k\|u_{\nu_a}\|_{L^\infty}\Big)\Big(\prod_{a=1}^k\|v_{\nu_a}\|_{L^\infty}\Big)\frac{1}{k!}\prod_{a=1}^k\|b_{\nu_a}\|_{L^\infty}\\
&=\sum_{k=0}^{\infty}\frac{1}{k!}\sum_{(\nu_1,\ldots,\nu_k)\in\Xi^k}\prod_{a=1}^k\Big(\|u_{\nu_a}\|_{L^\infty}\|v_{\nu_a}\|_{L^\infty}\|b_{\nu_a}\|_{L^\infty}\Big)=\sum_{k=0}^{\infty}\frac{1}{k!}\Big(\sum_{\nu\in\Xi}\|u_\nu\|_{L^\infty}\|v_\nu\|_{L^\infty}\|b_\nu\|_{L^\infty}\Big)^k\\
&\qquad\qquad\qquad\qquad\qquad\overset{\eqref{residual-phase-coefficient-bound}}{\leq}\sum_{k=0}^{\infty}\frac{1}{k!}\Big(\frac{\kappa}{16\sqrt{d-2}}\Big)^k=e^{\frac{\kappa}{16\sqrt{d-2}}}.
\end{split}
\label{convex-phase-exponential-bound}
\end{equation}
In particular, the series $\sum_{\gamma\in \{0\}\cup\Xi'} U_\gamma(x)V_\gamma(x')C_\gamma(w)$ converges absolutely and uniformly. For each fixed $s\in[0,1]$, the functions $(1-s)b_\nu(\cdot-(\ell+1)h_je)+sb_\nu$ have $L^\infty$ norms at most $\|b_\nu\|_{L^\infty}$, so \cref{separated-expansions} expands the exponential of the averaged phase below with the majorant \eqref{convex-phase-exponential-bound}, which does not depend on $s$. The integral in $s$ may therefore be taken term by term, and together with \eqref{separated-phase-series} and the definition of $C_\gamma$ this justifies the following chain of identities
\begin{equation}
\begin{split}
&\int_0^1e^{i\lambda((1-s)\Delta_{h_j,e}\phi(x,x',w-(\ell+1)h_je)+s\Delta_{h_j,e}\phi(x,x',w))}ds\\
&\qquad=e^{i\theta_0}\int_0^1\sum_{k=0}^{\infty}\frac{i^k}{k!}\Big(\sum_{\nu\in\Xi}u_\nu(x)v_\nu(x')\big((1-s)b_\nu(w-(\ell+1)h_je)+s b_\nu(w)\big)\Big)^kds\\
&\qquad=\sum_{k=0}^{\infty}\sum_{(\nu_1,\ldots,\nu_k)\in\Xi^k}\Big(\prod_{a=1}^ku_{\nu_a}(x)\Big)\Big(\prod_{a=1}^kv_{\nu_a}(x')\Big)\frac{e^{i\theta_0}i^k}{k!}\int_0^1\prod_{a=1}^k\Big((1-s)b_{\nu_a}(w-(\ell+1)h_je)+s b_{\nu_a}(w)\Big)ds\\
&\qquad=\sum_{\gamma\in\{0\}\cup\Xi'}U_\gamma(x)V_\gamma(x')C_\gamma(w).
\end{split}
\label{averaged-phase-exponential-expansion}
\end{equation}
The last piece of $c_{Q,Q',j,e,\ell}$ is therefore $\rho_{j,e}(x,x')$, which is the subject of the next paragraph.

\smallskip

\spacedlowsmallcaps{The expansion of the annular cutoff.}
On the fixed pair $(Q,Q')\in\mathcal P_{j,e}$, we have 
$$\mathbb1_Q(x)\mathbb1_{Q'}(x')\rho_{j,e}(x,x')=\mathbb1_Q(x)\mathbb1_{Q'}(x')H\Big(\frac{x-x'}{r_j}\Big),\qquad H(y):=\rho\left(\frac{y}{2}\right)-\rho(y).$$
Since $\supp H\subset\overline{B(0,4)}$, its restriction to $[-8,8]^{d-2}$ extends smoothly with period $16$ in each coordinate. The Fourier series of such a periodic expansion is
\begin{equation}
H(y)=\sum_{m\in\mathbb Z^{d-2}}c_me^{\frac{\pi i}{8}\langle m,y\rangle},\quad \text{where}\quad c_m:=16^{-(d-2)}\int_{[-8,8]^{d-2}}H(y)e^{-\frac{\pi i}{8}\langle m,y\rangle}dy,\quad\text{for every }y\in[-8,8]^{d-2}.
\label{fourierrho}
\end{equation}
For $x\in Q$ and $x'\in Q'$, the cube geometry gives
$$|x-x'|\leq|x-c|+|c-c'|+|c'-x'|\leq\Big(\frac1{128}+5+\frac1{128}\Big)r_j<8r_j.$$
We may therefore substitute $y=\frac{x-x'}{r_j}$ into \eqref{fourierrho}, obtaining
\begin{equation}
\begin{split}
\mathbb1_Q(x)\mathbb1_{Q'}(x')\rho_{j,e}(x,x')
&=\mathbb1_Q(x)\mathbb1_{Q'}(x')\sum_{m\in\mathbb Z^{d-2}}c_me^{\frac{\pi i\langle m,x-x'\rangle}{8r_j}}\\
&=\sum_{m\in\mathbb Z^{d-2}}\Big(c_m\mathbb1_Q(x)e^{\frac{\pi i\langle m,x\rangle}{8r_j}}\Big)\Big(\mathbb1_{Q'}(x')e^{-\frac{\pi i\langle m,x'\rangle}{8r_j}}\Big).
\end{split}
\label{annular-separated-expansion}
\end{equation}
We now estimate the coefficients and verify absolute and uniform convergence of the above series. 
By the definition of $c_m$, and since $H$ is compactly supported in $[-8,8]^{d-2}$, a standard computation gives
$$\Big(1+\frac{\pi^2|m|^2}{64}\Big)^{d-1}c_m=16^{-(d-2)}\int(1-\Delta_y)^{d-1}H(y)e^{-\frac{\pi i\langle m,y\rangle}{8}}\,dy,$$
where $\Delta_y$ is the Laplacian operator with respect to the variable $y$.
Taking absolute values, we obtain
$$|c_m|\leq\frac{16^{-(d-2)}\|(1-\Delta_y)^{d-1}H\|_{L^1([-8,8]^{d-2})}}{(1+\frac{\pi^2|m|^2}{64})^{d-1}}\leq C_d(1+|m|^2)^{-(d-1)},$$
where the constant $C_d$ depends only on the constant $d$ and the profile $\rho$ chosen. To sum these estimates, observe that if $m\in \Z^{d-2}$ with $|m|<2^{k+1}$ implies $|m_a|<2^{k+1}$ for every coordinate $a=1,\ldots,d-2$. Hence
$$\mathrm{Card}(\{m\in\mathbb Z^{d-2}:2^k\leq|m|<2^{k+1}\})\leq(2^{k+2}+1)^{d-2}\leq5^{d-2}2^{k(d-2)}.$$
This implies in particular that
\begin{equation*}
\begin{split}
\sum_{m\in\mathbb Z^{d-2}}|c_m|
&=|c_0|+\sum_{k=0}^{\infty}\sum_{2^k\leq|m|<2^{k+1}}|c_m|\leq|c_0|+C_d5^{d-2}\sum_{k=0}^{\infty}2^{k(d-2)}2^{-2k(d-1)}=|c_0|+C_d5^{d-2}\sum_{k=0}^{\infty}2^{-dk}<\infty.
\nonumber
\end{split}
\end{equation*}
Each term in \eqref{fourierrho} and \eqref{annular-separated-expansion} is bounded in absolute value by $|c_m|$. The summability just proved therefore gives absolute and uniform convergence of both series. Moreover,
\begin{equation*}
\begin{split}
&\sum_{m\in\mathbb Z^{d-2}}\|c_m\mathbb1_Qe^{\frac{\pi i\langle m,\cdot\rangle}{8r_j}}\|_{L^\infty}\|\mathbb1_{Q'}e^{-\frac{\pi i\langle m,\cdot\rangle}{8r_j}}\|_{L^\infty}=\sum_{m\in\mathbb Z^{d-2}}|c_m|\leq|c_0|+\frac{5^{d-2}C_d}{1-2^{-d}}.
\end{split}
\nonumber
\end{equation*}
Thus \eqref{annular-separated-expansion} represents the localized annular cutoff as a sum of products of functions of $x$ and $x'$, with a sum of supremum norm products bounded only in terms of $d$ and the fixed cutoff $\rho$, independently of $j,Q,Q'$.

\medskip
\spacedlowsmallcaps{The decomposition of $c_{Q,Q',j,e,\ell}$ into $c_{Q,Q',j,e,\ell}^{\mathfrak g}$ and $c_{Q,Q',j,e,\ell}^{\mathfrak t}$.}
By \eqref{defcdime} and \eqref{rappresentazione-con-cui-lavoriamo}, we have
\begin{equation}
\begin{split}
c_{Q,Q',j,e,\ell}(x,x',w)
&=i\lambda\mathbb1_Q(x)\mathbb1_{Q'}(x')\rho_{j,e}(x,x')\cdot\Big(\Delta_{h_j,e}\phi(x,x',w)-\Delta_{h_j,e}\phi(x,x',w-(\ell+1)h_je)\Big)\cdot \\
&\qquad\qquad\qquad\cdot \prod_{k=0}^{\ell+1}q_j(x,x',w-kh_je)\cdot \int_0^1e^{i\lambda((1-s)\Delta_{h_j,e}\phi(x,x',w-(\ell+1)h_je)+s\Delta_{h_j,e}\phi(x,x',w))}ds.
\end{split}
\label{remainder-coefficient-before-annular-expansion}
\end{equation}
Subtracting \eqref{separated-phase-series} at $w-(\ell+1)h_je$ from the same identity at $w$ gives
\begin{equation*}
\begin{split}
&\lambda\Big(\Delta_{h_j,e}\phi(x,x',w)-\Delta_{h_j,e}\phi(x,x',w-(\ell+1)h_je)\Big)=\sum_{\nu\in\Xi}u_\nu(x)v_\nu(x')\Big(b_\nu(w)-b_\nu(w-(\ell+1)h_je)\Big).
\end{split}
\nonumber
\end{equation*}
Since $\|b_\nu-b_\nu(\cdot-(\ell+1)h_je)\|_{L^\infty}\leq2\|b_\nu\|_{L^\infty}$, the bounds \eqref{residual-phase-coefficient-bound}, \eqref{translated-reciprocal-size} and \eqref{convex-phase-exponential-bound} show that this series, \eqref{translated-reciprocal-products} and \eqref{averaged-phase-exponential-expansion} are separated expansions of finite size, and so is \eqref{annular-separated-expansion}, by the summability of the $|c_m|$. The rule for products in \cref{separated-expansions} therefore allows us to multiply them. Substituting into \eqref{remainder-coefficient-before-annular-expansion}, we obtain
\begin{equation}
\begin{split}
&c_{Q,Q',j,e,\ell}(x,x',w)=\sum_{m\in\mathbb Z^{d-2}}\sum_{\nu\in\Xi}\sum_{\alpha_0,\ldots,\alpha_{\ell+1}\in\Xi''}\sum_{\gamma\in\{0\}\cup\Xi'}
\Bigg(\Big(c_m\mathbb1_Q(x)e^{\frac{\pi i\langle m,x\rangle}{8r_j}}u_\nu(x)U_\gamma(x)\prod_{k=0}^{\ell+1}\widetilde u_{\alpha_k}(x)\Big)\cdot\\
&\cdot\Big(\mathbb1_{Q'}(x')e^{-\frac{\pi i\langle m,x'\rangle}{8r_j}}v_\nu(x')V_\gamma(x')\prod_{k=0}^{\ell+1}\widetilde v_{\alpha_k}(x')\Big)\cdot \Big(i\big(b_\nu(w)-b_\nu(w-(\ell+1)h_je)\big)C_\gamma(w)\prod_{k=0}^{\ell+1}\widetilde b_{\alpha_k}(w-kh_je)\Big)\Bigg),
\label{expansionofc}
\end{split}
\end{equation}
where this series converges absolutely and uniformly. Let $\beta=(m,\nu,\alpha_0,\ldots,\alpha_{\ell+1},\gamma)$, where $m\in\mathbb Z^{d-2}$, $\nu\in\Xi$, $\alpha_0,\ldots,\alpha_{\ell+1}\in\Xi''$ and $\gamma\in\{0\}\cup\Xi'$, and define
\begin{equation*}
\begin{split}
\widehat U_\beta(x)&:=c_m\mathbb1_Q(x)e^{\frac{\pi i\langle m,x\rangle}{8r_j}}u_\nu(x)U_\gamma(x)\prod_{k=0}^{\ell+1}\widetilde u_{\alpha_k}(x),\\
\widehat V_\beta(x')&:=\mathbb1_{Q'}(x')e^{-\frac{\pi i\langle m,x'\rangle}{8r_j}}v_\nu(x')V_\gamma(x')\prod_{k=0}^{\ell+1}\widetilde v_{\alpha_k}(x'),\\
\widehat B_\beta(w)&:=i\big(b_\nu(w)-b_\nu(w-(\ell+1)h_je)\big)C_\gamma(w)\prod_{k=0}^{\ell+1}\widetilde b_{\alpha_k}(w-kh_je).
\end{split}
\nonumber
\end{equation*}
Thanks to the definition of the coefficients above and thanks to \eqref{expansionofc}, we have
\begin{equation}
c_{Q,Q',j,e,\ell}(x,x',w)
=\sum_\beta\widehat U_\beta(x)\widehat V_\beta(x')\widehat B_\beta(w).
\label{remainder-separated-expansion}
\end{equation}
The same bounds give $\sum_\beta\|\widehat U_\beta\|_{L^\infty}\|\widehat V_\beta\|_{L^\infty}\|\widehat B_\beta\|_{L^\infty}<\infty$, and $\widehat U_\beta$ and $\widehat V_\beta$ are supported in $Q$ and $Q'$, respectively.

We now split $c_{Q,Q',j,e,\ell}$ into $c_{Q,Q',j,e,\ell}^\mathfrak g$ and $c_{Q,Q',j,e,\ell}^\mathfrak t$. For every $\mu,j,e,\ell$, define
\begin{equation}
D_{\mu,j,e,\ell}(w):=\frac{g_\mu(w+h_je)-g_\mu(w)-g_\mu(w-\ell h_je)+g_\mu(w-(\ell+1)h_je)}{h_j}.
\label{definition-normalized-increment}
\end{equation}
The definitions of the coefficients above \eqref{separated-phase-series} gives the identity
$$b_{\mu,a}(w)-b_{\mu,a}(w-(\ell+1)h_je)=\kappa D_{\mu,j,e,\ell}(w),\qquad\text{for every } a=1,2,3.$$
Thus $\widehat B_\beta=0$ when the component $\nu$ of $\beta$ is $1$ or $2$, while for $\nu=(\mu,a)$ the definition of $\widehat B_\beta$ gives
\begin{equation}
\widehat B_\beta(w)=i\kappa D_{\mu,j,e,\ell}(w)C_\gamma(w)\prod_{k=0}^{\ell+1}\widetilde b_{\alpha_k}(w-kh_je).
\label{existing-coefficient-with-increment}
\end{equation}

Before estimating the coefficients, we specify their decomposition and estimates. We split $\widehat B_\beta$ directly, leaving $\widehat U_\beta$ and $\widehat V_\beta$ unchanged. The convergence of the resulting series will follow from the summability already established for \eqref{remainder-separated-expansion}.
Choose numbers $\tau_{\mu,j,e,\ell}\geq0$, independently of $Q,Q'$ and $z$. For $\beta=(m,(\mu,a),\alpha_0,\ldots,\alpha_{\ell+1},\gamma)$, set
$$\tau_\beta:=\tau_{\mu,j,e,\ell}\Big\|\kappa C_\gamma(\cdot)\prod_{k=0}^{\ell+1}\widetilde b_{\alpha_k}(\cdot-kh_je)\Big\|_{L^\infty},$$
and define $\widehat B_\beta^{\mathfrak g}(w):=\widehat B_\beta(w)\mathbb1_{\{|\widehat B_\beta(w)|\leq\tau_\beta\}}$ and $\widehat B_\beta^{\mathfrak t}(w):=\widehat B_\beta(w)\mathbb1_{\{|\widehat B_\beta(w)|>\tau_\beta\}}$. We then define
$$c_{Q,Q',j,e,\ell}^{\mathfrak g}(x,x',w):=\sum_\beta\widehat U_\beta(x)\widehat V_\beta(x')\widehat B_\beta^{\mathfrak g}(w),\qquad\text{and}\qquad c_{Q,Q',j,e,\ell}^{\mathfrak t}(x,x',w):=\sum_\beta\widehat U_\beta(x)\widehat V_\beta(x')\widehat B_\beta^{\mathfrak t}(w)$$
Since $|\widehat B_\beta^{\mathfrak g}|,|\widehat B_\beta^{\mathfrak t}|\leq|\widehat B_\beta|$, both series converge absolutely and uniformly by \eqref{remainder-separated-expansion}. Moreover, by definition
$$c_{Q,Q',j,e,\ell}=c_{Q,Q',j,e,\ell}^{\mathfrak g}+c_{Q,Q',j,e,\ell}^{\mathfrak t}.$$
Thus the decomposition is defined for every choice of the levels. We now prove \eqref{bounded-coefficient-condition} for the series defining $c_{Q,Q',j,e,\ell}^{\mathfrak g}$ and \eqref{integrable-coefficient-condition} for $c_{Q,Q',j,e,\ell}^{\mathfrak t}$, specifying $A_j$ and $B_j$ in terms of the levels $\tau_{\mu,j,e,\ell}$.

Fix $\beta=(m,(\mu,a),\alpha_0,\ldots,\alpha_{\ell+1},\gamma)$. We first estimate $\widehat B_\beta^{\mathfrak g}$ directly from its definition. For every $w$,
$$|\widehat B_\beta^{\mathfrak g}(w)|=|\widehat B_\beta(w)|\mathbb1_{\{|\widehat B_\beta(w)|\leq\tau_\beta\}}\leq\tau_\beta=\tau_{\mu,j,e,\ell}\Big\|\kappa C_\gamma(\cdot)\prod_{k=0}^{\ell+1}\widetilde b_{\alpha_k}(\cdot-kh_je)\Big\|_{L^\infty}.$$

To estimate $\widehat B_\beta^{\mathfrak t}$, we first show that it vanishes almost everywhere on the set where $|D_{\mu,j,e,\ell}|\leq\tau_{\mu,j,e,\ell}$. In fact, at almost every $w$ in this set, \eqref{existing-coefficient-with-increment} gives
\begin{equation*}
\begin{split}
|\widehat B_\beta(w)|=|D_{\mu,j,e,\ell}(w)|\Big|\kappa C_\gamma(w)\prod_{k=0}^{\ell+1}\widetilde b_{\alpha_k}(w-kh_je)\Big|
\leq\tau_{\mu,j,e,\ell}\Big\|\kappa C_\gamma(\cdot)\prod_{k=0}^{\ell+1}\widetilde b_{\alpha_k}(\cdot-kh_je)\Big\|_{L^\infty}
=\tau_\beta.
\nonumber
\end{split}
\end{equation*}
Thus, up to a null set we have
$$\{w:|\widehat B_\beta(w)|>\tau_\beta\}\subseteq\{w:|D_{\mu,j,e,\ell}(w)|>\tau_{\mu,j,e,\ell}\}.$$
Therefore, by \eqref{existing-coefficient-with-increment} once more, we obtain, for almost every $w$,
\begin{equation*}
\begin{split}
|\widehat B_\beta^{\mathfrak t}(w)|&=|\widehat B_\beta(w)|\mathbb1_{\{|\widehat B_\beta(w)|>\tau_\beta\}}\leq|\widehat B_\beta(w)|\mathbb1_{\{|D_{\mu,j,e,\ell}(w)|>\tau_{\mu,j,e,\ell}\}}\\
&=|D_{\mu,j,e,\ell}(w)|\mathbb1_{\{|D_{\mu,j,e,\ell}(w)|>\tau_{\mu,j,e,\ell}\}}
\Big|\kappa C_\gamma(w)\prod_{k=0}^{\ell+1}\widetilde b_{\alpha_k}(w-kh_je)\Big|\\
&\leq|D_{\mu,j,e,\ell}(w)|\mathbb1_{\{|D_{\mu,j,e,\ell}(w)|>\tau_{\mu,j,e,\ell}\}}
\Big\|\kappa C_\gamma(\cdot)\prod_{k=0}^{\ell+1}\widetilde b_{\alpha_k}(\cdot-kh_je)\Big\|_{L^\infty}.
\nonumber
\end{split}
\end{equation*}

If the norm on the right is zero, then $\widehat B_\beta=0$ and there is no contribution. It follows that
\begin{equation}
\begin{split}
\|\widehat B_\beta^{\mathfrak t}\|_{L^1}
&\leq\int |D_{\mu,j,e,\ell}(w)|\mathbb1_{\{|D_{\mu,j,e,\ell}(w)|>\tau_{\mu,j,e,\ell}\}}
\Big\|\kappa C_\gamma(\cdot)\prod_{k=0}^{\ell+1}\widetilde b_{\alpha_k}(\cdot-kh_je)\Big\|_{L^\infty}dw\\
&=\Big\|\kappa C_\gamma(\cdot)\prod_{k=0}^{\ell+1}\widetilde b_{\alpha_k}(\cdot-kh_je)\Big\|_{L^\infty}
\big\|D_{\mu,j,e,\ell}\mathbb1_{\{|D_{\mu,j,e,\ell}|>\tau_{\mu,j,e,\ell}\}}\big\|_{L^1}.
\end{split}
\label{direct-truncation-tail-bound}
\end{equation}
Summing up, we have shown that 
$$\|\widehat B_\beta^{\mathfrak g}\|_{L^\infty}
\leq\tau_{\mu,j,e,\ell}
\Big\|\kappa C_\gamma(\cdot)\prod_{k=0}^{\ell+1}\widetilde b_{\alpha_k}(\cdot-kh_je)\Big\|_{L^\infty}\,\,\,\text{and}\,\,\,\|\widehat B_\beta^{\mathfrak t}\|_{L^1}
\leq\big\|D_{\mu,j,e,\ell}\mathbb1_{\{|D_{\mu,j,e,\ell}|>\tau_{\mu,j,e,\ell}\}}\big\|_{L^1}
\Big\|\kappa C_\gamma(\cdot)\prod_{k=0}^{\ell+1}\widetilde b_{\alpha_k}(\cdot-kh_je)\Big\|_{L^\infty}.$$
Let us now claim that 
\begin{equation}
\begin{split}
&\sum_{a=1}^3\sum_{m\in\mathbb Z^{d-2}}\sum_{\gamma\in\{0\}\cup\Xi'}\sum_{\alpha_0,\ldots,\alpha_{\ell+1}\in\Xi''}
\|\widehat U_\beta\|_{L^\infty}\|\widehat V_\beta\|_{L^\infty}
\Big\|\kappa C_\gamma(\cdot)\prod_{k=0}^{\ell+1}\widetilde b_{\alpha_k}(\cdot-kh_je)\Big\|_{L^\infty}
\leq C_d\Lip(f_\mu),
\end{split}
\label{coefficient-sum-needed-for-truncation}
\end{equation}
where $\beta=(m,(\mu,a),\alpha_0,\ldots,\alpha_{\ell+1},\gamma)$ and $C_d$ is independent of $j,e,\ell,Q,Q'$.

If \eqref{coefficient-sum-needed-for-truncation} holds, we would have that 
\begin{equation}
    \begin{split}
        \sum_\beta\|\widehat U_\beta\|_{L^\infty}\|\widehat V_\beta\|_{L^\infty}\|\widehat B_\beta^{\mathfrak g}\|_{L^\infty}\leq   \sum_\beta\tau_{\mu,j,e,\ell}\|\widehat U_\beta\|_{L^\infty}\|\widehat V_\beta\|_{L^\infty}\|\Big\|\kappa C_\gamma(\cdot)\prod_{k=0}^{\ell+1}\widetilde b_{\alpha_k}(\cdot-kh_je)\Big\|_{L^\infty}\leq C_d\sum_\mu\Lip(f_\mu)\tau_{\mu,j,e,\ell}.
        \label{stimasommacappelliconf}
    \end{split}
\end{equation}
On the other hand, let us further notice that by Tonelli's theorem we have
\begin{equation}
\begin{split}
\operatorname*{ess\,sup}_{(x,x')\in Q\times Q'}\int|c_{Q,Q',j,e,\ell}^{\mathfrak t}&(x,x',w)|dw
\leq\sum_\beta\|\widehat U_\beta\|_{L^\infty}\|\widehat V_\beta\|_{L^\infty}\|\widehat B_\beta^{\mathfrak t}\|_{L^1}\\
&\leq \sum_\beta\|\widehat U_\beta\|_{L^\infty}\|\widehat V_\beta\|_{L^\infty}\Big\|\kappa C_\gamma(\cdot)\prod_{k=0}^{\ell+1}\widetilde b_{\alpha_k}(\cdot-kh_je)\Big\|_{L^\infty}
\big\|D_{\mu,j,e,\ell}\mathbb1_{\{|D_{\mu,j,e,\ell}|>\tau_{\mu,j,e,\ell}\}}\big\|_{L^1}\\
&\leq C_d\sum_\mu\Lip(f_\mu)\big\|D_{\mu,j,e,\ell}\mathbb1_{\{|D_{\mu,j,e,\ell}|>\tau_{\mu,j,e,\ell}\}}\big\|_{L^1}.
\label{ctstimax}
\end{split}
\end{equation}
These are the two estimates from which we will define $A_j$ and $B_j$, as we shall see after we prove \eqref{coefficient-sum-needed-for-truncation}, keeping $\mu$ fixed and estimating the other factors before summing over $a,m,\gamma,\alpha_0,\ldots,\alpha_{\ell+1}$.

Let us begin by observing that thanks to the definition of  $u_{\mu,a}$ and $v_{\mu,a}$ above \eqref{separated-phase-series}, we have 
\begin{equation}
\sum_{a=1}^3\|u_{\mu,a}\|_{L^\infty}\|v_{\mu,a}\|_{L^\infty}\leq\Big(\frac1{128}+5+\frac1{128}\Big)\Lip(f_\mu)\leq6\Lip(f_\mu).
   \label{sommanormeLinfty} 
\end{equation}
We prove the estimate first for a single term and then sum over its indices. Fix $\mu$ and write
$$
\beta=(m,(\mu,a),\alpha_0,\ldots,\alpha_{\ell+1},\gamma),
$$
with $a\in\{1,2,3\}$. We estimate only the factors multiplying $D_{\mu,j,e,\ell}$ in \eqref{existing-coefficient-with-increment}.

The definitions of $\widehat U_\beta$ and $\widehat V_\beta$, we have
\begin{equation*}
\begin{split}
\|\widehat U_\beta\|_{L^\infty}\leq |c_m|\|u_{\mu,a}\|_{L^\infty}\|U_\gamma\|_{L^\infty}\prod_{k=0}^{\ell+1}\|\widetilde u_{\alpha_k}\|_{L^\infty},\qquad\text{and}\qquad\|\widehat V_\beta\|_{L^\infty}\leq \|v_{\mu,a}\|_{L^\infty}\|V_\gamma\|_{L^\infty}\prod_{k=0}^{\ell+1}\|\widetilde v_{\alpha_k}\|_{L^\infty}.
\nonumber
\end{split}
\end{equation*}
Since translations preserve the $L^\infty$ norm, we also have
$$\Big\|\kappa C_\gamma(\cdot)\prod_{k=0}^{\ell+1}\widetilde b_{\alpha_k}(\cdot-kh_je)\Big\|_{L^\infty}\leq\kappa\|C_\gamma\|_{L^\infty}\prod_{k=0}^{\ell+1}\|\widetilde b_{\alpha_k}\|_{L^\infty}.$$
Multiplying these three inequalities gives
\begin{equation*}
\begin{split}
&\|\widehat U_\beta\|_{L^\infty}\|\widehat V_\beta\|_{L^\infty}\Big\|\kappa C_\gamma(\cdot)\prod_{k=0}^{\ell+1}\widetilde b_{\alpha_k}(\cdot-kh_je)\Big\|_{L^\infty}\\
&\qquad\qquad\qquad\leq\kappa |c_m|
\Big(\|u_{\mu,a}\|_{L^\infty}\|v_{\mu,a}\|_{L^\infty}\Big)
\Big(\|U_\gamma\|_{L^\infty}\|V_\gamma\|_{L^\infty}\|C_\gamma\|_{L^\infty}\Big)\prod_{k=0}^{\ell+1}
\Big(\|\widetilde u_{\alpha_k}\|_{L^\infty}\|\widetilde v_{\alpha_k}\|_{L^\infty}\|\widetilde b_{\alpha_k}\|_{L^\infty}\Big).
\end{split}
\nonumber
\end{equation*}
By \eqref{translated-reciprocal-size}, the sum over $\alpha_0,\ldots,\alpha_{\ell+1}\in\Xi''$ of the last products is at most $(\frac{2\pi\sqrt{d-2}}{\kappa})^{\ell+2}$.
For fixed $\mu,m,a$, write $\beta=(m,(\mu,a),\alpha_0,\ldots,\alpha_{\ell+1},\gamma)$. The estimate for a single term and the preceding computation of the sums over $\alpha_0,\ldots,\alpha_{\ell+1}$ give
\begin{equation*}
\begin{split}
&\sum_{\gamma\in\{0\}\cup\Xi'}\sum_{\alpha_0,\ldots,\alpha_{\ell+1}\in\Xi''}
\|\widehat U_\beta\|_{L^\infty}\|\widehat V_\beta\|_{L^\infty}
\Big\|\kappa C_\gamma(\cdot)\prod_{k=0}^{\ell+1}\widetilde b_{\alpha_k}(\cdot-kh_je)\Big\|_{L^\infty}\\
&\quad\leq\kappa|c_m|\|u_{\mu,a}\|_{L^\infty}\|v_{\mu,a}\|_{L^\infty}
\sum_{\gamma\in\{0\}\cup\Xi'}\Bigg(
\|U_\gamma\|_{L^\infty}\|V_\gamma\|_{L^\infty}\|C_\gamma\|_{L^\infty}
\sum_{\alpha_0,\ldots,\alpha_{\ell+1}\in\Xi''}
\prod_{k=0}^{\ell+1}
\Big(\|\widetilde u_{\alpha_k}\|_{L^\infty}\|\widetilde v_{\alpha_k}\|_{L^\infty}\|\widetilde b_{\alpha_k}\|_{L^\infty}\Big)
\Bigg)\\
&\qquad\qquad\quad\leq\kappa|c_m|\|u_{\mu,a}\|_{L^\infty}\|v_{\mu,a}\|_{L^\infty}
\Big(\frac{2\pi\sqrt{d-2}}{\kappa}\Big)^{\ell+2}
\sum_{\gamma\in\{0\}\cup\Xi'}\|U_\gamma\|_{L^\infty}\|V_\gamma\|_{L^\infty}\|C_\gamma\|_{L^\infty}\\
&\qquad\qquad\qquad\quad\overset{\eqref{convex-phase-exponential-bound}}{\leq}
\kappa|c_m|\|u_{\mu,a}\|_{L^\infty}\|v_{\mu,a}\|_{L^\infty}
\Big(\frac{2\pi\sqrt{d-2}}{\kappa}\Big)^{\ell+2}e^{\frac{\kappa}{16\sqrt{d-2}}}.
\end{split}
\nonumber
\end{equation*}
We next sum this inequality over $m$, keeping $\mu,a$ fixed:
\begin{equation}
\begin{split}
&\sum_{m\in\mathbb Z^{d-2}}\sum_{\gamma\in\{0\}\cup\Xi'}\sum_{\alpha_0,\ldots,\alpha_{\ell+1}\in\Xi''}
\|\widehat U_\beta\|_{L^\infty}\|\widehat V_\beta\|_{L^\infty}
\Big\|\kappa C_\gamma(\cdot)\prod_{k=0}^{\ell+1}\widetilde b_{\alpha_k}(\cdot-kh_je)\Big\|_{L^\infty}\\
&\quad\leq\sum_{m\in\mathbb Z^{d-2}}
\kappa|c_m|\|u_{\mu,a}\|_{L^\infty}\|v_{\mu,a}\|_{L^\infty}
\Big(\frac{2\pi\sqrt{d-2}}{\kappa}\Big)^{\ell+2}e^{\frac{\kappa}{16\sqrt{d-2}}}\\
&\quad=\kappa e^{\frac{\kappa}{16\sqrt{d-2}}}
\Big(\frac{2\pi\sqrt{d-2}}{\kappa}\Big)^{\ell+2}
\Big(\sum_{m\in\mathbb Z^{d-2}}|c_m|\Big)
\|u_{\mu,a}\|_{L^\infty}\|v_{\mu,a}\|_{L^\infty}.
\label{stimaintermediacoicappucci}
\end{split}
\end{equation}
Summing the preceding estimate over $a=1,2,3$, by \eqref{sommanormeLinfty} and \eqref{stimaintermediacoicappucci}, we therefore obtain
\begin{equation*}
\begin{split}
\sum_{a=1}^3\sum_{m\in\mathbb Z^{d-2}}&\sum_{\gamma\in\{0\}\cup\Xi'}\sum_{\alpha_0,\ldots,\alpha_{\ell+1}\in\Xi''}
\|\widehat U_\beta\|_{L^\infty}\|\widehat V_\beta\|_{L^\infty}
\Big\|\kappa C_\gamma(\cdot)\prod_{k=0}^{\ell+1}\widetilde b_{\alpha_k}(\cdot-kh_je)\Big\|_{L^\infty}\\
&\qquad\qquad\qquad\qquad\qquad\qquad\leq6\kappa e^{\frac{\kappa}{16\sqrt{d-2}}}
\Big(\frac{2\pi\sqrt{d-2}}{\kappa}\Big)^{\ell+2}
\Big(\sum_{m\in\mathbb Z^{d-2}}|c_m|\Big)\Lip(f_\mu).
\end{split}
\nonumber
\end{equation*}
This proves \eqref{coefficient-sum-needed-for-truncation} and therefore also \eqref{stimasommacappelliconf} and \eqref{ctstimax}. We can now define the sequences announced above. For every scale $j$ for which at least one of the families $\mathcal P_{j,e}$  is nonempty, set
\begin{equation}
A_j:=C_d\max_{\substack{e\in\mathscr E^\pm\\0\leq\ell\leq d-2}}2^\ell\sum_\mu\Lip(f_\mu)\tau_{\mu,j,e,\ell},
\label{definition-Aj}
\end{equation}
and
\begin{equation}
B_j:=C_d\|\eta_0\|_{L^\infty}^2\|\eta_1\|_{L^\infty}^2\Big(\frac{r_j}{64\sqrt{d-2}}\Big)^{d-2}\max_{\substack{e\in\mathscr E^\pm\\0\leq\ell\leq d-2}}2^\ell\sum_\mu\Lip(f_\mu)\big\|D_{\mu,j,e,\ell}\mathbb1_{\{|D_{\mu,j,e,\ell}|>\tau_{\mu,j,e,\ell}\}}\big\|_{L^1}.
\label{definition-Bj}
\end{equation}
If the family $\mathcal P_{j,e}$ introduced on page \pageref{familyPje} is empty for every $e\in\mathscr E^\pm$, we set $A_j=B_j=0$. 
Indeed, for every $e\in\mathscr E^\pm$ and $0\leq\ell\leq d-2$, \eqref{stimasommacappelliconf} and \eqref{definition-Aj} give
$$\sum_\beta\|\widehat U_\beta\|_{L^\infty}\|\widehat V_\beta\|_{L^\infty}\|\widehat B_\beta^{\mathfrak g}\|_{L^\infty}\leq C_d\sum_\mu\Lip(f_\mu)\tau_{\mu,j,e,\ell}\leq2^{-\ell}A_j.$$
Thus the series defining $c_{Q,Q',j,e,\ell}^{\mathfrak g}$ satisfies \eqref{bounded-coefficient-condition}. Moreover, the cubes have side-length $\frac{r_j}{64\sqrt{d-2}}$ and hence \eqref{ctstimax} and \eqref{definition-Bj} imply
\begin{equation*}
\begin{split}
&2^\ell\|\eta_0\|_{L^\infty}^2\|\eta_1\|_{L^\infty}^2\Leb^{d-2}(Q)^{\frac{1}{2}}\Leb^{d-2}(Q')^{\frac{1}{2}}\operatorname*{ess\,sup}_{(x,x')\in Q\times Q'}\int|c_{Q,Q',j,e,\ell}^{\mathfrak t}(x,x',w)|\,dw\leq B_j.
\nonumber
\end{split}
\end{equation*}
Therefore \eqref{integrable-coefficient-condition} also holds. The abstract argument following \eqref{remainder-coefficient-splitting} now applies and gives \eqref{required-coefficient-estimate}, \eqref{required-remainder-representation}, and hence there holds
$\|\mathcal K_{e,j,z}^{Q,Q',\ell}\|_{2\to2}\leq A_jN_\lambda+B_j$.
Combining this estimate with \eqref{cube-pair-summation} and \eqref{required-principal-estimate} gives again \eqref{operator-norm-before-absorption}. It remains to choose the levels $\tau_{\mu,j,e,\ell}$ so that
$C\sum_{j\geq0}A_j\leq\frac12$
and to estimate $\sum_{j\geq0}B_j$. The fact that such a choice is possible depends on the specific structure of the functions $f_\mu$ and $g_\mu$ in the expansion \eqref{target-local-phase}, and this will be done in Section~\ref{sss:regularity}.

\subsubsection{The representation of \texorpdfstring{$\Phi_\delta$}{Phi}}\label{sss:representation}

We now show that the phase constructed in Section~\ref{sss:normal} has the form required in Section~\ref{sss:model}. Let $S:=D^2H(z_0)$. The specific expression of $H$ gives
$$D^2H(z_0)[\zeta,\zeta]=\frac{|\zeta|^2}{H(z_0)}-\frac{\langle z_0,\zeta\rangle^2}{H(z_0)^3}\geq\frac{|\zeta|^2}{H(z_0)^3}\qquad\text{for every }\zeta\in \R^{d-1}.$$
In particular, $S$ is positive definite, $\|S\|\leq1$, and $\|S^{-1}\|\leq H(z_0)^3$. These bounds are uniform because $z_0$ belongs to a bounded set determined by $d$ and $C_0$.

By \eqref{strutturaphidelta} and the fundamental theorem of calculus,
\begin{equation}\begin{split}
\Phi_\delta(x,y)&=-\int_0^1\int_0^1D^2H(z_0+\delta(tU(x)-sW(y)))[U(x),W(y)]\,dt\,ds\\
&=-\langle U(x),SW(y)\rangle-\mathcal R_{\delta,z_0}(U(x),W(y)),
\end{split}\label{phase-and-smooth-remainder}\end{equation}
where
$$\mathcal R_{\delta,z_0}(z,z'):=\delta\int_0^1\int_0^1\int_0^1D^3H(z_0+\theta\delta(tz-sz'))[tz-sz',z,z']d\theta dtds.$$
With $L$ as fixed at the beginning of Section~\ref{sss:normal}, define
$\mathcal S_0:=(-2L,2L)^{2(d-1)}$, and $\mathcal S_1:=(-4L,4L)^{2(d-1)}$.
The rescaled cubes are contained in $(-\frac{1}{2},\frac{1}{2})^{d-2}$. Since $P_\delta(0)=Q_\delta(0)=0$ and their Lipschitz constants are at most $\varepsilon_*\leq1$, we have
$$|P_\delta(x)|\leq|x|,\qquad |Q_\delta(y)|\leq|y|.$$
Thus $(U(x),W(y))\in(-L,L)^{2(d-1)}\subset\mathcal S_0$ on these cubes.
For every $(z,z')\in\mathcal S_1$ and $s,t,\theta\in[0,1]$, there holds 
$$|z_0+\theta\delta(tz-sz')-z_0|\leq|z|+|z'|\leq8L\sqrt{d-1}.$$
All derivatives of $H$ are bounded on this fixed ball. Differentiating the formula for $\mathcal R_{\delta,z_0}$ preserves the factor $\delta$ preceding the integral. Each derivative falling on $D^3H$ introduces another factor bounded by $\delta$, and the other factors are polynomials of degree at most three in $z,z'$. Consequently,
$$\|\mathcal R_{\delta,z_0}\|_{C^{2d+2}(\mathcal S_1)}\leq C\delta.$$
The constant is independent of $\delta$ and of the chosen pair of cubes. Linearity also gives
$\mathcal R_{\delta,z_0}(z,0)=\mathcal R_{\delta,z_0}(0,z')=0$.

Apply \cref{fourier-expansions} in the two groups of variables $z,z'\in\mathbb R^{d-1}$, with $N=d+1$ and coefficient weight $(1+|k|)^3$. These choices are admissible because $2d+2>2(d-1)+3$. On $\mathcal S_0$ we obtain
$$\mathcal R_{\delta,z_0}(z,z')=\sum_{k\in\mathbb Z^{2(d-1)}}c_kE_{k,1}(z)E_{k,2}(z'),\qquad \sum_k(1+|k|)^3|c_k|\leq C\delta.$$
Moreover, we have $|E_{k,1}|=|E_{k,2}|=1$ and $\Lip(E_{k,1})+\Lip(E_{k,2})\leq C(1+|k|)$. 
Define
$$A_k(x):=c_kE_{k,1}(U(x)),\qquad B_k(y):=E_{k,2}(W(y)).$$
The maps $U,W$ have Lipschitz constants at most $2$, so $\Lip(A_k)\leq C|c_k|(1+|k|)$ and $\Lip(B_k)\leq C(1+|k|)$.

The series converges absolutely and uniformly because $|A_k(x)B_k(y)|\leq|c_k|$. Since $\mathcal R_{\delta,z_0}$ is real-valued, we may take real parts term by term. Set
$$F_{k,1}:=\mathrm{Re}\,A_k,\qquad G_{k,1}:=\mathrm{Re}\,B_k,\qquad F_{k,2}:=-\mathrm{Im}\,A_k,\qquad G_{k,2}:=\mathrm{Im}\,B_k.$$
The identity $\mathrm{Re}(A_kB_k)=F_{k,1}G_{k,1}+F_{k,2}G_{k,2}$ gives
\begin{equation}
\mathcal R_{\delta,z_0}(U(x),W(y))=\sum_k\sum_{\nu=1}^2F_{k,\nu}(x)G_{k,\nu}(y),\qquad \sum_k\sum_{\nu=1}^2\Lip(F_{k,\nu})\Lip(G_{k,\nu})\leq C\delta.
\label{smooth-remainder-expansion}
\end{equation}
Since $W(0)=0$ and $\mathcal R_{\delta,z_0}(U(x),0)=0$, subtracting this expansion at $y=0$ yields
$$\mathcal R_{\delta,z_0}(U(x),W(y))=\sum_k\sum_{\nu=1}^2F_{k,\nu}(x)\big(G_{k,\nu}(y)-G_{k,\nu}(0)\big).$$
Every factor depending on $y$ now vanishes at zero.

We now make a change of variable in the second argument of 
\eqref{phase-and-smooth-remainder}. Let $B:=\pi S\pi^T$, where $\pi$ is the projection on the first $d-2$ coordinates. Applying the bounds for $S$ to $(\xi,0)\in\mathbb R^{d-1}$ gives
$$\frac{|\xi|^2}{H(z_0)^3}\leq\langle S(\xi,0),(\xi,0)\rangle=\langle B\xi,\xi\rangle\leq|\xi|^2.$$
Thus $B$ is positive definite and $\|B^{-1}\|\leq H(z_0)^3$. This bound is uniform because $z_0$ belongs to a bounded set determined by $d$ and $C_0$.
Make the change of variables $w=-By$ and put $\widetilde Q(w):=Q_\delta(-B^{-1}w)$. 
 With $\pi$ denoting the projection onto the first $d-2$ coordinates, we have $B\xi=\pi(S(\xi,0))$. Writing $U(x)=(x,0)+P_\delta(x)(0,1)$ and $W(-B^{-1}w)=-(B^{-1}w,0)+\widetilde Q(w)(0,1)$, bilinearity gives
\begin{equation}\begin{split}
-\langle U(x),SW(-B^{-1}w)\rangle=&\langle(x,0),S(B^{-1}w,0)\rangle-\widetilde Q(w)\langle(x,0),S(0,1)\rangle\\
&+P_\delta(x)\langle(0,1),S(B^{-1}w,0)\rangle-P_\delta(x)\widetilde Q(w)\langle(0,1),S(0,1)\rangle.
\end{split}\end{equation}
By the definition of $B$ and the symmetry of $S$, respectively, we have 
\begin{equation}\begin{split}
\langle(x,0),S(B^{-1}w,0)\rangle&=\langle x,\pi(S(B^{-1}w,0))\rangle=\langle x,BB^{-1}w\rangle=\langle x,w\rangle,\\
\langle(0,1),S(B^{-1}w,0)\rangle&=\langle S(0,1),(B^{-1}w,0)\rangle=\langle\pi(S(0,1)),B^{-1}w\rangle.
\end{split}\end{equation}
Moreover, $\langle(x,0),S(0,1)\rangle=\langle x,\pi(S(0,1))\rangle$. Substituting these identities and the series for the remainder into \eqref{phase-and-smooth-remainder}, we obtain
\begin{equation}
\begin{split}
\Phi_\delta(x,-B^{-1}w)=&\langle x,w\rangle-\underbrace{\langle x,\pi(S(0,1))\rangle\widetilde Q(w)}_{I_a(x)I_b(w)}+\underbrace{P_\delta(x)\langle\pi(S(0,1)),B^{-1}w\rangle}_{II_a(x)II_b(w)}\\
&-\underbrace{\langle(0,1),S(0,1)\rangle P_\delta(x)\widetilde Q(w)}_{III_a(x)III_b(w)}-\underbrace{\sum_k\sum_{\nu=1}^2F_{k,\nu}(x)(G_{k,\nu}(-B^{-1}w)-G_{k,\nu}(0))}_{IV_a(x)IV_b(w)}.
\end{split}
\label{normalized-phase-expansion}
\end{equation}

We check now that $\Phi_\delta$ verifies the second item in \eqref{target-local-phase}. Recall that by construction
$$\Lip(P_\delta)+\Lip(Q_\delta)\leq\varepsilon_*,\qquad \Lip(\widetilde Q)\leq\|B^{-1}\|\Lip(Q_\delta).$$
Let us now check items (I) to (IV) separately. Since $I_a$ and $II_b$ are linear with gradients $\pi(S(0,1))$ and $B^{-T}\pi(S(0,1))$, whose norms are at most $\|S\|$ and $\|S\|\|B^{-1}\|$ respectively, while $\|S\|\leq1$ and $\|B^{-1}\|\leq H(z_0)^3\leq A^3$ with $A$ as in Section~\ref{sss:mixed}, we have
$$\Lip(I_a)\Lip(I_b)\leq \|S\|\|B^{-1}\|\Lip(Q_\delta)\leq A^3\varepsilon_*,$$
$$\Lip(II_a)\Lip(II_b)\leq \Lip(P_\delta)\|S\| \|B^{-1}\|\leq A^3\varepsilon_*,$$
$$\Lip(III_a)\Lip(III_b)\leq \|S\|\Lip(P_\delta)\|B^{-1}\|\Lip(Q_\delta) \leq A^3\varepsilon_*^2.$$
Finally, we estimate the contribution given by the fourth term. 
\begin{equation}\begin{split}
\sum_k\sum_{\nu=1}^2\Lip(F_{k,\nu})\Lip(G_{k,\nu}(-B^{-1}\cdot)-G_{k,\nu}(0))\leq\|B^{-1}\|\sum_k\sum_{\nu=1}^2\Lip(F_{k,\nu})\Lip(G_{k,\nu})\leq C\delta.
\end{split}\end{equation}
Adding the four bounds, the sum of the products of the Lipschitz constants of the factors following $\langle x,w\rangle$ in \eqref{normalized-phase-expansion} is at most $C(d,C_0)(\varepsilon_*+\varepsilon_*^2+\delta)$, which is smaller than the number $\varepsilon$ in \eqref{target-local-phase} once $\varepsilon_*$, and then $\delta$, are chosen sufficiently small. Using the factors in \eqref{normalized-phase-expansion}, define
\begin{equation*}
\begin{split}
&f_1(x):=-I_a(x),\qquad g_1(w):=I_b(w)=\widetilde Q(w),\\
&f_2(x):=II_a(x),\qquad g_2(w):=II_b(w),\\
&f_3(x):=-III_a(x),\qquad g_3(w):=III_b(w)=\widetilde Q(w).
\end{split}
\nonumber
\end{equation*}
For every $k\in\mathbb Z^{2(d-1)}$ and $\nu\in\{1,2\}$, also define
\begin{equation}
    f_{(k,\nu)}(x):=-F_{k,\nu}(x),\qquad g_{(k,\nu)}(w):=G_{k,\nu}(-B^{-1}w)-G_{k,\nu}(0).
    \label{definizionefkmu}
\end{equation}
The index $\mu$ therefore runs over $1,2,3$ and all pairs $(k,\nu)$. With these definitions, the sum $\sum_\mu f_\mu(x)g_\mu(w)$ is exactly the sum of the terms following $\langle x,w\rangle$ in \eqref{normalized-phase-expansion} and with these definitions \eqref{target-local-phase} is satisfied.

Let us now summarize what we have obtained. Set
$$\eta_0(x):=\chi_{0,a}(x_0+\delta x),\qquad \eta_1(w):=\chi_{1,b}(y_0-\delta B^{-1}w),$$
and
$$\phi(x,w):=\langle x,w\rangle+\sum_\mu f_\mu(x)g_\mu(w).$$
The factors $f_\mu,g_\mu$ are defined on the whole of $\mathbb R^{d-2}$, so that $\phi$ is a global phase of the form \eqref{target-local-phase}, and by \eqref{normalized-phase-expansion} it agrees with $\Phi_\delta(x,-B^{-1}w)$ whenever $x\in\supp\eta_0$ and $w\in\supp\eta_1$, since there $(U(x),W(-B^{-1}w))\in\mathcal S_0$. This global phase is the one used in the translated family \eqref{translated-operator-family} and in the supremum defining $N_\lambda$, and the increment estimates of Section~\ref{sss:regularity} apply to its factors on all of $\mathbb R^{d-2}$, so it never has to agree with $\Phi_\delta(x,-B^{-1}w)$ outside the supports of the cutoffs. Thanks to the discussion above, the change of variables $w=-By$ gives
\begin{equation*}
    \begin{split}
  \widetilde{V}_\lambda^{a,b}f(x)&=\chi_{0,a}(x_0+\delta x)\int_I e^{i\lambda\delta^2\Phi_\delta(x,y)}\chi_{1,b}(y_0+\delta y)f(y)dy\\
  &=|\det B|^{-1}\chi_{0,a}(x_0+\delta x)\int_{\mathbb R^{d-2}} e^{i\lambda\delta^2\Phi_\delta(x,-B^{-1}w)}\chi_{1,b}(y_0-\delta B^{-1}w)f(-B^{-1}w)dw \\
  &=|\det B|^{-1}\eta_0(x)\int_{\mathbb R^{d-2}} e^{i\lambda\delta^2\phi(x,w)}\eta_1(w)f(-B^{-1}w)dw\\
  &=|\det B|^{-1}T^{a,b}_{\lambda\delta^2}(f(-B^{-1}\cdot))(x),
  \nonumber
    \end{split}
\end{equation*}
where $T^{a,b}_{\lambda\delta^2}$ is the operator of Section~\ref{sss:model} with the cutoffs $\eta_0,\eta_1$, the phase $\phi$ and $\lambda\delta^2$ in place of $\lambda$. The integration extends over $\mathbb R^{d-2}$ because $\chi_{1,b}(y_0+\delta\cdot)$ is supported in $(-\frac{1}{2},\frac{1}{2})^{d-2}\subset I$. Consequently,
\begin{equation}
\begin{split}
\frac{\|\widetilde V_\lambda^{a,b}f\|_2}{\|f\|_2}&=|\det B|^{-1}\frac{\|T_{\lambda\delta^2}^{a,b}(f(-B^{-1}\cdot))\|_2}{|\det B|^{-\frac{1}{2}}\|f(-B^{-1}\cdot)\|_2}=|\det B|^{-\frac{1}{2}}\frac{\|T_{\lambda\delta^2}^{a,b}(f(-B^{-1}\cdot))\|_2}{\|f(-B^{-1}\cdot)\|_2}.
\end{split}
\end{equation}
We thus conclude, thanks to \eqref{identitadiV}, that 
$$\|V^{a,b}_\lambda\|_{2\to 2}=\delta^{d-2}\|\widetilde V_\lambda^{a,b}\|_{2\to2}=|\det B|^{-\frac{1}{2}}\delta^{d-2}\|T_{\lambda\delta^2}^{a,b}\|_{2\to2}.
$$
Apply \eqref{operator-norm-after-absorption} to $T_{\lambda\delta^2}^{a,b}$, replacing $\lambda$ by $\lambda\delta^2$ in the definitions of $A_j$ and $B_j$. If $\lambda\delta^2\geq1$ and $C\sum_jA_j\leq\frac{1}{2}$, the preceding identity for the operator norms gives
\begin{equation}
\begin{split}
\|V_\lambda^{a,b}\|_{2\to2}^2&\leq\frac{2C}{|\det B|}\Big(\Big(M^{d-2}+\frac1M\Big)\frac{\delta^{2(d-2)}}{(\lambda\delta^2)^{d-2}}+\delta^{2(d-2)}\sum_{j\geq0}B_j\Big)\\
&=\frac{2C}{|\det B|}\Big(\Big(M^{d-2}+\frac1M\Big)\lambda^{-(d-2)}+\delta^{2(d-2)}\sum_{j\geq0}B_j\Big).
\end{split}
\label{localized-operator-final-reduction}
\end{equation}

\subsubsection{The regularity estimates}\label{sss:regularity}

We now choose the levels $\tau_{\mu,j,e,\ell}$ in the decomposition already constructed in Section~\ref{sss:model}. We will verify \eqref{sommaA_j} and estimate $\sum_jB_j$, so that \eqref{localized-operator-final-reduction} gives the conclusions of the proposition. First, if $1\leq\lambda<\delta^{-2}$, Cauchy-Schwarz gives
$$\|V_\lambda\|_{2\to2}^2\leq\|\chi_0\|_{L^2}^2\|\chi_1\|_{L^2}^2\leq\delta^{-2(d-2)}\|\chi_0\|_{L^2}^2\|\chi_1\|_{L^2}^2\lambda^{-(d-2)}.$$
Since $\delta$ has already been fixed, this proves the required estimates in this range. We henceforth assume $\lambda\delta^2\geq1$ and fix one pair $(a,b)$. In the definitions of Section~\ref{sss:model}, we replace $\lambda$ by $\lambda\delta^2$, so that in particular
$$r_j=\frac{2^jM}{\lambda\delta^2},\qquad h_j=\frac{\kappa}{2^jM}.$$
All estimates below are uniform in $(a,b)$ with constants independent of the vector $z\in\mathbb R^{d-2}$ appearing in the factor $\eta_1(w-z)$ of $T_{\lambda,z}$.

\smallskip

\spacedlowsmallcaps{Estimates for the derivative increments.}
In this paragraph we prove that the following estimate holds
\begin{equation}
\sum_\mu\Lip(f_\mu)\|D_{\mu,j,e,\ell}\|_{L^\infty}\leq C\delta h_j+C\Big\|\int_0^1|\nabla\widetilde Q(w+th_je)-\nabla\widetilde Q(w+(t-\ell-1)h_je)|\,dt\Big\|_{L^\infty}.
\label{weighted-derivative-increment-bound}
\end{equation}
By \eqref{definition-normalized-increment} and the fundamental theorem of calculus, we have
\begin{equation*}
\begin{split}
D_{\mu,j,e,\ell}(w)
&=\frac{g_\mu(w+h_je)-g_\mu(w)}{h_j}-\frac{g_\mu(w-\ell h_je)-g_\mu(w-(\ell+1)h_je)}{h_j}\\
&=\int_0^1\Big(\partial_eg_\mu(w+th_je)-\partial_eg_\mu(w+(t-\ell-1)h_je)\Big)\,dt.
\end{split}
\nonumber
\end{equation*}
For the terms I and III in \eqref{normalized-phase-expansion}, we have $g_1=g_3=\widetilde Q$. Since $|e|=1$, it follows that
\begin{equation}
    |D_{1,j,e,\ell}(w)|=|D_{3,j,e,\ell}(w)|\leq\int_0^1|\nabla\widetilde Q(w+th_je)-\nabla\widetilde Q(w+(t-\ell-1)h_je)|\,dt.
    \label{stimaprimitermini}
\end{equation}
For the term $II$, we have
$\partial_eg_2(w)=\langle\pi(S(0,1)),B^{-1}e\rangle$ and since in this case $\partial_eg_2(w)$ is constant, we have  $D_{2,j,e,\ell}=0$.

Finally, for the last type of coefficient which is type IV, the definitions preceding \eqref{smooth-remainder-expansion} give
$$G_{k,1}=\operatorname{Re}B_k,\qquad G_{k,2}=\operatorname{Im}B_k,\qquad B_k=E_{k,2}\circ W.$$
Consequently, $g_{(k,\nu)}$ is the real or imaginary part of $B_k\circ(-B^{-1})-B_k(0)$, and
$$B_k(-B^{-1}w)=E_{k,2}(W(-B^{-1}w))=E_{k,2}(-B^{-1}w,Q_\delta(-B^{-1}w))=E_{k,2}(-B^{-1}w,\widetilde Q(w)).$$
We estimate its derivative increments by differentiating this last expression. 
By the definition $\widetilde Q(w)=Q_\delta(-B^{-1}w)$,
$$|\widetilde Q(u)-\widetilde Q(v)|\leq\Lip(Q_\delta)|B^{-1}(u-v)|\leq\Lip(Q_\delta)\|B^{-1}\||u-v|,$$
and hence $\Lip(\widetilde Q)\leq\Lip(Q_\delta)\|B^{-1}\|$. Recall from Section~\ref{sss:representation} that $\|B^{-1}\|\leq H(z_0)^3$, while $H(z_0)=\Psi(x_0,y_0)\leq A$ by Section~\ref{sss:mixed}. Since $\Lip(Q_\delta)\leq1$, we have
$$\|B^{-1}\|+\Lip(\widetilde Q)\leq\|B^{-1}\|\big(1+\Lip(Q_\delta)\big)\leq2\|B^{-1}\|\leq2H(z_0)^3\leq2A^3.$$
Moreover, the Fourier exponentials supplied by \cref{fourier-expansions} satisfy
$$\|DE_{k,2}\|_{L^\infty}\leq C(1+|k|),\qquad \|D^2E_{k,2}\|_{L^\infty}\leq C(1+|k|)^2,$$
where $C$ depends only on $d,L$. Fix $C_2:=\max\{1,4CA^6\}$.
The chain rule implies in particular
$$\partial_e(B_k\circ(-B^{-1}))(w)=DE_{k,2}(W(-B^{-1}w))[(-B^{-1}e,\partial_e\widetilde Q(w))].$$
This implies in particular that 
\begin{equation}
\begin{split}
&\partial_e(B_k\circ(-B^{-1}))(u)-\partial_e(B_k\circ(-B^{-1}))(v)\\
&\qquad=(DE_{k,2}(W(-B^{-1}u))-DE_{k,2}(W(-B^{-1}v)))[(-B^{-1}e,\partial_e\widetilde Q(u))]+DE_{k,2}(W(-B^{-1}v))\big[(0,\partial_e\widetilde Q(u)-\partial_e\widetilde Q(v))\big].
\end{split}
\label{identitadellebk}
\end{equation}
The mean value inequality gives
$$\|DE_{k,2}(W(-B^{-1}u))-DE_{k,2}(W(-B^{-1}v))\|\leq\|D^2E_{k,2}\|_{L^\infty}|W(-B^{-1}u)-W(-B^{-1}v)|.$$
Since $W(-B^{-1}w)=(-B^{-1}w,\widetilde Q(w))$, we also have
$$|W(-B^{-1}u)-W(-B^{-1}v)|\leq\big(\|B^{-1}\|+\Lip(\widetilde Q)\big)|u-v|,\qquad |(-B^{-1}e,\partial_e\widetilde Q(u))|\leq\|B^{-1}\|+\Lip(\widetilde Q).$$
Taking real or imaginary parts in \eqref{identitadellebk} and using these estimates, we obtain, for $\mu=(k,\nu)$,
\begin{equation}
\begin{split}
|\partial_eg_\mu(u)-\partial_eg_\mu(v)|
&\leq\|D^2E_{k,2}\|_{L^\infty}\big(\|B^{-1}\|+\Lip(\widetilde Q)\big)^2|u-v|+\|DE_{k,2}\|_{L^\infty}|\partial_e\widetilde Q(u)-\partial_e\widetilde Q(v)|\\
&\leq4CA^6(1+|k|)^2|u-v|+C(1+|k|)|\nabla\widetilde Q(u)-\nabla\widetilde Q(v)|.
\end{split}
\label{fourier-derivative-increment}
\end{equation}

We now apply \eqref{fourier-derivative-increment} to the integral expression for $D_{\mu,j,e,\ell}$ established at the beginning of this paragraph. Since
$|(w+th_je)-(w+(t-\ell-1)h_je)|=(\ell+1)h_j$,
we have
\begin{equation}
\begin{split}
|D_{\mu,j,e,\ell}(w)|
&\leq\int_0^1|\partial_eg_\mu(w+th_je)-\partial_eg_\mu(w+(t-\ell-1)h_je)|\,dt\\
&\leq C_2(\ell+1)(1+|k|)^2h_j+C_2(1+|k|)\int_0^1|\nabla\widetilde Q(w+th_je)-\nabla\widetilde Q(w+(t-\ell-1)h_je)|\,dt.
\end{split}
\label{fourier-remainder-increment-bound}
\end{equation}
It remains to sum these estimates with the weights $\Lip(f_\mu)$ and prove \eqref{weighted-derivative-increment-bound}. For I and III, the definitions in \eqref{normalized-phase-expansion} give
\begin{equation*}
\begin{split}
\Lip(f_1)+\Lip(f_3)=|\pi(S(0,1))|+|\langle(0,1),S(0,1)\rangle|\Lip(P_\delta)\leq\|S\|\big(1+\Lip(P_\delta)\big)\leq2.
\nonumber
\end{split}
\end{equation*}
For the coefficients of type II there is nothing to do. For the coefficients of type IV, recall from the definitions preceding \eqref{smooth-remainder-expansion} that
$$f_{(k,\nu)}=-F_{k,\nu},\qquad F_{k,1}=\operatorname{Re}A_k,\qquad F_{k,2}=-\operatorname{Im}A_k,\qquad A_k=c_kE_{k,1}\circ U.$$
Using $\Lip(U)\leq2$ and $\Lip(E_{k,1})\leq C(1+|k|)$, we obtain
$$\sum_{\nu=1}^2\Lip(F_{k,\nu})\leq2\Lip(A_k)\leq2|c_k|\Lip(E_{k,1})\Lip(U)\leq C|c_k|(1+|k|).$$
The estimate $\sum_k(1+|k|)^3|c_k|\leq C\delta$, proved in Section~\ref{sss:representation}, therefore yields
\begin{equation}
\begin{split}
\sum_k\sum_{\nu=1}^2\Lip(F_{k,\nu})(1+|k|)
&\leq\sum_k\sum_{\nu=1}^2\Lip(F_{k,\nu})(1+|k|)^2\leq C\sum_k|c_k|(1+|k|)^3\leq C\delta.
\end{split}
\label{weighted-fourier-factor-bound}
\end{equation}
Combining the bounds for I, II, III with \eqref{fourier-remainder-increment-bound}, we conclude that
\begin{equation*}
\begin{split}
&\qquad\qquad\qquad\sum_\mu\Lip(f_\mu)\|D_{\mu,j,e,\ell}\|_{L^\infty}
=\big(\Lip(f_1)+\Lip(f_3)\big)\|D_{1,j,e,\ell}\|_{L^\infty}+\sum_k\sum_{\nu=1}^2\Lip(F_{k,\nu})\|D_{(k,\nu),j,e,\ell}\|_{L^\infty}\\
&\leq C_2(d-1)h_j\sum_k\sum_{\nu=1}^2\Lip(F_{k,\nu})(1+|k|)^2\\
&\qquad\qquad\qquad+\Big(2+C_2\sum_k\sum_{\nu=1}^2\Lip(F_{k,\nu})(1+|k|)\Big)\Big\|\int_0^1|\nabla\widetilde Q(w+th_je)-\nabla\widetilde Q(w+(t-\ell-1)h_je)|\,dt\Big\|_{L^\infty}\\
&\leq C\delta h_j+C\Big\|\int_0^1|\nabla\widetilde Q(w+th_je)-\nabla\widetilde Q(w+(t-\ell-1)h_je)|\,dt\Big\|_{L^\infty},
\end{split}
\nonumber
\end{equation*}
where the last inequality uses \eqref{weighted-fourier-factor-bound} and $\delta\leq1$. This proves \eqref{weighted-derivative-increment-bound}, with a constant independent of $j,e,\ell$ and of the pair $(a,b)$. This proves our claim.

\spacedlowsmallcaps{Dini regularity.}
Assume that $\int_0^1\omega_Q(t)\frac{dt}{t}<\infty$. In this case we let $\tau_{\mu,j,e,\ell}:=\|D_{\mu,j,e,\ell}\|_{L^\infty}$
and notice that $B_j=0$ for every $j\geq0$, and that, once $M$ is chosen sufficiently large, we also have
$$C\sum_{j\geq0}A_j\leq\frac12,$$
where $C$ is the constant in \eqref{sommaA_j}. We will then apply \eqref{localized-operator-final-reduction} to prove $\|V_\lambda\|_{2\to2}^2\leq C\lambda^{-(d-2)}$. 
It is immediate to see that $B_j=0$, indeed 
$D_{\mu,j,e,\ell}\mathbb1_{\{|D_{\mu,j,e,\ell}|>\tau_{\mu,j,e,\ell}\}}=0$ almost everywhere, so \eqref{definition-Bj} gives $B_j=0$. Moreover, \eqref{definition-Aj} gives
$$A_j\leq C_d\max_{\substack{e\in\mathscr E^\pm\\0\leq\ell\leq d-2}}2^\ell\sum_\mu\Lip(f_\mu)\|D_{\mu,j,e,\ell}\|_{L^\infty}.$$

To estimate the right-hand side, we apply \eqref{weighted-derivative-increment-bound}. By the definitions of $\widetilde Q$ and $Q_\delta$, we have
$$\nabla\widetilde Q(w)=-B^{-T}\nabla Q_\delta(-B^{-1}w)=-B^{-T}\nabla Q(y_0-\delta B^{-1}w).$$
Consequently, using $\|B^{-1}\|\leq A^3$, we obtain
\begin{equation*}
\begin{split}
|\nabla\widetilde Q(u)-\nabla\widetilde Q(v)|
&\leq\|B^{-1}\|\,|\nabla Q(y_0-\delta B^{-1}u)-\nabla Q(y_0-\delta B^{-1}v)|\leq A^3\omega_Q(\delta A^3|u-v|).
\end{split}
\nonumber
\end{equation*}
Since $|(w+th_je)-(w+(t-\ell-1)h_je)|=(\ell+1)h_j$, substitution into \eqref{weighted-derivative-increment-bound} yields
\begin{equation*}
\begin{split}
\sum_\mu\Lip(f_\mu)\|D_{\mu,j,e,\ell}\|_{L^\infty}
&\leq C\delta h_j+CA^3\int_0^1\omega_Q(\delta A^3(\ell+1)h_j)dt\leq C(\delta h_j+\omega_Q(\delta A^3(d-1)h_j)).
\end{split}
\nonumber
\end{equation*}
Using $2^\ell\leq2^{d-2}$ in the preceding bound for $A_j$, therefore for all couples $(a,b)$, uniformly we have
$$A_j\leq C(\delta h_j+\omega_Q(\delta A^3(d-1)h_j)),$$

It remains to verify \eqref{sommaA_j}. 
Using $h_j=\frac{\kappa}{2^jM}$ and summing over the intervals $[\delta A^3(d-1)h_j,2\delta A^3(d-1)h_j]$, whose interiors are disjoint, we obtain
\begin{equation}
\begin{split}
\sum_{j\geq0}A_j
&\leq C\sum_{j\geq0}\Big(\frac{\delta\kappa}{2^jM}+\omega_Q\Big(\frac{\delta A^3(d-1)\kappa}{2^jM}\Big)\Big)\leq\frac{2C\delta\kappa}{M}+\frac{C}{\log2}\int_0^{\frac{2\delta A^3(d-1)\kappa}{M}}\omega_Q(u)\frac{du}{u}.
\end{split}
\label{dini-coefficient-sum}
\end{equation}
The right-hand side tends to zero as $M\to\infty$. We therefore fix $M\geq1$, independently of $\lambda$ and $(a,b)$, so that $\frac{2\delta A^3(d-1)\kappa}{M}\leq1$ and \eqref{sommaA_j} holds. Since $B_j=0$, \eqref{localized-operator-final-reduction} now gives
$$\|V_\lambda^{a,b}\|_{2\to2}^2\leq\frac{2C}{|\det B|}\Big(M^{d-2}+\frac1M\Big)\lambda^{-(d-2)}\leq C\lambda^{-(d-2)},$$
where we used $|\det B|^{-1}\leq\|B^{-1}\|^{d-2}\leq A^{3(d-2)}$ and the fact that $M$ has been fixed. Returning to the finite decomposition in Section~\ref{sss:normal}, we conclude that
$$\|V_\lambda\|_{2\to2}^2\leq\Big(\sum_{a,b}\|V_\lambda^{a,b}\|_{2\to2}\Big)^2\leq C\lambda^{-(d-2)}.$$
Suppose now that $\int_0^1\omega_P(t)\,\frac{dt}{t}<\infty$. By the definition of the adjoint, for every $f\in L^2$ we have
\begin{equation*}
\begin{split}
\overline{V_\lambda^*(\overline f)}(y)=\chi_1(y)\int_I e^{i\lambda\sqrt{1+|-v+y-x|^2+(Q(y)-P(x))^2}}\chi_0(x)f(x)dx.
\nonumber
\end{split}
\end{equation*}
The last expression has the form of the operator in the statement, with
$(P,Q,v,\chi_0,\chi_1)$ replaced by $(Q,P,-v,\chi_1,\chi_0)$. 
Since this substitution preserves the hypotheses, the estimate proved above applies and gives
$$\big\|\overline{V_\lambda^*(\overline f)}\big\|_{L^2}^2\leq C\lambda^{-(d-2)}\|f\|_{L^2}^2.$$
Taking the supremum over $f\ne0$, and using the invariance of the $L^2$ norm under complex conjugation, we conclude that
$$\|V_\lambda\|_{2\to2}^2=\|V_\lambda^*\|_{2\to2}^2=\sup_{f\ne0}\frac{\|\overline{V_\lambda^*(\overline f)}\|_{L^2}^2}{\|f\|_{L^2}^2}\leq C\lambda^{-(d-2)}.$$
This proves assertion (i) of the thesis.

\smallskip
\spacedlowsmallcaps{Sobolev regularity and the sum of the increments.}
Assume now that $\nabla P,\nabla Q\in W^{s,p}(\mathbb R^{d-2})$, where $0<s<1$ and $1\leq p<\infty$. We first estimate the Sobolev seminorm of $\widetilde{Q}$. 
Using $u=y_0-\delta B^{-1}w$ and $v=y_0-\delta B^{-1}w'$, we obtain
\begin{equation}
\begin{split}
[\nabla\widetilde Q]_{W^{s,p}}^p
&=\int_{\mathbb R^{d-2}}\int_{\mathbb R^{d-2}}\frac{|B^{-T}(\nabla Q(y_0-\delta B^{-1}w)-\nabla Q(y_0-\delta B^{-1}w'))|^p}{|w-w'|^{d-2+sp}}\,dw\,dw'\\
&\leq|\det B|^2\|B^{-1}\|^{p+d-2+sp}\delta^{sp-(d-2)}\int_{\mathbb R^{d-2}}\int_{\mathbb R^{d-2}}\frac{|\nabla Q(u)-\nabla Q(v)|^p}{|u-v|^{d-2+sp}}dudv\leq C\delta^{sp-(d-2)}[\nabla Q]_{W^{s,p}}^p.
\end{split}
\label{rescaled-gradient-sobolev-bound}
\end{equation}
The same change of variables gives $\|\nabla\widetilde Q\|_{L^p}\leq C\delta^{-\frac{d-2}{p}}\|\nabla Q\|_{L^p}$, and in particular we have that $\nabla\widetilde Q\in W^{s,p}(\mathbb R^{d-2})$. 
In the following we shall prove that, for $0<s<1$ and $1\leq p<\infty$,
\begin{equation}
\begin{split}
&\sum_{j\geq0}h_j^{-sp}\max_{\substack{e\in\mathscr E^\pm\\0\leq\ell\leq d-2}}\Big\|\int_0^1|\nabla\widetilde Q(\cdot+th_je)-\nabla\widetilde Q(\cdot+(t-\ell-1)h_je)|\,dt\Big\|_{L^p}^p\leq C[\nabla\widetilde Q]_{W^{s,p}}^p,
\end{split}
\label{dyadic-sobolev-increment-sum}
\end{equation}
where $C$ depends only on $d,s,p$. Let $F:\mathbb R^{d-2}\to\mathbb R^{d-2}$ be a measurable map with $[F]_{W^{s,p}}<\infty$. For every $v,z\in\mathbb R^{d-2}$, we have
\begin{equation}
\begin{split}
\|F(\cdot+v)-F\|_{L^p}^p
&=\int_{\mathbb R^{d-2}}|F(w+v)-F(w+z)+F(w+z)-F(w)|^p\,dw\\
&\leq2^{p-1}(\|F(\cdot+v-z)-F\|_{L^p}^p+\|F(\cdot+z)-F\|_{L^p}^p).
\label{differenzaLptriang}
\end{split}
\end{equation}
Now, let us observe that for $v\neq0$, by setting $u:=v-z$, we have $|u-\frac{v}{2}|=|\frac{v}{2}-z|=|z-\frac{v}{2}|$. Thus, the map $v\mapsto v-z$ maps $B(\frac{v}{2},\frac{|v|}{4})$ onto itself and it is an isometry. Consequently,
\begin{equation*}
\begin{split}
\int_{B(\frac{v}{2},\frac{|v|}{4})}\|F(\cdot+v-z)-F\|_{L^p}^pdz=\int_{B(\frac{v}{2},\frac{|v|}{4})}\|F(\cdot+u)-F\|_{L^p}^pdu.
\nonumber
\end{split}
\end{equation*}
Using the preceding inequality and this identity, we obtain
\begin{equation*}
\begin{split}
\|F(\cdot+v)-F\|_{L^p}^p
&=\fint_{B(\frac{v}{2},\frac{|v|}{4})}\|F(\cdot+v)-F\|_{L^p}^p\,dz\leq2^{p-1}\fint_{B(\frac{v}{2},\frac{|v|}{4})}\Big(\|F(\cdot+v-z)-F\|_{L^p}^p+\|F(\cdot+z)-F\|_{L^p}^p\Big)\,dz\\
&\qquad\qquad\qquad=2^{p}\fint_{B(\frac{v}{2},\frac{|v|}{4})}\|F(\cdot+z)-F\|_{L^p}^p\,dz.
\nonumber
\end{split}
\end{equation*}
Since $B(\frac{v}{2},\frac{|v|}{4})\subset\{z:\frac{|v|}{4}\leq|z|\leq\frac{3|v|}{4}\}$, we obtain
\begin{equation}
\begin{split}
|v|^{-sp}\|F(\cdot+v)-F\|_{L^p}^p
&\leq 2^p|v|^{-sp}\fint_{B(\frac{v}{2},\frac{|v|}{4})}\|F(\cdot+z)-F\|_{L^p}^p\,dz\leq C\int_{\frac{|v|}{4}\leq|z|\leq\frac{3|v|}{4}}\frac{\|F(\cdot+z)-F\|_{L^p}^p}{|z|^{d-2+sp}}dz,
\end{split}
\label{directional-increment-annular-bound}
\end{equation}
where the constant depends only on $s,p$ and $d$. 
Fix $e\in\mathscr E^\pm$ and $0\leq\ell\leq d-2$, and apply \eqref{directional-increment-annular-bound} with $v=(\ell+1)h_je$. Since
$h_j^{-sp}=(\ell+1)^{sp}|v|^{-sp}\leq(d-1)^{sp}|v|^{-sp}$, it follows that
$$h_j^{-sp}\|F(\cdot+(\ell+1)h_je)-F\|_{L^p}^p\leq C\int_{\frac{(\ell+1)h_j}{4}\leq|z|\leq\frac{3(\ell+1)h_j}{4}}\frac{\|F(\cdot+z)-F\|_{L^p}^p}{|z|^{d-2+sp}}\,dz.$$
Since $h_{j+1}=\frac{h_j}{2}$, for each $z\ne0$ there are at most three indices $j$ satisfying $\frac{4|z|}{3}\leq(\ell+1)h_j\leq4|z|$. Hence Tonelli's theorem gives
\begin{equation}
\begin{split}
&\qquad\quad\sum_{j\geq0}h_j^{-sp}\|F(\cdot+(\ell+1)h_je)-F\|_{L^p}^p\leq C\sum_{j\geq0}\int_{\mathbb R^{d-2}}\mathbb1_{\{\frac{(\ell+1)h_j}{4}\leq|z|\leq\frac{3(\ell+1)h_j}{4}\}}\frac{\|F(\cdot+z)-F\|_{L^p}^p}{|z|^{d-2+sp}}\,dz\\
&=C\int_{\mathbb R^{d-2}}\Big(\sum_{j\geq0}\mathbb1_{\{\frac{(\ell+1)h_j}{4}\leq|z|\leq\frac{3(\ell+1)h_j}{4}\}}\Big)\frac{\|F(\cdot+z)-F\|_{L^p}^p}{|z|^{d-2+sp}}\,dz\leq3C\int_{\mathbb R^{d-2}}\frac{\|F(\cdot+z)-F\|_{L^p}^p}{|z|^{d-2+sp}}\,dz=3C[F]_{W^{s,p}}^p.
\end{split}
\label{directional-dyadic-sobolev-bound}
\end{equation}
We now apply this estimate to $F=\nabla\widetilde Q$. Jensen's inequality gives
\begin{equation*}
\begin{split}
&\Big\|\int_0^1|\nabla\widetilde Q(\cdot+th_je)-\nabla\widetilde Q(\cdot+(t-\ell-1)h_je)|\,dt\Big\|_{L^p}^p\leq\int_0^1\int_{\mathbb R^{d-2}}|\nabla\widetilde Q(w+th_je)-\nabla\widetilde Q(w+(t-\ell-1)h_je)|^p\,dw\,dt\\
&\qquad\qquad\qquad=\int_0^1\|\nabla\widetilde Q(\cdot+(\ell+1)h_je)-\nabla\widetilde Q\|_{L^p}^p\,dt=\|\nabla\widetilde Q(\cdot+(\ell+1)h_je)-\nabla\widetilde Q\|_{L^p}^p.
\nonumber
\end{split}
\end{equation*}
Combining the preceding inequality with \eqref{directional-dyadic-sobolev-bound}, we obtain
\begin{equation*}
\begin{split}
&\sum_{j\geq0}h_j^{-sp}\max_{\substack{e\in\mathscr E^\pm\\0\leq\ell\leq d-2}}\Big\|\int_0^1|\nabla\widetilde Q(\cdot+th_je)-\nabla\widetilde Q(\cdot+(t-\ell-1)h_je)|\,dt\Big\|_{L^p}^p\\&\qquad\leq\sum_{j\geq0}h_j^{-sp}\max_{\substack{e\in\mathscr E^\pm\\0\leq\ell\leq d-2}}\|\nabla\widetilde Q(\cdot+(\ell+1)h_je)-\nabla\widetilde Q\|_{L^p}^p\\
&\qquad\leq\sum_{e\in\mathscr E^\pm}\sum_{\ell=0}^{d-2}\sum_{j\geq0}h_j^{-sp}\|\nabla\widetilde Q(\cdot+(\ell+1)h_je)-\nabla\widetilde Q\|_{L^p}^p\\&\qquad\leq6C(d-2)(d-1)[\nabla\widetilde Q]_{W^{s,p}}^p.
\nonumber
\end{split}
\end{equation*}
This proves \eqref{dyadic-sobolev-increment-sum}.

\smallskip

\spacedlowsmallcaps{Choice of the thresholds $\tau_{\mu,j,e,\ell}$.}
We shall choose the thresholds so that
\begin{equation}
\sum_{\substack{j\geq0\\2^jM\leq\lambda\delta^2\sqrt{d-2}}}A_j\leq\frac{1}{2C},
\label{stimaridottaAj}
\end{equation}
and then estimate the sum of $B_j$ over the same indices.
Recall that, immediately after \eqref{initial-operator-decomposition}, we observed that the annular kernels vanish unless $r_j\leq|x-x'|\leq4r_j$ and $x,x'\in\supp\eta_0$. Here $\eta_0(x)=\chi_{0,a}(x_0+\delta x)$, and the support condition imposed in Section~\ref{sss:normal} gives $\supp\eta_0\subset(-\frac{1}{2},\frac{1}{2})^{d-2}$. Consequently, $\eta_0(x)\overline{\eta_0(x')}\rho_{j,e}(x,x')\neq0$ implies that 
$$\frac{2^jM}{\lambda\delta^2}=r_j\leq|x-x'|\leq\sqrt{d-2},$$
and in particular $2^jM\leq\lambda\delta^2\sqrt{d-2}$.
Since every kernel $\mathcal K_{e,j,z}^{Q,Q',\ell}$ contains the factor $\eta_0(x)\overline{\eta_0(x')}\rho_{j,e}(x,x')$, we infer that 
$$\mathcal K_{e,j,z}^{Q,Q',\ell}=0\qquad\text{whenever }\qquad 2^jM>\lambda\delta^2\sqrt{d-2}.$$
This is why it is sufficient to estimate \eqref{stimaridottaAj}.

We now fix the constants. Let $C$ be the constant in \eqref{sommaA_j}, and let $C_d$ be the constant in \eqref{definition-Aj}. By \eqref{weighted-fourier-factor-bound} and $\Lip(f_1)+\Lip(f_3)\leq2$, we can fix $C_*\geq1$, uniformly in $(a,b)$, such that
\begin{equation*}
\begin{split}
2^{d-1}CC_dC_2(d-1)\sum_k\sum_{\nu=1}^2\Lip(F_{k,\nu})(1+|k|)^2&\leq C_*\delta,\\
2^{d-1}CC_d\Big(\Lip(f_1)+\Lip(f_3)+C_2\sum_k\sum_{\nu=1}^2\Lip(F_{k,\nu})(1+|k|)\Big)&\leq C_*.
\end{split}
\nonumber
\end{equation*}
Set $M:=\max\{1,8C_*\delta\kappa\}$. If $\lambda\delta^2\sqrt{d-2}<M$, all the annular operators vanish, so \eqref{initial-operator-decomposition} becomes
$$T_{\lambda\delta^2}^{a,b}(T_{\lambda\delta^2}^{a,b})^*=\mathcal K_{\lambda\delta^2,0}^{\mathrm{near}}.$$
Consequently, the identity proved in Section~\ref{sss:representation} and the bound obtained from \eqref{near-row-bound} and \eqref{near-column-bound}, applied with $\lambda$ replaced by $\lambda\delta^2$, give
\begin{equation*}
\begin{split}
\|V_\lambda^{a,b}\|_{2\to2}^2
&=\frac{\delta^{2(d-2)}}{|\det B|}\|T_{\lambda\delta^2}^{a,b}\|_{2\to2}^2=\frac{\delta^{2(d-2)}}{|\det B|}\|\mathcal K_{\lambda\delta^2,0}^{\mathrm{near}}\|_{2\to2}\leq\frac{C_d}{|\det B|}\|\eta_0\|_{L^\infty}^2\|\eta_1\|_{L^2}^2M^{d-2}\frac{\delta^{2(d-2)}}{(\lambda\delta^2)^{d-2}}\\
&=\frac{C_d}{|\det B|}\|\eta_0\|_{L^\infty}^2\|\eta_1\|_{L^2}^2M^{d-2}\lambda^{-(d-2)}\leq C\lambda^{-(d-2)}.
\end{split}
\nonumber
\end{equation*}
The last inequality uses $|\det B|^{-1}\leq A^{3(d-2)}$, the uniform bounds for the cutoffs, and the choice of $M$.
This proves the required bound in this case.

We thus assume $\lambda\delta^2\sqrt{d-2}\geq M$ and define
$J:=\max\{j\geq0:2^jM\leq\lambda\delta^2\sqrt{d-2}\}$.
We now fix
\begin{equation}
t_j:=\begin{cases}
\displaystyle\frac{1}{4C_*}\frac{2^{\frac{j(d-2-sp)}{p-1}}}{\sum_{i=0}^J2^{\frac{i(d-2-sp)}{p-1}}},&p>1,\quad 0\leq j\leq J,\\
0,&p=1\ \text{or}\ j>J,
\end{cases}
\label{sobolev-scale-levels}
\end{equation}
and finally define
\begin{equation}
\tau_{\mu,j,e,\ell}:=\begin{cases}
2t_j,&\mu\in\{1,3\},\\
0,&\mu=2,\\
2C_2((d-1)(1+|k|)^2h_j+(1+|k|)t_j),&\mu=(k,\nu).
\end{cases}
\label{sobolev-truncation-levels}
\end{equation}
These thresholds are independent of $Q,Q'$ and of the vector $z$ in $\eta_1(w-z)$. Now that all definitions are in place, we verify the required bound for $\sum_{j=0}^JA_j$.

Recall that in \eqref{definizionefkmu} we defined $f_{(k,\nu)}(x):=-F_{k,\nu}(x)$. Consequently, $\Lip(f_{(k,\nu)})=\Lip(F_{k,\nu})$ and hence, substituting this in \eqref{sobolev-truncation-levels} gives
\begin{equation*}
\begin{split}
&\sum_\mu\Lip(f_\mu)\tau_{\mu,j,e,\ell}
=2t_j\big(\Lip(f_1)+\Lip(f_3)\big)+2C_2\sum_k\sum_{\nu=1}^2\Lip(F_{k,\nu})\Big((d-1)(1+|k|)^2h_j+(1+|k|)t_j\Big)\\
&=2C_2(d-1)h_j\sum_k\sum_{\nu=1}^2\Lip(F_{k,\nu})(1+|k|)^2+2t_j\Big(\Lip(f_1)+\Lip(f_3)+C_2\sum_k\sum_{\nu=1}^2\Lip(F_{k,\nu})(1+|k|)\Big)\\
&\qquad\qquad\qquad\quad\leq \frac{C_*}{2^{d-2}CC_d}(\delta h_j+t_j),
\nonumber
\end{split}
\end{equation*}
where the last inequality comes from the choice of the constant $C_*$. By definition of $A_j$ notice that 
$$CA_j=CC_d\max_{\substack{e\in \mathscr E^\pm\\0\leq \ell\leq d-2}}2^\ell\sum_\mu\Lip(f_\mu)\tau_{\mu,j,e,\ell}\leq C_*(\delta h_j+t_j).$$

We now shall prove that the choice in 
\eqref{sobolev-truncation-levels} satisfies
\begin{equation*}
C\sum_{j=0}^JA_j\leq\frac12\qquad\text{and}\qquad\sum_{j=0}^JB_j\leq C[\nabla\widetilde Q]_{W^{s,p}}^p\begin{cases}
(\lambda\delta^2)^{-(d-2)},&\text{if }sp>d-2,\\
(\lambda\delta^2)^{-(d-2)}(\log(2+\lambda\delta^2))^{p-1},&\text{if }sp=d-2,\\
(\lambda\delta^2)^{-sp},&\text{if }sp<d-2.
\end{cases}
\nonumber
\end{equation*}
By \eqref{sobolev-scale-levels}, $\sum_{j=0}^Jt_j=\frac{1}{4C_*}$ when $p>1$, whereas this sum is zero when $p=1$. 
Using the estimate $CA_j\leq C_*(\delta h_j+t_j)$ just proved, together with $h_j=\frac{\kappa}{2^jM}$ and the choice $M\geq8C_*\delta\kappa$, we obtain
\begin{equation*}
\begin{split}
C\sum_{j=0}^JA_j\leq C_*\delta\sum_{j=0}^Jh_j+C_*\sum_{j=0}^Jt_j\leq\frac{C_*\delta\kappa}{M}\sum_{j=0}^{\infty}2^{-j}+\frac14=\frac{2C_*\delta\kappa}{M}+\frac14\leq\frac12.
\end{split}
\nonumber
\end{equation*}
We next estimate $B_j$. For every $p>1$, we will show that 
\begin{equation}
B_j\leq Cr_j^{d-2}t_j^{1-p}\max_{\substack{e\in\mathscr E^\pm\\0\leq\ell\leq d-2}}\Big\|\int_0^1|\nabla\widetilde Q(\cdot+th_je)-\nabla\widetilde Q(\cdot+(t-\ell-1)h_je)|\,dt\Big\|_{L^p}^p.
\label{sobolev-tail-coefficient-bound}
\end{equation}
This is the estimate to which we will apply \eqref{dyadic-sobolev-increment-sum}. 
By \eqref{definition-Bj}, we must estimate
$\|D_{\mu,j,e,\ell}\mathbb1_{\{|D_{\mu,j,e,\ell}|>\tau_{\mu,j,e,\ell}\}}\|_{L^1}$. 
Let us first observe that if $\mu\in\{1,3\}$, we have $g_\mu=\widetilde Q$. Thus \eqref{stimaprimitermini} shows that
$$|D_{\mu,j,e,\ell}(w)|\leq\int_0^1|\nabla\widetilde Q(w+th_je)-\nabla\widetilde Q(w+(t-\ell-1)h_je)|\,dt.$$
Since $\tau_{\mu,j,e,\ell}=2t_j$ by \eqref{sobolev-truncation-levels}, on $\{|D_{\mu,j,e,\ell}|>\tau_{\mu,j,e,\ell}\}$ we have
$$2t_j<|D_{\mu,j,e,\ell}(w)|\leq\int_0^1|\nabla\widetilde Q(w+th_je)-\nabla\widetilde Q(w+(t-\ell-1)h_je)|\,dt.$$
Recall, as always, that the term with $\mu=2$ contributes zero because $D_{2,j,e,\ell}=0$.
For $\mu=(k,\nu)$, recall from \eqref{definizionefkmu} that $g_\mu(w)=G_{k,\nu}(-B^{-1}w)-G_{k,\nu}(0)$. For the corresponding function $D_{\mu,j,e,\ell}$ defined in \eqref{definition-normalized-increment}, we shall prove that
\begin{equation}
\begin{split}
2C_2(1+|k|)t_j
&\leq\tau_{\mu,j,e,\ell}<|D_{\mu,j,e,\ell}(w)|\leq2C_2(1+|k|)\int_0^1|\nabla\widetilde Q(w+th_je)-\nabla\widetilde Q(w+(t-\ell-1)h_je)|\,dt
\label{ineq:Qtau}
\end{split}
\end{equation}
whenever $|D_{\mu,j,e,\ell}(w)|>\tau_{\mu,j,e,\ell}$.
Indeed, thanks to the definition in \eqref{sobolev-truncation-levels}, for a $w$ satisfying $|D_{\mu,j,e,\ell}(w)|>\tau_{\mu,j,e,\ell}$, we have
$$C_2(d-1)(1+|k|)^2h_j\leq\frac{\tau_{\mu,j,e,\ell}}2<\frac{|D_{\mu,j,e,\ell}(w)|}2.$$
Combining this inequality with \eqref{fourier-remainder-increment-bound}, we obtain
\begin{equation*}
\begin{split}
\frac{|D_{\mu,j,e,\ell}(w)|}2\leq C_2(1+|k|)\int_0^1|\nabla\widetilde Q(w+th_je)-\nabla\widetilde Q(w+(t-\ell-1)h_je)|\,dt.
\nonumber
\end{split}
\end{equation*}
This concludes the proof of \eqref{ineq:Qtau}.
Dividing the chain by $2C_2(1+|k|)>0$, we also obtain
$$t_j<\int_0^1|\nabla\widetilde Q(w+th_je)-\nabla\widetilde Q(w+(t-\ell-1)h_je)|\,dt.$$

Suppose now that $p>1$. The numbers $t_j$ are positive for $0\leq j\leq J$, and $u\mathbb1_{\{u>t_j\}}\leq t_j^{1-p}u^p$ for $u\geq0$. Applying this inequality to the integrals in the preceding bounds, integrating in $w$ and using $\Lip(f_{(k,\nu)})=\Lip(F_{k,\nu})$ from \eqref{definizionefkmu}, we obtain
\begin{equation*}
\begin{split}
&|D_{\mu,j,e,\ell}(w)|\mathbb1_{\{|D_{\mu,j,e,\ell}(w)|>\tau_{\mu,j,e,\ell}\}}\\&\qquad\leq t_j^{1-p}\Big(\int_0^1|\nabla\widetilde Q(w+th_je)-\nabla\widetilde Q(w+(t-\ell-1)h_je)|\,dt\Big)^p\begin{cases}
1,&\mu\in\{1,3\},\\
0,&\mu=2,\\
2C_2(1+|k|),&\mu=(k,\nu).
\end{cases}
\nonumber
\end{split}
\end{equation*}
The coefficients are those of the pointwise bounds proved above. Since $\Lip(f_{(k,\nu)})=\Lip(F_{k,\nu})$ by \eqref{definizionefkmu}, integration and summation give
\begin{equation}
\begin{split}
&\qquad\qquad\qquad\sum_\mu\Lip(f_\mu)\|D_{\mu,j,e,\ell}\mathbb1_{\{|D_{\mu,j,e,\ell}|>\tau_{\mu,j,e,\ell}\}}\|_{L^1}
=\sum_\mu\Lip(f_\mu)\int_{\mathbb R^{d-2}}|D_{\mu,j,e,\ell}(w)|\mathbb1_{\{|D_{\mu,j,e,\ell}(w)|>\tau_{\mu,j,e,\ell}\}}\,dw\\
&\leq t_j^{1-p}\Big(\Lip(f_1)+\Lip(f_3)+2C_2\sum_k\sum_{\nu=1}^2\Lip(F_{k,\nu})(1+|k|)\Big)\Big\|\int_0^1|\nabla\widetilde Q(\cdot+th_je)-\nabla\widetilde Q(\cdot+(t-\ell-1)h_je)|\,dt\Big\|_{L^p}^p\\
&\qquad\qquad\qquad\qquad\leq Ct_j^{1-p}\Big\|\int_0^1|\nabla\widetilde Q(\cdot+th_je)-\nabla\widetilde Q(\cdot+(t-\ell-1)h_je)|\,dt\Big\|_{L^p}^p.
\end{split}
\label{sobolev-tail-integral-bound}
\end{equation}
For the last inequality, we used $\Lip(f_1)+\Lip(f_3)\leq2$, proved before \eqref{weighted-fourier-factor-bound}, and that same estimate to obtain
$$\Lip(f_1)+\Lip(f_3)+2C_2\sum_k\sum_{\nu=1}^2\Lip(F_{k,\nu})(1+|k|)\leq2+C\delta\leq C.$$
Substituting \eqref{sobolev-tail-integral-bound} into \eqref{definition-Bj}, and using $2^\ell\leq2^{d-2}$ and the uniform bounds for the cutoffs, proves \eqref{sobolev-tail-coefficient-bound}. 
Let  $p>1$ and recall that $r_j=\frac{2^jM}{\lambda\delta^2}$, $h_j=\frac{\kappa}{2^jM}$ and the definition of $t_j$, we obtain
\begin{equation}
\begin{split}
&\sum_{j=0}^JB_j
\overset{\eqref{sobolev-tail-coefficient-bound}}{\leq}C\sum_{j=0}^J(r_j^{d-2}h_j^{sp}t_j^{1-p})h_j^{-sp}\max_{\substack{e\in\mathscr E^\pm\\0\leq\ell\leq d-2}}\Big\|\int_0^1|\nabla\widetilde Q(\cdot+th_je)-\nabla\widetilde Q(\cdot+(t-\ell-1)h_je)|\,dt\Big\|_{L^p}^p\\
&\qquad\overset{\eqref{sobolev-scale-levels}}{=}C\sum_{j=0}^J\Big(\frac{2^jM}{\lambda\delta^2}\Big)^{d-2}\Big(\frac{\kappa}{2^jM}\Big)^{sp}t_j^{1-p}h_j^{-sp}\max_{\substack{e\in\mathscr E^\pm\\0\leq\ell\leq d-2}}\Big\|\int_0^1|\nabla\widetilde Q(\cdot+th_je)-\nabla\widetilde Q(\cdot+(t-\ell-1)h_je)|\,dt\Big\|_{L^p}^p\\
&=C\kappa^{sp}(\lambda\delta^2)^{-(d-2)}M^{d-2-sp}\sum_{j=0}^J2^{(d-2-sp)j}t_j^{1-p}h_j^{-sp}\max_{\substack{e\in\mathscr E^\pm\\0\leq\ell\leq d-2}}\Big\|\int_0^1|\nabla\widetilde Q(\cdot+th_je)-\nabla\widetilde Q(\cdot+(t-\ell-1)h_je)|\,dt\Big\|_{L^p}^p
\end{split}
\end{equation}
By the choice of $t_j$ in \eqref{sobolev-scale-levels}, for every $0\leq j\leq J$ we have
$$2^{(d-2-sp)j}t_j^{1-p}=(4C_*)^{p-1}\Big(\sum_{i=0}^{J}2^{\frac{i(d-2-sp)}{p-1}}\Big)^{p-1},$$
and therefore
\begin{equation*}
    \begin{split}
        &\qquad\qquad\sum_{j=0}^J2^{(d-2-sp)j}t_j^{1-p}h_j^{-sp}\max_{\substack{e\in\mathscr E^\pm\\0\leq\ell\leq d-2}}\Big\|\int_0^1|\nabla\widetilde Q(\cdot+th_je)-\nabla\widetilde Q(\cdot+(t-\ell-1)h_je)|\,dt\Big\|_{L^p}^p\\
        &=(4C_*)^{p-1}\Big(\sum_{i=0}^{J}2^{\frac{i(d-2-sp)}{p-1}}\Big)^{p-1}\sum_{j=0}^Jh_j^{-sp}\max_{\substack{e\in\mathscr E^\pm\\0\leq\ell\leq d-2}}\Big\|\int_0^1|\nabla\widetilde Q(\cdot+th_je)-\nabla\widetilde Q(\cdot+(t-\ell-1)h_je)|\,dt\Big\|_{L^p}^p\\
&\qquad\qquad\qquad\qquad\qquad\overset{\eqref{dyadic-sobolev-increment-sum}}{\leq} C(4C_*)^{p-1}\Big(\sum_{i=0}^{J}2^{\frac{i(d-2-sp)}{p-1}}\Big)^{p-1}[\nabla\widetilde Q]_{W^{s,p}}^p.
        \nonumber
    \end{split}
\end{equation*}
The three claimed bounds for $p>1$ now follow by evaluating the geometric sum. If $sp>d-2$, then $M^{d-2-sp}\leq1$ and
$$\sum_{i=0}^J2^{\frac{i(d-2-sp)}{p-1}}\leq\sum_{i=0}^{\infty}2^{\frac{i(d-2-sp)}{p-1}}=\frac1{1-2^{\frac{d-2-sp}{p-1}}}.$$
If $sp=d-2$, then
$$M^{d-2-sp}\Big(\sum_{i=0}^J2^{\frac{i(d-2-sp)}{p-1}}\Big)^{p-1}=(J+1)^{p-1}\leq\Big(1+\log_2\Big(\frac{\lambda\delta^2\sqrt{d-2}}{M}\Big)\Big)^{p-1}\leq C(\log(2+\lambda\delta^2))^{p-1}.$$
Finally, if $sp<d-2$, then
$$\sum_{i=0}^J2^{\frac{i(d-2-sp)}{p-1}}=\frac{2^{\frac{(J+1)(d-2-sp)}{p-1}}-1}{2^{\frac{d-2-sp}{p-1}}-1}\leq C2^{\frac{J(d-2-sp)}{p-1}},$$
and hence
$$M^{d-2-sp}\Big(\sum_{i=0}^J2^{\frac{i(d-2-sp)}{p-1}}\Big)^{p-1}\leq C(2^JM)^{d-2-sp}\leq C(\lambda\delta^2)^{d-2-sp}.$$
It remains to prove that, for $p=1$, we have
$$\sum_{j=0}^JB_j\leq C(\lambda\delta^2)^{-s}[\nabla\widetilde Q]_{W^{s,1}}.$$
In this case $t_j=0$ by \eqref{sobolev-scale-levels}. We use \eqref{stimaprimitermini} for $\mu\in\{1,3\}$ and the upper bound for $|D_{\mu,j,e,\ell}|$ in \eqref{ineq:Qtau} for $\mu=(k,\nu)$. The latter was proved using only $t_j\geq0$ and therefore remains valid for this choice. More precisely, its last inequality gives
$$|D_{\mu,j,e,\ell}(w)|\mathbb1_{\{|D_{\mu,j,e,\ell}(w)|>\tau_{\mu,j,e,\ell}\}}\leq2C_2(1+|k|)\int_0^1|\nabla\widetilde Q(w+th_je)-\nabla\widetilde Q(w+(t-\ell-1)h_je)|\,dt,\qquad \mu=(k,\nu).$$
Indeed, \eqref{ineq:Qtau} supplies this bound on the set selected by the indicator, and the left-hand side is zero outside that set.

We apply these bounds inside the integrals defining the $L^1$ norms in \eqref{definition-Bj}. Since $D_{2,j,e,\ell}=0$ and $\Lip(f_{(k,\nu)})=\Lip(F_{k,\nu})$ by \eqref{definizionefkmu}, we obtain
\begin{equation*}
\begin{split}
&\qquad\qquad\sum_\mu\Lip(f_\mu)\|D_{\mu,j,e,\ell}\mathbb1_{\{|D_{\mu,j,e,\ell}|>\tau_{\mu,j,e,\ell}\}}\|_{L^1} =\sum_\mu\Lip(f_\mu)\int_{\mathbb R^{d-2}}|D_{\mu,j,e,\ell}(w)|\mathbb1_{\{|D_{\mu,j,e,\ell}(w)|>\tau_{\mu,j,e,\ell}\}}\,dw\\
&\leq\Big(\Lip(f_1)+\Lip(f_3)+2C_2\sum_k\sum_{\nu=1}^2\Lip(F_{k,\nu})(1+|k|)\Big)\int_{\mathbb R^{d-2}}\int_0^1|\nabla\widetilde Q(w+th_je)-\nabla\widetilde Q(w+(t-\ell-1)h_je)|\,dt\,dw\\
&\qquad\qquad\qquad\qquad\qquad\qquad\overset{\eqref{weighted-fourier-factor-bound}}{\leq}C\Big\|\int_0^1|\nabla\widetilde Q(\cdot+th_je)-\nabla\widetilde Q(\cdot+(t-\ell-1)h_je)|\,dt\Big\|_{L^1},
\nonumber
\end{split}
\end{equation*}
where we also used $\Lip(f_1)+\Lip(f_3)\leq2$, established immediately before \eqref{weighted-fourier-factor-bound}.

We substitute this estimate into \eqref{definition-Bj} and sum over $j$, obtaining 
\begin{equation*}
\begin{split}
\sum_{j=0}^JB_j
&\overset{\eqref{definition-Bj}}{\leq}C\sum_{j=0}^J(r_j^{d-2}h_j^s)h_j^{-s}\max_{\substack{e\in\mathscr E^\pm\\0\leq\ell\leq d-2}}\Big\|\int_0^1|\nabla\widetilde Q(\cdot+th_je)-\nabla\widetilde Q(\cdot+(t-\ell-1)h_je)|\,dt\Big\|_{L^1}\\
&\leq C\max_{0\leq j\leq J}(r_j^{d-2}h_j^s)\sum_{j=0}^Jh_j^{-s}\max_{\substack{e\in\mathscr E^\pm\\0\leq\ell\leq d-2}}\Big\|\int_0^1|\nabla\widetilde Q(\cdot+th_je)-\nabla\widetilde Q(\cdot+(t-\ell-1)h_je)|\,dt\Big\|_{L^1}\\
&\overset{\eqref{dyadic-sobolev-increment-sum}}{\leq}C\max_{0\leq j\leq J}(r_j^{d-2}h_j^s)[\nabla\widetilde Q]_{W^{s,1}}=C\kappa^s(\lambda\delta^2)^{-(d-2)}\max_{0\leq j\leq J}(2^jM)^{d-2-s}[\nabla\widetilde Q]_{W^{s,1}}.
\end{split}
\nonumber
\end{equation*}
The first inequality includes the uniform cutoff bounds and $2^\ell\leq2^{d-2}$. In the two equalities we used $r_j=\frac{2^jM}{\lambda\delta^2}$, $h_j=\frac{\kappa}{2^jM}$.
Since $p=1$ and $0<s<1\leq d-2$, the exponent $d-2-s$ is positive. Therefore the maximum is attained at $j=J$, and the definition of $J$ gives
\begin{equation*}
\begin{split}
\sum_{j=0}^JB_j
&\leq C\kappa^s(\lambda\delta^2)^{-(d-2)}(2^JM)^{d-2-s}[\nabla\widetilde Q]_{W^{s,1}}\leq C\kappa^s(\lambda\delta^2)^{-(d-2)}(\lambda\delta^2\sqrt{d-2})^{d-2-s}[\nabla\widetilde Q]_{W^{s,1}}\leq C(\lambda\delta^2)^{-s}[\nabla\widetilde Q]_{W^{s,1}}.
\end{split}
\nonumber
\end{equation*}
This completes the case $p=1$. Combining it with the three estimates already proved for $p>1$, we obtain, for every $1\leq p<\infty$,
\begin{equation}
\sum_{j=0}^JB_j\leq C[\nabla\widetilde Q]_{W^{s,p}}^p\begin{cases}
(\lambda\delta^2)^{-(d-2)},&sp>d-2,\\
(\lambda\delta^2)^{-(d-2)}(\log(2+\lambda\delta^2))^{p-1},&sp=d-2,\\
(\lambda\delta^2)^{-sp},&sp<d-2.
\end{cases}
\label{stimeB_jfinali}
\end{equation}
We can now conclude the operator estimate for all the cases considered above.

\smallskip
\spacedlowsmallcaps{Conclusion.}
We now combine \eqref{stimaridottaAj} with the estimates for $\sum_{j=0}^JB_j$ to prove assertion (ii). As established immediately after \eqref{stimaridottaAj}, the annular operators vanish for $j>J$. Thus the argument giving \eqref{operator-norm-before-absorption}, applied only to the remaining operators, yields
\begin{equation*}
\begin{split}
\frac12N_{\lambda\delta^2}
&\overset{\eqref{stimaridottaAj}}{\leq}\Big(1-C\sum_{j=0}^JA_j\Big)N_{\lambda\delta^2}\leq C\Big(M^{d-2}+\frac1M\Big)(\lambda\delta^2)^{-(d-2)}+C\sum_{j=0}^JB_j.
\end{split}
\nonumber
\end{equation*}
Recall the operator norm identity proved immediately before \eqref{localized-operator-final-reduction}. Using that identity and \eqref{translated-operator-family}, we obtain
\begin{equation*}
\begin{split}
\|V_\lambda^{a,b}\|_{2\to2}^2
=\frac{\delta^{2(d-2)}}{|\det B|}\|T_{\lambda\delta^2}^{a,b}\|_{2\to2}^2
&\overset{\eqref{translated-operator-family}}{\leq}\frac{\delta^{2(d-2)}}{|\det B|}N_{\lambda\delta^2}\leq\frac{2C}{|\det B|}\Big(\Big(M^{d-2}+\frac1M\Big)\frac{\delta^{2(d-2)}}{(\lambda\delta^2)^{d-2}}+\delta^{2(d-2)}\sum_{j=0}^JB_j\Big)\\
&\leq C\lambda^{-(d-2)}+C\delta^{2(d-2)}\sum_{j=0}^JB_j.
\end{split}
\nonumber
\end{equation*}
Recall that $M$ has been fixed, and the bound for $B^{-1}$ proved in Section~\ref{sss:representation} before the change of variables $w=-By$ gives $|\det B|^{-1}\leq A^{3(d-2)}$ and hence the constant in the last line is uniform in $(a,b)$.

We estimate the second term separately in the three cases. If $sp>d-2$, \eqref{stimeB_jfinali} and \eqref{rescaled-gradient-sobolev-bound} give
\begin{equation*}
\begin{split}
\delta^{2(d-2)}\sum_{j=0}^JB_j
&\leq C\delta^{2(d-2)}(\lambda\delta^2)^{-(d-2)}[\nabla\widetilde Q]_{W^{s,p}}^p\overset{\eqref{rescaled-gradient-sobolev-bound}}{\leq}C\delta^{2(d-2)}(\lambda\delta^2)^{-(d-2)}\delta^{sp-(d-2)}[\nabla Q]_{W^{s,p}}^p\\
&=C\lambda^{-(d-2)}\delta^{sp-(d-2)}[\nabla Q]_{W^{s,p}}^p
\leq C\lambda^{-(d-2)}[\nabla Q]_{W^{s,p}}^p,
\end{split}
\nonumber
\end{equation*}
where the last inequality uses $0<\delta\leq1$ and $sp-(d-2)>0$.

If $sp=d-2$, then $\delta^{sp-(d-2)}=1$, and the same substitution gives
\begin{equation*}
\begin{split}
\delta^{2(d-2)}\sum_{j=0}^JB_j
&\leq C\delta^{2(d-2)}(\lambda\delta^2)^{-(d-2)}(\log(2+\lambda\delta^2))^{p-1}[\nabla\widetilde Q]_{W^{s,p}}^p\overset{\eqref{rescaled-gradient-sobolev-bound}}{\leq}C\lambda^{-(d-2)}(\log(2+\lambda\delta^2))^{p-1}[\nabla Q]_{W^{s,p}}^p\\
&\leq C\lambda^{-(d-2)}(\log(2+\lambda))^{p-1}[\nabla Q]_{W^{s,p}}^p,
\end{split}
\nonumber
\end{equation*}
since $\lambda\delta^2\leq\lambda$ and $p-1\geq0$.

Finally, if $sp<d-2$, including the case $p=1$, we obtain
\begin{equation*}
\begin{split}
\delta^{2(d-2)}\sum_{j=0}^JB_j
&\leq C\delta^{2(d-2)}(\lambda\delta^2)^{-sp}[\nabla\widetilde Q]_{W^{s,p}}^p\overset{\eqref{rescaled-gradient-sobolev-bound}}{\leq}C\delta^{2(d-2)}(\lambda\delta^2)^{-sp}\delta^{sp-(d-2)}[\nabla Q]_{W^{s,p}}^p\\
&=C\lambda^{-sp}\delta^{d-2-sp}[\nabla Q]_{W^{s,p}}^p
\leq C\lambda^{-sp}[\nabla Q]_{W^{s,p}}^p.
\end{split}
\nonumber
\end{equation*}
Here the last inequality follows from $d-2-sp>0$ and $0<\delta\leq1$.

Substituting these three estimates into the bound for $\|V_\lambda^{a,b}\|_{2\to2}^2$, we conclude that
\begin{equation}
\|V_\lambda^{a,b}\|_{2\to2}^2\leq\begin{cases}
C\lambda^{-(d-2)}\big(1+[\nabla Q]_{W^{s,p}}^p\big),&sp>d-2,\\
C\lambda^{-(d-2)}\big(1+(\log(2+\lambda))^{p-1}[\nabla Q]_{W^{s,p}}^p\big),&sp=d-2,\\
C\lambda^{-(d-2)}+C\lambda^{-sp}[\nabla Q]_{W^{s,p}}^p,&sp<d-2.
\end{cases}
\label{localized-regularity-estimates}
\end{equation}
The ranges $\lambda\delta^2<1$ and $\lambda\delta^2\sqrt{d-2}<M$ were already treated at the beginning of Section~\ref{sss:regularity} and before the definition of $J$, respectively, with the bound $C\lambda^{-(d-2)}$. Thus \eqref{localized-regularity-estimates} holds for every $\lambda\geq1$.

Returning to the finite decomposition of $V_\lambda$ established in Section~\ref{sss:normal}, the uniform bounds \eqref{localized-regularity-estimates} give the same three estimates for $\|V_\lambda\|_{2\to2}^2$, with a constant enlarged by a factor depending only on the fixed decomposition. Since $[\nabla Q]_{W^{s,p}}^p\leq[\nabla P]_{W^{s,p}}^p+[\nabla Q]_{W^{s,p}}^p$, these estimates prove assertion (ii), with constants depending only on $d,C_0,\chi_0,\chi_1,s,p$. Assertion (i) was proved in the paragraph on Dini regularity, and the proof of \cref{estimate-in-d-minus-two-variables} is complete.

\subsection{The spherical integral} 

\begin{lemma}\label{spherical-integral-expansion}
Let $d\geq2$ and $\vartheta_0\in C^\infty(\mathbb S^{d-1})$. There exist complex-valued functions $c_+,c_-$ and $\mathcal E$, defined for $w\in\mathbb R^d\setminus\{0\}$ and $\lambda\geq1$ and smooth in $w$, such that
\begin{equation}
\int_{\mathbb S^{d-1}}e^{-i\lambda\langle\omega,w\rangle}\vartheta_0(\omega)\,d\omega=\lambda^{-\frac{d-1}{2}}(e^{i\lambda|w|}c_+(w,\lambda)+e^{-i\lambda|w|}c_-(w,\lambda))+\mathcal E(w,\lambda)
\label{spherical-integral-formula}
\end{equation}
for every $w\ne0$ and $\lambda\geq1$. 
For these same functions, every compact set $K\subset\mathbb R^d\setminus\{0\}$ and every integer $L\geq0$ satisfy
$$\sup_{\lambda\geq1}\sup_{w\in K}\sum_{|\alpha|\leq L}|\partial_w^\alpha c_\pm(w,\lambda)|\leq C_{K,\vartheta_0,L}.$$
Moreover, for every such $K$ and every integer $N\geq1$,
$$\sup_{\lambda\geq1}\sup_{w\in K}\lambda^N|\mathcal E(w,\lambda)|\leq C_{K,\vartheta_0,N}.$$
\end{lemma}

\begin{proof} Let us fix $w\neq 0$ and let  $\chi\in C^\infty(\mathbb R)$ such that $0\leq\chi\leq1$, $\chi(t)=0$ for $t\leq\frac{3}{4}$, and $\chi(t)=1$ for $t\geq\frac{7}{8}$. The factors $\chi(\frac{\langle\omega,w\rangle}{|w|})$ and $\chi(-\frac{\langle\omega,w\rangle}{|w|})$ are supported near $\frac{w}{|w|}$ and $-\frac{w}{|w|}$, respectively. Let us define
$$c_\pm(w,\lambda):=\lambda^{\frac{d-1}{2}}e^{\mp i\lambda|w|}\int_{\mathbb S^{d-1}}e^{-i\lambda\langle\omega,w\rangle}\chi\Big(\mp\frac{\langle\omega,w\rangle}{|w|}\Big)\vartheta_0(\omega)\,d\omega,$$
and
$$\mathcal E(w,\lambda):=\int_{\mathbb S^{d-1}}e^{-i\lambda\langle\omega,w\rangle}\Big(1-\chi\Big(\frac{\langle\omega,w\rangle}{|w|}\Big)-\chi\Big(-\frac{\langle\omega,w\rangle}{|w|}\Big)\Big)\vartheta_0(\omega)\,d\omega,$$
and let us notice with these definitions, \eqref{spherical-integral-formula} follows immediately. The three functions have now been fixed, independently of $K,L,N$. We now check that the claimed properties hold. 

Before moving to the proof of the bounds, let us recall that the derivative on the sphere of $\omega\mapsto\langle\omega,w\rangle$ vanishes precisely at $\omega=\pm \frac{w}{|w|}$, indeed
$$\nabla_{\mathbb S^{d-1}}\langle\omega,w\rangle=w-\langle\omega,w\rangle\omega,\qquad |\nabla_{\mathbb S^{d-1}}\langle\omega,w\rangle|^2=|w|^2-\langle\omega,w\rangle^2.$$
We separate neighbourhoods of these two points. On their complement we integrate by parts, and near each of them we use the method of stationary phase in the quantitative form of \cite[Section 7.7]{HormanderI}.

\smallskip
\spacedlowsmallcaps{Estimate for $\mathcal E(w,\lambda)$.}
Fix $w_0\in K$ and choose $0<r<\frac{|w_0|}{4}$ sufficiently small. On the open set $\R^d\setminus\{0\}$, we can choose a smooth, matrix-valued function $O:\R^d\setminus\{0\}\to\R^{d\times d}$, such that $O(w)$ is an orthogonal transformation for every $w\ne0$ and $O(w)e_d=\frac{w}{|w|}$ for every $w\in B(w_0,2r)$. Substituting $\omega=O(w)(\sqrt{1-t^2}\,\nu,t)$, with $t\in(-1,1)$ and $\nu\in\mathbb S^{d-2}$, in the definition of $\mathcal E$ by $O(w)^Tw=|w|e_d$ and by the change of variables, we have
\begin{equation*}
\begin{split}
\mathcal E(w,\lambda)
&=\int_{-1}^1\int_{\mathbb S^{d-2}}e^{-i\lambda|w|t}(1-\chi(t)-\chi(-t))\vartheta_0\big(O(w)(\sqrt{1-t^2}\,\nu,t)\big)(1-t^2)^{\frac{d-3}{2}}\,d\nu\,dt
=\int_{\mathbb R}e^{-i\lambda|w|t}a(t,w)\,dt.
\end{split}
\nonumber
\end{equation*}
where
$$a(t,w):=(1-\chi(t)-\chi(-t))(1-t^2)^{\frac{d-3}{2}}\int_{\mathbb S^{d-2}}\vartheta_0\big(O(w)(\sqrt{1-t^2}\,\nu,t)\big)\,d\nu$$
for $|t|<1$, and $a(t,w):=0$ otherwise. 
By the choice of $\chi$, the function $a(\cdot,w)$ is smooth and supported in $[-\frac{7}{8},\frac{7}{8}]$. Thus integration by parts has no boundary term and gives
\begin{equation*}
\begin{split}
\int_{\mathbb R}e^{-i\lambda|w|t}a(t,w)\,dt
&=-\frac1{i\lambda|w|}\int_{\mathbb R}\partial_t(e^{-i\lambda|w|t})a(t,w)\,dt=\frac1{i\lambda|w|}\int_{\mathbb R}e^{-i\lambda|w|t}\partial_ta(t,w)\,dt.
\end{split}
\nonumber
\end{equation*}
By induction, differentiating $N$ times, we obtain
\begin{equation*}
\begin{split}
|\mathcal E(w,\lambda)|
&=\Big|(i\lambda|w|)^{-N}\int_{\mathbb R}e^{-i\lambda|w|t}\partial_t^Na(t,w)\,dt\Big|\leq(\lambda|w|)^{-N}\int_{-\frac{7}{8}}^{\frac{7}{8}}|\partial_t^Na(t,w)|\,dt\leq C_{K,\vartheta_0,N}\lambda^{-N}.
\end{split}
\nonumber
\end{equation*}
The last inequality follows because $|w|$ is bounded away from zero on $K$, while the derivatives of $a$ are bounded on compact sets.

In addition, if we assume that there exists 
$\varepsilon\in(0,1)$, such that 
$$|\langle\omega,w\rangle|\leq(1-\varepsilon)|w|\qquad\text{for every }\omega\in\supp\vartheta_0,$$
then the same integration-by-parts trick we used above implies the following estimate
\begin{equation}
\Big|\int_{\mathbb S^{d-1}}e^{-i\lambda\langle\omega,w\rangle}\vartheta_0(\omega)\,d\omega\Big|\leq C_{d,\varepsilon,\vartheta_0,N}(\lambda|w|)^{-N}.
\label{spherical-integration-by-parts}
\end{equation}
Indeed the assumption implies that $\vartheta_0\big(O(w)(\sqrt{1-t^2}\,\nu,t)\big)=0$ whenever $1-\varepsilon<|t|<1$, we have
\begin{equation*}
\begin{split}
\Big|\int_{\mathbb S^{d-1}}e^{-i\lambda\langle\omega,w\rangle}\vartheta_0(\omega)\,d\omega\Big|
&\leq(\lambda|w|)^{-N}\int_{-1+\varepsilon}^{1-\varepsilon}\int_{\mathbb S^{d-2}}\Big|\partial_t^N\Big((1-t^2)^{\frac{d-3}{2}}\vartheta_0\big(O(w)(\sqrt{1-t^2}\,\nu,t)\big)\Big)\Big|\,d\nu\,dt\\
&\leq C_{d,\varepsilon,N}\|\vartheta_0\|_{C^N(\mathbb S^{d-1})}(\lambda|w|)^{-N},
\end{split}
\nonumber
\end{equation*}
with a constant independent of $w$ and $\lambda$.

\smallskip
\spacedlowsmallcaps{Estimates for $c_\pm$.}
Fix $w_0\in K$, and let $r$ and $O$ be as in the previous paragraph. Substituting $\omega=O(w)\eta$ in the definition of $c_\pm$ and using $\langle O(w)\eta,w\rangle=|w|\eta_d$, we find that $c_\pm(w,\lambda)=\lambda^{\frac{d-1}{2}}F_\pm(\lambda|w|,w)$ for $w\in B(w_0,r)$, where
$$
F_\pm(\mu,w):=\int_{\mathbb S^{d-1}}e^{-i\mu(\eta_d\pm1)}\chi(\mp\eta_d)\vartheta_0(O(w)\eta)d\eta.
$$
The integrand vanishes unless $\mp\eta_d\geq\frac{3}{4}$, and on this cap we use the coordinates $\eta=(y,\mp\sqrt{1-|y|^2})$. The phase becomes $\mu f_\pm(y)$, with $f_\pm(y)=\pm(\sqrt{1-|y|^2}-1)$, and the amplitude becomes
$$
a_\pm(y,w):=\frac{\chi(\sqrt{1-|y|^2})\vartheta_0(O(w)(y,\mp\sqrt{1-|y|^2}))}{\sqrt{1-|y|^2}},
$$
a smooth function supported in $|y|\leq\frac{\sqrt7}{4}$, all of whose derivatives in $(y,w)$ are bounded for $w\in B(w_0,r)$. The function $f_\pm$ vanishes to second order at $y=0$, which is its only critical point, its Hessian there is $\mp$ the identity, and $\frac{|y|}{|\nabla f_\pm(y)|}=\sqrt{1-|y|^2}\leq1$. Differentiating under the integral sign, we obtain
$$
\partial_\mu^j\partial_w^\beta F_\pm(\mu,w)=\int_{\mathbb R^{d-1}}e^{i\mu f_\pm(y)}(if_\pm(y))^j\partial_w^\beta a_\pm(y,w)dy,
$$
and the amplitude on the right vanishes to order $2j$ at the origin. We apply \cite[Theorem 7.7.5]{HormanderI} with $k=j+d$. Its terms of order less than $j$ involve only derivatives of order less than $2j$ of the amplitude at the origin, so they vanish, while the remaining terms and the remainder are bounded by a constant multiple of $\mu^{-\frac{d-1}{2}-j}$ when $\mu\geq1$. The same bound holds trivially when $\frac{3|w_0|}{4}\leq\mu\leq1$, and therefore
$$
|\partial_\mu^j\partial_w^\beta F_\pm(\mu,w)|\leq C_{j,\beta}\mu^{-\frac{d-1}{2}-j}
$$
for every $\mu\geq\frac{3|w_0|}{4}$ and every $w\in B(w_0,r)$, with $C_{j,\beta}$ depending also on $d$, $\vartheta_0$, $\chi$, $O$, $w_0$ and $r$.

Now $\lambda|w|\geq\frac{3|w_0|}{4}$ whenever $\lambda\geq1$ and $w\in B(w_0,r)$. Every derivative in $w$ of $\lambda^{\frac{d-1}{2}}F_\pm(\lambda|w|,w)$ either falls on the second variable of $F_\pm$ or produces the factor $\lambda$ times a derivative of $|w|$, so for $|\alpha|\leq L$ the derivative $\partial_w^\alpha c_\pm(w,\lambda)$ is a finite sum of terms bounded by a constant multiple of
$$
\lambda^{\frac{d-1}{2}+j}|\partial_\mu^j\partial_w^\beta F_\pm(\lambda|w|,w)|\leq C_{j,\beta}\lambda^{\frac{d-1}{2}+j}(\lambda|w|)^{-\frac{d-1}{2}-j}=C_{j,\beta}|w|^{-\frac{d-1}{2}-j},
$$
where $j+|\beta|\leq|\alpha|$. The right-hand side is bounded on $B(w_0,r)$ independently of $\lambda\geq1$, and covering $K$ by finitely many such balls concludes the proof.
\end{proof}

\begin{lemma}\label{fractional-slicing}
Let us fix $d\geq3$, $m\geq1$, $0<s<1$, and $1\leq p<\infty$. Let $g\in W^{s,p}(\mathbb R^{d-1},\mathbb R^m)$. Then
$$\int_{\mathbb R}[g(\sigma,\cdot)]_{W^{s,p}(\mathbb R^{d-2})}^p\,d\sigma\leq C[g]_{W^{s,p}(\mathbb R^{d-1})}^p,$$
where $C$ depends only on $d,s,p$.
\end{lemma}

\begin{proof}
Write points of $\mathbb R^{d-1}$ as $(\sigma,z)\in\mathbb R\times\mathbb R^{d-2}$ and set $D(v):=\|g(\cdot+v)-g\|_{L^p(\mathbb R^{d-1})}$ for $v\in\mathbb R^{d-1}$. Substituting $z'=z+h$ in the slice seminorms and applying Tonelli's theorem, we see that
$$
\int_{\mathbb R}[g(\sigma,\cdot)]_{W^{s,p}(\mathbb R^{d-2})}^pd\sigma=\int_{\mathbb R^{d-2}}\frac{D((0,h))^p}{|h|^{d-2+sp}}dh.
$$
Fix $h\neq0$, and let $b_h$ be the ball of radius $\frac{|h|}{4}$ centred at $\frac{(0,h)}{2}$. For $v\in b_h$, the inequality \eqref{differenzaLptriang} yields $D((0,h))^p\leq2^{p-1}(D(v)^p+D((0,h)-v)^p)$, and since $v\mapsto(0,h)-v$ maps $b_h$ onto itself, averaging over $b_h$ leads to
$$
D((0,h))^p\leq\frac{2^p}{|b_h|}\int_{b_h}D(v)^pdv=\frac{2^p4^{d-1}}{|B(0,1)||h|^{d-1}}\int_{b_h}D(v)^pdv.
$$
We insert this bound into the previous identity and apply Tonelli's theorem once more, which shows that the left-hand side of the lemma is at most
$$
\frac{2^p4^{d-1}}{|B(0,1)|}\int_{\mathbb R^{d-1}}D(v)^p\int_{\{h:v\in b_h\}}\frac{dh}{|h|^{2d-3+sp}}dv.
$$
If $v\in b_h$, then $\frac{|h|}{4}<|v|<\frac{3|h|}{4}$, so the inner integral is at most the integral of $|h|^{-(2d-3+sp)}$ over the annulus $\frac{4|v|}{3}<|h|<4|v|$ in $\mathbb R^{d-2}$, which equals $C_{d,s,p}|v|^{-(d-1+sp)}$. By \eqref{normaWsp}, the remaining integral $\int_{\mathbb R^{d-1}}D(v)^p|v|^{-(d-1+sp)}dv$ is equal to $[g]_{W^{s,p}(\mathbb R^{d-1})}^p$, which proves the lemma for every $1\leq p<\infty$, the case $p=1$ included.
\end{proof}

\begin{proposizione}\label{spherical-graph-estimate}
Let $d\geq3$, let $U_0\Subset U_1\Subset U\subset\mathbb R^{d-1}$ be bounded open cubes, and fix $c_0>0$. There exists $\varepsilon_0>0$, depending only on $d,U_0,U_1,U,c_0$, such that the following holds.
Let $A\in C^1(U)$ satisfy $\|\nabla A\|_{L^\infty(U)}\leq\varepsilon_0$, set $X(u):=(u,A(u))$, and let $\vartheta_0\in C^\infty(\mathbb S^{d-1})$ be nonnegative. Assume that
$$\langle\omega,\partial_1X(u)\rangle\geq c_0\qquad\text{for every }u\in U_1,\ \omega\in\supp\vartheta_0.$$
For every $b\in L^2(U)$ supported in $U_0$ and every $R\geq2$, the following estimates hold.
If $\nabla A$ is Dini continuous, then
$$\int_{\mathbb S^{d-1}}\vartheta_0(\omega)\Big|\int_Ue^{-iR\langle\omega,X(u)\rangle}b(u)\,du\Big|^2\,d\omega\leq CR^{-d+\frac{3}{2}}\|b\|_{L^2(U)}^2.$$
If $\nabla A\in W^{s,p}(U)$, where $0<s<1$ and $1\leq p<\infty$, then
\begin{equation}
\int_{\mathbb S^{d-1}}\vartheta_0(\omega)\Big|\int_Ue^{-iR\langle\omega,X(u)\rangle}b(u)\,du\Big|^2\,d\omega\leq C\|b\|_{L^2(U)}^2\begin{cases}
R^{-d+\frac{3}{2}},&sp>d-2,\\
R^{-d+\frac{3}{2}}(\log(2+R))^{\frac{p}{2}},&sp=d-2,\\
R^{-\frac{d-1+sp}{2}}(\log(2+R))^{\frac{1}{2}},&sp<d-2.
\end{cases}
\label{spherical-graph-bound}
\end{equation}
In each estimate, $C$ depends only on $d,U_0,U_1,U,c_0,\vartheta_0$ and, respectively, on the Dini modulus of $\nabla A$ or on $s,p,[\nabla A]_{W^{s,p}(U)}$, but not on $R$ or $b$.
\end{proposizione}

\begin{proof} We divide the proof into several parts.

\smallskip
\spacedlowsmallcaps{Expansion of the square and reduction to two contributions.}
Write $u=(\sigma,z)$ and $u'=(\sigma',z')$, with $\sigma,\sigma'\in\mathbb R$ and $z,z'\in\mathbb R^{d-2}$, and extend $b$ by zero outside $U$. For $u,u'\in U$, define
$$K_R(u,u'):=\int_{\mathbb S^{d-1}}e^{-iR\langle\omega,X(u)-X(u')\rangle}\vartheta_0(\omega)\,d\omega,$$
and extend $K_R$ by zero outside $U\times U$ and set 
$$\mathcal B_R(\sigma,\sigma'):=\int_{\mathbb R^{d-2}}\int_{\mathbb R^{d-2}}K_R((\sigma,z),(\sigma',z'))b(\sigma,z)\overline{b(\sigma',z')}\,dz\,dz'.$$
Since $\|b\|_{L^1}\leq|U_0|^{\frac{1}{2}}\|b\|_{L^2}$, we may expand the square and apply Fubini's theorem. Moreover, setting $C_U:=\max\{1,\operatorname{diam}(U)\}$, the support of $b$ implies that $\mathcal B_R(\sigma,\sigma')=0$ whenever $|\sigma-\sigma'|>C_U$. Thus
\begin{equation}
\begin{split}
\int_{\mathbb S^{d-1}}\vartheta_0(\omega)\Big|\int_Ue^{-iR\langle\omega,X(u)\rangle}b(u)\,du\Big|^2\,d\omega
&=\int_{\mathbb R}\int_{\mathbb R}\mathcal B_R(\sigma,\sigma')\,d\sigma\,d\sigma'\\
&=\int_{|\sigma-\sigma'|\leq R^{-1}}\mathcal B_R(\sigma,\sigma')\,d\sigma\,d\sigma'+\int_{R^{-1}<|\sigma-\sigma'|\leq C_U}\mathcal B_R(\sigma,\sigma')\,d\sigma\,d\sigma'\\
&\leq\underbrace{\int_{|\sigma-\sigma'|\leq R^{-1}}|\mathcal B_R(\sigma,\sigma')|\,d\sigma\,d\sigma'}_{(I)}+\underbrace{\int_{R^{-1}<|\sigma-\sigma'|\leq C_U}|\mathcal B_R(\sigma,\sigma')|\,d\sigma\,d\sigma'}_{(II)}.
\end{split}
\label{graph-expanded-square}
\end{equation}
We estimate the first integral by Plancherel's theorem. To treat the second, we will localize $z,z'$ to cubes of side comparable to $|\sigma-\sigma'|$ and apply \cref{estimate-in-d-minus-two-variables} to the rescaled operators. That application requires an extension of $A$ with a small gradient and controlled regularity. We construct this extension after estimating the first integral.

We will fix $\varepsilon_0\leq\frac{c_0}{2}$, depending only on $d,U_0,U_1,U,c_0$, just before \cref{estimate-in-d-minus-two-variables} is applied, and the same choice is used throughout the proof.

\smallskip
\spacedlowsmallcaps{The contribution of $|\sigma-\sigma'|\leq R^{-1}$.}
We now prove that
$$\int_{|\sigma-\sigma'|\leq R^{-1}}|\mathcal B_R(\sigma,\sigma')|\,d\sigma\,d\sigma'\leq C_{\mathrm{near}}R^{-d+1}\lVert b\rVert_{L^2}^2\qquad\text{where}\qquad
C_{\mathrm{near}}:=\frac{8(2\pi)^{d-2}\|\vartheta_0\|_{L^\infty}}{c_0}.$$
Since $X(u)=(u_1,\ldots,u_{d-1},A(u))$, differentiating with respect to $u_1$ gives $\partial_1X(u)=(1,0,\ldots,0,\partial_1A(u))$ 
$$\langle\omega,\partial_1X(u)\rangle=\omega_1+\omega_d\partial_1A(u).$$
For $u\in U_1$ and $\omega\in\supp\vartheta_0$, the hypothesis $\langle\omega,\partial_1X(u)\rangle\geq c_0$, together with $|\omega_d|\leq1$, therefore yields
$$\omega_1=\langle\omega,\partial_1X(u)\rangle-\omega_d\partial_1A(u)\geq c_0-|\omega_d|\,|\partial_1A(u)|\geq c_0-\|\nabla A\|_{L^\infty(U)}\geq c_0-\varepsilon_0\geq \frac{c_0}{2},$$
where the last inequality uses $\varepsilon_0\leq \frac{c_0}{2}$. 
This shows in particular that $\supp\vartheta_0$ is contained in the set $\{\omega\in \mathbb S^{d-1}:\omega_1\geq \frac{c_0}{2}\}$. Let us write any element of $\mathbb S^{d-1}$ as $\omega=(\omega_1,\bar\omega,\omega_d)$. 
We apply the area formula to the parametrization $(\bar\omega,\omega_d)\mapsto((1-|\bar\omega|^2-\omega_d^2)^{\frac{1}{2}},\bar\omega,\omega_d)$, defined for $|\bar\omega|^2+\omega_d^2<1$. An easy computation shows  $J\omega_1=(1-|\bar\omega|^2-\omega_d^2)^{-\frac{1}{2}}=\omega_1^{-1}$, and therefore Fubini's theorem yields
\begin{equation*}
\begin{split}
&\qquad\qquad\qquad\qquad\qquad\mathcal B_R(\sigma,\sigma')
=\int_{\mathbb R^{d-2}}\int_{\mathbb R^{d-2}}K_R((\sigma,z),(\sigma',z'))b(\sigma,z)\overline{b(\sigma',z')}\,dz'\,dz\\
&\qquad\qquad\qquad\qquad=\int_{\mathbb R^{d-2}}\int_{\mathbb R^{d-2}}\int_{\mathbb S^{d-1}}e^{-iR\langle\omega,X(\sigma,z)-X(\sigma',z')\rangle}\vartheta_0(\omega)b(\sigma,z)\overline{b(\sigma',z')}\,d\omega\,dz'\,dz\\
&\qquad=\int_{\mathbb D}\frac{\vartheta_0(\omega)}{\omega_1}\int_{\mathbb R^{d-2}}\int_{\mathbb R^{d-2}}e^{-iR(\omega_1(\sigma-\sigma')+\langle\bar\omega,z-z'\rangle+\omega_d(A(\sigma,z)-A(\sigma',z')))} b(\sigma,z)\overline{b(\sigma',z')}dz'dzd\bar\omega d\omega_d\\
&=\int_{\mathbb D}\frac{e^{-iR\omega_1(\sigma-\sigma')}\vartheta_0(\omega)}{\omega_1}\Big(\int_{\mathbb R^{d-2}}e^{-iR\langle\bar\omega,z\rangle}b(\sigma,z)e^{-iR\omega_dA(\sigma,z)}dz\Big)\overline{\Big(\int_{\mathbb R^{d-2}}e^{-iR\langle\bar\omega,z'\rangle}b(\sigma',z')e^{-iR\omega_dA(\sigma',z')}dz'\Big)}d\bar\omega d\omega_d,
\nonumber
\end{split}
\end{equation*}  
where $\mathbb D:=\{(\bar\omega,\omega_d)\in \R^{d-1}:|\bar\omega|^2+\omega_d^2<1\}$. Therefore, we have 
\begin{equation*}
\begin{split}
\mathcal B_R(\sigma,\sigma')=\int_{\mathbb D}e^{-iR\omega_1(\sigma-\sigma')}\frac{\vartheta_0(\omega)}{\omega_1}\mathscr F(b(\sigma,\cdot)e^{-iR\omega_dA(\sigma,\cdot)})(R\bar\omega)\overline{\mathscr F(b(\sigma',\cdot)e^{-iR\omega_dA(\sigma',\cdot)}})(R\bar\omega)\,d\bar\omega\,d\omega_d,
\nonumber
\end{split}
\end{equation*}
where we recall that $\mathscr Ff$ stands for the Fourier transform of $f$, as defined at the beginning of the paper. 

\begin{equation}
    \begin{split}
        \int_{\mathbb R^{d-2}}\Big|\mathscr F(b(\sigma,\cdot)e^{-iR\omega_dA(\sigma,\cdot)})(R\bar\omega)\Big|^2\,d\bar\omega&=\int_{\mathbb R^{d-2}}\Big|\int_{\mathbb R^{d-2}}e^{-i\langle z,R\bar\omega\rangle}b(\sigma,z)e^{-iR\omega_dA(\sigma,z)}\,dz\Big|^2\,d\bar\omega\\
        &=R^{-(d-2)}\int_{\mathbb R^{d-2}}\Big|\int_{\mathbb R^{d-2}}e^{-i\langle z,\zeta\rangle}b(\sigma,z)e^{-iR\omega_dA(\sigma,z)}\,dz\Big|^2\,d\zeta,
    \end{split}
\end{equation}
where in the last identity we used the change of variables $R\bar \omega=\zeta$. Therefore, Plancherel's identity yields
\begin{equation*}
\begin{split}
\int_{\mathbb R^{d-2}}|\mathscr F(b(\sigma,\cdot)e^{-iR\omega_dA(\sigma,\cdot)})(R\bar\omega)|^2\,d\bar\omega
&=R^{-(d-2)}\int_{\mathbb R^{d-2}}|\mathscr F(b(\sigma,\cdot)e^{-iR\omega_dA(\sigma,\cdot)})(\zeta)|^2\,d\zeta\\
&=(2\pi)^{d-2}R^{-(d-2)}\|b(\sigma,\cdot)e^{-iR\omega_dA(\sigma,\cdot)}\|_{L^2}^2=(2\pi)^{d-2}R^{-(d-2)}\|b(\sigma,\cdot)\|_{L^2}^2.
\nonumber
\end{split}
\end{equation*}
Thanks to the above discussion, we know that $\frac{\vartheta_0(\omega)}{\omega_1}\leq\frac{2\|\vartheta_0\|_{L^\infty}}{c_0}$ and hence, applying Cauchy-Schwarz in $\bar\omega$, we consequently obtain
\begin{equation*}
\begin{split}
&|\mathcal B_R(\sigma,\sigma')|\\
&\quad\leq\frac{2\|\vartheta_0\|_{L^\infty}}{c_0}\int_{-1}^1\Big(\int_{\mathbb R^{d-2}}\Big|\mathscr F(b(\sigma,\cdot)e^{-iR\omega_dA(\sigma,\cdot)})(R\bar\omega)\Big|^2\,d\bar\omega\Big)^{\frac{1}{2}}\Big(\int_{\mathbb R^{d-2}}\Big|\mathscr F(b(\sigma',\cdot)e^{-iR\omega_dA(\sigma',\cdot)})(R\bar\omega)\Big|^2\,d\bar\omega\Big)^{\frac{1}{2}}\,d\omega_d\\
&\quad=\frac{2(2\pi)^{d-2}\|\vartheta_0\|_{L^\infty}}{c_0R^{d-2}}\int_{-1}^1\|b(\sigma,\cdot)\|_{L^2}\|b(\sigma',\cdot)\|_{L^2}\,d\omega_d=\frac{4(2\pi)^{d-2}\|\vartheta_0\|_{L^\infty}}{c_0R^{d-2}}\|b(\sigma,\cdot)\|_{L^2}\|b(\sigma',\cdot)\|_{L^2}.
\nonumber
\end{split}
\end{equation*}
Integrating the preceding estimate and using $2ab\leq a^2+b^2$, we obtain
\begin{equation}
\begin{split}
(I)&\leq\frac{2(2\pi)^{d-2}\|\vartheta_0\|_{L^\infty}}{c_0R^{d-2}}\int_{|\sigma-\sigma'|\leq R^{-1}}\Big(\|b(\sigma,\cdot)\|_{L^2}^2+\|b(\sigma',\cdot)\|_{L^2}^2\Big)\,d\sigma\,d\sigma'\\
&=\frac{8(2\pi)^{d-2}\|\vartheta_0\|_{L^\infty}}{c_0R^{d-1}}\|b\|_{L^2}^2
=C_{\mathrm{near}}R^{-d+1}\|b\|_{L^2}^2.
\label{stima(I)}
\end{split}
\end{equation}

\smallskip

\spacedlowsmallcaps{Extension of $A$.}
Before estimating $(II)$, we extend $A$ to $\mathbb R^{d-1}$, preserving its values near $U_1$ and controlling the regularity of its gradient. Choose $\zeta\in C_c^\infty(U)$, with $0\leq\zeta\leq1$ and $\zeta=1$ in a neighborhood of $U_1$, and set
$$A_U:=\frac1{|U|}\int_U A(u)\,du,\qquad
\widetilde A(u):=\begin{cases}
A_U+\zeta(u)(A(u)-A_U),&u\in U,\\
A_U,&u\notin U.
\end{cases}$$
Then $\widetilde A\in C^1(\mathbb R^{d-1})$ and $\nabla\widetilde A$ has compact support. Since $U$ is convex, $|A(u)-A_U|\leq\operatorname{diam}(U)\|\nabla A\|_{L^\infty(U)}$, and hence
$$\|\nabla\widetilde A\|_{L^\infty}\leq(1+\operatorname{diam}(U)\|\nabla\zeta\|_{L^\infty})\|\nabla A\|_{L^\infty(U)}
\leq(1+\operatorname{diam}(U)\|\nabla\zeta\|_{L^\infty})\varepsilon_0<1.$$
Thus, by choosing $\varepsilon_0$ sufficiently small, we retain the smallness required to apply \cref{estimate-in-d-minus-two-variables}.

Moreover, the standard estimates for multiplication by a smooth cutoff and extension by zero, together with $W^{1,p}(U)\subset W^{s,p}(U)$ and Poincar\'e's inequality, yield
\begin{equation}
    [\nabla\widetilde A]_{W^{s,p}(\mathbb R^{d-1})}^p+\|\nabla\widetilde A\|_{L^p(\mathbb R^{d-1})}^p
\leq C([\nabla A]_{W^{s,p}(U)}^p+\|\nabla A\|_{L^p(U)}^p).
\label{stimaestensioneWsp}
\end{equation}
Here $C$ depends only on $d,s,p,U,U_1$, see \cite[Proposition 2.2, Lemmas 5.1, 5.3]{Hitchhiker}.
In addition, the identity $\nabla\widetilde A=\zeta\nabla A+(A-A_U)\nabla\zeta$ also gives
\begin{equation}
\omega_{\nabla\widetilde A}(t)\leq\omega_{\nabla A}(t)+(2\|\nabla\zeta\|_{L^\infty}+\operatorname{diam}(U)\|D^2\zeta\|_{L^\infty})t\|\nabla A\|_{L^\infty(U)},\qquad\text{for all } 0<t\leq1,
\label{dinicontinuoustilde}
\end{equation}
In particular, Dini continuity of the gradient is preserved.

\smallskip
\spacedlowsmallcaps{The cube decomposition for $|\sigma-\sigma'|>R^{-1}$.}
Fix $\sigma,\sigma'$ with $R^{-1}<|\sigma-\sigma'|\leq\max\{1,\operatorname{diam}(U)\}$ and set $r:=|\sigma-\sigma'|$ and $I:=(-\frac{1}{2},\frac{3}{2})^{d-2}$. Fix once and for all a function $\psi\in C_c^\infty(I)$ with $0\leq\psi\leq1$ and $\sum_{k\in\mathbb Z^{d-2}}\psi(x-k)=1$ for every $x\in\mathbb R^{d-2}$. Such a function exists. Dividing a nonnegative smooth function supported in $(-\frac{1}{2},\frac{3}{2})$ and positive on $[0,1]$ by the sum of its integer translates produces a function of one variable with the same properties, and the product of $d-2$ copies of it, one in each coordinate, has the required properties. The functions $\psi_k(z):=\psi(\frac{z}{r}-k)$, $k\in\mathbb Z^{d-2}$, then form a smooth partition of unity subordinate to the open covering $\{r(k+I)\}_{k\in\mathbb Z^{d-2}}$ of $\mathbb R^{d-2}$. These cubes have overlap at most $3^{d-2}$. For $k,l\in\mathbb Z^{d-2}$, set
$$\mathcal B_R^{k,l}(\sigma,\sigma'):=\int_{\mathbb R^{d-2}}\int_{\mathbb R^{d-2}}K_R((\sigma,z),(\sigma',z'))\psi_k(z)\psi_l(z')b(\sigma,z)\overline{b(\sigma',z')}\,dz\,dz'.$$
Let us define  $\newC\label{C:graph-diameter}:=\max\{1,\operatorname{diam}(U)\}$ and $\newC\label{C:graph-adm}:=\frac{4}{c_0}+2\sqrt{d-2}$. We split $\mathcal B_R$ in two parts: the contributions given by the cubes that are close to each other and the cubes that are far away
\begin{equation}
\mathcal B_R(\sigma,\sigma')=\sum_{|k-l|\leq \oldC{C:graph-adm}}\mathcal B_R^{k,l}(\sigma,\sigma')+\sum_{|k-l|>\oldC{C:graph-adm}}\mathcal B_R^{k,l}(\sigma,\sigma').
\label{graph-cube-decomposition}
\end{equation}
Recall that $X(\sigma,z)=(\sigma,z,A(\sigma,z))$, thus
\begin{equation}
|X(\sigma,z)-X(\sigma',z')|^2
=r^2+|z-z'|^2+|A(\sigma,z)-A(\sigma',z')|^2
\geq r^2+|z-z'|^2\geq r^2\Big(1+\frac{16}{c_0^2}\Big).
\label{pdcan}
\end{equation}
where the last inequality comes from $|z-z'|>\frac{4r}{c_0}$. It follows that $\frac{|\sigma-\sigma'|}{|X(\sigma,z)-X(\sigma',z')|}
<\frac{c_0}4$  and since $\omega_1\geq \frac{c_0}{2}$ on $\supp\vartheta_0$, we have 
$$\Big|\omega\mp\frac{X(\sigma,z)-X(\sigma',z')}{|X(\sigma,z)-X(\sigma',z')|}\Big|
\geq\omega_1-\frac {|\sigma-\sigma'|}{|X(\sigma,z)-X(\sigma',z')|}
>\frac{c_0}4.$$
Writing $\omega=(\omega_1,\ldots,\omega_d)$, we have
\begin{equation}
    2\mp2\frac{\langle\omega,X(\sigma,z)-X(\sigma',z')\rangle}{|X(\sigma,z)-X(\sigma',z')|}=\Big|\omega\mp\frac{X(\sigma,z)-X(\sigma',z')}{|X(\sigma,z)-X(\sigma',z')|}\Big|^2\geq \Big(\omega_1\mp\frac{\sigma-\sigma'}{|X(\sigma,z)-X(\sigma',z')|}\Big)^2\geq\frac{c_0^2}{16}.
\end{equation}
Rearranging the preceding inequality gives
\begin{equation}
    |\langle\omega,X(\sigma,z)-X(\sigma',z')\rangle|
\leq\left(1-\frac{c_0^2}{32}\right)|X(\sigma,z)-X(\sigma',z')|,
\qquad\text{for every }\omega\in\supp\vartheta_0.
\label{stimadiffdiX}
\end{equation}
For $|k-l|>\oldC{C:graph-adm}$, $z\in r(k+I)$ and $z'\in r(l+I)$, we prove
\begin{equation}
    |K_R((\sigma,z),(\sigma',z'))|
\leq\oldC{C:graph-KR}R^{-d}r^{-d}(1+|k-l|)^{-d},
\label{stimaKRdiscretizzata}
\end{equation}
with a constant $\newC\label{C:graph-KR}$ depending only on $d,c_0,\vartheta_0$.
By definition of $K_R$, and \eqref{stimadiffdiX}, applying \eqref{spherical-integration-by-parts} with the choice $N=d$ and $\varepsilon=\min\{\frac{c_0^2}{32},\frac{1}{2}\}$, we have 
$$|K_R((\sigma,z),(\sigma',z'))|
\leq\newC R^{-d}|X(\sigma,z)-X(\sigma',z')|^{-d}\overset{\eqref{pdcan}}{\leq} \newC(1+r^{-2}|z-z'|^2)^{-\frac{d}{2}} R^{-d}r^{-d}.$$
Therefore, in order to prove \eqref{stimaKRdiscretizzata}, we just need to prove $1+|k-l|
\leq2(1+\sqrt{d-2})
(1+r^{-2}|z-z'|^2)^{\frac{1}{2}}$. Notice that thanks to the choice of $z$ and $z'$ that belong to $r(k+I)$ and $r(l+I)$ respectively, we know that $|k-\frac zr|\leq2\sqrt{d-2}$ and similarly $|\frac{z'}r-l|\leq2\sqrt{d-2}$. In particular, we have
\begin{equation}
    1+|k-l|\leq1+4\sqrt{d-2}+\frac{|z-z'|}{r}\leq 4(1+\sqrt{d-2})
(1+r^{-2}|z-z'|^2)^{\frac{1}{2}},
\end{equation}

Now, we show that there exists a constant  $\newC\label{C:graph-far-sum}>0$, depending only on
$d,c_0,\vartheta_0$, such that 
\begin{equation}
\begin{split}
\int_{R^{-1}<|\sigma-\sigma'|\leq\oldC{C:graph-diameter}}
\sum_{|k-l|>\oldC{C:graph-adm}}
|\mathcal B_R^{k,l}(\sigma,\sigma')|\,d\sigma\,d\sigma'\leq\oldC{C:graph-far-sum}R^{-d+1}\|b\|_{L^2(U)}^2.
\end{split}
\label{graph-distant-sum-bound}
\end{equation}
Thanks to the definition of $\mathcal B_R^{k,l}(\sigma,\sigma')$, we have that 
\begin{equation}
    \begin{split}
        |\mathcal B_R^{k,l}(\sigma,\sigma')|\leq&\newC R^{-d}r^{-d}(1+|k-l|)^{-d}
\int_{r(k+I)}|b(\sigma,z)|\,dz
\int_{r(l+I)}|b(\sigma',z')|\,dz'\\
\leq&\oldC{C:graph-cotto}R^{-d}r^{-2}(1+|k-l|)^{-d}
\|b(\sigma,\cdot)\|_{L^2(r(k+I))}
\|b(\sigma',\cdot)\|_{L^2(r(l+I))},
\label{Cotto}
    \end{split}
\end{equation}
where in the last inequality we applied Cauchy-Schwarz to both the integrals of $|b|$ on $r(k+I)$ and  $r(l+I)$ respectively, and the constant $\newC\label{C:graph-cotto}$ depends only on $d,c_0,\vartheta_0$.
Now, let us observe that since each point belongs to at most $3^{d-2}$ cubes $r(k+I)$, we have
$$\sum_k\|b(\sigma,\cdot)\|_{L^2(r(k+I))}^2=\sum_k\int_{r(k+I)}|b(\sigma,z)|^2\,dz=\int_{\mathbb R^{d-2}}|b(\sigma,z)|^2\sum_k\mathbf1_{r(k+I)}(z)\,dz\leq3^{d-2}\|b(\sigma,\cdot)\|_{L^2}^2.$$
Applying Cauchy-Schwarz to the sum over $k$, with $m$ fixed, yields
\begin{equation}
\begin{split}
&\sum_k\|b(\sigma,\cdot)\|_{L^2(r(k+I))}\|b(\sigma',\cdot)\|_{L^2(r(k-m+I))}
\leq\Big(\sum_k\|b(\sigma,\cdot)\|_{L^2(r(k+I))}^2\Big)^{\frac{1}{2}}\Big(\sum_k\|b(\sigma',\cdot)\|_{L^2(r(k-m+I))}^2\Big)^{\frac{1}{2}}\\
&\qquad\qquad\leq3^{\frac{d-2}{2}}\|b(\sigma,\cdot)\|_{L^2}\,3^{\frac{d-2}{2}}\|b(\sigma',\cdot)\|_{L^2}=3^{d-2}\|b(\sigma,\cdot)\|_{L^2}\|b(\sigma',\cdot)\|_{L^2}.
\label{stimaCSb}
\end{split}
\end{equation}
Therefore, we see that, with a change of indices in the sum, we have 
\begin{equation}
    \begin{split}
        &\sum_{|k-l|>\oldC{C:graph-adm}}(1+|k-l|)^{-d}\|b(\sigma,\cdot)\|_{L^2(r(k+I))}\|b(\sigma',\cdot)\|_{L^2(r(l+I))}\\&\qquad=\sum_{|m|>\oldC{C:graph-adm}}(1+|m|)^{-d}\sum_k\|b(\sigma,\cdot)\|_{L^2(r(k+I))}\|b(\sigma',\cdot)\|_{L^2(r(k-m+I))}\\
        &\qquad\leq 3^{d-2}\Big(\sum_{|m|>\oldC{C:graph-adm}}(1+|m|)^{-d}\Big)\|b(\sigma,\cdot)\|_{L^2}\|b(\sigma',\cdot)\|_{L^2}.
        \label{stimolasemola}
    \end{split}
\end{equation}
For $q=0,1,\ldots$, the annulus $2^q\leq|m|<2^{q+1}$ lies in $[-2^{q+1},2^{q+1}]^{d-2}$ and contains at most $(2^{q+2}+1)^{d-2}$ lattice points. Since each weight on this annulus is at most $(1+|m|)^{-d}\leq 2^{-qd}$, we infer that 
$$\sum_{m\in\mathbb Z^{d-2}}(1+|m|)^{-d}\leq1+ \sum_{q=0}^\infty (2^{q+2}+1)^{d-2}2^{-qd}\leq 1+5^{d-2}\sum_{q=0}^{\infty}2^{-2q}= 1+\frac43 5^{d-2}.$$
Thanks to this, we infer that 
\begin{equation}
\begin{split}
&\sum_{|k-l|>\oldC{C:graph-adm}}|\mathcal B_R^{k,l}(\sigma,\sigma')|\\&\qquad\overset{\eqref{Cotto}}{\leq} \oldC{C:graph-cotto}R^{-d}r^{-2}\sum_{|k-l|>\oldC{C:graph-adm}}(1+|k-l|)^{-d}
\|b(\sigma,\cdot)\|_{L^2(r(k+I))}
\|b(\sigma',\cdot)\|_{L^2(r(l+I))}\\
&\overset{\eqref{stimolasemola}}{\leq}\oldC{C:graph-cotto}3^{d-2}R^{-d}r^{-2}\Big(\sum_{|m|>\oldC{C:graph-adm}}(1+|m|)^{-d}\Big)\|b(\sigma,\cdot)\|_{L^2}\|b(\sigma',\cdot)\|_{L^2}\\&\qquad\leq \oldC{C:graph-stimadem}R^{-d}r^{-2}\|b(\sigma,\cdot)\|_{L^2}\|b(\sigma',\cdot)\|_{L^2}.
\label{stimadem}
\end{split}
\end{equation}
The constant $\newC\label{C:graph-stimadem}$ depends only on $d,c_0,\vartheta_0$. In order to conclude the proof of \eqref{graph-distant-sum-bound}, we integrate \eqref{stimadem} to obtain
\begin{equation}
\begin{split}
&\qquad\int_{R^{-1}<|\sigma-\sigma'|\leq\oldC{C:graph-diameter}}
\sum_{|k-l|>\oldC{C:graph-adm}}
|\mathcal B_R^{k,l}(\sigma,\sigma')|\,d\sigma\,d\sigma'\leq \oldC{C:graph-stimadem}R^{-d}
\int_{R^{-1}<|\sigma-\sigma'|\leq\oldC{C:graph-diameter}}
\frac{\|b(\sigma,\cdot)\|_{L^2}\|b(\sigma',\cdot)\|_{L^2}}
{|\sigma-\sigma'|^2}\,d\sigma\,d\sigma'\\
\leq&\frac{\oldC{C:graph-stimadem}}{2}R^{-d}
\int_{R^{-1}<|\sigma-\sigma'|\leq\oldC{C:graph-diameter}}
\frac{\|b(\sigma,\cdot)\|_{L^2}^2+\|b(\sigma',\cdot)\|_{L^2}^2}
{|\sigma-\sigma'|^2}\,d\sigma\,d\sigma'=2\oldC{C:graph-stimadem}R^{-d}\|b\|_{L^2(U)}^2
\int_{R^{-1}}^{\oldC{C:graph-diameter}}t^{-2}\,dt\leq2\oldC{C:graph-stimadem}R^{-d+1}\|b\|_{L^2(U)}^2,
\end{split}
\end{equation}
where the first inequality in the second line is a consequence of Young's inequality. 

\smallskip
We now estimate the contribution of the first sum in \eqref{graph-cube-decomposition}, using \cref{spherical-integral-expansion} and \cref{estimate-in-d-minus-two-variables}. We will show that 
\begin{equation*}
\begin{split}
&\int_{R^{-1}<|\sigma-\sigma'|\leq\oldC{C:graph-diameter}}
\sum_{|k-l|\leq\oldC{C:graph-adm}}|\mathcal B_R^{k,l}(\sigma,\sigma')|\,d\sigma\,d\sigma'\\&\qquad\leq\newC\|b\|_{L^2(U)}^2
\begin{cases}
R^{-d+\frac{3}{2}},&\text{in the Dini case, or if }sp>d-2,\\
R^{-d+\frac{3}{2}}(\log(2+R))^{\frac{p}{2}},&sp=d-2,\\
R^{-\frac{d-1+sp}{2}}(\log(2+R))^{\frac{1}{2}},&sp<d-2.
\end{cases}
\nonumber
\end{split}
\end{equation*}
The constant depends only on $d,U,U_1,c_0,\vartheta_0,\psi$, on the Dini modulus in the Dini case, and on $s,p$ and $[\nabla A]_{W^{s,p}(U)}^p+\|\nabla A\|_{L^p(U)}^p$ in the Sobolev cases. We first reduce each contribution to the norm of an oscillatory operator, prove the estimate needed for this reduction, and then apply \cref{estimate-in-d-minus-two-variables} and integrate. For each pair $k,l$ under consideration, let $T_R^{k,l}:L^2(r(l+I))\to L^2(r(k+I))$ be the operator
$$T_R^{k,l}F(z):=\psi_k(z)\int_{r(l+I)}e^{iR\sqrt{r^2+|z-z'|^2+(\widetilde A(\sigma,z)-\widetilde A(\sigma',z'))^2}}\psi_l(z')F(z')\,dz'.$$
We first prove that there is a constant $\newC\label{C:graph-spherical}$, depending only on $d,c_0,\vartheta_0,\psi$, such that
\begin{equation}
\begin{split}
|\mathcal B_R^{k,l}(\sigma,\sigma')|
&\leq\oldC{C:graph-spherical}
((Rr)^{-\frac{d-1}{2}}\|T_R^{k,l}\|_{2\to2}+r^{d-2}(Rr)^{-d})\|b(\sigma,\cdot)\|_{L^2(r(k+I))}
\|b(\sigma',\cdot)\|_{L^2(r(l+I))}.
\end{split}
\label{graph-rescaled-operator-bound}
\end{equation}
On the support of the two factors of $b$, we have $\widetilde A=A$. Since $\vartheta_0$ is real, exchanging $(\sigma,z)$ with $(\sigma',z')$ turns $K_R$ into its complex conjugate, and hence $|\mathcal B_R^{k,l}(\sigma,\sigma')|=|\mathcal B_R^{l,k}(\sigma',\sigma)|$. It is therefore enough to prove \eqref{graph-rescaled-operator-bound} when $\sigma>\sigma'$, which we assume from now on. Then $r=\sigma-\sigma'$ and by construction
\begin{equation}
\frac{X(\sigma,z)-X(\sigma',z')}{r}
=\Big(1,\frac{z-z'}r,\frac{\widetilde A(\sigma,z)-\widetilde A(\sigma',z')}r\Big).
\label{graph-rescaled-chord}
\end{equation}
Applying \cref{spherical-integral-expansion} with parameter $Rr>1$ and $N=d$, we infer that
\begin{equation}
\begin{split}
K_R((\sigma,z),(\sigma',z'))
&=(Rr)^{-\frac{d-1}{2}}\sum_\pm e^{\pm iR|X(\sigma,z)-X(\sigma',z')|}
c_\pm\Big(\frac{X(\sigma,z)-X(\sigma',z')}{r},Rr\Big)+\mathcal E\Big(\frac{X(\sigma,z)-X(\sigma',z')}{r},Rr\Big).
\end{split}
\label{graph-spherical-decomposition}
\end{equation}
The bounds for these coefficients and the remainder are uniform on the cubes under consideration. Indeed, for $z\in r(k+I)$ and $z'\in r(l+I)$, we infer that
$$\frac{|z-z'|}{r}
\leq|k-l|+\Big|\frac zr-k-\frac{z'}r+l\Big|
\leq\oldC{C:graph-adm}+2\sqrt{d-2}.$$
Moreover, $\|\nabla\widetilde A\|_{L^\infty}\leq1$ implies
$$\frac{|\widetilde A(\sigma,z)-\widetilde A(\sigma',z')|}{r}
\leq \Big(1+\frac{|z-z'|^2}{r^2}\Big)^{\frac{1}{2}}
\leq1+\oldC{C:graph-adm}+2\sqrt{d-2}.$$
Thus the last $d-1$ coordinates in \eqref{graph-rescaled-chord} belong to $[-\newC\label{C:graph-size},\oldC{C:graph-size}]^{d-1}$, where $\oldC{C:graph-size}:=1+\oldC{C:graph-adm}+4\sqrt{d-2}$ while the first coordinate is one.

To shorten the formulas, we write $\gimel_r(\sigma,\sigma',z,z'):=\frac{X(\sigma,z)-X(\sigma',z')}{r}$. We first estimate the remainder $\mathcal E$ using \cref{spherical-integral-expansion}, and then prove the bound for the term containing $c_+$ by separating its variables through a Fourier expansion. 
By \cref{spherical-integral-expansion}, with $N=d$, the quantity $2^{d-2}(Rr)^d|\mathcal E(w,Rr)|$ is bounded uniformly for $w$ in the compact set above and $Rr\geq1$. 
Since $0\leq\psi_k,\psi_l\leq1$, notice that
\begin{equation*}
\begin{split}
&\Big|\int_{r(k+I)}\int_{r(l+I)}\mathcal E(\gimel_r(\sigma,\sigma',z,z'),Rr)\psi_k(z)\psi_l(z')b(\sigma,z)\overline{b(\sigma',z')}\,dz'\,dz\Big|\\&\qquad\leq2^{-(d-2)}C(Rr)^{-d}\int_{r(k+I)}|b(\sigma,z)|\,dz\int_{r(l+I)}|b(\sigma',z')|\,dz'\\&\leq\oldC{C:graph-spherical}r^{d-2}(Rr)^{-d}\|b(\sigma,\cdot)\|_{L^2(r(k+I))}\|b(\sigma',\cdot)\|_{L^2(r(l+I))},
\nonumber
\end{split}
\end{equation*}
where in the last inequality we used the Cauchy-Schwarz inequality to estimate the integrals of $b$ by the $L^2$ norms. It remains to prove
\begin{equation}
\begin{split}
&\Big|\int_{r(k+I)}\int_{r(l+I)}e^{iRr|\gimel_r(\sigma,\sigma',z,z')|}c_+(\gimel_r(\sigma,\sigma',z,z'),Rr)\psi_k(z)\psi_l(z')b(\sigma,z)\overline{b(\sigma',z')}\,dz'\,dz\Big|\\
&\qquad\qquad\qquad\qquad\qquad\qquad\qquad\qquad\qquad\qquad\leq\oldC{C:graph-spherical}\|T_R^{k,l}\|_{2\to2}\|b(\sigma,\cdot)\|_{L^2(r(k+I))}\|b(\sigma',\cdot)\|_{L^2(r(l+I))}.
\label{simtasaoisnfaosingasoindf}
\end{split}
\end{equation}
The proof for $c_-$ is identical, because the operator with the negative phase is $F\mapsto\overline{T_R^{k,l}(\overline F)}$ and has the same norm. We therefore prove the estimate only for $c_+$.

To separate the dependence on $z$ and $z'$, we expand $c_+$ in its last $d-1$ variables. We first multiply $c_+((1,\cdot),Rr)$ by a fixed smooth cutoff $\eta$ equal to $1$ on $[-\oldC{C:graph-size},\oldC{C:graph-size}]^{d-1}$ and supported in $(-2\oldC{C:graph-size},2\oldC{C:graph-size})^{d-1}$. The product vanishes near the boundary of the larger box, so it extends smoothly with period $4\oldC{C:graph-size}$ in each coordinate.

We apply \cref{spherical-integral-expansion} with $L=2d$ on the compact set $\{1\}\times[-2\oldC{C:graph-size},2\oldC{C:graph-size}]^{d-1}$. Restricting to derivatives in the last $d-1$ coordinates gives
$$\sup_{\lambda\geq1}\sum_{|\alpha|\leq2d}\|\partial_u^\alpha[c_+((1,u),\lambda)]\|_{L^\infty((-2\oldC{C:graph-size},2\oldC{C:graph-size})^{d-1})}\leq C,$$
where $\alpha\in\mathbb N_0^{d-1}$. Since $\eta$ is fixed, the product rule then gives
$$\sup_{\lambda\geq1}\sum_{|\alpha|\leq2d}\|\partial_u^\alpha[\eta(u)c_+((1,u),\lambda)]\|_{L^\infty((-2\oldC{C:graph-size},2\oldC{C:graph-size})^{d-1})}\leq\widetilde C.$$
We apply \cref{fourier-expansions} to the function $\eta(u)c_+((1,u),Rr)$, with $q=d-1$, $N=d$ and the rectangles and $\mathcal R_0=[-\oldC{C:graph-size},\oldC{C:graph-size}]^{d-1}$ and $\mathcal R_1=(-2\oldC{C:graph-size},2\oldC{C:graph-size})^{d-1}$. The lemma supplies coefficients $a_{+,\ell}(Rr)$ such that 
$$\eta(u)c_+((1,u),Rr)=\sum_{\ell\in\mathbb Z^{d-1}}a_{+,\ell}(Rr)E_\ell(u),\qquad \text{for }u\in\mathcal R_0,\text{ and where } E_\ell(u)=e^{i(\langle\xi_\ell,u\rangle+\theta_\ell)},\, \xi_\ell\in\mathbb R^{d-1},\,\theta_\ell\in\mathbb R$$
In addition, thanks to the conclusion of \cref{fourier-expansions} and since $2d>(d-1)+1$, its coefficient estimate with weight $1+|\ell|$ gives
\begin{equation}
    \sum_{\ell\in\mathbb Z^{d-1}}(1+|\ell|)|a_{+,\ell}(Rr)|\leq C\sum_{|\alpha|\leq2d}\|\partial_u^\alpha[\eta(u)c_+((1,u),Rr)]\|_{L^\infty(\mathcal R_1)}\leq C,
    \label{stimacoefficentia(rR)}
\end{equation}
uniformly in $Rr\geq1$. Since $\eta=1$ on $\mathcal R_0$, the expansion equals $c_+((1,u),Rr)$ there.

We now use this expansion to prove \eqref{graph-rescaled-operator-bound}. Recall that $\sigma>\sigma'$, so that $r=\sigma-\sigma'$ and the first coordinate of $\gimel_r(\sigma,\sigma',z,z')$ is one. Write $\xi_\ell=(\xi_{\ell,z},\xi_{\ell,d})\in\mathbb R^{d-2}\times\mathbb R$. On the support of the two factors of $b$, the corresponding exponential separates as
$$E_\ell\Big(\frac{z-z'}r,\frac{\widetilde A(\sigma,z)-\widetilde A(\sigma',z')}r\Big)=e^{i\theta_\ell}e^{\frac{i}{r}(\langle\xi_{\ell,z},z\rangle+\xi_{\ell,d}\widetilde A(\sigma,z))}e^{-\frac{i}{r}(\langle\xi_{\ell,z},z'\rangle+\xi_{\ell,d}\widetilde A(\sigma',z'))}.$$
We identify the integral corresponding to each Fourier coefficient directly from the definition of $T_R^{k,l}$. On the support of the two factors of $b$, the identities $X(\sigma,z)=(\sigma,z,A(\sigma,z))$, $|\sigma-\sigma'|=r$ and $\widetilde A=A$ give
$$Rr|\gimel_r(\sigma,\sigma',z,z')|=R\sqrt{r^2+|z-z'|^2+(\widetilde A(\sigma,z)-\widetilde A(\sigma',z'))^2}.$$
Moreover, we have
\begin{equation*}
\begin{split}
\langle\xi_\ell,(\frac{z-z'}r,\frac{\widetilde A(\sigma,z)-\widetilde A(\sigma',z')}r)\rangle+\theta_\ell
&=\theta_\ell+\frac{\langle\xi_{\ell,z},z\rangle+\xi_{\ell,d}\widetilde A(\sigma,z)}r-\frac{\langle\xi_{\ell,z},z'\rangle+\xi_{\ell,d}\widetilde A(\sigma',z')}r.
\nonumber
\end{split}
\end{equation*}
This implies in particular that 
\begin{equation*}
\begin{split}
&T_R^{k,l}(e^{-\frac{i}{r}(\langle\xi_{\ell,z},\cdot\rangle+\xi_{\ell,d}\widetilde A(\sigma',\cdot))}\overline{b(\sigma',\cdot)})(z)\\&\qquad=\psi_k(z)\int_{r(l+I)}e^{iR\sqrt{r^2+|z-z'|^2+(\widetilde A(\sigma,z)-\widetilde A(\sigma',z'))^2}}\psi_l(z')e^{-\frac{i}{r}(\langle\xi_{\ell,z},z'\rangle+\xi_{\ell,d}\widetilde A(\sigma',z'))}\overline{b(\sigma',z')}\,dz'.
\nonumber
\end{split}
\end{equation*}
Multiplying this equality by $e^{i\theta_\ell}e^{\frac{i}{r}(\langle\xi_{\ell,z},z\rangle+\xi_{\ell,d}\widetilde A(\sigma,z))}b(\sigma,z)$ and integrating in $z$, we deduce that 
\begin{equation*}
\begin{split}
&\int_{r(k+I)}\int_{r(l+I)}e^{iRr|\gimel_r(\sigma,\sigma',z,z')|}E_\ell\Big(\frac{z-z'}r,\frac{\widetilde A(\sigma,z)-\widetilde A(\sigma',z')}r\Big)\psi_k(z)\psi_l(z')b(\sigma,z)\overline{b(\sigma',z')}\,dz'\,dz\\
&\qquad\qquad\qquad\qquad\qquad\qquad\qquad\quad=e^{i\theta_\ell}\int_{r(k+I)}e^{\frac{i}{r}(\langle\xi_{\ell,z},z\rangle+\xi_{\ell,d}\widetilde A(\sigma,z))}b(\sigma,z)T_R^{k,l}(e^{-\frac{i}{r}(\langle\xi_{\ell,z},\cdot\rangle+\xi_{\ell,d}\widetilde A(\sigma',\cdot))}\overline{b(\sigma',\cdot)})(z)\,dz.
\nonumber
\end{split}
\end{equation*}
We apply Cauchy-Schwarz in $z$ to the preceding integral and then use the operator norm of $T_R^{k,l}$
\begin{equation*}
\begin{split}
&\Big|\int_{r(k+I)}e^{\frac{i}{r}(\langle\xi_{\ell,z},z\rangle+\xi_{\ell,d}\widetilde A(\sigma,z))}b(\sigma,z)T_R^{k,l}(e^{-\frac{i}{r}(\langle\xi_{\ell,z},\cdot\rangle+\xi_{\ell,d}\widetilde A(\sigma',\cdot))}\overline{b(\sigma',\cdot)})(z)dz\Big|\\&\qquad\leq\|T_R^{k,l}\|_{2\to2}\|b(\sigma,\cdot)\|_{L^2(r(k+I))}\|b(\sigma',\cdot)\|_{L^2(r(l+I))}.
\nonumber
\end{split}
\end{equation*}
We now insert the Fourier expansion into the integral containing $c_+$ and the absolute convergence allows us to integrate term by term, thus
\begin{equation*}
\begin{split}
&\Big|\int_{r(k+I)}\int_{r(l+I)}e^{iRr|\gimel_r(\sigma,\sigma',z,z')|}c_+(\gimel_r(\sigma,\sigma',z,z'),Rr)\psi_k(z)\psi_l(z')b(\sigma,z)\overline{b(\sigma',z')}dz'dz\Big|\\
&\qquad\leq \Big|\int_{r(k+I)}\int_{r(l+I)}e^{iRr|\gimel_r(\sigma,\sigma',z,z')|}\Big(\sum_{\ell\in\mathbb Z^{d-1}}a_{+,\ell}(Rr)E_\ell\Big(\frac{z-z'}r,\frac{\widetilde A(\sigma,z)-\widetilde A(\sigma',z')}r\Big)\Big)\\&\qquad\qquad\times\psi_k(z)\psi_l(z')b(\sigma,z)\overline{b(\sigma',z')}dz'dz\Big|\\
&\qquad=\Big|\sum_{\ell\in \mathbb Z^{d-1}}a_{+,\ell}(Rr)e^{i\theta_\ell}\int_{r(k+I)}e^{\frac{i}{r}(\langle\xi_{\ell,z},z\rangle+\xi_{\ell,d}\widetilde A(\sigma,z))}b(\sigma,z)\\&\qquad\qquad\times T_R^{k,l}(e^{-\frac{i}{r}(\langle\xi_{\ell,z},\cdot\rangle+\xi_{\ell,d}\widetilde A(\sigma',\cdot))}\overline{b(\sigma',\cdot)})(z)\,dz\Big|\\
&\qquad\qquad\leq \Big(\sum_{\ell\in \mathbb Z^{d-1}}|a_{+,\ell}(Rr)|\Big)\|T_R^{k,l}\|_{2\to2}\|b(\sigma,\cdot)\|_{L^2(r(k+I))}\|b(\sigma',\cdot)\|_{L^2(r(l+I))}.
\nonumber
\end{split}
\end{equation*}
Thanks to \eqref{stimacoefficentia(rR)}, this implies the   required estimate \eqref{simtasaoisnfaosingasoindf} by choosing $\oldC{C:graph-spherical}$ so that 
$$\sum_{\ell\in\mathbb Z^{d-1}}(|a_{+,\ell}(Rr)|+|a_{-,\ell}(Rr)|)\leq\oldC{C:graph-spherical},\qquad 2^{d-2}(Rr)^d|\mathcal E(w,Rr)|\leq\oldC{C:graph-spherical}.$$
We now return to $\mathcal B_R^{k,l}$. Its definition and \eqref{graph-spherical-decomposition} give
\begin{equation*}
\begin{split}
|\mathcal B_R^{k,l}(\sigma,\sigma')|
&\leq(Rr)^{-\frac{d-1}{2}}\sum_\pm|\int_{r(k+I)}\int_{r(l+I)}e^{\pm iRr\Big|\gimel_r(\sigma,\sigma',z,z')|}c_\pm(\gimel_r(\sigma,\sigma',z,z'),Rr)\psi_k(z)\psi_l(z')b(\sigma,z)\overline{b(\sigma',z')}\,dz'\,dz\Big|\\
&\qquad\qquad\qquad\qquad\qquad\qquad+\Big|\int_{r(k+I)}\int_{r(l+I)}\mathcal E(\gimel_r(\sigma,\sigma',z,z'),Rr)\psi_k(z)\psi_l(z')b(\sigma,z)\overline{b(\sigma',z')}\,dz'\,dz\Big|\\
&\leq \oldC{C:graph-spherical}((Rr)^{-\frac{d-1}{2}}\|T_R^{k,l}\|_{2\to2}+r^{d-2}(Rr)^{-d})\|b(\sigma,\cdot)\|_{L^2(r(k+I))}\|b(\sigma',\cdot)\|_{L^2(r(l+I))}.
\nonumber
\end{split}
\end{equation*}
Notice that the constant $\oldC{C:graph-spherical}$ depends only on $d,c_0,\vartheta_0$ and the fixed cutoffs. This proves \eqref{graph-rescaled-operator-bound}.

We have thus reduced the problem to estimate $\|T_R^{k,l}\|_{2\to2}$ and in order to do this, we will apply \cref{estimate-in-d-minus-two-variables}. We now put $T_R^{k,l}$ in the form required by \cref{estimate-in-d-minus-two-variables}. Since $\psi_k(z)=\psi(\frac{z}{r}-k)$, the substitutions $z=r(k+x)$ and $z'=r(l+y)$ give
\begin{equation*}
\begin{split}
T_R^{k,l}F(r(k+x))
=r^{d-2}\psi(x)\int_I e^{iRr\big(1+|k-l+x-y|^2+\big(\frac{\widetilde A(\sigma,r(k+x))-\widetilde A(\sigma',r(l+y))}{r}\big)^2\big)^{\frac{1}{2}}}\psi(y)F(r(l+y))\,dy,
\nonumber
\end{split}
\end{equation*}
and let us set  $\lambda:=Rr$, $v:=k-l$ and finally $P(x):=r^{-1}\widetilde A(\sigma,r(k+x))$ and $ Q(y):=r^{-1}\widetilde A(\sigma',r(l+y))$ and with these notations, we have 
\begin{equation}
    T_R^{k,l}F(r(k+x))=r^{d-2}\psi(x)\int_I e^{i\lambda\sqrt{1+|v+x-y|^2+(P(x)-Q(y))^2}}\psi(y)F(r(l+y))\,dy=:r^{d-2}\mathcal T_\lambda(F(r(l+\cdot)))(x).
    \label{changeofcoordinatesoperator}
\end{equation}
The operator $\mathcal T_\lambda$ is exactly the operator in \cref{estimate-in-d-minus-two-variables}, with $\chi_0=\chi_1=\psi$. We now verify the hypothesis of \cref{estimate-in-d-minus-two-variables} for the quantities just defined. First, $\lambda=Rr>1$ and $P,Q\in C^1(\mathbb R^{d-2})$. We can now fix $\varepsilon_0\leq\frac{c_0}{2}$, depending only on $d,U_0,U_1,U,c_0$, so small that
$2(1+\operatorname{diam}(U)\|\nabla\zeta\|_{L^\infty})\varepsilon_0\leq\varepsilon_*(d,2\oldC{C:graph-size})$. Since $\varepsilon_*\leq1$, this choice implies the bound $\|\nabla\widetilde A\|_{L^\infty}<1$ used in the construction of $\widetilde A$, as well as $\|\nabla P\|_{L^\infty}+\|\nabla Q\|_{L^\infty}\leq2\|\nabla\widetilde A\|_{L^\infty}\leq\varepsilon_*(d,2\oldC{C:graph-size})$, which is the first condition in \eqref{stimenecessarioes}. For every $x,y\in I$, it also yields
\begin{equation*}
\begin{split}
|v|+|P(x)-Q(y)|
&=|k-l|+r^{-1}|\widetilde A(\sigma,r(k+x))-\widetilde A(\sigma',r(l+y))|\leq|k-l|+\|\nabla\widetilde A\|_{L^\infty}(1+|k-l+x-y|^2)^{\frac{1}{2}}\\
&\leq1+2|k-l|+|x-y|\leq1+2\oldC{C:graph-adm}+2\sqrt{d-2}\leq2\oldC{C:graph-size}.
\nonumber
\end{split}
\end{equation*}
Here we used $|k-l|\leq\oldC{C:graph-adm}$, the bound $|x-y|\leq2\sqrt{d-2}$ for $x,y\in I$, and the definition $\oldC{C:graph-size}$. Taking the supremum over $x,y\in I$ therefore verifies the second condition in \eqref{stimenecessarioes} with $C_0=2\oldC{C:graph-size}$. 
Now we verify separately the hypothesis of items (i) and (ii) of \cref{estimate-in-d-minus-two-variables}.

\spacedlowsmallcaps{The Dini regularity case.} Suppose that $\nabla A$ is Dini continuous. The estimate in \eqref{dinicontinuoustilde}, shows that $\nabla\widetilde{A}$ is Dini continuous. Moreover, for every $x,x'\in\mathbb R^{d-2}$, there holds
$$|\nabla P(x)-\nabla P(x')|=|\nabla_z\widetilde A(\sigma,r(k+x))-\nabla_z\widetilde A(\sigma,r(k+x'))|\leq\omega_{\nabla \widetilde A}(r|x-x'|),$$
and the same estimate holds for $\nabla Q$. Thus, both gradients are controlled by the common modulus $t\mapsto\omega_{\nabla\widetilde A}(rt)$. Since $r\leq\oldC{C:graph-diameter}$ and $\omega_{\nabla\widetilde A}$ is nondecreasing, this modulus is at most $\omega_{\nabla\widetilde A}(\oldC{C:graph-diameter}t)$, which no longer depends on $r$. Integrating \eqref{dinicontinuoustilde} on $(0,1)$ and using $\omega_{\nabla\widetilde A}\leq2\|\nabla\widetilde A\|_{L^\infty}\leq2$ on $[1,\oldC{C:graph-diameter}]$, we obtain, with $C_\zeta:=2\|\nabla\zeta\|_{L^\infty}+\operatorname{diam}(U)\|D^2\zeta\|_{L^\infty}$,
$$
\begin{aligned}
\int_0^1\omega_{\nabla\widetilde A}(\oldC{C:graph-diameter}t)\frac{dt}{t}&=\int_0^{\oldC{C:graph-diameter}}\omega_{\nabla\widetilde A}(t)\frac{dt}{t}\\
&\leq\int_0^1\omega_{\nabla A}(t)\frac{dt}{t}+C_\zeta\|\nabla A\|_{L^\infty(U)}+2\log\oldC{C:graph-diameter}<\infty.
\end{aligned}
$$
Item (i) of \cref{estimate-in-d-minus-two-variables}, applied with this modulus, therefore gives
$\|\mathcal T_\lambda\|_{2\to2}^2\leq\newC\label{C:graph-op-D}\lambda^{-(d-2)}$, with a constant independent of $\sigma,\sigma',r,k,l$. To obtain back the estimate for $T_R^{k,l}$, we use the identity \eqref{changeofcoordinatesoperator}. Notice that thanks to \eqref{changeofcoordinatesoperator}, we infer in particular that 
$r^{\frac{d-2}{2}}(T_R^{k,l}F)(r(k+x))=r^{d-2}\mathcal T_\lambda(r^{\frac{d-2}{2}}F(r(l+\cdot)))(x)$ and hence a change of variables shows that 
$$\|T_R^{k,l}\|_{2\to2}=r^{d-2}\|\mathcal T_\lambda\|_{2\to2}.$$
Using $\lambda=Rr$, we conclude that
$\|T_R^{k,l}\|_{2\to2}^2=r^{2(d-2)}\|\mathcal T_\lambda\|_{2\to2}^2\leq\oldC{C:graph-op-D}R^{-(d-2)}r^{d-2}$.

\smallskip
\spacedlowsmallcaps{Estimates in the fractional Sobolev regularity case.}
Suppose now that $\nabla A\in W^{s,p}(U)$, where $0<s<1$ and $1\leq p<\infty$. To apply item (ii) of \cref{estimate-in-d-minus-two-variables}, we first verify that $\nabla P,\nabla Q\in W^{s,p}(\mathbb R^{d-2})$ for almost every $\sigma,\sigma'$, and then compute their seminorms. The estimate \eqref{stimaestensioneWsp} for the extension $\widetilde A$ shows that $\nabla_z\widetilde A\in W^{s,p}(\mathbb R^{d-1},\mathbb R^{d-2})$, since taking components does not increase either the fractional seminorm or the $L^p$ norm. We may therefore apply \cref{fractional-slicing} to $\nabla_z\widetilde A$, which bounds the integral of the slice seminorms by the global seminorm $[\nabla_z\widetilde A]_{W^{s,p}(\mathbb R^{d-1})}^p$. For the $L^p$ norms, Fubini's theorem gives
$\int_{\mathbb R}\|\nabla_z\widetilde A(\sigma,\cdot)\|_{L^p(\mathbb R^{d-2})}^p\,d\sigma=\|\nabla_z\widetilde A\|_{L^p(\mathbb R^{d-1})}^p$, and hence 
combining these two relations with the estimate for the extension, we obtain
\begin{equation}
\begin{split}
\int_{\mathbb R}([\nabla_z\widetilde A(\sigma,\cdot)]_{W^{s,p}(\mathbb R^{d-2})}^p+\|\nabla_z\widetilde A(\sigma,\cdot)\|_{L^p(\mathbb R^{d-2})}^p)\,d\sigma
&\leq\oldC{C:graph-slices}([\nabla A]_{W^{s,p}(U)}^p+\|\nabla A\|_{L^p(U)}^p).
\label{graph-slice-bound}
\end{split}
\end{equation}
The constant $\newC\label{C:graph-slices}$ depends only on $d,s,p,U,U_1$.
Using the definition of $P$ and substituting $z=r(k+x)$ and $z'=r(k+x')$ in the double integral defining its seminorm, we obtain
\begin{equation*}
\begin{split}
[\nabla P]_{W^{s,p}(\mathbb R^{d-2})}^p
&=\int_{\mathbb R^{d-2}}\int_{\mathbb R^{d-2}}\frac{|\nabla_z\widetilde A(\sigma,r(k+x))-\nabla_z\widetilde A(\sigma,r(k+x'))|^p}{|x-x'|^{d-2+sp}}dx'dx=r^{sp-d+2}[\nabla_z\widetilde A(\sigma,\cdot)]_{W^{s,p}(\mathbb R^{d-2})}^p.
\nonumber
\end{split}
\end{equation*}
and that 
$\|\nabla P\|_{L^p(\mathbb R^{d-2})}^p=r^{-(d-2)}\|\nabla_z\widetilde A(\sigma,\cdot)\|_{L^p(\mathbb R^{d-2})}^p$. The same calculation applies to $Q$. 
Thus, the regularity hypothesis of item (ii) holds for almost every $\sigma,\sigma'$, and thanks to the above formula, we have
\begin{equation}
\begin{split}
[\nabla P]_{W^{s,p}(\mathbb R^{d-2})}^p+[\nabla Q]_{W^{s,p}(\mathbb R^{d-2})}^p
&=r^{sp-d+2}([\nabla_z\widetilde A(\sigma,\cdot)]_{W^{s,p}(\mathbb R^{d-2})}^p+[\nabla_z\widetilde A(\sigma',\cdot)]_{W^{s,p}(\mathbb R^{d-2})}^p).
\label{idbablaqbablaP}
\end{split}
\end{equation}
We now apply item (ii) of \cref{estimate-in-d-minus-two-variables} to $\mathcal T_\lambda$ and compute explicitly the factors needed to return to $T_R^{k,l}$. The proposition gives
\begin{equation*}
\begin{split}
\|\mathcal T_\lambda\|_{2\to2}^2
&\leq C\lambda^{-(d-2)}+C([\nabla P]_{W^{s,p}(\mathbb R^{d-2})}^p+[\nabla Q]_{W^{s,p}(\mathbb R^{d-2})}^p)
\begin{cases}
\lambda^{-(d-2)},&sp>d-2,\\
\lambda^{-(d-2)}(\log(2+\lambda))^{p-1},&sp=d-2,\\
\lambda^{-sp},&sp<d-2.
\end{cases}
\nonumber
\end{split}
\end{equation*}
The identity $\|T_R^{k,l}\|_{2\to2}^2=r^{2(d-2)}\|\mathcal T_\lambda\|_{2\to2}^2$ and \eqref{idbablaqbablaP}, and a few algebraic computations imply that
\begin{equation}
\begin{split}
\|T_R^{k,l}\|_{2\to2}^2
&\leq\oldC{C:graph-op-S}R^{-(d-2)}r^{d-2}\\&\qquad+\oldC{C:graph-op-S}([\nabla_z\widetilde A(\sigma,\cdot)]_{W^{s,p}(\mathbb R^{d-2})}^p+[\nabla_z\widetilde A(\sigma',\cdot)]_{W^{s,p}(\mathbb R^{d-2})}^p)\\&\qquad\qquad\times
\begin{cases}
R^{-(d-2)}r^{sp},&sp>d-2,\\
R^{-(d-2)}r^{sp}(\log(2+Rr))^{p-1},&sp=d-2,\\
R^{-sp}r^{d-2},&sp<d-2.
\end{cases}
\label{graph-sobolev-operator-bound}
\end{split}
\end{equation}
The constant $\newC\label{C:graph-op-S}$ depends only on $d,\oldC{C:graph-size},\psi,s,p$, and is therefore independent of $\sigma,\sigma',r,k,l$. Indeed, the smallness and size hypotheses of item (ii) of \cref{estimate-in-d-minus-two-variables} were verified uniformly with $C_0=2\oldC{C:graph-size}$, and the dependence on the gradient seminorms of $P,Q$ is explicit in its estimate. The slice seminorms on the right-hand side are controlled only in integral by the estimate above.

\smallskip
\spacedlowsmallcaps{Conclusion.}
We now combine the estimates for the two contributions in \eqref{graph-expanded-square}. Recall that
\begin{equation*}
\begin{split}
&\int_{\mathbb S^{d-1}}\vartheta_0(\omega)\Big|\int_Ue^{-iR\langle\omega,X(u)\rangle}b(u)\,du\Big|^2\,d\omega
\\&\qquad\leq\int_{|\sigma-\sigma'|\leq R^{-1}}|\mathcal B_R(\sigma,\sigma')|\,d\sigma\,d\sigma'\\&\qquad\quad+\int_{R^{-1}<|\sigma-\sigma'|\leq\oldC{C:graph-diameter}}|\mathcal B_R(\sigma,\sigma')|\,d\sigma\,d\sigma'=(I)+(II)\\
&\qquad\overset{\substack{\eqref{stima(I)},\eqref{graph-cube-decomposition}\\\eqref{graph-distant-sum-bound}}}{\leq} (C_{\mathrm{near}}+\oldC{C:graph-far-sum})R^{-d+1}\|b\|_{L^2(U)}^2\\&\qquad\quad+\int_{R^{-1}<|\sigma-\sigma'|\leq\oldC{C:graph-diameter}}\sum_{|k-l|\leq\oldC{C:graph-adm}}|\mathcal B_R^{k,l}(\sigma,\sigma')|\,d\sigma\,d\sigma'.
\nonumber
\end{split}
\end{equation*}
We estimate the last integral separately under the Dini and fractional Sobolev assumptions. If $\sigma>\sigma'$, we set $r=\sigma-\sigma'$ and use $\|T_R^{k,l}\|_{2\to2}\leq\oldC{C:graph-op-D}^{\frac{1}{2}}R^{-\frac{d-2}{2}}r^{\frac{d-2}{2}}$. Thus,
\begin{equation}
\begin{split}
&\int_{R^{-1}<|\sigma-\sigma'|\leq\oldC{C:graph-diameter}}\sum_{|k-l|\leq\oldC{C:graph-adm}}|\mathcal B_R^{k,l}(\sigma,\sigma')|\,d\sigma\,d\sigma'=2\int_{R^{-1}<\sigma-\sigma'\leq\oldC{C:graph-diameter}}\sum_{|k-l|\leq\oldC{C:graph-adm}}|\mathcal B_R^{k,l}(\sigma,\sigma')|\,d\sigma\,d\sigma'\\
&\qquad\overset{\eqref{graph-rescaled-operator-bound}}{\leq}2\oldC{C:graph-spherical}\int_{R^{-1}<\sigma-\sigma'\leq\oldC{C:graph-diameter}}\sum_{|k-l|\leq\oldC{C:graph-adm}}\big((Rr)^{-\frac{d-1}{2}}\|T_R^{k,l}\|_{2\to2}+r^{d-2}(Rr)^{-d}\big)\|b(\sigma,\cdot)\|_{L^2(r(k+I))}\|b(\sigma',\cdot)\|_{L^2(r(l+I))}\,d\sigma\,d\sigma'\\
&\qquad\leq2\oldC{C:graph-spherical}\int_{R^{-1}<\sigma-\sigma'\leq\oldC{C:graph-diameter}}\big(\oldC{C:graph-op-D}^{\frac{1}{2}}R^{-d+\frac{3}{2}}r^{-\frac{1}{2}}+R^{-d}r^{-2}\big)\sum_{|k-l|\leq\oldC{C:graph-adm}}\|b(\sigma,\cdot)\|_{L^2(r(k+I))}\|b(\sigma',\cdot)\|_{L^2(r(l+I))}\,d\sigma\,d\sigma'.
\label{graph-dini-pair-reduction}
\end{split}
\end{equation}
To integrate the bound just obtained, we first apply \eqref{stimaCSb} for each fixed difference $m=k-l$ and sum over the at most $(2\lceil\oldC{C:graph-adm}\rceil+1)^{d-2}$ admissible differences. Thus, the elementary inequality $2ab\leq a^2+b^2$ gives
\begin{equation}
\begin{split}
&\int_{R^{-1}<|\sigma-\sigma'|\leq\oldC{C:graph-diameter}}\sum_{|k-l|\leq\oldC{C:graph-adm}}|\mathcal B_R^{k,l}(\sigma,\sigma')|\,d\sigma\,d\sigma'\\
&\qquad\overset{\eqref{stimaCSb}}{\leq}2\oldC{C:graph-spherical}3^{d-2}(2\lceil\oldC{C:graph-adm}\rceil+1)^{d-2}\int_{R^{-1}<\sigma-\sigma'\leq\oldC{C:graph-diameter}}\big(\oldC{C:graph-op-D}^{\frac{1}{2}}R^{-d+\frac{3}{2}}r^{-\frac{1}{2}}+R^{-d}r^{-2}\big)\|b(\sigma,\cdot)\|_{L^2}\|b(\sigma',\cdot)\|_{L^2}\,d\sigma\,d\sigma'\\
&\leq\oldC{C:graph-spherical}3^{d-2}(2\lceil\oldC{C:graph-adm}\rceil+1)^{d-2}\int_{R^{-1}<\sigma-\sigma'\leq\oldC{C:graph-diameter}}\big(\oldC{C:graph-op-D}^{\frac{1}{2}}R^{-d+\frac{3}{2}}r^{-\frac{1}{2}}+R^{-d}r^{-2}\big)\big(\|b(\sigma,\cdot)\|_{L^2}^2+\|b(\sigma',\cdot)\|_{L^2}^2\big)\,d\sigma\,d\sigma'.
\label{graph-dini-pair-summation}
\end{split}
\end{equation}
We now substitute $r=\sigma-\sigma'$ and apply Fubini's theorem, integrating first in $\sigma'$ in the term containing $\|b(\sigma,\cdot)\|_{L^2}^2$ and first in $\sigma$ in the other term. Computing the inner integrals gives
\begin{equation}
\begin{split}
&\int_{R^{-1}<|\sigma-\sigma'|\leq\oldC{C:graph-diameter}}\sum_{|k-l|\leq\oldC{C:graph-adm}}|\mathcal B_R^{k,l}(\sigma,\sigma')|\,d\sigma\,d\sigma'\\
&\qquad\leq\oldC{C:graph-spherical}3^{d-2}(2\lceil\oldC{C:graph-adm}\rceil+1)^{d-2}\int_{\mathbb R}\|b(\sigma,\cdot)\|_{L^2}^2\int_{\sigma-\oldC{C:graph-diameter}}^{\sigma-R^{-1}}\Big(\frac{\oldC{C:graph-op-D}^{\frac{1}{2}}R^{-d+\frac{3}{2}}}{(\sigma-\sigma')^{\frac{1}{2}}}+\frac{R^{-d}}{(\sigma-\sigma')^2}\Big)d\sigma'd\sigma\\
&\qquad\qquad\qquad\qquad\qquad\qquad+\oldC{C:graph-spherical}3^{d-2}(2\lceil\oldC{C:graph-adm}\rceil+1)^{d-2}\int_{\mathbb R}\|b(\sigma',\cdot)\|_{L^2}^2\int_{\sigma'+R^{-1}}^{\sigma'+\oldC{C:graph-diameter}}\Big(\frac{\oldC{C:graph-op-D}^{\frac{1}{2}}R^{-d+\frac{3}{2}}}{(\sigma-\sigma')^{\frac{1}{2}}}+\frac{R^{-d}}{(\sigma-\sigma')^2}\Big)d\sigma d\sigma'\\
&\qquad=2\oldC{C:graph-spherical}3^{d-2}(2\lceil\oldC{C:graph-adm}\rceil+1)^{d-2}(2\oldC{C:graph-op-D}^{\frac{1}{2}}R^{-d+\frac{3}{2}}(\oldC{C:graph-diameter}^{\frac{1}{2}}-R^{-\frac{1}{2}})+R^{-d}(R-\oldC{C:graph-diameter}^{-1}))\|b\|_{L^2(U)}^2\\
&\qquad\leq2\oldC{C:graph-spherical}3^{d-2}(2\lceil\oldC{C:graph-adm}\rceil+1)^{d-2}(2\oldC{C:graph-op-D}^{\frac{1}{2}}\oldC{C:graph-diameter}^{\frac{1}{2}}R^{-d+\frac{3}{2}}+R^{-d+1})\|b\|_{L^2(U)}^2.
\label{graph-dini-scalar-bound}
\end{split}
\end{equation}
Substituting the last estimate into the bound for $(I)+(II)$ obtained above, and using $R^{-d+1}\leq R^{-d+\frac{3}{2}}$, we have therefore proved that, under the Dini assumption on $\nabla A$,
$$\int_{\mathbb S^{d-1}}\vartheta_0(\omega)\Big|\int_Ue^{-iR\langle\omega,X(u)\rangle}b(u)\,du\Big|^2\,d\omega\overset{\eqref{graph-expanded-square}}{\leq}(I)+(II)\leq C R^{-d+\frac{3}{2}}\|b\|_{L^2(U)}^2,$$
where
$$C:=C_{\mathrm{near}}+\oldC{C:graph-far-sum}+2\oldC{C:graph-spherical}3^{d-2}(2\lceil\oldC{C:graph-adm}\rceil+1)^{d-2}\big(2\oldC{C:graph-op-D}^{\frac{1}{2}}\oldC{C:graph-diameter}^{\frac{1}{2}}+1\big).$$
Since all the constants in this expression have already been fixed, $C$ depends only on $d,U_0,U_1,U,c_0,\vartheta_0$ and the Dini modulus of $\nabla A$, and is independent of $R$ and $b$. This proves the Dini assertion of the proposition.

Let us move to the Sobolev case. 
Suppose now that $\nabla A\in W^{s,p}(U)$. We first reduce the estimate for the admissible pairs to the $L^2$ norm of a one-dimensional kernel, and then compute this norm in the three ranges of $sp$. Define
$$\mathfrak{k}_R(t):=
\begin{cases}
R^{-d+\frac{3}{2}}|t|^{\frac{sp-d+1}{2}},&sp>d-2,\\
R^{-d+\frac{3}{2}}|t|^{-\frac{1}{2}}(\log(2+R|t|))^{\frac{p-1}{2}},&sp=d-2,\\
R^{-\frac{d-1+sp}{2}}|t|^{-\frac{1}{2}},&sp<d-2,
\end{cases}
\qquad \text{for every }R^{-1}<|t|\leq\oldC{C:graph-diameter},$$
and set $\mathfrak{k}_R(t)=0$ elsewhere. 
On $R^{-1}<\sigma-\sigma'\leq\oldC{C:graph-diameter}$, recall that $r=\sigma-\sigma'$ and multiplying \eqref{graph-sobolev-operator-bound} by $(Rr)^{-(d-1)}$ gives
\begin{equation*}
\begin{split}
(Rr)^{-(d-1)}\|T_R^{k,l}\|_{2\to2}^2
&\overset{\eqref{graph-sobolev-operator-bound}}{\leq}
\oldC{C:graph-op-S}(
R^{-2d+3}r^{-1}+\mathfrak{k}_R(r)^2
([\nabla_z\widetilde A(\sigma,\cdot)]_{W^{s,p}(\mathbb R^{d-2})}^p
+[\nabla_z\widetilde A(\sigma',\cdot)]_{W^{s,p}(\mathbb R^{d-2})}^p)).
\nonumber
\end{split}
\end{equation*}
Since $\sqrt{a+b}\leq\sqrt a+\sqrt b$, we have in particular that 
\begin{equation*}
\begin{split}
(Rr)^{-\frac{d-1}{2}}\|T_R^{k,l}\|_{2\to2}&\leq\oldC{C:graph-op-S}^{\frac{1}{2}}(R^{-d+\frac{3}{2}}r^{-\frac{1}{2}}+\mathfrak{k}_R(r)([\nabla_z\widetilde A(\sigma,\cdot)]_{W^{s,p}(\mathbb R^{d-2})}^p+[\nabla_z\widetilde A(\sigma',\cdot)]_{W^{s,p}(\mathbb R^{d-2})}^p)^{\frac{1}{2}}).
\nonumber
\end{split}
\end{equation*}
We repeat the symmetry argument and the operator substitution of \eqref{graph-dini-pair-reduction} and the summation over admissible differences of \eqref{graph-dini-pair-summation}, now using \eqref{graph-sobolev-operator-bound} in place of the Dini operator estimate. The scalar terms are integrated as in \eqref{graph-dini-scalar-bound}, with $\oldC{C:graph-op-S}^{\frac{1}{2}}$ in place of $\oldC{C:graph-op-D}^{\frac{1}{2}}$. For the additional term containing $\mathfrak{k}_R$, its symmetry and support allow us to write the integral over $\mathbb R^2$. Consequently,
\begin{equation*}
\begin{split}
&\int_{R^{-1}<|\sigma-\sigma'|\leq\oldC{C:graph-diameter}}\sum_{|k-l|\leq\oldC{C:graph-adm}}|\mathcal B_R^{k,l}(\sigma,\sigma')|\,d\sigma\,d\sigma'\\&\overset{\eqref{graph-dini-pair-reduction},\eqref{graph-dini-pair-summation},\eqref{graph-dini-scalar-bound}}{\leq}\oldC{C:graph-spherical}3^{d-2}(2\oldC{C:graph-adm}+1)^{d-2}\big(4\oldC{C:graph-op-S}^{\frac{1}{2}}\oldC{C:graph-diameter}^{\frac{1}{2}}R^{-d+\frac{3}{2}}+2R^{-d+1}\big)\|b\|_{L^2(U)}^2\\
&\qquad+\oldC{C:graph-op-S}^{\frac{1}{2}}\oldC{C:graph-spherical}3^{d-2}(2\oldC{C:graph-adm}+1)^{d-2}\\&\qquad\quad\times\int_{\mathbb R^2}\mathfrak{k}_R(\sigma-\sigma')\big([\nabla_z\widetilde A(\sigma,\cdot)]_{W^{s,p}(\mathbb R^{d-2})}^p+[\nabla_z\widetilde A(\sigma',\cdot)]_{W^{s,p}(\mathbb R^{d-2})}^p\big)^{\frac{1}{2}}\\&\qquad\qquad\times\|b(\sigma,\cdot)\|_{L^2}\|b(\sigma',\cdot)\|_{L^2}d\sigma\,d\sigma'.
\nonumber
\end{split}
\end{equation*}
We now estimate the last integral by Cauchy-Schwarz. Fubini's theorem implies that
$$\int_{\mathbb R^2}\mathfrak{k}_R(\sigma-\sigma')^2
[\nabla_z\widetilde A(\sigma,\cdot)]_{W^{s,p}(\mathbb R^{d-2})}^p
\,d\sigma\,d\sigma'
=\|\mathfrak{k}_R\|_{L^2(\mathbb R)}^2
\int_{\mathbb R}[\nabla_z\widetilde A(\sigma,\cdot)]_{W^{s,p}(\mathbb R^{d-2})}^p\,d\sigma,$$
and the same identity holds for the term containing the slice at $\sigma'$. Moreover, we have
$$
\int_{\mathbb R^2}\|b(\sigma,\cdot)\|_{L^2}^2
\|b(\sigma',\cdot)\|_{L^2}^2\,d\sigma\,d\sigma'
=\|b\|_{L^2(U)}^4
$$
and consequently, we have
\begin{equation*}
\begin{split}
&\int_{\mathbb R^2}\mathfrak{k}_R(\sigma-\sigma')
\big([\nabla_z\widetilde A(\sigma,\cdot)]_{W^{s,p}(\mathbb R^{d-2})}^p
+[\nabla_z\widetilde A(\sigma',\cdot)]_{W^{s,p}(\mathbb R^{d-2})}^p\big)^{\frac{1}{2}}
\|b(\sigma,\cdot)\|_{L^2}\|b(\sigma',\cdot)\|_{L^2}
\,d\sigma\,d\sigma'\\
&\qquad\leq \sqrt2\|\mathfrak{k}_R\|_{L^2(\mathbb R)}
\Big(\int_{\mathbb R}
[\nabla_z\widetilde A(\sigma,\cdot)]_{W^{s,p}(\mathbb R^{d-2})}^p\,d\sigma\Big)^{\frac{1}{2}}
\|b\|_{L^2(U)}^2\\&\qquad\overset{\eqref{graph-slice-bound}}{\leq}
\big(2\oldC{C:graph-slices}\big)^{\frac{1}{2}}
\|\mathfrak{k}_R\|_{L^2(\mathbb R)}
\big([\nabla A]_{W^{s,p}(U)}^p+\|\nabla A\|_{L^p(U)}^p\big)^{\frac{1}{2}}
\|b\|_{L^2(U)}^2,
\nonumber
\end{split}
\end{equation*}
where the first inequality in the second line above comes from the Cauchy-Schwarz inequality. Substituting the last bound into the estimate for the admissible pairs, we obtain
\begin{equation*}
\begin{split}
&\int_{\mathbb S^{d-1}}\vartheta_0(\omega)
\Big|\int_Ue^{-iR\langle\omega,X(u)\rangle}b(u)\,du\Big|^2\,d\omega\overset{\eqref{graph-expanded-square}}{\leq}(I)+(II)\\
&\qquad\overset{\substack{
\eqref{stima(I)},\,\eqref{graph-cube-decomposition}\\
\eqref{graph-distant-sum-bound}}}{\leq}
(C_{\mathrm{near}}+\oldC{C:graph-far-sum})
R^{-d+1}\|b\|_{L^2(U)}^2\\&\qquad\quad+
\int_{R^{-1}<|\sigma-\sigma'|\leq\oldC{C:graph-diameter}}
\sum_{|k-l|\leq\oldC{C:graph-adm}}
|\mathcal B_R^{k,l}(\sigma,\sigma')|\,d\sigma\,d\sigma'\\&\qquad\leq C(R^{-d+\frac{3}{2}}
+\|\mathfrak{k}_R\|_{L^2(\mathbb R)})\|b\|_{L^2(U)}^2.
\nonumber
\end{split}
\end{equation*}
The last inequality follows from the two preceding estimates and
$R^{-d+1}\leq R^{-d+\frac{3}{2}}$.
Here $C$ depends only on $d,U_0,U_1,U,c_0,\vartheta_0,s,p$
and $[\nabla A]_{W^{s,p}(U)}$, since
$\|\nabla A\|_{L^p(U)}^p\leq |U|\varepsilon_0^p$.
By the definition of $\mathfrak{k}_R$, we have 
\begin{equation*}
\begin{split}
\|\mathfrak{k}_R\|_{L^2(\mathbb R)}^2
&=2R^{-2d+3}
\begin{cases}
\displaystyle
\frac{\oldC{C:graph-diameter}^{sp-d+2}-R^{-(sp-d+2)}}{sp-d+2},
&sp>d-2,\\[7pt]
\displaystyle
\int_{R^{-1}}^{\oldC{C:graph-diameter}}
\frac{(\log(2+Rt))^{p-1}}t\,dt,
&sp=d-2,\\[7pt]
\displaystyle
R^{d-2-sp}\log(R\oldC{C:graph-diameter}),
&sp<d-2.
\end{cases}
\nonumber
\end{split}
\end{equation*}
In the critical case, $p\geq1$ gives
$$\int_{R^{-1}}^{\oldC{C:graph-diameter}}
\frac{(\log(2+Rt))^{p-1}}t\,dt
\leq
(\log(2+R\oldC{C:graph-diameter}))^{p-1}
\log(R\oldC{C:graph-diameter})
\leq
(\log(2+R\oldC{C:graph-diameter}))^p.$$
Since $\oldC{C:graph-diameter}\geq1$ is fixed and $R\geq2$,
$\log(2+R\oldC{C:graph-diameter})
\leq\log\oldC{C:graph-diameter}+\log(2+R)
\leq C\log(2+R)$.
Taking square roots therefore yields
$$\|\mathfrak{k}_R\|_{L^2(\mathbb R)}
\leq C
\begin{cases}
R^{-d+\frac{3}{2}},&sp>d-2,\\
R^{-d+\frac{3}{2}}(\log(2+R))^{\frac{p}{2}},&sp=d-2,\\
R^{-\frac{d-1+sp}{2}}(\log(2+R))^{\frac{1}{2}},&sp<d-2.
\end{cases}$$
Each of these bounds also absorbs $R^{-d+\frac{3}{2}}$, indeed 
$\log(2+R)\geq1$, and, when $sp<d-2$, we have
$R^{-d+\frac{3}{2}}
=R^{-\frac{d-1+sp}{2}}R^{-\frac{d-2-sp}{2}}
\leq R^{-\frac{d-1+sp}{2}}$. 
Substituting into the preceding reduction gives
\begin{equation*}
\begin{split}
&\int_{\mathbb S^{d-1}}\vartheta_0(\omega)
\Big|\int_Ue^{-iR\langle\omega,X(u)\rangle}b(u)\,du\Big|^2\,d\omega\leq C\|b\|_{L^2(U)}^2
\begin{cases}
R^{-d+\frac{3}{2}},&sp>d-2,\\
R^{-d+\frac{3}{2}}(\log(2+R))^{\frac{p}{2}},&sp=d-2,\\
R^{-\frac{d-1+sp}{2}}(\log(2+R))^{\frac{1}{2}},&sp<d-2.
\end{cases}
\nonumber
\end{split}
\end{equation*}
The constant $C$ has the dependence stated above and is independent of $R$ and $b$.
Together with the Dini estimate, this concludes the proof of the proposition.
\end{proof}

\section*{Use of AI tools}

The authors used Chat GPT and Clause to check the various versions of the paper, create diagrams, and for help with some calculations. The authors reviewed every change and take full responsibility for the content of the paper.

\printbibliography[title=References, heading=bibintoc]

\end{document}